\documentclass[12pt,a4paper]{article} %ГришаТамасян
\usepackage[cp866]{inputenc}
\usepackage[english,russian]{babel}
\usepackage{amsmath}
\usepackage{amsfonts}
\usepackage{longtable}
\usepackage{amssymb}
\usepackage{graphics}

\newcommand{\bc}{\begin{center}}
\newcommand{\br}{\begin{right}}
\newcommand{\ec}{\end{center}}
\newcommand{\be}{\begin{equation}}
\newcommand{\ee}{\end{equation}}

\newcommand{\pr}{\parallel}
\newcommand{\vl}{\mid}

\newcommand{\rar}{\rightarrow}

\newcommand{\p}{\partial}

\newcommand{\bint}{\mbox{int} \,}

\newcommand{\bco}{\mbox{co} \,}

\newcommand{\blim}{\mbox{lim} \,}

\newcommand{\meq}{\geq}

\newcommand{\al}{\alpha}

\newcommand{\eps}{\varepsilon}

\newtheorem{thm}{Теорема}[subsection]

\newtheorem{rem}{Замечание}[subsection]

\newtheorem{defi}{Определение}[subsection]
\newtheorem{cor}{Следствие}[subsection]
\newtheorem{ex}{Пример}[subsection]

\newtheorem{lem}{Лемма}[subsection]

\numberwithin{equation}{section}

\newcounter{secnum}

\makeindex \makeatletter
\renewcommand{\@oddhead}{\hfil --- \textrm{\thepage} ---\hfil}
\renewcommand{\@evenhead}{\hfil --- \textrm{\thepage} ---\hfil}
\makeatother
\newsavebox{\fmbox}

\begin{document}
\large{
\setcounter{page}{1}

\vspace{20cm}

\hspace{10cm}

{\bc  ПРУДНИКОВ Игорь Михайлович \ec}

\vspace{1cm}

{\bf\bc  РЕШЕНИЕ ПРОБЛЕМЫ АЛЕКСАНДРОВА А.Д. О ФУНКЦИЯХ,
ПРЕДСТАВИМЫХ В ВИДЕ РАЗНОСТИ ВЫПУКЛЫХ ФУНКЦИЙ \ec}

\vspace{12cm}

\vspace{3cm} {\bc  Смоленск,\,\,  2025 \ec}

\thispagestyle{empty}

\nopagebreak
}

%\pagestyle{plain} % нумерация вкл.
%\large{ \tableofcontents }
%\setcounter{page}{1}

\newpage

\large{\tableofcontents}
\newpage

{\bf  СПИСОК СОКРАЩЕНИЙ И ОБОЗНАЧЕНИЙ}

\vspace{1cm}

п.в. - почти всюду,

п.о. - положительно однородная (функция),

$\mathbb{R}^n - $ n-мерное евклидовое действительное пространство,

ПРВ - представимые в виде разности выпуклых,

DC - difference of convex (разность выпуклых),

$\pr x \pr - $ норма вектора $x$ в $ \mathbb{R}^n $,

$S^{n-1}_1(0)=\{q \in R^n \vl \pr q \pr = 1 \}$,

$B^{n}_1(0)=\{q \in R^n \vl \pr q \pr \leq 1 \}$,

$N^+ - $ множество целых положительных чисел,

$N - $ натуральное множество чисел,

$\bar{A} - $ замыкание множества $A,$

$\bco  A -$ выпуклая оболочка множества $A,$

$bd \, A - $ граница множества $A,$

$\mbox{int} A -$ внутренность множества $A$,

$\rho_H (A,B) -$ расстояние между множествами $A,B$ в метрике
Хаусдорфа,

$\nabla f(x) = f' (x) - $ производная функции $f(\cdot)$ в точке
$x$,

$ \mbox{dom}f -  $ область определения функции $f(\cdot)$,

$\Gamma_f= \{ (y,x) \in  \mathbb{R}^n | y = f(x), \,\,\, x \in
\mbox{dom}f \} -$ график функции $f(\cdot)$,

$\Pi-$ плоскость в $\mathbb{R}^n$,

$(a,b)-$ скалярное произведение векторов $a$ и $b$,

$K = \mbox{con} \{ a_1, a_2, \dots, a_k \} -$ конус $K$, равный
конической оболочке векторов $  a_1, a_2, \dots, a_k$,

\newpage

\bc {\bf   СПИСОК ОБОЗНАЧЕНИЙ КРИВЫХ, ИСПОЛЬЗУЕМЫХ В МОНОГРАФИИ}
\ec

%\vspace{1cm}

$Pr_{\Pi}(Q) -$ проекция множества $Q $ на  плоскость $\Pi$,

$\Re(D)-$ класс непрерывных кривых $r(\cdot)$, ограничивающих в
$D$ выпуклые компактные множества,

$\hat{\Re} - $ класс непрерывных кривых на поверхности единичного
шара с центром в нуле $B_1^3(0)=\{q \in \mathbb{R}^3 |\,\,\, \| q
\| \leq 1 \} $, получающихся в результате сечения единичной сферы
$S_1^2(0)=\{q \in \mathbb{R}^3 |\,\,\, \| q \| =1 \} $
произвольными плоскостями $\hat \Pi  $,

$\Re -$ класс непрерывных кривых $r(\cdot)$ на $ S_1^2(0)$, каждая
координата $r_i(\cdot)$ которой имеет не более трех участков
монотонности,

$\wp -$ множество непрерывных кривых $r(\cdot)$ на сфере
$S_1^{n-1}(0) $, проекции которых на координатные плоскости $\Pi$
есть кривые, ограничивающие выпуклые компактные множества,

$\wp(\Pi)-$  класс непрерывных кривых на сфере $S^{n-1}_1(0)$,
проекции которых на плоскость $ \Pi \in \mathbb{R}^3$, проходящую
через начало координат, ограничивают на $\Pi$ звездные множества,

$\hat \wp (\Pi)-$ класс непрерывных кривых $r(t)$ на сфере
$S^2_1(0)$, проекции которых на некоторую плоскость $ \Pi \in
\mathbb{R}^3$, проходящую через начало координат, ограничивают на
$\Pi$ выпуклые компактные множества,

$\wp(D)-$ класс непрерывных кривых на плоскости $XOY$,
ограничивающих звездные области в множестве $D$, замыкание которых
компактно,

$\varrho(D)-$ класс непрерывных кривых на плоскости  $XOY$,
ограничивающих выпуклые компактные множества в выпуклом множестве
$D$,

$\tilde \wp(D)-$ класс непрерывных кривых в  $D \subset
\mathbb{R}^n$, проекции которых на некоторую плоскость $\Pi$, $0
\in \Pi$, ограничивают звездные множества,

$\hat \wp(D)-$ класс непрерывных кривых в $D \subset
\mathbb{R}^n$, проекции которых на одну из координатных плоскостей
ограничивают звездные множества.

\newpage
\vspace{1cm}

{\bf АННОТАЦИЯ}

\vspace{1cm}

В монографии дается полное решение проблемы академика Александрова
А.Д. об условиях представимости функции в виде разности выпуклых.
Выпуклые функции нашли широкое применение в оптимизации из-за их
хороших свойств, в геометрии при построении внутренней геометрии
выпуклых поверхностей и в теории игр (теорема Нэша).

Так производная по направлению есть монотонно возрастающая функция
вдоль этого направления. Также локальная точка минимума выпуклой
функции является глобальной точкой минимума..Кроме того, выпуклая
функция является липшицевой функцией, которая, как известно, почти
всюду дифференцируемая, но также, как доказал Александров А.Д.
выпуклая функция почти всюду дважды дифференцируемая в области
определения.

Разность выпуклых - это следующий шаг. Функции в виде разности
выпуклых также находят широкое применение в оптимизации, геометрии
и теории коалиционных и бескоалиционных игр.  .

Для одномерной функции условия представимости в виде разности
выпуклых были получены более 100 лет назад. Функции, представимые
в виде разности выпуклых, называют ПРВ функциями или DC функциями
(difference of convex). Для многомерной функции от число
переменных два и более таких условий не было, хотя, конечно, ими
интересовались давно и интересуются поныне согласно публикациям.
Впервые эту проблему четко сформулировал и опубликовал академик
Александров А.Д. в 40-ых годах 20 века для построения внутренней
геометрии. В монографии даются необходимые и достаточные условия
представимости в виде разности выпуклых для многомерной функции.
Для этой цели рассматриваются звездные множества на плоскости и
кривые, ограничивающие их.   Приводится геометрическая
интерпретация этих условий путем введения поворота кривой через
поворот касательной к этой кривой.

В монографии собраны и даются ссылки на все наиболее важные
результаты по данной тематике, которых не так уж много. Поэтому
читатель может легко вникнуть в проблематику проблемы
представления функции в виде разности выпуклых, прочитав
монографию.

\vspace{0.5cm}

\section{ПРЕДИСЛОВИЕ}

\vspace{0.5cm}

Тематикой представления и условия для такого представления функции
в виде разности выпуклых автор стал интересоваться будучи
аспирантом Ленинградского государственного университета (ныне
Санкт-Петербургского университета). Условия представления функции
одной переменной в виде разности выпуклых известны давно, более
100 лет назад. Определение производной было введено Исааком
Ньютоном и Готфридом Лейбницем в 60-70-ых годах 17 века. Понятие
вариации функции (производной) было введено Жозефом Лагранжем в
1755 году, а выпуклые функции ввел И. Йенсен в своей работе 1906
года. Тогда же, в начале 20 века, были получены условия
представимости в виде разности выпуклых для функций одной
переменной. Конечно, ученых интересовали не только одномерные
функции (функции одной переменной), но и многомерные функции
(функции от многих переменных). С тех пор ответ на поставленный
вопрос не был найден. Проблему о представлении функции в виде
разности выпуклых впервые четко сформулировал и опубликовал
академик Александров А.Д. еще в 40-ых годах 20-ого века
\cite{aleksandrov1} для построения внутренней геометрии выпуклых
поверхностей и поверхностей, являющихся графиками функций,
представимых разностью выпуклых. Эта тематика тесно примыкает к
оптимизации \cite{strekal}. Были введены в оптимизации
\cite{demvas}, \cite{demrub} так называемые квазидифференцируемые
функции (КВД функции), у которых производная по направлению, как
функция от направления $g$, есть разность двух выпуклых
положительно-однородных первого порядка функций. КВД функции
получили широкое распространение в негладкой оптимизации.

Задача, поставленная Александровым А.Д., оказалась довольно
сложной. Многие математики разных специальностей интересовались и
интересуются условиями и алгоритмами представления функции в виде
разности выпуклых.  Вот уже более 100 лет проблема оставалась
открытой. В данной монографии автор решает эту проблему и приводит
необходимые и достаточные условия представимости произвольной
липшицевой функции в виде разности выпуклых на выпуклом компактном
множестве в $\mathbb{R}^n$, а также приводит алгоритмы
представления функции в виде разности выпуклых в некоторых частных
случаях, наиболее часто встречающихся на практике.

Книга состоит из статей, опубликованных автором в разные годы,
поэтому определения и термины могут повторяться в разных главах
монографии. С другой стороны, это удобно для читателя, для
которого не надо искать необходимые определения по всей
монографии. Помимо рассмотрения более общего многомерного случая
автор поставил перед собой задачу расширения класса кривых с
помощью которых можно сформулировать не только необходимые, но и
достаточные условия представимости функции в виде разности
выпуклых. Из доказательства теорем имеем предельного типа
алгоритмы такого представления.

\newpage
\vspace{1cm}

\section{О ПРЕДСТАВИМОСТИ ПОЛОЖИТЕЛЬНО ОДНОРОДНОЙ
ФУНКЦИИ ПЕРВОЙ СТЕПЕНИ ДВУХ ПЕРЕМЕННЫХ В ВИДЕ РАЗНОСТИ ВЫПУКЛЫХ
ФУНКЦИЙ}

\vspace{0.5cm}

В этой главе приведены необходимые и достаточные условия
представимости произвольной положительно однородной функции первой
степени двух переменных в виде разности выпуклых функций. Дана
также геометрическая интерпретация этих условий. Приведен алгоритм
такого представления, результатом которого есть последовательность
равномерно сходящихся на произвольном компакте, внутренности
которого принадлежит начало координат, выпуклых положительно
однородных первой степени многогранных функций. Объясняется связь
поставленной задачи с оптимизацией.

\noindent

\vspace{0.5cm}
\subsection{Введение}
\vspace{0.5cm}

 Задача об условиях представимости функции в виде
разности выпуклых интересна как для геометров, так и для
специалистов других специальностей
\cite{aleksandrov0}-\cite{Hartman}. Представление функции в виде
разности выпуклых нашло применение в оптимизации. Методы
оптимизации функций, представленных в виде разности выпуклых,
активно развиваются в \cite{strekal}.

Академик А.Д. Александров развил геометрию выпуклых поверхностей
\cite{aleksandrov0}. Следующим шагом были поверхности, являющиеся
графиками функций, представимых разностью выпуклых, так называемые
ПРВ функции. В западной литературе их называют DC функциями.

Дадим определение выпуклой функции.

\begin{defi} Функция $\varphi (\cdot)$ называется выпуклой в области
определения $D$, если для любых точек $x_1 , x_2 \in D$ и для
любых $\al_1, \al_2 \geq 0$ таких, что \mbox{$ \al_1 + \al_2 = 1
$} , выполняется неравенство
$$
\varphi (\al_1 x_1 + \al_2 x_2) \leq \al_1 \varphi(x_1) + \al_2
\varphi(x_2).
$$
\end{defi}
Из определения ясно, что область определения $D$ выпуклой функции
$-$ выпуклое множество, так как вместе с точками $x_1, x_2 \in D$
множеству $D$ принадлежит также точка $\al_1 x_1 + \al_2 x_2$, где
$\al_1 + \al_2 = 1, \,\, \al_1, \al_2  \geq 0.$

По определению функция $f(\cdot)$ называется {\em представимой в
виде разности выпуклых} (ПРВ функция), если верно представление
$$
f(\cdot)= f_1(\cdot) - f_2(\cdot),
$$
где $ f_1(\cdot), f_2(\cdot)- $ выпуклые функции.

К проблеме о представимости функции в виде разности выпуклых автор
пришел, решая задачу о нахождении условия, когда функция является
квазидифференцируемой в точке \cite{demvas}. Первый результат на
эту тему автор получил будучи аспирантом Ленинградского
государственного университета. Этот и другие результаты вошли в
кандидатскую (PhD) диссертацию автора \cite{lupikovdissertation}.
Интересно, повторение этого же результата (конечно, без ссылки)
автор нашел в сборнике тезисов международной конференции
"Конструктивный анализ и смежные вопросы", посвященный памяти
профессора В.Ф. Демьянова, проходившей 22-27 мая 2017 года в
Санкт-Петербурге \cite{pallaschke}. Свою кандидатскую диссертацию
автор защищал еще в 80-ых годах 20 века, то есть за 30 лет до
публикаций в сборнике.

Оптимизация функций многих переменных начала свое развитие с
выпуклых функций. Сюда относится линейное программирование
\cite{karpelsadovskii}. Также, как и в геометрии, следующим шагом
была теория оптимизации функций, представимых разностью выпуклых
функций, так называемых ПРВ (DC) функций, которые в свою очередь
принадлежат более широкому множеству квазидифференцируемых функций
\cite{demrub}.

Ясно, что ПРВ функция обязательно дифференцируема по направлениям,
так как дифференцируемы по направлениям выпуклые функции
\cite{pshenichnyi}.

Введем функцию $h(\cdot): \mathbb{R}^n \rar \mathbb{R} :$
$$
h(g) = f'(x_{0},g)= \frac{\p f(x_0)}{\p g} = \blim _{ \al \rar +0}
(f(x_{0}+\alpha g ) - f(x_{0}))/\alpha ,
$$
которая есть производная по направлению $g \in \mathbb{R}^n $
функции $f(\cdot)$ в точке $x_0$. Доказывается \cite{pshenichnyi},
что $h(\cdot)$ есть положительно однородная (п.о.) функция первого
порядка, т.е.
$$
h(\lambda g) = \lambda h( g)
$$
для любого $\lambda > 0$.

По определению \cite{demvas} функция $f(\cdot)$ называется
квазидифференцируемой (КВД) в точке $x_0$ , если
$$
h(g) = h_1(g) - h_2(g),
$$
где $h_1(\cdot), h_2(\cdot) - $ выпуклые функции.

Согласно двойственности Минковского \cite{demrub} любой выпуклой
конечной п.о. функции соответствует выпуклое компактное множество
в $\mathbb{R}^n$, называемое субдифференциалом этой функции в
нуле. Обозначим субдифференциалы функций $h_1(\cdot), h_2(\cdot)$
в нуле через $ \p h_1(0), \p h_2(0)$ соответственно. Тогда верны
равенства \cite{demrub}
$$
h_1(g) = \max_{ v \in \p h_1(0)} (v, g),\,\,\,\,\,\,\,\, h_2(g) =
\max_{ v \in \p h_2(0)} (v, g) \,\,\, \forall g \in \mathbb{R}^n.
$$
Здесь в правой части равенств стоят скалярные произведения
векторов $v$ и $g$: $(v,g)$.  Таким образом, вопрос о
квазидифференцируемости функции $f(\cdot)$ в точке $x_0$ сводится
к вопросу о представимости функции $h(\cdot)$ в виде разности
выпуклых п.о. функций.

Это была как раз та первая задача, за которую наряду с другими
задачами в других разделах математики взялся совсем еще молодой
автор статьи. Были получены необходимые и достаточные условия
представимости произвольной липшицевой п.о. первой степени функции
от двух переменных в виде разности выпуклых. Результат на данную
тему вошел в диссертацию, защищенную автором
\cite{lupikovdissertation}.

Напомним о результатах на тему о представлении функции в виде
разности выпуклых, известных автору, когда он был аспирантом и
решал поставленную задачу.

Необходимые и достаточные условия представимости функции одной
переменной в виде разности выпуклых, т.е. условия. когда функция
является ПРВ функцией, хорошо известны. Эти условия могут быть
записаны в следующем виде.

Пусть $x \rightarrow f(x): [a,b] \rightarrow \mathbb{R}$ -
произвольная липшицевая функция, т.е.
$$
\vl f(x) - f(y) \vl \leq L \vl x - y \vl \,\,\, \forall \,  x,y
\in [a,b].
$$
Известно, что множество $N_f$, где функция $ f(\cdot)$
дифференцируемая, есть множество полной меры на [a,b]. Для того,
чтобы функция $f(\cdot)$ была представима в виде разности выпуклых
функций, необходимо и достаточно, чтобы выполнялось условие
$$
         \vee (f'; a,b) < \infty ,
$$
где производные вычисляются там, где они существуют. Символ $\vee$
означает вариацию функции $ f'$ на отрезке [a,b].

В статье \cite{aleksandrov1}  А.Д.Александров задает вопрос о
представимости функции в виде разности выпуклых, если она является
таковой для любой прямой в области определения. Ответ на этот
вопрос отрицательный (см. \cite{hiriarturruty},
\cite{veselyzajicek} ).

Автор пошел по пути поиска кривых в $ \mathbb{R}^n$, с помощью
которых можно сформулировать необходимые и достаточные условия
представимости функции в виде разности выпуклых. Ясно, что для
положительно однородных функций такие кривые надо искать на
$S^{n-1}_1(0). $ Кривых в $ \mathbb{R}^n$ бесконечно много.
Возникает вопрос, какие кривые выбирать? И второй вопрос, как
формулировать необходимые и достаточные условия представимости
функции в виде разности выпуклых?

Согласно терминологии А.Д.Александрова под {\em многогранной
кусочно-линейной функцией} с конечным числом граней, определенной
в $\mathbb{R}^2$, будем понимать  такую функцию, график которой
состоит из конечного числа частей плоскостей, которые называются
гранями.

Введем понятие {\em двугранного угла}. Будем понимать под
двугранным углом функцию, график которой состоит из полуплоскостей
с общей граничной прямой, называемой ребром двугранного угла.
Иначе говоря, двугранный угол $-$ это функция, равная максимуму
или минимуму двух линейных функций.

В \cite{aleksandrov1} академик А.Д. Александров доказал, что
многогранная функция с конечным числом граней и функция, первая
производная которой липшицевая, являются ПРВ функциями в области
их определения.

Мы будем использовать в дальнейшем метод доказательства теоремы,
что любая многогранная функция $f(\cdot): \mathbb{R}^2 \rar
\mathbb{R}$ с конечным числом граней  является ПРВ функцией.
Повторим доказательство этой теоремы, взятое из
\cite{aleksandrov1}.

Рассмотрим все выпуклые двугранные углы, части графиков которых
принадлежат графику функции $f(\cdot)$. Просуммируем все такие
выпуклые двугранные углы. В итоге получим выпуклую многогранную
функцию $f_{1}(\cdot):\mathbb{R}^2 \rightarrow \mathbb{R}$.
Доказывается \cite{aleksandrov1}, что разность
$$
 f_{1}(\cdot) - f(\cdot) = f_{2}(\cdot)
$$
есть также выпуклая многогранная функция.

Действительно, для доказательства достаточно показать, что все
двугранные углы, части графиков которых принадлежат графику
функции $f_{1}(\cdot) - f(\cdot), $ являются выпуклыми. Для этого
покажем, что любая точка, лежащая на проекции ребра произвольного
двугранного угла функции $f_{1}(\cdot) - f(\cdot),$ имеет малую
окрестность, где функция $f_{1}(\cdot) - f(\cdot)$ выпуклая.

Если берем точку, в малой окрестности которой функция $f(\cdot)$
линейная, то локальная выпуклость разности $f_{1}(\cdot) -
f(\cdot)$ очевидна. Пусть берем точку, лежащую на проекции на
плоскость ребра выпуклого двугранного угла, часть графика которого
принадлежит графику функции $f(\cdot)$. Поскольку согласно
алгоритму этот же двугранный угол входит в сумму выпуклых
двугранных углов, образующих функцию $f_{1}(\cdot)$, то опять
разность $f_{1}(\cdot) - f(\cdot)$ будет локально выпуклой в
окрестности рассматриваемой точки. Если же точка лежит на проекции
ребра вогнутого двугранного угла, часть графика которого
принадлежит графику функции $f(\cdot)$, то $-f(\cdot) \,\, -$
локально выпуклая в окрестности этой точки, а поэтому разность
$f_{1}(\cdot) - f(\cdot)$ снова локально выпуклая в той же
окрестности. Из локальной выпуклости всех двугранных углов функции
$f_{1}(\cdot) - f(\cdot)$ следует ее выпуклость в $\mathbb{R}^2$.
Итак, доказано, что многогранная функция от двух переменных
является ПРВ функцией.

Докажем, что дифференцируемая функция с липшицевой производной
также является ПРВ функций.

Для доказательства воспользуемся замечательным свойством выпуклых
функций. Это свойство замечательно тем, что оно является
характерным только для выпуклых функций. Оказывается
\cite{pshenichnyi}, производная по произвольному направлению $g
\in \mathbb{R}^n$ выпуклой функции $\varphi(\cdot)$ в точке
$x_0+\al g,\,\, \al>0, $
$$
\varphi'(x_{0}+ \al g,g) = \frac{\p \varphi(x_0 + \al g)}{\p g} =
\blim _{ \al \rar +0} (\varphi(x_{0}+\alpha g + \tau g) -
\varphi(x_{0} + \al g))/ \tau ,
$$
есть монотонно возрастающая функция по $\al$.

Пусть задана дифференцируемая функция $f(\cdot): \mathbb{R}^n \rar
\mathbb{R}$ с липшицевой производной, т.е. для любых $x, y \in
\mathbb{R}^2 $ верно неравенство
$$
 \|  f'(x) - f'(y) \| \leq L \| x - y \|,
$$
где $L- $ константа Липшица, а справа и слева от знака неравенства
стоят нормы разности векторов. Покажем, что ее можно представить в
виде разности выпуклых.

Возьмем в качестве одной из выпуклых функций функцию $L \| x \|
^2$. Вычислим и сравним производные по направлению $g$ функции
$$
f_1(x) = L \| x \|^2 - f(x)
$$
в точках $x_0$ и $x_0 + \al g$, $\al >0$, где $ \| x \|^2 =(x,x) -
$ норма вектора $x$.

В точке $x_0$ производная по направлению $g$ равна
$$
f'_1(x_0,g)= 2L (x_0 ,g) -  f'(x_0,g)
$$
В точке $x_0 + \al g$ производная по направлению $g$ равна
$$
f'_1(x_0 + \al g,g)= 2L (x_0 + \al g ,g) -  f'(x_0 + \al g,g).
$$
Здесь, как и ранее, в круглых скобках через запятую стоят
скалярные произведения векторов.

Найдем разность этих выражений и определим ее знак
$$
f'_1(x_0 + \al g,g) - f'_1(x_0,g) = 2L (x_0 + \al g ,g) - 2L (x_0
,g) -
$$
$$
- f'(x_0 + \al g,g) + f'(x_0,g)\meq 2L \al \| g \|^2 - L \al \| g
\|^2 \meq 0,
$$
что по свойству, указанному выше, следует, что $f_1(\cdot) - $
выпуклая функция. Таким образом, доказано, что функция $f(\cdot)$
представима в виде разности выпуклых функций $L \| x \|^2$ и $
f_1(x)$:
$$
f(x) = L \| x \|^2 - f_1(x).
$$
Заметим, что данное свойство доказано для функции, определенной в
$\mathbb{R}^n$.

Итак, для представления дифференцируемой функции $f(\cdot)$ в виде
разности выпуклых надо знать  константу Липшица ее производной.
Если функция $f(\cdot)$ дважды непрерывно дифференцируемая, то в
качестве константы $L$ можно взять верхнюю грань нормы матрицы
вторых смешанных производных, т.е.
$$
\sup_{x \in D} \| f''(x)\| = \sup_{x \in D} \| \nabla ^2 f(x) \| =
L,
$$
если матрица $ f''(\cdot) =  \nabla ^2 f(\cdot) $ существует  и ее
норма ограничена сверху во всей области определения $D$ функции
$f(\cdot).$

В данном параграфе будут записаны необходимые и достаточные
условия представимости положительно однородной (п.о.) функции  в
виде разности выпуклых. Начнем с положительно однородной функции
первой степени.

\vspace{1cm}

\subsection {Положительно однородная функция первой степени.
Двухмерный случай}

\vspace{0.5cm}

Пусть $\varphi (\cdot): \mathbb{R}^2 \rar \mathbb{R},\,\,-$ п.о.
функция, удовлетворяющая условию Липшица с константой $L$, т.е.
$$
\varphi(\lambda x ) = \lambda \varphi(x) \,\,\, \forall x \in
\mathbb{R}^n, \forall \lambda >0.
$$
$$
\vl \varphi(x) - \varphi(y) \vl \leq L \| x - y \| \,\,\, \forall
x,y  \in  \mathbb{R}^n.
$$
Обозначим через
$$
\Phi(t)=\varphi(\cos t, \sin t) = \varphi(r(t)), \,\,\, r(t)=(\cos
t, \sin t),
$$
где $t \in [0, 2\pi]- $  естественная параметризация единичной
окружности с центром начале координат \cite{pogorelov1}.

Функция $\Phi (\cdot) $ также будет удовлетворять условию Липшица
с той же константой $L$. Это следует из очевидных неравенств
$$
|\Phi(t_1)-\Phi(t_2)|= |\varphi(r(t_1)) - \varphi(r(t_2))| \leq
$$
$$
\leq L \| r(t_1) - r(t_2) \| = 2L \sin \frac{ | t_1 - t_2 |}{2}
\leq L |t_1 - t_2 |.
$$
Следовательно, функция $\Phi(\cdot)$ почти всюду (п.в.)
дифференцируема по $t$ на отрезке $[0, 2\pi]$
\cite{kolmogorovfomin}.

Обозначим, как это принято в литературе \cite{kolmogorovfomin},
через $\vee (\Phi'; 0,2\pi) $ вариацию производной $ \Phi'(\cdot)
$ на отрезке $[0, 2\pi],$  предполагая, что рассматриваются точки,
где производная существует.

\begin{thm} Для того, чтобы липшицевая п.о. первой степени функция
$\varphi (\cdot): \mathbb{R}^2 \rar \mathbb{R}$ была представима в
виде разности выпуклых функций (была ПРВ функцией), необходимо и
достаточно, чтобы
$$
 (\exists c(\varphi)>0):  \;\;\; \vee (\Phi'; 0,2\pi)<  c(\varphi),
$$
где $\Phi(t)=\varphi(r(t)), \,\,r(t)=(\cos t, \sin t), \,\, t \in
[0, 2\pi] $. \label{difconvgomogthm1}
\end{thm}
{\bf Доказательство.} {\bf Необходимость}. Пусть функция
$\varphi(\cdot)$ представима в виде разности выпуклых п.о. функций
$\varphi_1(\cdot),  \varphi_1(\cdot):$
$$
\varphi(q) = \varphi_1(q) -  \varphi_2(q) \,\,\, \forall
q\in\mathbb{R}^2.
$$
Определим функции
$$
\Phi_1(t)= \varphi_1(r(t)), \,\, \Phi_2(t)= \varphi_1(r(t)),
 \,\,\, r(t)=(\cos t, \sin t), \,\, \forall t \in [0, 2\pi].
$$
Очевидно, что
$$
\Phi(t)= \Phi_1(t) - \Phi_2(t) \,\,\,\forall t \in [0, 2\pi].
$$
По свойству вариации \cite{kolmogorovfomin}
$$
\vee (\Phi'; 0,2\pi) \leq  \vee (\Phi'_1; 0,2\pi) + \vee
(-\Phi'_2; 0,2\pi) = \vee (\Phi'_1; 0,2\pi) + \vee (\Phi'_2;
0,2\pi).
$$
Следовательно, для доказательства необходимости достаточно
доказать неравенство
$$
\vee (\Phi'_1; 0,2\pi)<  c(\varphi_1).
$$
Возьмем точки $\{ t_i \}, i \in 1:s, $ достаточно равномерно
расположенные на окружности $S^1_1(0) = \{ (x,y) \in \mathbb{R}^2
| x^2+y^2 = 1 \}, $ где функция $\Phi_1(\cdot)$ дифференцируема.
Это всегда можно сделать,  так как  $\Phi_1 (\cdot) $ липшицева на
окружности $S^1_1(0),$ поскольку выпуклая функция
$\varphi_1(\cdot) $ липшицева \cite{demvas}, константу Липшица
которой обозначим через $L_1$.

Вариация по определению \cite{kolmogorovfomin}  есть верхний
предел сумм \be \sup_{ \{ t_i \} }  \sum_1^s | \Phi'_1(t_i) -
\Phi'_1(t_{i+1}) |.
\label{difconvgomog1}
\ee
Посчитаем
производные
$$
 \Phi'_1(t_i) = (\varphi'_1(r(t_i)), p(t_i)),
$$
где $p(t_i)= r'(t_i), \,\, \| p(t_i) \|=1,$ и вычислим модуль
разности
$$
| \Phi'_1(t_{i+1}) - \Phi'_1(t_{i}) | = |(\varphi'_1(r(t_{i+1})),
p(t_{i+1})) - (\varphi'_1(r(t_i)), p(t_i)) | =
$$
$$
= | (\varphi'_1(r(t_{i+1})), p(t_{i+1})) - (\varphi'_1(r(t_i)),
p(t_{i+1})) +
$$
$$
+(\varphi'_1(r(t_i)), p(t_{i+1})) - (\varphi'_1(r(t_i)), p(t_i)) |
\leq
$$
$$
\leq | (\varphi'_1(r(t_{i+1})), p(t_{i+1})) - (\varphi'_1(r(t_i)),
p(t_{i+1})) | +
$$
$$
+ | (\varphi'_1(r(t_i)), p(t_{i+1})) - (\varphi'_1(r(t_i)),
p(t_i)) | \leq
$$
$$
\leq  \| \varphi'_1(r(t_{i+1})) - \varphi'_1(r(t_i)) \|  \, \|
p(t_{i+1}) \| + \| p(t_{i+1}) - p(t_i) \| \, \| \varphi'_1(r(t_i))
\| \leq
$$
\be \leq \| \varphi'_1(r(t_{i+1})) - \varphi'_1(r(t_i)) \| + L_1
\, \| p(t_{i+1}) - p(t_i) \|, \label{difconvgomog2} \ee так как
$\| p(t_{i+1}) \| = 1, \, \| \varphi'_1(r(t_i)) \| \leq L_1.$

Сумма
$$
\sum_i \| p(t_{i+1}) - p(t_i) \|
$$
не превосходит длины окружности $S^1_1(0)$, равное $2\pi.$

Докажем ограниченность сверху суммы
$$
\sum_i \| \varphi'_1(r(t_{i+1})) - \varphi'_1(r(t_i)) \|.
$$
Так как функция $\varphi_1 (\cdot) -$ п.о. выпуклая функция, то
согласно принципу двойственности Минковского
\cite{kutateladzerubinov} имеет место равенство
$$
 \varphi_1(r(t)) = \max_{v \in \partial \varphi_1(0)}
 (v, r(t))=(v(t),r(t)), \;\; v(t) \in \partial \varphi_1(0),
$$
где $\partial \varphi_1 (0)$ - субдифференциал функции
$\varphi_1(\cdot)$ в нуле.

Субдифференциал $-$ это выпуклая оболочка всех обобщенных
градиентов выпуклой функции $\varphi_1(\cdot)$ в нуле, который
является выпуклым компактным множеством. В точках $r(t) \in
S^1_1(0)$, где функция $\varphi_1(\cdot)$ дифференцируема,
граничный вектор $v(t) \in
\partial \varphi_1(0)$ в написанном выше равенстве единственный и
равен $ v(t)= \varphi'_1(r(t))= \nabla \varphi_1(r(t))$. По
свойству субдифференциала \cite{demvas} нормальный вектор к
границе множества $\partial \varphi_1(0)$ в точке $v(t)$ равен $
r(t).$

Тогда
$$
\sum_i \| \varphi'_1(r(t_{i+1})) - \varphi'_1(r(t_i)) \| = \sum_i
\| v(t_{i+1}) - v(t_i) \|.
$$
Соединив пары векторов $v(t_i),v(t_{i+1}), i \in 1:s,$ отрезком,
мы получим замкнутую ломаную, являющуюся границей замкнутого
выпуклого многоугольника, вписанного в  $\partial \varphi_1(0)$.

Поэтому супремум суммы
$$
\sup_{t_i} \sum_i \| \varphi'_1(r(t_{i+1})) - \varphi'_1(r(t_i))
\|
$$
не превосходит длины кривой $L_{\partial \varphi_1(0)}$,
ограничивающей замкнутое выпуклое множество $\partial
\varphi_1(0)$. Длину кривой $L_{\partial \varphi_1(0)}$ обозначим
через $P(L_{\partial \varphi_1(0)}).$

Окончательно, супремум суммы  (\ref{difconvgomog1} ) согласно
неравенствам (\ref{difconvgomog2}) не превосходит
$$
P(L_{\partial \varphi_1(0)}) + 2\pi L_1 .
$$
Поэтому в случае, когда $\varphi(\cdot) = \varphi_1(\cdot) -
\varphi_2(\cdot),$ можно записать
$$
\vee (\Phi'; 0,2\pi) \leq  \vee (\Phi'_1; 0,2\pi) + \vee (\Phi'_2;
0,2\pi) \leq
$$
$$
\leq P(L_{\partial \varphi_1(0)})+ P(L_{\partial \varphi_2(0)}) +
2\pi L_1 + 2\pi L_2  ,
$$
где $\partial \varphi_1(0), \partial \varphi_2(0), L_1, L_2 -   $
субдифференциалы в нуле и константы Липшица функций $
\varphi_1(\cdot)$  и $  \varphi_2(\cdot)$ соответственно. Здесь
$P(L_{\partial \varphi_1(0)}), P(L_{\partial \varphi_2(0)}) - $
длины кривых, ограничивающих субдифференциалы $\partial
\varphi_1(0), \partial \varphi_2(0)$. Необходимость доказана.

{\bf Достаточность}. Пусть условия теоремы \ref{difconvgomogthm1}
выполняются. Разобьем единичный круг $B_1^2(0)= \{ (x,y) \in
\mathbb{R}^2 | x^2+y^2 \leq 1 \} $ точками $ \{ r(t_i) \} \in
S^1_1(0), \, i \in 1:m, $  на $m$ равных секторов. Построим по
этому разбиению п.о. многогранную функцию $\varphi_m(\cdot)$
следующим образом.

Значения функции $\varphi_m(\cdot)$ на векторах $0, r(t_i),
r(t_{i+1})$ равны
$$
\varphi_m(0)=0, \, \varphi_m(r(t_i))= \varphi(r(t_i)), \,
\varphi_m(r(t_{i+1}))= \varphi(r(t_{i+1})).
$$
Во всех остальных точках $i-$ого  сектора, определяемого векторами
$0,r(t_i), r(t_{i+1}),$ функцию $\varphi_m(\cdot)$ определим
линейным образом, т.е. для вектора
$$
r=\lambda_1  r(t_i) +
\lambda_2 r(t_{i+1}), \lambda_1, \lambda_2 >0,
$$
верно равенство
$$
\varphi_m(r)= \lambda_1\varphi_m(r(t_i)) + \lambda_2
\varphi_m(r(t_{i+1})).
$$
Так сделаем для всех $m$ секторов, на которые мы разбили круг
$B_1^2(0)$. Линейную функцию, график которой совпадает в $i-$ ом
секторе с графиком функции $\varphi_m(\cdot)$, обозначим через
$\pi_i(\cdot).$

Введем некоторые определения и понятия, используемые в дальнейшем.
Будем понимать под {\em градиентом плоскости} градиент функции,
график которой есть эта плоскость. {\em Двугранный угол} по
определению есть непрерывная функция, график которой состоит из
двух полуплоскостей с общей прямой, называемой {\em ребром, } а
полуплоскости называются  {\em гранями двугранного угла}.

Первое, что мы докажем $-$ это равномерную липшицевость функций
$\varphi_m(\cdot)$ по $m$. Для этого требуется доказать, что
градиенты линейных функций $\pi_i(\cdot)$ ограничены сверху
константой, не зависящей от $i, m.$

Рассмотрим различные расположения на плоскости градиента $\nabla
\pi_i = \pi'_i$. Если градиент $\nabla \pi_i $ находится на
границе $i-$ ого сектора и расположен вдоль вектора $r(t_i)$, то
$$
(\nabla \pi_i, l_i)= \| \nabla \pi_i \| =\frac{\p \varphi(0)}{\p
l_i}\leq L,
$$
где $l_i= {r(t_i)}/{\| r(t_i) \|},$ $ {\p \varphi(0)}/{\p l_i}- $
производная по направлению  $ l_i $ функции $ \varphi(\cdot) $ в
точке $0$. Аналогичные рассуждения проводим, если градиент $\nabla
\pi_i $ находится на границе $i-$ ого сектора и расположен вдоль
вектора $r(t_{i+1})$.

Если вектор $\nabla \pi_i $ находится внутри $i-$ ого сектора, то
поскольку производные по направлению $ l_i $ функций $ \pi_i
(\cdot) $ и $ \varphi (\cdot)$ совпадают, то
$$
 \frac{\p \pi_i(0)}{\p l_i} =   (\nabla \pi_i , l_i) = \frac{\p \varphi(0)}{\p
l_i}\leq L.
$$
Перепишем скалярное произведение в виде произведения норм векторов
на косинус угла $\al_i $  между ними
$$
\| \nabla \pi_i \| \,\, \| l_i  \| \,\,\cos \al_i   =  \| \nabla
\pi_i \| \,\, \cos \al_i \leq L.
$$
Откуда
$$
\| \nabla \pi_i \| \leq \frac{L}{\cos \al_i} \leq \frac{L}{\cos
\frac{2\pi}{m}},
$$
так как при больших $m$ угол $\al_i \leq 2\pi / m$. Очевидно, что
правая часть неравенства равномерно ограничена для всех больших
$m$.

Если градиент $ \nabla \pi_i  $ не принадлежит $ i- $ ому сектору,
то проведем прямую $l$ с направляющим единичным вектором $l_i$,
образующую с вектором $q_i = \nabla \pi_i / \| \nabla \pi_i \| $
угол $\al_i$, не больший $2\pi / m$. Пусть прямая $l$ пересекает
лучи, на которых лежат векторы $r(t_i)$ и $r(t_{i+1})$, в точках
$A$ и $B$. Параметризуем прямую $l$ естественным образом.
Параметры точек $A$ и $B$ на прямой $l$ обозначим через $a$ и $b$.
Тогда
$$
\pi_i(b) - \pi_i(a)= \int_a^b \frac{\p \pi_i(\tau)}{\p l_i} d \tau
= \int_a^b \frac{\p \varphi(\tau)}{\p l_i} d \tau = \varphi(b) -
\varphi(a).
$$
Отсюда заключаем, что существует точка $s_i \in [a, b]$, для
которой
$$
\frac{\p \pi_i(s_i)}{\p l_i}= \| \nabla \pi_i \| \,\, \| l_i  \|
\,\,\cos \al_i = \| \nabla \pi_i \| \,\,\cos \al_i \leq \frac{\p
\varphi(s_i)}{\p l_i} \leq L.
$$
Проводим рассуждения, аналогичные  приведенным выше. Получим ту же
оценку для нормы $\| \nabla \pi_i \|.$ Тем самым мы доказали
равномерную липшицевость по $m$  функций $ \varphi_m (\cdot).  $

Будем представлять функцию $ \varphi_m (\cdot)$ в виде разности
выпуклых. Для этого воспользуемся алгоритмом академика А.Д.
Александрова для представления многогранной функции в виде
разности выпуклых, описанным в \cite{aleksandrov1}.

Введем определение {\em вариации двугранного угла}. Под вариацией
двугранного угла будем понимать максимальную вариацию производной
этой функции вдоль некоторой прямой. Нетрудно видеть, что
максимальная вариация производной  будет для прямой,
перпендикулярной проекции на плоскость $\mathbb{R}^2$ ребра
двугранного угла, которую обозначим через $l$. Поясним сказанное.

Пусть направляющий единичный вектор прямой $l$ есть вектор $ q \in
\mathbb{R}^2 $. В этом случае вариация  двугранного угла равна $\|
a_1 - a_2 \|, $ где $a_1, a_2 \in \mathbb{R}^2 - $ градиенты
плоскостей (полуплоскостей) $\pi_1, \pi_2$, образующих двугранный
угол. Действительно, поскольку производные  по направлению $q$
плоскостей $\pi_1, \pi_2$ совпадают, то
$$
       (\nabla \pi_1, q) = (\nabla \pi_2, q).
$$
Отсюда следует, что вектор $ a_1 - a_2  $ перпендикулярен прямой
$l$. Поэтому для всех прямых, не перпендикулярных прямой $l$,
вариация производной двугранного угла будем меньше $\| a_1 - a_2
\|$. Итак, вариация двугранного угла, записанная через производные
по направлению, равна
$$
 | (\nabla \pi_1, p) - (\nabla \pi_2, p) | =
 | (\nabla \pi_1 - \nabla \pi_2 , p) |,
$$
где $p \in \mathbb{R}^2 -$ единичный вектор, перпендикулярный
вектору $q$.

Сделаем оценку вариации двугранного угла через вариацию
производной функции $\Phi(\cdot)$. Рассмотрим секторы $i, i+1,$ по
которым построены плоскости $\pi_i, \pi_{i+1}$, являющиеся гранями
$i-$ ого двугранного угла. Путь сектор $i$ определяют векторы
$r(t_i), r(t_{t+1})$, а сектор $i+1 - $ векторы $ r(t_{t+1}),
r(t_{i+2}).$ Для определенности предположим, что
$$
(\nabla \pi_{i+1}, p(t_{i+1})) \meq (\nabla \pi_i, p(t_{i+1})),
$$
где $p(t_{i+1}) = r'(t_{i+1})$.  Очевидны равенства
$$
\pi_i(t_{i+1}) - \pi_i(t_{i}) = \Phi (t_{i+1}) - \Phi (t_{i}) =
\int_{t_i}^{t_{i+1}}  \Phi'(\tau)  d \tau = \int_{t_i}^{t_{i+1}}
\frac{d \pi_i(r(\tau)) }{d \tau} d \tau.
$$
Отсюда следует, что существует точка $\tau_1 \in [t_i, t_{i+1}], $
для которой
$$
\frac{d \pi_i(r(\tau_1))}{ d t } \meq \Phi'(\tau_1).
$$
Аналогично для отрезка $[t_{i+1}, t_{i+2} ]$:
$$
\pi_{i+1}(t_{i+2}) - \pi_{i+1}(t_{i+1}) = \Phi (t_{i+2}) - \Phi
(t_{i+1}) =
$$
$$
=\int_{t_{i+1}}^{t_{i+2}}  \Phi'(\tau)  d \tau =
\int_{t_{i+1}}^{t_{i+2}} \frac{d \pi_{i+1}(r(\tau)) }{d \tau} d
\tau.
$$
Отсюда заключаем, что существует точка $\tau_2 \in [t_{i+1},
t_{i+2} ], $ для которой
$$
\frac{d \pi_{i+1}(r(\tau_2))}{ d t } \leq \Phi'(\tau_2).
$$
Получим оценку сверху для вариации $i -$ ого двугранного угла с
учетом сделанного предположения.
$$
\| \nabla \pi_{i+1} - \nabla \pi_i \| = (\nabla \pi_{i+1},
p(t_{i+1})) - (\nabla \pi_i, p(t_{i+1}))=
$$
$$
=\frac{d \pi_{i+1}(r(\tau_2))}{ d t } - \frac{d
\pi_{i}(r(\tau_1))}{ d t } +
$$
$$
+ (\nabla \pi_{i+1}, p(t_{i+1}))-\frac{d \pi_{i+1}(r(\tau_2))}{ d
t } - (\nabla \pi_i, p(t_{i+1}))+\frac{d \pi_{i}(r(\tau_1))}{ d t
}.
$$
Производные  функций $\pi_i(\cdot), \pi_{i+1}(\cdot) $ по $t$
можно переписать в ином виде
$$
\frac{d \pi_{i+1}(r(\tau_2))}{ d t }= (\nabla \pi_{i+1},
p(\tau_2)); \,\, \frac{d \pi_{i}(r(\tau_1))}{ d t }= (\nabla
\pi_i, p(\tau_1)),
$$
где $ p(\tau_1)=r'(\tau_1), \, p(\tau_2)=r'(\tau_2).$ Отсюда
получаем требуемую оценку
$$
\| \nabla \pi_{i+1} - \nabla \pi_i \| \leq \Phi'(\tau_2) -
\Phi'(\tau_1) + \| \nabla \pi_{i+1} \| \, |\tau_2 - t_{i+1} | + \|
\nabla \pi_{i} \| \, |\ \tau_1 - t_{i} | \leq
$$
$$
\leq \int_{t_{i}}^{t_{i+2}}  \Phi'(\tau)  d \tau + \frac{L}{\cos
\frac{2\pi}{m}} \,\, | t_{i+2} - t_i |.
$$
Теперь опишем алгоритм представления функции в виде разности
выпуклых. Рассмотрим все выпуклые двугранные углы, части графиков
которых принадлежат графику функции $\varphi_m(\cdot)$.
Просуммируем все такие выпуклые двугранные углы. В итоге получим
п.о. выпуклую многогранную функцию, которую обозначим
$\varphi_{1m}(\cdot)$. Разность
$$
\varphi_{2m}(\cdot) = \varphi_{1m}(\cdot) - \varphi_{m}(\cdot),
$$
как было ранее доказано, является снова п.о. выпуклой
многогранной функцией. Также мы доказали, что функции
$\varphi_{1m}(\cdot)$, $\varphi_{2m}(\cdot)$ являются равномерно
липшицевыми по $m$  с нулевым значением в начальной точке. Поэтому
из последовательности функций $\varphi_{1m}(\cdot)$,
$\varphi_{2m}(\cdot)$ можно выбрать равномерно сходящиеся на
единичном круге $B_1^2(0)$. Обозначим их равномерные пределы по
$m$ через
$$
\lim_m \varphi_{1m}(x,y) = \varphi_1 (x,y), \,\,\,\, \lim_m
\varphi_{2m}(x,y) = \varphi_2 (x,y) \,\,\, \forall (x,y) \in
B_1^2(0).
$$
Переходя в равенстве
$$
\varphi_{m}(x,y) = \varphi_{1m}(x,y) - \varphi_{2m}(x,y) \,\,\,
\forall (x,y) \in B_1^2(0)
$$
к пределу по $m \rar \infty$, получим требуемое представление
$$
\varphi(x,y) = \varphi_{1}(x,y) - \varphi_{2}(x,y) \,\,\, \forall
(x,y) \in B_1^2(0),
$$
где $\varphi_{1}(\cdot)$, $\varphi_{2}(\cdot) -$ п.о. выпуклые
функции, как равномерный поточечный предел последовательности
выпуклых функций. Достаточность доказана, а вместе с ней и
теорема.

\begin{rem}
Доказанная теорема дает конструктивный путь представления п.о.
функции двух переменных в виде разности выпуклых.
\end{rem}

\newpage
\vspace{1cm}

\section{УСЛОВИЯ ПРЕДСТАВИМОСТИ ПОЛОЖИТЕЛЬНО ОДНОРОДНОЙ
ФУНКЦИИ $M -$ ОЙ СТЕПЕНИ ДВУХ ПЕРЕМЕННЫХ В ВИДЕ РАЗНОСТИ ВЫПУКЛЫХ
ФУНКЦИЙ}

\vspace{1cm}

В данной главе приведены необходимые и достаточные условия
представимости произвольной положительно однородной функции $m -$
ого порядка двух переменных в виде разности выпуклых функций. Дана
также геометрическая интерпретация этих условий. Приведен алгоритм
такого представления, результатом которого есть последовательность
равномерно сходящихся на произвольном компакте, внутренности
которого принадлежит начало координат, выпуклых положительно
однородных $m - $ ой степени функций.

\noindent

\vspace{0.5cm}

\subsection{Введение}

\vspace{0.5cm}

История вопроса о представлении функции в виде разности выпуклых
берет свое начало в работах \cite{aleksandrov1},
\cite{aleksandrov2}. На тему о представлении функции в виде
разности выпуклых было написано много работ \cite{aleksandrov0} -
\cite{lupikovdissertation} авторами разных специальностей.
Функции, представимые в виде разности выпуклых (ПРВ функции) нашли
применение в оптимизации \cite{strekal}. ПРВ функции в западной
литературе называют DC (difference of convex) функциями.

Важный результат после работ академика А.Д. Александрова был
получен автором в его кандидатской диссертации
\cite{lupikovdissertation}, где приведены необходимые и
достаточные условия представимости положительно однородной первой
степени функции двух переменных в виде разности выпуклых. Этот
результат напрямую связан с оптимизацией, так как эти условия
являются необходимыми и достаточными условиями
квазидифференцируемости функции в точке \cite{demvas}. Метод
доказательства был обобщен позднее автором для функций от
произвольного количества переменных. Результат, подобный
результату автора \cite{lupikovdissertation}, можно найти без
ссылки на первоисточник  в статье Pallaschke D \cite{pallaschke}.
Возможно, Pallaschke D не знал о работах автора

Вышла статья \cite{ginchev}, где речь идет об условиях
представимости функции в виде разности выпуклых в бесконечномерных
пространствах.

В \cite{ginchev} доказана теорема, дающая необходимые и
достаточные условия представимости функции $f(\cdot)$ в линейном
бесконечномерном пространстве $X$. Основной результат этой статьи
следующий.

Пусть $X -$ линейное пространство, $\Omega -$ выпуклое множество в
$X$ и $f(\cdot) : \Omega \rar R -$ произвольная функция. $f$
является разностью двух выпуклых функций, если и только если
найдутся (конечное или бесконечное) множество индексов $I$ и
множество выпуклых функций $h_i : \Omega \rar R, i \in I,$ таких,
что сумма $\sum_{i \in I} h_i(x) $  существует и конечна в
$\Omega$ и для любой пары точек $a, b \in \Omega$  существует
множество $J \subset I$ такое, что $f + \sum_{i \in J}  h_i$
является выпуклой функцией на сегменте $[a; b]$.

Непонятен основной результат указанной статьи. Почему нельзя в
формулировке теоремы говорить об одной выпуклой функции $f_1
(\cdot)$, для которой $f + f_1$ выпуклая на $\Omega$, вместо
семейства функций $\{ h_i (\cdot) \}$? Зачем выбирать из семейства
функций $\{ h_i (\cdot) \}$ подмножество $\{ h_j (\cdot) \}$, если
сумма $\sum_{i \in I} h_i(x) $ конечная и выпуклая на всем
$\Omega$?

Все это говорит о том, что вопрос о нахождении необходимых и
достаточных условий представимости функции в виде разности
выпуклых - это довольно сложная задача. В данном разделе задача
решается для п.о. функции степени $m$ от двух переменных, где $m
-$ натуральное число.

Напомним полученный ранее результат.
\begin{thm} Для того, чтобы липшицевая п.о. степени 1 функция $\varphi (\cdot):
\mathbb{R}^2 \rar \mathbb{R}$ была представима в виде разности
выпуклых функций (была ПРВ функцией), необходимо и достаточно,
чтобы
$$
 (\exists c(\varphi)>0):  \;\;\; \vee (\Phi'; 0,2\pi)<  c(\varphi),
$$
где $\Phi(t)=\varphi(r(t)), \, r(t)=(\cos t, \sin t), \,  t \in
[0, 2\pi] $. \label{difconvgomogMDegreethm1}
\end{thm}
Здесь, как и прежде, используется общепринятое обозначение
вариации $\vee  $  производной $\Phi'(\cdot)$ на отрезке $[0,
2\pi]. $

\vspace{0,5cm}

\subsection {Положительно однородная функция $m -$ ого порядка от
двух переменных}

\vspace{0.5cm}

Пусть задана липшицевая п.о. функция степени $m: q \rar
\varphi(q): \mathbb{R}^2 \rar \mathbb{R}$, $m >0$, т.е. $
\varphi(\lambda q )=  \lambda^m \varphi(q), \lambda \in
\mathbb{R}, \lambda >0, m - $ натуральное число. Обозначим через
$S_1^1(0)= \{ z \in \mathbb{R}^2 | \| z \| =1 \} - $ единичная
окружность с центром в начале координат. Пусть $ r(t)=(\cos t,
\sin t), t \in [0, 2\pi]- $ естественная параметризация единичной
окружности. Введем функцию
$$
\Phi(t)=\varphi(r(t)) = \varphi(\cos t, \sin t), \,\,\forall t \in
[0, 2\pi].
$$
Далее попытаемся свести наш случай к однородному случаю первой
степени. А именно:
\begin{enumerate}
\item для любой выпуклой локально липшицевой функции $z \rar f(z):
\mathbb{R}^2 \rar \mathbb{R}, f(0,0)=0,$ мы построим п.о. степени
1 выпуклую функцию $z \rar \psi(z): \mathbb{R}^2 \rar \mathbb{R},$
принимающую на $S_1^1(0)$ те же значения, что и функция
$f(\cdot)$;

\item обратно, для любой выпуклой п.о. степени 1 функции
$\psi(\cdot),$  принимающей на $S_1^1(0)$ положительные значения,
построим выпуклую п.о. степени $m$ функцию $\varphi(\cdot) $,
принимающую на $S_1^1(0)$ те же значения, что и функция
$\psi(\cdot).$
\end{enumerate}

Итак, пусть задана произвольная локально липшицевая выпуклая
функция $z \rar f(z): \mathbb{R}^2 \rar \mathbb{R}, f(0,0)=0,$ у
которой начало координат $-$ точка минимума.

\begin{lem}  Пусть $(x,y) \rightarrow f(x,y): \mathbb{R}^2
\rightarrow \mathbb{R}-$ непрерывная выпуклая функция и $r(t) =
(\cos t, \sin t), t \in [0,2 \pi].$ Тогда существует константа
$c(f)>0,$ что
\begin{equation}
   \vee (\Phi' ; 0, 2 \pi) \leq  c(f),
\label{difconvGomogMDegree1} \end{equation} где $\Phi(t) =
f(r(t)), t \in [0,2 \pi].$ \label{difconvGomogMDegreelem1}
\end{lem}

{\bf  Доказательство.} На начальном этапе будем считать, что
$f(\cdot,\cdot)$ дважды непрерывно дифференцируемая функция,
которая принимает неотрицательные значения и начало координат $-$
ее точка минимума. Таким образом точка минимума  $0=(0,0)$
принадлежит внутренности шара $ B^2_1(0)= \{ z \in  \mathbb{R}^2 |
\| z \| \leq 1 \} .$

Построим для функции $f(\cdot,\cdot)$ п.о. степени 1 функцию
$\psi(\cdot),$  которая на $r(\cdot)$ принимает значения, равные
$f(r(\cdot)).$ Покажем, что $\psi(\cdot)$ - выпуклая.

Рассмотрим функцию
$$
f_{\varepsilon}(x,y)=f(x,y)+\varepsilon ( \mid \mid x \mid \mid^2
+ \mid \mid y \mid \mid^2  ), \,\,\, \varepsilon>0.
$$
Разобьем отрезок $[ 0, 2\pi] $ точками  $\{t_i\} , i \in 1:J ,$ на
равные отрезки. Построим плоскости $\pi_i$ в $\mathbb{R}^3,$
проходящие соответственно через точки $(0,0,0),
(r(t_i),f_{\varepsilon}(r(t_i))),
(r(t_{i+1}),f_{\varepsilon}(r(t_{i+1})), i \in 1:J $. Части
плоскостей $\pi_i ,i \in 1:J$ , определенных в секторах,
образуемых векторами $(0,0), r(t_i), r(t_{i+1})$, определяют
график п.о. степени 1 многогранной функцию $(\psi_{\varepsilon})_J
(r(\cdot)).$ Будем понимать под двугранным углом функцию, график
которой состоит из полуплоскостей с общей граничной прямой,
включающих плоскости $\pi_i,$ построенные в соседних секторах.
Покажем, что все двугранные углы функции $(\psi_{\varepsilon})_J
(r(\cdot)),$ образуемые плоскостями $\pi_i, i \in J,$ построенными
по соседним секторам, $-$ выпуклые.

Под градиентом плоскости $\pi_i$ будем понимать градиент линейной
функции, график которой совпадает с плоскостью $\pi_i$. Обозначим
градиенты плоскостей $\pi_i$ и $\pi_{i+1}$ через $\nabla \pi_i$ и
$\nabla \pi_{i+1}$ соответственно.  Воспользуемся теоремой о
средней точке, согласно которой существует такая точка $t_m \in
[t_i, t_{i+1}],$ что
$$
           \partial f_{\varepsilon}(r(t_m))/ \partial e_i =
           (\nabla \pi_i, e_i),
$$
где
$$
e_i=(r(t_{i+1})-r(t_i))/ \mid \mid r(t_{i+1})-r(t_i) \mid \mid .
$$
Аналогично для плоскости $\pi_{i+1}$ и некоторой точки $t_c \in
[t_{i+1},  t_{i+2}]$ имеем
$$
           \partial f_{\varepsilon}(r(t_c))/ \partial e_{i+1} =
           (\nabla \pi_{i+1}, e_{i+1}),
$$
где
$$
e_{i+1}=(r(t_{i+2})-r(t_{i+1}))/ \mid \mid r(t_{i+2})-r(t_{i+1})
\mid \mid .
$$
Функция $f_{\varepsilon}(\cdot)$ сильно выпуклая, так как ее
матрица вторых частных производных положительно определенная.
Любая выпуклая функция имеет неубывающую производную по
направлению вдоль произвольного луча. Но для сильно выпуклой
функции производная по касательному направлению к кривой вида
$r(x_0,\tau, g)= x_0+\tau g +o_{\varepsilon}(\tau)$, $g \in
\mathbb{R}^n$,$ \tau >0$ в малой окрестности точки $x_0$ есть
возрастающая функция вдоль этой кривой. Поэтому для достаточно
большом $J$ и равномерном разбиении кривой $r(\cdot)$ точками
$t_i$ имеем
$$
\partial f_{\varepsilon}(r(t_m))/ \partial e_i < \partial f_{\varepsilon}(r(t_c))/ \partial
e_{i+1},
$$
или
$$
(\nabla \pi_i, e_i) < (\nabla \pi_{i+1}, e_{i+1}).
$$
Учтем также, что разность $\nabla \pi_{i+1} - \nabla \pi_i$
перпендикулярна вектору $r(t_{i+1}).$ Отсюда и из неравенства выше
следует, что двугранный угол $\pi_i, \pi_{i+1}$ - выпуклый. При $J
\rightarrow \infty$
$$ (\psi_{\varepsilon})_J(\cdot) \Rightarrow (\psi_{\varepsilon})(\cdot).$$
Так как точечный предел для выпуклых функций равносилен
равномерному пределу, то $\psi_{\varepsilon}(\cdot)$ - выпуклая
функция. Также $\psi_{\varepsilon}(\cdot) \Rightarrow \psi
(\cdot)$ при $\varepsilon \rightarrow +0,$ т.е. $\psi(\cdot)-$
выпуклая, что и требовалось доказать.

Очевидно, что градиенты линейных функций, графики которых есть
$\pi_i, i \in J,$ ограничены константой, зависящей только от самой
функции $ f(\cdot,\cdot).$ Верно равенство
$$
    \psi (r(t)) = f(r(t)) \;\; \forall t \in [0, 2\pi]].
$$
Ясно, что $\psi (\cdot,\cdot)$ строится однозначно по функции
$f(\cdot,\cdot)$ и выбранной кривой $r(\cdot).$ Из сказанного выше
следует, что функция $\psi (\cdot,\cdot)$ есть липшицевая с
константой $L(f)$.

Пусть
$$
 \Psi(t) = \psi (r(t)) \;\; \forall t \in [0, 2\pi]].
$$
Поскольку
$$
  \vee (\Phi'; 0, 2\pi)  =  \vee( \Psi ' ; 0, 2\pi) ,
$$
то из доказанной Теоремы \ref{difconvgomogMDegreethm1}   следует,
что
$$
\vee (\Phi'; 0,  2 \pi) \leq  c(f).
$$
Если функция $ f(\cdot,\cdot)$ не является  дважды непрерывно
дифференцируемой, то ее можно приблизить выпуклой дважды
непрерывно дифференцируемой функцией $ \tilde{f}(\cdot,\cdot)$ и
построить соответствующую ей функцию $\tilde{\psi} (\cdot,\cdot)$
так, чтобы значения  функций $\psi (\cdot,\cdot)$, $\tilde{\psi}
(\cdot,\cdot)$ и их производных там, где они существуют, как
угодно мало отличались друг от друга. Но тогда аналогичное будет
верно для функций $\Psi(\cdot)$, $\tilde{\Psi}(\cdot)$,
построенных по $\psi (\cdot,\cdot), \tilde{\psi} (\cdot,\cdot)$
соответственно, и их производных. Значит написанное выше
неравенство для вариации производных функции $\Psi(\cdot)$ верно
для общего случая. Лемма \ref{difconvGomogMDegreelem1} доказана.
$\Box$
\begin{lem}
Пусть задана произвольная ПРВ функция от двух переменных $(x,y)
\rar f(x,y): \mathbb{R}^2 \rar  \mathbb{R}. $ Тогда существует
константа $c(f)$, для которой  верно неравенство
\begin{equation}
   \vee (\Phi' ; 0, 2 \pi) \leq  c(f),
\label{difconvGomogMDegree0} \end{equation} где $\Phi(t) =
f(r(t)), t \in [0,2 \pi].$ \label{difconvGomogMDegreelem2}
\end{lem}

{\bf  Доказательство.} По условию
$$
f(x,y)=f_1(x,y) - f_2(x,y)   \,\,\,\, \forall (x,y) \in
\mathbb{R}^2,
$$
где $f_1(\cdot, \cdot), f_2(\cdot, \cdot) - $ выпуклые функции. По
лемме \ref{difconvGomogMDegreelem1}
$$
\vee (\Phi_1' ; 0, 2 \pi) \leq  c_1(f), \,\,\, \vee (\Phi_2' ; 0,
2 \pi) \leq  c_2(f),
$$
где
$$
\Phi_1(t) = f_1(r(t)), \,\,\, \Phi_2(t) = f_2(r(t)),\,\,\, \forall
t \in [0,2 \pi].
$$
На основании неравенства для вариации суммы функций
\cite{kolmogorovfomin} имеем
$$
\vee (\Phi' ; 0, 2 \pi) \leq \vee (\Phi_1' ; 0, 2 \pi)+ \vee
(\Phi_2' ; 0, 2 \pi) \leq c_1 + c_2.
$$
Лемма доказана. $\Box $

Рассмотрим теперь произвольную выпуклую п.о. степени 1 функцию $z
\rar \psi(z): \mathbb{R}^2 \rar \mathbb{R},$ принимающую на
$S_1^1(0)$ положительные значения, и по ней определим п.о. степени
$m$ функцию $z \rar \varphi(z): \mathbb{R}^2 \rar \mathbb{R}$
следующим образом
$$
\varphi(z)=\psi(z) \,\, \forall z \in S^1_1(0), \,\,
\varphi(\lambda z) = \lambda^m \varphi(z) \,\, \forall \lambda>0.
$$
Покажем, что $\varphi(\cdot) - $ выпуклая. Нетрудно показать, что
если $q -$ точка дифференцируемости функции $\varphi(\cdot)$ (или
$\psi(\cdot) $ ), то $\lambda q, \lambda >0, $ также точка
дифференцируемости функции $\varphi(\cdot)$ (или $\psi(\cdot) $ ).
Причем верны равенства
$$
\nabla \varphi ( \lambda q  ) = \varphi'(\lambda q)= \lambda^{m-1}
\nabla \varphi (  q  ) = \lambda^{m-1}\varphi'(q),
$$
$$
\nabla \psi ( \lambda q )= \psi'(\lambda q)= \nabla \psi (  q
)=\psi'(q) \,\,\, \forall \lambda >0.
$$
Без ограничения общности будем считать, что функция $\psi(\cdot)
-$ гладкая функция. В противном случае мы перешли бы к
последовательности выпуклых, гладких п.о. первой степени функций
$\{\psi_n (\cdot)\}$, равномерно сходящихся к $\psi (\cdot)$ на
$B_1^2(0)$, по которым мы построим последовательность $\{\varphi_n
(\cdot)\}$ выпуклых, гладких п.о. функций $m -$ ого порядка, также
равномерно сходящихся к $\varphi (\cdot)$ на $B_1^2(0)$.

Для доказательства выпуклости функции $\varphi(\cdot)$
воспользуемся замечательным свойством выпуклой функции, пользуясь
которым академик А.Д. Александров доказал почти всюду дважды
дифференцируемость выпуклой функции \cite{alexandrovSecondDeriv}.

\begin{thm}\cite{pshenichnyi} Для того чтобы функция $z \rar \theta(z):
\mathbb{R}^n \rar \mathbb{R}$ была выпуклой, необходимо и
достаточно, чтобы она была дифференцируема по направлениям и для
любой точки $z \in \mathbb{R}^n $ и любого направления $p \in
\mathbb{R}^n $ функция $\alpha \rar h( \alpha ): \mathbb{R}^+ \rar
\mathbb{R}:$
$$
h(\alpha) = \frac{\p \theta (z + \alpha p)}{\p p }
$$
была неубывающей по $\alpha >0.$
\end{thm}

Возьмем произвольную точку $z \in S^1_1(0)$ и направление $p \in
\mathbb{R}^2. $ Из п.о. функции $\varphi(\cdot)$ следует, что
достаточно рассмотреть случай, когда вектор ${p}$ перпендикулярен
вектору $z$. Обозначим через $\xi $ точку на луче $z + \alpha p,
\alpha >0,$  а через $\eta_{\xi} - $ точку пересечения прямой,
проходящей через начало координат и точку $\xi $ с окружностью
$S^1_1(0)$. Обозначим также через $\varphi'(\eta_{\xi},
\tau_{\xi}), $ производную функции $\varphi(\cdot) $ в точке $
\eta_{\xi} $ по направлению $ \tau_{\xi}$,  где $ \tau_{\xi} -$
касательная к единичной окружности $S^1_1(0)$ в точке $\eta_{\xi}
$, сонаправленная с вектором $p$ (см. рис. 1). Без ограничения
общности будем считать, что $\varphi'(z , p) \meq 0. $ В противном
случае мы  возьмем направление $-p$. Заметим, что
$$
\varphi'(\xi, \tau_{\xi}) = \| \xi \|^{m-1} \varphi'(\eta_{\xi},
\tau_{\xi}) = \| \xi \|^{m-1} \psi'(\eta_{\xi}, \tau_{\xi}),
$$
$$
\varphi'(\xi, \eta_{\xi})= \| \xi \|^{m-1} \varphi'(\eta_{\xi},
\eta_{\xi}) = m \| \xi \|^{m-1} \varphi(\eta_{\xi})= m \| \xi
\|^{m-1} \psi(\eta_{\xi}).
$$
Тогда
$$
\nabla \varphi (\xi)= \varphi'(\xi)=\varphi'(\xi,
\tau_{\xi})\tau_{\xi}+\varphi'(\xi, \eta_{\xi})\eta_{\xi}=
$$
$$
=\| \xi \|^{m-1} (\varphi'(\eta_{\xi}, \tau_{\xi})\tau_{\xi}+m
\varphi(\eta_{\xi})\eta_{\xi})
$$
и
$$
(\varphi'(\xi), p) = \| \xi \|^{m-1} ((\varphi'(\eta_{\xi},
\tau_{\xi})\tau_{\xi},p)+m (\varphi(\eta_{\xi})\eta_{\xi},p))=
$$
$$
=\| \xi \|^{m-1} ((\psi'(\eta_{\xi}, \tau_{\xi})\tau_{\xi},p)+m
(\psi(\eta_{\xi})\eta_{\xi},p))=
$$
\be =\| \xi \|^{m-1} ((\psi'(\xi),p)+(m-1)(\psi( \eta_{\xi}
)\eta_{\xi},p)) \label{difconvGomogMDegree2} \ee

Возьмем две произвольные точки $\xi_1$ и $\xi_2$ на луче $z+\alpha
p, \al>0:$ $ \xi_1=z+\alpha_1 p , \,\, \xi_2=z+\alpha_2 p,$
$\alpha_2
> \alpha_1.$ Из формулы (\ref{difconvGomogMDegree2}) имеем
$$
(\varphi'(\xi_1), p) = \| \xi_1 \|^{m-1}
((\psi'(\xi_1),p)+(m-1)(\psi( \eta_{\xi_1} )\eta_{\xi_1},p)),
$$
$$
(\varphi'(\xi_2), p) = \| \xi_2 \|^{m-1}
((\psi'(\xi_2),p)+(m-1)(\psi( \eta_{\xi_2} )\eta_{\xi_2},p)).
$$
Пусть $\| \xi_2 \| > \| \xi_1 \| $. Сравним два числа $
(\varphi'(\xi_1), p) $ и $ (\varphi'(\xi_2), p) $. Поскольку
$\psi(\cdot)- $ выпуклая, то \be (\psi'(\xi_2), p) \meq
\psi'(\xi_1), p) \label{difconvGomogMDegree3} \ee Кроме того, так
как функции $ \varphi(\cdot), \psi(\cdot) - $ липшицевые и их
значения равны на единичной окружности, то $ \varphi'(z, p)=
\psi'(z, p) \meq 0$, откуда следует, что $\psi(\xi_2) \meq
\psi(\xi_1). $ Также для $ \eta_{\xi_1}, \eta_{\xi_2},  $ близких
к $z$, $\psi(\eta_{\xi_2}) \meq \psi(\eta_{\xi_1})>0, $ а поэтому
\be (\psi(\eta_{\xi_2}) \eta_{\xi_2},p)  \meq (
\psi(\eta_{\xi_1})\eta_{\xi_1},p). \label{difconvGomogMDegree4}\ee
Из (\ref{difconvGomogMDegree3}) и (\ref{difconvGomogMDegree4})
получим
$$
(\varphi'(\xi_1), p) \leq (\varphi'(\xi_2), p).
$$
Следовательно, функция $ \varphi(\cdot) -$  выпуклая в окрестности
точки $z$. Из локальной выпуклости функции $\varphi(\cdot) $
следует ее глобальная выпуклость.

Итак, пункты 1) и 2) выполнимы. Выясним теперь условия, при
которых п.о. степени $m$ функция $\varphi(\cdot)  $ представима в
виде разности выпуклых функций $\varphi_1(\cdot),
\varphi_2(\cdot). $

Пусть
$$
\varphi(q) = \varphi_1(q) - \varphi_2(q) \,\,\,\forall q \in
\mathbb{R}^2,
$$
где $\varphi_1(\cdot) , \varphi_2(\cdot) -  $ выпуклые функции.

Построим, как это мы делали ранее, п.о. степени 1 выпуклые функции
$\psi_1(\cdot), \psi_2(\cdot), $ соответствующие выпуклым функциям
$\varphi_1(\cdot) , \varphi_2(\cdot). $ Тогда функция $
\psi(\cdot), $ определяемая равенством
$$
\psi(q) = \psi_1(q) - \psi_2(q) \,\,\,\forall q \in \mathbb{R}^2,
$$
есть ПРВ функция. Поскольку
$$
\vee(\Psi';0, 2\pi) = \vee(\Phi';0, 2\pi),
$$
где
$$
\Phi(t)=\varphi(r(t)), \,\,\, \Psi(t)=\psi(r(t)) \,\,\, \forall t
\in [0,2\pi],
$$
то выполняется неравенство
$$
\vee(\Phi';0, 2\pi) < c(\varphi).
$$
Без ограничения общности считаем, что $\psi_1(\cdot),
\psi_2(\cdot) $ принимают на $S_1^1(0) $ положительные значения.
Построим теперь по выпуклым функциям $\psi_1(\cdot),
\psi_2(\cdot)$ п.о. степени $m$ функции $\tilde{\varphi}_1(\cdot)
, \tilde{\varphi}_2(\cdot). $ По доказанному ранее
$\tilde{\varphi}_1(\cdot) , \tilde{\varphi}_2(\cdot)- $ выпуклые
функции. Очевидно, что
$$
 \varphi(q) =  \tilde{\varphi}_2(q)- \tilde{\varphi}_1(q)
 \,\,\,\forall q \in \mathbb{R}^2.
$$
Пусть теперь выполняется неравенство
$$
\vee(\Phi';0, 2\pi) < c(\varphi).
$$
Покажем, что $\varphi(\cdot)- $ ПРВ функция. Определим п.о.
степени 1 функцию  $q \rar \psi(q): \mathbb{R}^2 \rar \mathbb{R},
$ принимающую на $S_1^1(0) $ те же значения, что и функция
$\varphi(\cdot) $ на $S_1^1(0) $. Так как
$$
\vee(\Phi';0, 2\pi) = \vee(\Psi';0, 2\pi)< c(\psi),
$$
то $ \psi(\cdot) - $ ПРВ функция \cite{proudconvex1},
\cite{lupikovdissertation}, т.е.
$$
\psi(q) =  {\psi}_2(q)- {\psi}_1(q)
 \,\,\,\forall q \in \mathbb{R}^2,
$$
где $ \psi_i(\cdot), i=1,2, -  $ выпуклые п..о. степени 1 функции.
Очевидно, что функции $ \psi_i(\cdot), i=1,2, $ можно выбрать
такими, чтобы они принимали положительные значения на $S_1^1(0) $.
По функциям $ \psi_i(\cdot), i=1,2, $ построим п.о. степени $m$
функции $ \varphi_i(\cdot), i=1,2. $ По доказанному ранее $
\varphi_i(\cdot), i=1,2, - $ выпуклые. Кроме того, очевидно,
выполняется равенство
$$
\varphi(q) =  {\varphi}_2(q)- {\varphi}_1(q)
 \,\,\,\forall q \in \mathbb{R}^2,
$$
т.е. $\varphi(\cdot)-  $ ПРВ функция.

Таким образом доказана следующая теорема.
\begin{thm}
П.о. степени $m$ липшицевая функция $q \rar \varphi(q):
\mathbb{R}^2 \rar \mathbb{R} $ является ПРВ функцией тогда и
только тогда, когда
$$
\vee(\Phi';0, 2\pi) < c(\varphi),
$$
где $\Phi(t)=\varphi(r(t)), \, r(t)=(\cos t, \sin t), \,\, \forall
t \in [0,2\pi].$ \label{difconvgomogMDegreethm2}
\end{thm}
\begin{cor} Если п.о. степени $m$ функция $\varphi(\cdot) $ есть
ПРВ функция, то она представима в виде разности выпуклых п.о.
степени $m$ функций.
\end{cor}

\vspace{0.5cm}

\subsection{Класс кривых, ограничивающих выпуклые компактные
множества на плоскости}

\vspace{0.5cm}

Обозначим через $D-$ произвольное выпуклое открытое ограниченное
множество на плоскости $ \mathbb{R}^2 $, так что замыкание его $-$
компакт с $\mbox{int} D \neq \emptyset $ и $0 \in \mbox{int} D.$
Пусть $\Re(D)- $ множество кривых $r(\cdot)$, ограничивающих в $D$
выпуклые компактные множества. Параметризуем $r(\cdot)$
естественным образом, т.е. $t-$ расстояние вдоль кривой $r(\cdot)$
от фиксированной точки на кривой до точки $r(t)$.
Параметризованную кривую обозначим через $r(t), t \in [0,T_r].$
Здесь $T=T(r)-$ длина кривой $r(\cdot).$

Пусть на $D$ задана произвольная выпуклая функция $(x,y) \rar
f(x,y): \mathbb{R}^2 \rar \mathbb{R} $ с константой Липшица $
L(D). $ Обозначим через
$$
F(t)=f(r(t)) \,\,\, \forall t \in [0,T].
$$
Функция $F(\cdot) $ липшицевая с константой Липшица $L(D).$
Действительно, для любых $t_1, t_2 \in [0, T] $
$$
| F(t_1) - F(t_2) | = | f(r(t_2)) - f(r(t_1)) | \leq L(D) \|
r(t_2) - r(t_1) \| \leq L(D) | t_2 - t_1 |.
$$
Следовательно \cite{kolmogorovfomin}, $F(\cdot) - $ почти всюду
(п.в.) дифференцируемая на $ [0,T]. $

Справедлива следующая лемма.

\begin{lem}
Для любой выпуклой функции $f(\cdot): \mathbb{R}^2 \rar
\mathbb{R}$ и любой кривой $r(\cdot) \in \Re(D)$ существует
константа $c_1(f,D) >0 $ такая, что верно неравенство
$$
\vee(F';0,T) < c_1(f,D) \,\,\,\,\, \forall r(\cdot) \in \Re(D),
$$
где $F(t)=f(r(t)) \,\, \forall t \in [0, 2\pi].$
\label{difconvGomogMDegreelem3}
\end{lem}
Здесь, как и ранее, производные берутся там, где они существуют.
Предварительно докажем такую лемму.
\begin{lem}
Для любой кривой $r(\cdot) \in \Re(D)$ и любой выпуклой п.о.
степени 1 функции $(x,y) \rar \psi(x,y): \mathbb{R}^2 \rar
\mathbb{R} $ существует константа $c_2(\psi,D) >0 $ такая, что
верно неравенство
$$
\vee(\Psi';0,T) < c_2(\psi,D) \,\,\,\,\,\,\, \forall r(\cdot) \in
\Re(D),
$$
где $\Psi(t)=\psi(r(t)) \,\,\, \forall t \in [0,T].$
\label{difconvGomogMDegreelem4}
\end{lem}
{\bf Доказательство леммы \ref{difconvGomogMDegreelem4}} Без
ограничения общности будем считать, что $\psi(\cdot)-$ гладкая на
$\mathbb{R}^2 \backslash 0.$ Пусть
$$
\psi(r(t))=\max_{v \in \p \psi(0) } (v, r(t))=(v(t), r(t)),
\,\,\,\, v(t) \in \p \psi(0),
$$
где $  \p \psi(0) -  $ субдифференциал функции $\psi(\cdot)  $ в
нуле \cite{demvas}. Обозначим через $ d(D) - $ диаметр множества
$D$. Верны соотношения
$$
| \psi(r(t_1))-\psi(r(t_2)) | = |(v(t_1),r(t_1))-
(v(t_2),r(t_2))|=
$$
$$
=|(v(t_1)-v(t_2),r(t_1))+(v(t_2),r(t_1))-(v(t_2),r(t_2)) | \leq
$$
$$
\leq \|v(t_1)-v(t_2)\| \|r(t_1)\|+\|r(t_1))-r(t_2)\| \, \|v(t_2)\|
\leq
$$
$$
\leq \|v(t_1)-v(t_2)\| \,\, d(D) + L(D) \, |t_1 - t_2|.
$$
Отсюда следует, что \be \vee(\Psi';0,T) \leq P(\p \psi(0) ) \,
d(D) + L(D) \, T, \label{difconvGomogMDegree5}\ee где $ P(\p
\psi(0) )- $ длина кривой, ограничивающей выпуклое компактное
множество $\p \psi(0)$. Поскольку справа от знака неравенства
(\ref{difconvGomogMDegree5}) стоит конечная величина, зависящая
только от множества $D$, то лемма \ref{difconvGomogMDegreelem4}
доказана. $\Box$.

{\bf Доказательство леммы \ref{difconvGomogMDegreelem3}.} Строим
по функции $f(\cdot) $ также, как это делалось ранее, п.о. степени
1 функцию $(x,y) \rar \psi(x,y): \mathbb{R}^2 \rar \mathbb{R} $,
принимающую на $r(\cdot)$ те же значения, что и функция $ f(\cdot)
$. Повторяя рассуждения, проведенные ранее, показываем, что
$\psi(\cdot)- $ выпуклая функция. Пусть
$$
\Psi(t)=\psi(r(t)) \,\,\,\, \forall t \in [0,T].
$$
Поскольку
$$
\vee(\Psi';0,T)= \vee(F';0,T),
$$
то из леммы \ref{difconvGomogMDegreelem4} следует, что
$$
\vee(F';0,T) < c_2(f,D).
$$
Лемма \ref{difconvGomogMDegreelem3} доказана. $\Box  $

\begin{thm}
Для того чтобы п.о. степени $m$ липшицевая функция $(x,y) \rar
\varphi(x,y): \mathbb{R}^2 \rar \mathbb{R} $ была ПРВ функцией,
необходимо и достаточно, чтобы для любой кривой $r(\cdot) \in
\Re(D) $ нашлась константа $c_3(\varphi,D)>0$, для которой \be
\vee(\Phi';0,T) < c_3(\varphi, D) \,\,\,\,\,\, \forall r(\cdot)
\in \Re(D), \label{difconvGomogMDegree6} \ee где
$\Phi(t)=\varphi(r(t)) \,\,\, \forall t \in [0,T].$
\label{difconvgomogMDegreethm3}
\end{thm}
{\bf Доказательство.} Поскольку в класс $\Re(D) $ входит
окружность $S_\rho^1(0)$ радиуса $\rho>0  $  с центром в начале
координат и все производные на этой окружности там, где они
существуют, связаны с производными на единичной окружности
$S^1_1(0)$ одним и тем же коэффициентом пропорциональности $
\rho^{m-1} $, то по теореме \ref{difconvgomogMDegreethm2}
выполнение неравенства (\ref{difconvGomogMDegree6}) достаточно для
представимости функции $ \varphi(\cdot) $ в виде разности выпуклых
функций.

Докажем необходимость. Пусть
$$
\varphi(q)=\varphi_1(q)-\varphi_2(q) \,\,\,\,  \forall q \in
\mathbb{R}^2,
$$
где $ \varphi_i(\cdot), i=1,2,- $ выпуклые. По лемме
\ref{difconvGomogMDegreelem3} для любой кривой $r(\cdot) \in
\Re(D) $
$$
\vee(\Phi_i';0,T) < c_i(\varphi_i, D), \,\,\,\, i=1,2,
$$
где
$$
\Phi_i(t)=\varphi_i(r(t)) \,\,\,\,\, \forall t \in [0,T], i=1,2.
$$
Тогда
$$
\vee(\Phi';0,T)\leq
\vee(\Phi_1';0,T)+\vee(\Phi_2';0,T)<c_1(\varphi_1,
D)+c_2(\varphi_2, D)=c_3(\varphi, D).
$$
Теорема доказана. $\Box$

\vspace{0.5cm}

\subsection{Геометрическая интерпретация теоремы для положительно однородной функции
$m-$ ой степени двух переменных}

\vspace{0.5cm}

Перефразируем теорему \ref{difconvgomogMDegreethm3}, придав ей
более геометрический характер. Введем понятие поворота кривой
$r(\cdot)$ на графике $\Gamma_{\varphi } = \{(x,y,z) \in
\mathbb{R}^3 \mid z = \varphi(x , y)\}.$

Рассмотрим на $\Gamma_{\varphi }$ кривую $ R(t)=(r(t),{\varphi
}(r(t))), t \in [0, T(r)],$ где $r(\cdot) \in \Re(D).$   Так как
функция ${\varphi }(\cdot,\cdot)$ есть липшицевая, то п.в. на
$[0,T(r)]$ существует производная $R'(\cdot),$ которую обозначим
через $\tau(\cdot)=R'(\cdot),$  а множество точек, где она
существует, $-$ через $N_{\varphi }$.

\begin{defi}  Поворотом кривой $R(\cdot)$ на многообразии $\Gamma_f$
назовем величину
$$
sup_{ \{t_i\} \subset N_{\varphi }} \,\, \sum_i \Vert \tau(t_i)/
\Vert \tau (t_i) \Vert -  \tau(t_{i-1})/ \Vert \tau (t_{i-1})
\Vert \Vert = O_{\varphi }.
$$
\end{defi}

Таким образом, поворот $O_{\varphi }$ кривой $R(\cdot)$ есть
верхняя грань суммы углов между касательными $\tau(t)$  для $t \in
[0,T(r)].$ Нетрудно видеть, что для плоской гладкой кривой,
параметризованной естественным образом, величина $ O_\varphi$
равна интегралу
$$
\int^{T(r)}_0 \mid k(s) \mid ds,
$$
где $k(s)$ - кривизна рассматриваемой кривой $ r(\cdot)$ в точке
$s \in [0,T(r)],$ т.е. совпадает с обычным определением поворота
кривой в точке \cite{pogorelov1} .

\begin{thm} Для того, чтобы произвольная липшицевая п.о. степени
$m$ функция $z \rightarrow {\varphi }(z) :\mathbb{R}^2 \rightarrow
\mathbb{R}$ была ПРВ функцией, необходимо и достаточно, чтобы для
всех $r(\cdot) \in \Re(D)$ существовала константа $c_4({\varphi
})>0$ такая, что поворот кривой $R(\cdot)$ на $\Gamma_{\varphi }$
ограничен сверху константой $c_4({\varphi })>0,$ т.е.
\begin{equation} O_{\varphi } \leq c_4({\varphi}) \,\,\,\,\,
\forall r(\cdot) \in \Re(D). \label{difconvgomog4}
\end{equation}
\label{difconvgomogMDegreethm4}
\end{thm}

{\bf Доказательство. }{\bf Необходимость}. Пусть $ {\varphi
}(\cdot,\cdot)$ есть ПРВ функция. Покажем, что тогда справедливо
неравенство (\ref{difconvgomog4}). Воспользуемся неравенством,
вытекающим из неравенства треугольника,
$$
\Vert \tau(t_i) / \Vert \tau(t_i) \Vert  - \tau (t_{i-1}) / \Vert
\tau(t_{i-1} \Vert \Vert \leq
$$
$$
\leq \Vert r'(t_i) / \sqrt{ 1+\varphi'^2_t (r(t_i))}  -
r'(t_{i-1}) / \sqrt{ 1+\varphi'^2_t (r(t_{i-1}))} \Vert +
$$
$$
\mid {\varphi }'_t(r(t_i)) /  \sqrt{ 1+{\varphi }'^2_t (r(t_i))}-
{\varphi }'_t(r(t_{i-1})) / \sqrt{ 1+{\varphi }'^2_t (r(t_{i-1}))}
\mid .
$$
Так как $1 \leq  \sqrt{1+{\varphi }'^2_t (r(t_i))} \leq
\sqrt{1+L^2}$ для всех $ t_i \in [0, T(r)]$ , то очевидно,
существует такое $c_5 >1,$ для которого верно неравенство

\begin{equation} \Vert r'(t_i) / \sqrt{1+{\varphi }'^2_t (r(t_i))}-
r'(t_{i-1}) / \sqrt{1+{\varphi }'^2_t (r(t_{ i-1}))} \Vert \leq
c_5 \Vert r'(t_i) - r'(t_{i-1}) \Vert. \label{difconvgomog5}
\end{equation}

Из свойств функции $\theta(x)= x / \sqrt{ 1+x^2}$ следует
неравенство
$$
\mid {\varphi }'_t (r(t_i)) / \sqrt{ 1+{\varphi
}'^2_t (r(t_i))} - {\varphi }'_t (r(t_{i-1})) / \sqrt{1+{\varphi
}'^2_t (r(t_{i-1}))} \mid \leq
$$
\begin{equation}
\leq \mid {\varphi }'_t (r(t_i)) - {\varphi }'_t (r(t_{i-1}))
\mid. \label{difconvgomog6}
\end{equation}

Из (\ref{difconvgomog5}) и (\ref{difconvgomog6}) имеем
$$
 \sup_{\{t_i \} \in N_{\varphi } } \,\, \sum_i \,
\Vert \tau(t_i) / \Vert \tau(t_i) \Vert - \tau (t_{i-1}) / \Vert
\tau(t_{i-1}) \Vert \Vert \leq
$$
\begin{equation}
\leq c_5 (\vee (  r' ; 0 ,T(r)) + \vee (\Phi' ; 0,T(r)) ).
\label{difconvgomog7}
\end{equation}

Так как по условию $ \varphi(\cdot,\cdot)-$  ПРВ функция, то
согласно теореме \ref{difconvgomogMDegreethm3}
$$
\vee (\Phi'; 0,T(r)) \leq  c_3(\varphi,D),
$$
откуда с учетом (\ref{difconvgomog7}) и ограниченности вариации
$\vee ( r' ; 0 ,T(r)) $ следует неравенство (\ref{difconvgomog4}).
Необходимость доказана.

{\bf Достаточность}.  Пусть справедливо неравенство
(\ref{difconvgomog4}). Покажем, что ${\varphi }(\cdot,\cdot)$ -
ПРВ функция. Воспользуемся неравенством
$$
\Vert \tau(t_i) / \Vert \tau( t_i) \Vert  - \tau(t_{i-1}) /  \Vert
\tau (t_{i-1}) \Vert \Vert \geq \mid  {\varphi }'_t (r(t_i)) /
\sqrt{1+{\varphi }'_t (r(t_i))} -
$$
\begin{equation} - {\varphi }'_t (r(t_{i-1})) / \sqrt {1+{\varphi }'^2_t
(r(t_{i-1}))} \label{difconvgomog8} \end{equation} Из свойств
функции $\theta(x) = x / \sqrt{1+x^2}$ и из $\Vert {\varphi }'(z)
\Vert \leq L$ для всех $z \in D,$ где производная существует,
следует существование константы $ c_6 > 0,$ для которой
$$
\mid {\varphi }'_t (r(t_i)) / \sqrt{1+{\varphi }'^2_t (r(t_i))} -
{\varphi }'_t (r(t_{i-1})) / \sqrt{1+{\varphi }'^2_t(r(t_{i-1}))}
\geq
$$
$$
\geq c_6 \mid {\varphi }'_t (r(t_i)) - {\varphi }'_t (r(t_{i-1}))
\mid,
$$
откуда с учетом (\ref{difconvgomog8}) имеем
$$
c_4(\varphi,D) \geq \sup_{ \{ t_i \} \subset N_{\varphi}} \sum_i
\Vert \tau(t_i) / \Vert \tau( t_i) \Vert  - \tau(t_{i-1}) / \Vert
\tau (t_{i-1}) \Vert \Vert \geq
$$
$$
\geq c_6 \vee (\Phi';  0,T(r)) .
$$
Из теоремы \ref{difconvgomogMDegreethm3} следует, что ${\varphi
}(\cdot)$ - ПРВ функция. Достаточность доказана.  $\Box$

\newpage
\vspace{1cm}

\section{УСЛОВИЯ ПРЕДСТАВИМОСТИ
ПОЛОЖИТЕЛЬНО ОДНОРОДНОЙ ФУНКЦИИ ТРЕХ ПЕРЕМЕННЫХ В ВИДЕ РАЗНОСТИ
ВЫПУКЛЫХ ФУНКЦИЙ}

\vspace{1cm}

В данном параграфе приведены необходимые и достаточные условия
представимости произвольной положительно однородной функции
первого порядка трех переменных в виде разности выпуклых функций.
Дана также геометрическая интерпретация этих условий. Приведен
алгоритм такого представления, результатом которого есть
последовательность равномерно сходящихся на произвольном компакте,
внутренности которого принадлежит начало координат,  выпуклых
положительно однородных первого порядка функций.

\noindent

\vspace{0.5cm}

\subsection{Введение}

\vspace{0.5cm}

Обозначим через $\mathbb{R}^n - n$-мерное евклидово пространство
со скалярным произведением $(a,b)$ векторов $a$ и $b$. Пусть
задана липшицевая, дифференцируемая по направлениям функция
$f:\mathbb{R}^n \rar \mathbb{R}$. Обозначим через $ f'(x,q) $
производную по направлению $q \in \mathbb{R}^n$ функции $f(\cdot)$
в точке $x \in \mathbb{R}^n$.

В теории оптимизации важную роль играют квазидифференцируемые
(КВД) функции \cite{demvas}.
\begin{defi}
Функция $x \rar f(x)$ называется КВД в точке  $x \in
\mathbb{R}^n$, если существует пара выпуклых компактных множеств
$\underline{\p} f(x) $ и $\overline{\p} f(x)$, называемых
соответственно  субдифференциалом и  суппердифференциалом, таких,
что верно равенство \be \varphi(q)=f'(x,q)=\max_{v \in
\underline{\p} f(x)} (v,q) + \min_{w \in \overline{\p} f(x)} (w,q)
\,\,\, \forall q \in \mathbb{R}^n. \label{DiffconvthreevarGomog1}
\ee
\end{defi}
Легко видно, что (\ref{DiffconvthreevarGomog1}) есть разложение
функции $\varphi(\cdot)$ в виде разности выпуклых п.о. функций:
$$
\max_{v \in \underline{\p} f(x)} (v,q) \,\,\,\,\,\, \mbox{ и}
\,\,\,\, \max_{w \in -\overline{\p} f(x)} (w,q).
$$

Если для дифференцируемых функций $ f(\cdot)$ необходимое условие
экстремума в точке $x^*$ записывается в виде $ \nabla f(x^*) =
f'(x^*)=0 $, то для КВД функции условие минимума есть $
\underline{\p} f(x^*) \supset - \overline{\p} f(x^*) $, а условие
максимума:  $ \underline{\p} f(x^*) \subset - \overline{\p}
f(x^*). $

Для двумерного случая $q \in \mathbb{R}^2 $ необходимое и
достаточное условие представимости п.о. липшицевой функции $q \rar
\varphi(q) $ в виде разности выпуклых п.о. функций получены в
\cite{lupikovdissertation}. Сформулируем эти условия.

Обозначим через $r(t)=(\cos t, \sin t), t \in [0,2\pi],$
окружность $S_1^1(0)=\{q \in \mathbb{R}^2 |\,\,\, \| q \| =1 \} $
единичного радиуса  с центром в начале координат,
параметризованную естественным образом. Пусть $
\Phi(t)=\varphi(r(t)), t \in [0, 2\pi]. $ Нетрудно показать, что
функция $ \Phi(\cdot) -$ липшицевая, а поэтому почти всюду (п.в.)
дифференцируемая на $ [0, 2\pi]. $ По определению положим
$$
\vee(\Phi';0, 2\pi) = \sup \lim_{\begin{array}{c}
                                 n \rar \infty,  \\
                                 \{ t_i \} \in [0, 2\pi]
                                 \end{array}}
\sum_{i=0}^{n} | \Phi'(t_i) - \Phi'(t_{i+1} |,
$$
где берутся такие точки $\{ t_i \}$, где производные $\Phi'(t_i)$
существуют.

\begin{thm} {\em \cite{lupikovdissertation}} Для того, чтобы липшицевая
п.о. первой степени функция $\varphi(\cdot): \mathbb{R}^2 \rar
\mathbb{R}$ была представима в виде разности выпуклых функций
(была ПРВ функцией), необходимо и достаточно, чтобы
$$
 (\exists c(\varphi)>0):  \;\;\; \vee (\Phi'; 0,2\pi)<  c(\varphi),
$$
где $\Phi(t)=\varphi(r(t)), \,\, r(t)=(\cos t, \sin t), \,\, t \in
[0, 2\pi] $. \label{difconvgomogthm1}
\end{thm}

\vspace{0.5cm}

\subsection{Положительно однородные функции первой степени от трех
переменных}

\vspace{0.5cm}

Далее будем рассматривать случай $n=3$. Вначале дадим необходимое
и достаточное условие представимости функции $ q \rar \varphi(q)$
в виде разности выпуклых функций.

Обозначим через $\hat{\Re} $ класс кривых на поверхности
единичного шара с центром в нуле $B_1^3(0)=\{q \in \mathbb{R}^3
|\,\,\, \| q \| \leq 1 \} $, получающихся в результате сечения
единичной сферы $S_1^2(0)=\{q \in \mathbb{R}^3 |\,\,\, \| q \| =1
\} $ произвольными плоскостями $\hat \Pi  $. Очевидно, что класс
кривых $\hat{\Re} $ состоит из окружностей на поверхности шара
$B_1^3(0)$.

Возьмем любую кривую $\hat r \in \hat {\Re} $ и параметризуем ее
естественным образом. Тогда
$$
  \| r(t_2) - r(t_1) \| \leq \mid t_2 - t_1 \mid.
$$
Отрезок значений параметра $t$ обозначим через $ [0, T(\hat r)]$ и
положим
$$
\hat{\Phi}(t)=\varphi(\hat r(t)) \,\,\, \forall t \in [0, T(\hat
r)]
$$
Заметим, что при условии липшицевости  функции $\varphi(\cdot)$ с
константой Липшица $L$ функция $ \hat{\Phi}(\cdot) $ также
липшицевая с константой $L$:
$$
  \mid  \hat{\Phi}(t_2)-\hat{\Phi}(t_1) \mid = \mid \varphi(\hat r(t_2)) -
  \varphi(\hat r(t_1)) \mid
$$
$$
  \leq  L \| \hat r(t_2) - \hat r(t_1) \| \leq
  \mid t_2 - t_1 \mid.
$$
Откуда следует, что вектор-функция $\hat r(\cdot) $ почти всюду
(п.в.) на отрезке $[0, T]$ имеет касательную $\hat r'(\cdot)$, а
функция $ \hat{\Phi}(\cdot) $  п.в. имеет производную $
\hat{\Phi}'(\cdot)  $.

\begin{lem}
Если п.о. первой степени липшицевая функция $q \rar
\varphi(q):\mathbb{R}^3 \rar \mathbb{R} $ представима в виде
разности выпуклых функций, то для всех $\hat r \in \hat {\Re} $
существует такая константа $ \hat c$, что
$$
\vee(\hat{\Phi}'; 0,T(\hat r)) \leq \hat c \,\,\,\,\,\, \forall
\hat r \in \hat {\Re}.
$$
\end{lem}
{\bf Доказательство.} По условию леммы верно равенство
$$
\varphi(q)=\varphi_1(q)-\varphi_2(q) \,\,\,\, \forall q \in
\mathbb{R}^3,
$$
где $ \varphi_i(\cdot), i=1,2, - $ выпуклые п.о. функции с
константой Липшица $L_i, i=1,2,$ соответственно. Возьмем
произвольную кривую $ \hat r(\cdot) \in \hat \Re. $

Пусть
$$
\hat \Phi_i(t)=\varphi_i (\hat r(t)), \,\,\, i=1,2, \,\,\,\forall
t \in [0,T(\hat r)].
$$
Функции $\hat \Phi_i(\cdot)-$ липшицевые, поэтому они п.в.
дифференцируемы на отрезке $[0,T(\hat r)].  $ Так как
\cite{kolmogorovfomin}
$$
\vee(\hat \Phi';0,T(\hat r)) \leq \vee(\hat \Phi_1';0,T(\hat r))+
\vee(\hat \Phi_2';0,T(\hat r)),
$$
то достаточно доказать, что для некоторой константы $\hat c_1$ \be
\vee(\hat \Phi_1';0,T(\hat r)) \leq \hat c_1
\label{DiffconvthreevarGomog2} \ee для всех $\hat r(\cdot) \in
\hat \Re.$

Допустим, что кривая $ \hat r(\cdot)$ принадлежит плоскости $ \hat
\Pi. $ Проведем перпендикулярно $\hat \Pi$ вектор $e$, начальная
точка которого есть нуль, а конечная точка принадлежит $\hat \Pi$.
Обозначим через $\Pi$ плоскость, параллельную плоскости $\hat \Pi$
и проходящую через начало координат. Введем в $\mathbb{R}^3$
декартову систему координат, две оси которой принадлежат плоскости
$\Pi$ , а другая ось параллельна и сонаправлена с вектором $e$.
Введем также функцию $\psi_1(\cdot):\mathbb{R}^2 \rar \mathbb{R}$,
определенную на плоскости $\Pi$. Далее, трехмерный вектор $\tilde
r \in \Pi,$  записанный в системе координат плоскости $\Pi$,
обозначим через $r$, т.е. $r \in \mathbb{R}^2$.

По определению положим для $r \in \Pi  $
$$
\psi_1(r)=\varphi_1(e+\tilde r) - \varphi_1(e)=\varphi_1(\hat
r)-\varphi_1(e),
$$
где $\hat r=\tilde r+e.$ Покажем, что $\psi_1(\cdot) - $ выпуклая
липшицевая функция с константой Липшица $L_1$.

Действительно, для любых $\al_1, \al_2 \meq 0, \al_1+\al_2=1,$ и
любых $r_1, r_2 \in \Pi$ имеем
$$
\psi_1(\al_1 r_1 + \al_2 r_2) = \varphi_1(\al_1 e + \al_1 \tilde
r_1+ \al_2 e+ \al_2 \tilde r_2) -\varphi_1(e) \leq
$$
$$
\leq \al_1 \varphi_1(e + \tilde r_1) + \al_2 \varphi_1(e + \tilde
r_2)- \varphi_1(e)=
$$
$$
=\al_1 \varphi_1(\hat r_1) + \al_2 \varphi_1(\hat r_2)-
\varphi_1(e) =\al_1 \psi_1( r_1) +\al_2 \psi_1(r_2),
$$
где $ \hat r_i=e+\tilde r_i, i=1,2. $ Кроме того,
$$
|\psi_1(r_1)-\psi_1(r_2)| = |\varphi_1(e+\tilde r_1) -
\varphi_1(e+\tilde r_2)|\leq L_1 \| \tilde r_1 -\tilde r_2 \| =
L_1 \| r_1 - r_2 \|.
$$
Определим теперь п.о. первой степени функцию двух переменных
$\tilde \psi_1(\cdot): \mathbb{R}^2 \rar \mathbb{R}$. Обозначим
через $r(t), t \in [0,T(\hat r),]$  окружность в плоскости $\Pi,
r(t) \in \mathbb{R}^2,$ являющуюся проекцией окружности $\hat
r(t),t \in [0,T(\hat r)],$ на плоскость $\Pi$. По определению
значения функции $\tilde \psi_1(\cdot)$ на окружности $ r(\cdot) $
равны значениям функции $\psi_1(\cdot)$ на той же окружности. Вне
окружности $r(\cdot)$ функция $\tilde \psi_1(\cdot)$
распространяется по свойству положительной однородности. Покажем,
что $\tilde \psi_1(\cdot)- $ выпуклая функция.

Разобьем круг с окружностью $r(\cdot)$ на секторы. Построим по
данному разбиению п.о. первой степени многогранную функцию $\tilde
\psi_{1m}(\cdot)$ следующим образом. В каждом секторе $\tilde
\psi_{1m}(\cdot)$ линейна, а ее значения на граничных векторах
секторов, принадлежащих окружности $r(\cdot)$, равны значениям
функции $\tilde \psi_{1}(\cdot)$. Докажем, что $\tilde
\psi_{1m}(\cdot) -$ выпуклая.

Пусть $r(t_i), t_i \in [0,T(\hat r)], i \in 1:m, -$  точки на
окружности $r(\cdot)$, получающиеся в результате деления круга на
секторы. В дальнейшем мы будем точки отождествляем с векторами, у
которых начальная точка $-$ начало координат, а конечная $-$ сама
точка. Рассмотрим конусы $K_i = \mbox{con} \{ e, \hat r(t_i), \hat
r(t_{i+1}) \}, i \in 1:m $, в каждом из которых определим линейную
функцию, значения которой равны значениям функции
$\varphi_1(\cdot) $ на векторах $e, \hat r(t_i), \hat r(t_{i+1}
).$ Функцию, равную максимуму (минимуму) линейных функций,
построенных по векторам из соседних конусов $K_i$ и $K_{i+1}$,
имеющих общую часть плоскости $\sigma_i$, и совпадающую в  $K_i$ и
$K_{i+1}$  с этими линейными функциями, в назовем {\em двугранным
углом}. Градиенты этих линейных функций обозначим через $a_i$ и
$a_{i+1}$. Поскольку $\sigma_i=K_i \cap K_{i+1}$, то проекции
векторов $a_i$ и $a_{i+1}$ на $\sigma_i$ равны. Следовательно,
вектор $a_i - a_{i+1} $ перпендикулярен $\sigma_i$. Все
рассмотренные двугранные углы выпуклые, так как они построены по
выпуклой п.о. первой степени функции $\varphi_1(\cdot)$.

Заметим также, что двумерный вектор $ a_{i\Pi},$  равный
составляющей вектора $a_i$ на плоскости $\Pi$, есть градиент
линейной функции, построенной по сектору, образованному векторами
$ r(t_i), r(t_{i+1}).$ так как плоскость $\Pi$ перпендикулярна
вектору $e$, то $a_{i\Pi} - a_{(i+1)\Pi} = a_i - a_{i+1}.$  По
свойству выпуклых функций \cite{pshenichnyi} все двугранные углы
функции $\tilde \psi_{1m}(\cdot)-$  выпуклые. При $m \rar \infty$
функции $\tilde \psi_{1m}(\cdot)$ равномерно на $B_1^2(0)$
стремятся к $\tilde \psi_{1}(\cdot)$. Поэтому $\tilde
\psi_{1}(\cdot) -$ выпуклая. С другой стороны,
$$
| \tilde \psi_{1}(r(t)) | = | \psi_1(r(t)) | =| \varphi_1(e+\tilde
r(t)) - \varphi_1(e) | \leq L_1 \| \tilde r(t) \| = L_1 \| r(t)
\|.
$$
Поэтому функция $ \tilde \psi(\cdot)-$ липшицевая с константой
Липшица $L_1$. Так же, как в \cite{lupikovdissertation},
показываем, что для функции $\tilde \Phi_1(t)=\tilde \psi_1(r(t))
$ верно неравенство \be \vee(\tilde \Phi'_1;0,T(\hat r)) \leq c.
\label{DiffconvthreevarGomog3} \ee Но \be \vee(\hat
\Phi'_1;0,T(\hat r)) = \vee(\tilde \Phi'_1;0,T(\hat r)).
\label{DiffconvthreevarGomog4} \ee Из
(\ref{DiffconvthreevarGomog3}) и (\ref{DiffconvthreevarGomog4})
следует (\ref{DiffconvthreevarGomog2}). Лемма доказана. $\Box$

Определим на $ S_1^2(0)$ класс $\Re$ непрерывных кривых
$r(\cdot)$. Кривые параметризуем естественным образом. Считаем,
что $t \in [0,T(r)]$, где $T(r) - $  длина кривой $r(\cdot)$.
Предполагаем, что для всех кривых $r(\cdot) \in \Re $ существует
система координат такая, что для любой координаты $r_i(\cdot), i
\in 1:3,$ вектор функции $r(\cdot)$ отрезок $[0,T(r)]$ можно
разбить на не более чем на три отрезка $[0,T_{i1}], [T_{i1},
T_{i2}], [T_{i2},T_{i3}],$ на каждом из которых верны неравенства
$$
r_i'(t) \meq 0 \,\,\,\, (r_i'(t) \leq 0) \,\,\,\, \forall t \in
[0,T_{i1}],
$$
$$
r_i'(t) \leq 0 \,\,\,\, (r_i'(t) \meq 0) \,\,\,\, \forall t \in
[T_{i1}, T_{i2}],
$$
$$
r_i'(t) \meq 0 \,\,\,\, (r_i'(t) \leq 0) \,\,\,\, \forall t \in
[T_{i2},T_{i3}].
$$
Возможно, что $T_{i1}=0$. Нетрудно видеть, что $ T(r) - $ угол
конической поверхности, образованной  лучами с началом в точке
$0$, проходящими через $r(t), t \in [0, T(r)].$ Из определения
следует, что кривые $r(\cdot) \in \Re$ п.в. дифференцируемы по $t$
на отрезке $[0,T(r)]$.  Этот класс кривых был определен в
\cite{proudconvex1}.

Возьмем произвольную кривую $r(\cdot) \in \Re$. Определим функцию
$$
\Phi(t)=\varphi(r(t)) \,\,\,\, \forall t \in [0, T(r)].
$$
Покажем, что $\Phi(\cdot)-$ липшицевая на отрезке $[0, T(r)]$. Для
любых $t_1, t_2 \in [0, T(r)]$
$$
| \Phi(t_1) - \Phi(t_2) | = |\varphi(r(t_1)) - \varphi(r(t_2)|
\leq L \| r(t_1) - r(t_2) \|.
$$
Из очевидного неравенства $\| r(t_1) - r(t_2) \| \leq |t_1 - t_2
|$ имеем
$$
| \Phi(t_1) - \Phi(t_2) | \leq L |t_1 - t_2 |.
$$
Откуда следует, что $\Phi(\cdot)-$ п.в. дифференцируема на $[0,
T(r)]$, а вектор-функция $r(\cdot)$ п.в. имеем касательную
$r'(\cdot)$. Также, как и ранее, определяем вариацию
$\vee(\Phi';0,{T(r)}) $ функции $\Phi'(\cdot)$ на отрезке $
[0,T(r)]$.

\begin{thm}
Для того чтобы п.о. первой степени липшицевая функция
$\varphi(\cdot):\mathbb{R}^3 \rar \mathbb{R} $ с константой
Липшица $L$ была ПРВ функцией, необходимо, чтобы для любой
$r(\cdot) \in \Re$ и любого ее подмножества нашлась константа
$C(\varphi)>0$ такая, что
$$
\vee(\Phi';0,{T(r)}) \leq C(\varphi) + 2L \vee(r';0,{T(r)})
\,\,\,\,\,\, \forall r(\cdot) \in \Re,
$$
где
$$
\vee(r';0,{T(r)}) = \sup_{{t_i}, t_i \in [0,T(r)]} \sum_i \|
r'(t_i) - r'(t_{i+1}) \|.
$$
\label{difconvgomogthm2}
\end{thm}
\begin{rem}
Так как вариация $\vee(r';0,{T(r)})$ может быть бесконечной, то мы
требуем, чтобы теорема выполнялась для любых подмножеств кривой
$r(\cdot) \in \Re $, что означает следующее: для любой системы
отрезков $[T_i, T_{i+1} ] \subset [0, T(r)] $, $i \in 1:m,$
существуют константы $C(\varphi), 2 L$, что
$$
 \sum_1^m \vee (\Phi'; T_i,T_{i+1})<  C(\varphi) + 2 L \sum_1^m\vee (r';
 T_i,T_{i+1}).
$$
Необходимость рассмотрения системы подмножеств для кривой $r \in
\Re$ объясняется тем, что вариация $\vee (r'; 0,T(r))$ может быть
неограниченной для всей кривой $r(\cdot)$.
\end{rem}

{\bf Доказательство. Необходимость.} Пусть функция
$\varphi(\cdot)$ представима в виде разности п.о. выпуклых функций
$\varphi_i(\cdot):\mathbb{R}^3 \rar \mathbb{R}, \, i=1,2$ с
константами Липшица $L_i, i=1,2,$ соответственно. Зафиксируем
произвольную кривую $r(\cdot) \in \Re$ и соответствующую ей
систему координат в $\mathbb{R}^3$, в которой выполняются
сделанные выше предположения насчет монотонности координат
$r_i(\cdot)$ вектор-функции $r(\cdot)$. Обозначим через $\Pi$
плоскость, проходящую через начало координат перпендикулярно оси
$OZ$. Параметризуем $ r(\cdot)$ естественным образом и через
$[0,T]$ обозначим отрезок значений параметра $t$, где $T=T(r)$.
Для функций $\varphi(\cdot), \varphi_i(\cdot), \, i=1,2,$
определим соответственно функции $\Phi(\cdot), \Phi_i(\cdot),
i=1,2,$ как это делали ранее.

Поскольку \cite{kolmogorovfomin}
$$
\vee(\Phi'; 0,T) \leq \vee(\Phi_1'; 0,T)+\vee(\Phi_2'; 0,T),
$$
то достаточно доказать неравенство

\be \vee(\Phi_1';0,{T}) \leq C_1 + L_1 \vee(r';0,{T}),
\label{DiffconvthreevarGomog5} \ee где $C_1=C_1(\varphi_1).$

Разобьем шар $B_1^3(0)$ на конусы. Для этого проведем вертикальные
плоскости $\Pi_i$, проходящие через $0=(0,0,0)$ и ось $OZ$ и
образующие между собой равные углы. Части шара $B_1^3(0)$,
расположенные между плоскостями $\Pi_i$, равномерно разобьем на
конусы, внутренности которых не пересекаются, каждый из которых
образован тройкой линейно независимых векторов, принадлежащих
вертикальным плоскостям. В каждом таком конусе
$K=\mbox{con}(q_1,q_2,q_3) $ определим линейную функцию, значения
которой равны значениям функции $ \varphi(\cdot)$ на векторах $0,
q_1, q_2, q_3. $ Функцию, равную в каждом конусе построенной
линейной функции, по числу конусов $m$ обозначим через
$\varphi_{1,m}: \mathbb{R}^3 \rar \mathbb{R}. $ Нетрудно видеть,
что выпуклая оболочка градиентов всех построенных линейных функций
есть субдифференциал в нуле \cite{demvas}, \cite{demrub},
\cite{pshenichnyi}.

Обозначим, как и ранее, $\Phi_{1m}(t)= \varphi_{1m}(r(t)), \,\, t
\in [0, T]. $ Так как при $m \rar \infty$ функции
$\varphi_{1,m}(\cdot)$ равномерно на $B_1^3(0)$ сходятся к функции
$\varphi_{1}(\cdot)$, то
$$
\lim_{m \rar \infty} \rho_H (\p \varphi_1(0), \p \varphi_{1m}(0))
= 0
$$
и
$$
\lim_{m \rar \infty} \vee(\Phi_{1m}'; 0,T) = \vee(\Phi_1'; 0,T),
$$
где $\rho_H-$ метрика Хаусдорфа \cite{pshenichnyi}. То для
доказательства (\ref{DiffconvthreevarGomog5}) достаточно доказать,
что всех $m$
$$
\vee(\Phi_{1m}';0,{T}) \leq C_{1} + L_{1m} \vee(r';0,{T}),
$$
где $L_{1m}- $ константа Липшица функции $\varphi_{1m}(\cdot).$
Покажем, что

\be \vee(\Phi_{1m}';0,{T}) \leq p(l_{1m}) + L_{1m} \vee(r';0,{T}),
\label{DiffconvthreevarGomog5a} \ee где $p(l_{1m}) -$ длина
ломаной $l_{1m}$, вершинами которой являются вершины многогранника
$\p \varphi_{1m}(\cdot) $ с нормалями $r(t)$, когда $t$ пробегает
от $0$ до $T(r)$.

Пусть
$$
 \varphi_{1m} (r(t)) = \max _{v \in
 \partial \varphi_{1m}(0)} (v, r(t))=(v_m(t),r(t)), \;\;
$$
где $\partial \varphi_{1m}(0)$ - субдифференциал функции
$\varphi_{1m}(\cdot)$ в нуле, $v_{1m}(t)-$ вершины выпуклого
многогранника $\partial \varphi_{1m}(0)$ с нормалью $r(t)$.

Очевидно, что $\Phi_{1m}(\cdot)$ п.в.  дифференцируемая на $[0,t]$
и ее точки дифференцируемости функции
$\Phi_{1m}(t)=\varphi_{1m}(r(t)) -$ это точки $t \in [0,T]$, для
которых $r(t)$ принадлежит внутренности нормального конуса к
множеству $\partial \varphi_{1m}(0)$ в точках $v_{1m}(t) \in
\partial \varphi_{1m}(0)$. Для этих $t$
$$
\Phi_{1m}'(t)=(v_{1m}'(t),r(t))+(v_{1m}(t),r'(t)).
$$
Но для $t$, когда $r(t)$ принадлежит внутренности нормального
конуса, построенного в вершине $v_{1m}(t)$, вектор $v_{1m}(t)$
постоянен, а следовательно, $(v_{1m}'(t),r(t))=0.$

Поскольку кривая $r(\cdot)$ параметризована естественным образом,
то $\Vert r'(t) \Vert =1$ для любых $t \in [0,T(r)]$. Нетрудно
проверить следующую цепочку неравенств
$$
\mid \Phi_{1m}' (r(t_1)) - \Phi_{1m}' (r(t_2)) \mid = \mid
(v_{1m}(t_1),r'(t_1)) - (v_{1m}(t_2),r'(t_2)) \mid =
$$
$$
\mid (v_{1m}(t_1)- v_{1m}(t_2),r'(t_1)) +(v_{1m}(t_2),r'(t_1))-
(v_{1m}(t_2),r'(t_2)) \mid \leq
$$
$$
\leq  \Vert v_{1m}(t_1) - v_{1m}(t_2) \Vert \; \Vert r'(t_1) \Vert
+ \Vert r'(t_1) - r'(t_2) \Vert \; \Vert v_{1m}(t_2) \Vert \leq
$$
\begin{equation}
\leq \Vert v_{1m}(t_1) - v_{1m}(t_2) \Vert  + L_{1m} \Vert r'(t_1)
- r'(t_2) \Vert. \label{DiffconvthreevarGomog6}
\end{equation}
Из (\ref{DiffconvthreevarGomog6}) следует, что
$$
   \vee (\Phi_{1m}'; 0,T(r)) \leq \mbox{длина кривой}\,\,\,
   v_{1m}(t)\,\,\,
   \mbox{для } \,\,\, t \in [0, T(r)]
   + L_{1m} \cdot \vee (r'; 0,T(r))=
$$
$$
   = p(l_{1m}) + L_{1m} \cdot \vee (r'; 0,T(r)),
$$
что и требовалось показать.

Так как  $L_{1m} \rar L_1 $ при $m \rar \infty $, то для
доказательства (\ref{DiffconvthreevarGomog5}) достаточно показать,
что $p(l_{1m})$ ограничена сверху одной и той же константой для
всех $m$ и кривых $ r(\cdot) \in \Re$.

Далее без ограничения общности будем рассматривать такие ломаные
$l_{1m}$, отрезки которых $ a_k $ являются гранями размерности 1
многогранника $\partial \varphi_{1m}(0)$. Разобьем ломаную
$l_{1m}$ на участки  $s_i$, где отрезки $ a_k $ образуют углы
$\al_k$, для которых $0 \leq \al_k \leq \pi / 4 $, с
горизонтальными прямыми, принадлежащими тем же граням, что и
отрезки $a_k$. Оценим сумму
$$
\sum_{a_k \in \bigcup_i s_i} | a_k |,
$$
где $| a_k | -  $ длина отрезка $a_k$. Длина горизонтальной
проекции отрезка $a_k$ равна $| a_k | \cos \, \al_k. $ Заметим,
что сужение (расширение) по горизонтали граней многогранника
$\partial \varphi_{1m}(0)$ по мере приближения к граням с
нормалями из плоскости $\Pi$, происходит за счет расширения
(сужения) по горизонтали соседних граней.

Покажем, что сумма

\be \sum_{a_k \in \bigcup_i s_i} | a_k | \, \cos \, \al_k,
\label{DiffconvthreevarGomog7} \ee ограничена сверху константой,
независящей от $m$ и $r(\cdot) \in \Re$.

Рассмотрим функцию $\tilde{\varphi}_{1m}(\cdot): \mathbb{R}^2 \rar
\mathbb{R},$ равную по определению
\begin{equation}
\tilde{\varphi}_{1m}(q) = \varphi_{1m}(q) \,\,\, \forall q \in
\Pi.
\label{DiffconvthreevarGomog7a}
\end{equation}
Нетрудно видеть, что $\tilde{\varphi}_{1m}(\cdot)$ есть п.о.
выпуклая функция от двух переменных, определенная на плоскости
$\Pi,$ имеющая ту же константу Липшица, что и функция
${\varphi}_{1m}(\cdot)$, т.е. $L_{1m}$. Поэтому длина ломаной,
ограничивающей субдифференциал $\partial \tilde{\varphi}_{1m}(0)$
функции $\tilde{\varphi}_{1m}(\cdot)$ в нуле, ограничена сверху
константой $C_1$, независящей ни от выбранной кривой $r(\cdot) \in
\Re $, ни от расположения плоскости $\Pi$.

Заметим, что градиенты функции $\tilde{\varphi}_{1m}(\cdot)$ есть
двумерные векторы, длины которых равны длинам проекции на
плоскость $\Pi$ градиентов функции ${\varphi}_{1m}(\cdot)$,
вычисленных для $q \in \Pi$. Поэтому сумма
(\ref{DiffconvthreevarGomog7}) ограничена сверху константой $C_1$,
независящей ни от выбранной кривой $r(\cdot) \in \Re $, ни от
плоскости $\Pi$, т.е.
$$
\sum_{a_k \in \bigcup_i s_i} | a_k | \, \cos \, \al_k \leq C_1.
$$
Так как $ \cos \, \al_k \meq 1/ \sqrt{2}$, то

\be \sum_{a_k \in \bigcup_i s_i} | a_k | \, \leq C_1 \, \sqrt{2}.
\label{DiffconvthreevarGomog8} \ee

Отрезки ломаной $l_{1m}$, которые не принадлежат $ \bigcup _i
s_i$, обозначим через $ b_{k}.$ Сумму длин этих отрезков можно
оценить следующим образом.

Разобьем шар $ B_1^3(0) $ на $j$ равных вертикальных полос. Для
каждой вертикальной полосы определим свою систему координат, ось
OZ которой принадлежит плоскости $ \Pi $ и перпендикулярна одной
из из вертикальных плоскостей $ \Pi_i $, образующей вертикальную
полосу. Число $j$ выбираем так, чтобы вертикальные плоскости из
рассматриваемой полосы и плоскость, перпендикулярная оси OZ,
образовывали между собой достаточно малые углы. Тогда сумма длин
отрезков $ b_{k}$ ломаной $l_{1m}$ из одной и той же вертикальной
полосы не больше, чем ( по аналогии с
(\ref{DiffconvthreevarGomog8}))
$$
C_1 \, (\sqrt{2} + \delta(j)),
$$
где $\delta(j) \rar 0$ при $j \rar \infty$.

Поскольку кривая $ r(\cdot) \in \Re,$ то согласно требованиям
монотонности координат $r_i(\cdot)$,  проекция $ r(\cdot)$ на
плоскость $\Pi$ при изменении $t$ от $0$ до $T$ вращается вокруг
оси $OZ$ по (против) часовой стрелке, а в каждой вертикальной
полосе $r(\cdot)$ бывает ровно один раз без повторений с
монотонным возрастанием (убыванием) проекции на ось $OZ$, то,
просуммировав по всем вертикальным полосам, получим \be \sum_{b_k
} | b_k | \, \leq C_1 (\sqrt{2} + \delta(j))j   .
\label{DiffconvthreevarGomog9} \ee

Из (\ref{DiffconvthreevarGomog8}) и (\ref{DiffconvthreevarGomog9})
получим, что длина ломаной $l_{1m}  $ ограничена сверху величиной,
независящей от $m$ и кривой $r(\cdot) \in \Re$, что и требовалось
доказать. Необходимость доказана.

\vspace{0.5cm}

\subsection{Множество кривых на сфере $S^2_1(0)$, ограничивающих
выпуклые компактные множества}

\vspace{0.5cm}

При доказательстве достаточности теоремы (\ref{difconvgomogthm2})
было показано, что поиск кривой $r(\cdot)$ ограничивается поиском
направления $g$ и вектора $x \in S_1^2(0)$, $g \perp x$. Функция $
\varphi_m(\cdot)$, построенная по произвольному разбиению шара
$B_1^3(0)$ на конусы, является ПРВ функцией, но растяжение
многогранника $G_m$, соответствующего этому разбиению, вдоль
направления $g$ неограниченно возрастает при $m \rar \infty,$ так
что выбрать выпуклые многогранники к которому сходятся в метрике
Хаусдорфа $G'_m$ при $m \rar \infty$, получающиеся после
растяжения $G_m$, в случае, когда $\varphi(\cdot)$ не ПРВ функции,
невозможно.

В качестве вектора $g$ можно взять тот, вдоль которого растяжение
многогранника $G_m$ согласно описанному выше алгоритму для
преобразования его в выпуклый многогранник $G'_m$ происходит
наибольшим образом при $m \rar \infty$.

Мы определим более узкий класс кривых, с помощью которого
сформулируем необходимые условия представимости функции в виде
разности выпуклых.

Будем понимать под {\em координатной плоскостью}  в $\mathbb{R}^3$
плоскость, образованную двумя координатными осями.

\begin{defi}
Выпуклым  множеством на сфере $S_1^2(0)$ назовем такое множество
$Q \subset S_1^2(0)$, проекция которого на любую координатную
плоскость есть выпуклое множество.
\end{defi}

Рассмотрим на сфере $S_1^2(0)$ множество кривых $\wp$,
ограничивающих на $S_1^2(0)$ {\em выпуклые компактные множества}.

В дальнейшем мы будем рассматривать выпуклые компактные множества
$Q \subset S_1^2(0)$. Покажем, что для всех $ r(\cdot) \in  \wp $
верно неравенство \be
  \vee(r'; 0, T) < C_2,
\label{DiffconvthreevarGomog12} \ee где $C_2 -$ некоторая
константа. По определению
$$
\vee(r';0,{T(r)}) = \sup_{{t_i}, t_i \in [0,T(r)]} \sum_i \|
r'(t_i) - r'(t_{i+1}) \|.
$$
Но
$$
\| r'(t_i) - r'(t_{i+1}) \| \leq | r_1'(t_i) - r_1'(t_{i+1}) | + |
r_2'(t_i) - r_2'(t_{i+1}) | + | r_3'(t_i) - r_3'(t_{i+1}) |,
$$
где $r_i(\cdot), i \in 1:3, -$ координаты вектор-функции
$r(\cdot).$ Так как кривая $r(\cdot)$ параметризована естественным
образом, то $\| r'(\cdot) \| =1, $ а следовательно, $| r'_i
(\cdot) | \leq 1 $ для всех $i \in 1:3$. Поэтому пары координат
вектор-функции $r'(\cdot)$ ограничивают на соответствующей
координатной плоскости выпуклое множество, диаметр которого
равномерно ограничен для всех кривых $r(\cdot) \in  \wp $. Отсюда
следует ограниченность вариации проекции $r'(\cdot)$  на
произвольную координатную плоскость. А поэтому верно неравенство
(\ref{DiffconvthreevarGomog12}) для всех $ r(\cdot) \in \wp $.

Итак, для множества кривых $\wp $ и  для ПРВ функции
$\varphi(\cdot)$  неравенство, фигурирующее в формулировке теоремы
\ref{difconvgomogthm2a}, с необходимостью выполняется. Поскольку
вариация $\vee(r';0,{T(r)})$ ограничена равномерно для всех кривых
из $\wp $ и кривую $r(\cdot)$ можно разбить на участки, каждая из
которых образует с какой-то координатной плоскостью угол, не
больший $\frac{\pi}{4}$, а также выпуклые множества $-$ это
звездные множества, то приходим к неравенству
$$
\vee(\Phi';0,T) < C_3 \,\,\,\, \forall r(\cdot) \in \wp.
$$
Из проведенных рассуждений следует теорема.
\begin{thm}
Для того чтобы п.о. первой степени липшицевая функция
$\varphi(\cdot):\mathbb{R}^3 \rar \mathbb{R} $ с константой
Липшица $L$ была ПРВ функцией, необходимо, чтобы для всех кривых
$r(\cdot) \in \wp$, параметризованных естественным образом с
параметром $t \in [0,T], \, T=T(r),$ нашлась константа
$C_3=C_3(\varphi)>0$ такая, что \be \vee(\Phi';0,T) \leq C_3
\,\,\,\,\,\, \forall r(\cdot) \in \wp.
\label{DiffconvthreevarGomog13} \ee \label{difconvgomogthm3}
\end{thm}

\vspace{0.5cm}

\subsection{Звездные множества. Необходимые и достаточные условия
представимости положительно однородной функции от трех переменных
в виде разности выпуклых функций}

\vspace{0.5cm}

Заметим, что класса кривых $ \hat {\Re} $ и даже более широкого
класса $ {\Re} $ недостаточно для представимости функции в виде
разности выпуклых. Это видно уже из того, что отрезки $a_k$  с
максимальной длиной в каждой вертикальной полосе, являющиеся
субдифференциалами Кларка вогнутых двугранных углов многогранной
функции $\varphi_{m}(q)$, построенной по $\varphi(q)$ путем
разбиения шара $B_1^3(0)$ на конусы,  могут не помещаться на одну
и ту же кривую из $ {\Re} $. Но сумма длин отрезков $a_k$
определяет растяжение по горизонтали многогранника $\p
\varphi_{m}(0)$ и может быть бесконечной при $m \rar \infty$.

Отметим, что классы кривых $\hat \Re $ и $\Re$ принадлежат более
широкому классу кривых $\wp(\Pi)$ на сфере $S^2_1(0)$,
параметризованных естественным образом $r(t)$, $0 \in [0, T(r)]$,
проекции которых на некоторую плоскость $ \Pi \in \mathbb{R}^3$,
проходящую через начало координат, ограничивают на $\Pi$ звездные
множества. Плоскость $\Pi$ выбирается произвольным образом.

Отметим, что из неравенства
$$
    \mid r(t_2) - r(t_1) \mid \leq \mid t_2  - t_1 \mid
$$
следует, что вектор-функция $r(\cdot)$ п.в. на $[0, T(r)]  $ имеет
касательную $r'(\cdot)$.

Дадим определение звездной области.

\begin{defi} Область ${D}$ на плоскости $\Pi $
называется звездной, если существует точка $O \in \mbox{int } {D}
$, для которой отрезок, соединяющий точку $O$ с любой точкой
области ${D}$, целиком принадлежит области $D$.
\end{defi}

Ясно, что в качестве точки $O$ для кривых на сфере $S_1^{n-1}(0)$
из множества $\wp(\Pi)$ для некоторой плоскости $\Pi$, проходящей
через начало координат, можно взять начало координат.

%%%%%%%%%%%%%%%%%%%%%%%%%%%%%%%%%%%%%
Возьмем произвольную плоскость $\Pi$, $0 \in \Pi$.
Воспользовавшись классом кривых $\wp(\Pi)$, можно дать необходимые
и достаточные условия представимости п.о. функции $\varphi(\cdot):
\mathbb{R}^3 \rar \mathbb{R}$  в виде разности выпуклых.

Для $r(\cdot) \in \wp(\Pi) $ выделим те ее участки $[T_i, T_{i+1}]
\subset [0,T(r)]$, где производная $r'(t), t \in [0, T(r)],$
составляет с плоскостью $\Pi$ угол $\delta $: $0 < \delta \leq
\frac{\pi }{4} $. Обозначим объединение таких участков через
$\Lambda(r)$.

\begin{thm}
Для того чтобы п.о. первой степени липшицевая функция
$\varphi(\cdot):\mathbb{R}^3 \rar \mathbb{R} $ с константой
Липшица $L$ была ПРВ функцией, необходимо и достаточно, чтобы для
любой  плоскости $\Pi$,$0 \in \Pi$, и  $r(\cdot) \in \wp(\Pi)$, а
также для любого подмножества параметров $[T_i, T_{i+1} ] \subset
\Lambda(r) $ нашлась константа $C(\varphi)>0$ такая, что
$$
 \sum_1^m \vee (\Phi'; T_i,T_{i+1})<  C(\varphi) + 2 L \sum_1^m\vee (r';
 T_i,T_{i+1}),
$$
где
$$
\vee(r';0,{T(r)}) = \sup_{{t_i}, t_i \in [0,T(r)]} \sum_i \|
r'(t_i) - r'(t_{i+1}) \|.
$$
\label{difconvgomogthm2a}
\end{thm}
\begin{rem}
Так как вариация $\vee(r';0,{T(r)})$ может быть бесконечной, то мы
требуем, чтобы теорема выполнялась для любых подмножеств множества
$\Lambda(r)$, что означает следующее: для любой системы отрезков
$[T_i, T_{i+1} ] \subset \Lambda(r) $, $i \in 1:m,$ существует
константа $C(\varphi)$, что
$$
 \sum_1^m \vee (\Phi'; T_i,T_{i+1})<  C(\varphi) + 2 L \sum_1^m\vee (r';
 T_i,T_{i+1}).
$$
Необходимость рассмотрения системы подмножеств множества
$\Lambda(r) $ объясняется тем, что вариация $\vee (r'; 0,T(r))$
может быть неограниченной для всей кривой $r(\cdot)$.
\end{rem}

%%%%%%%%%%%%%%%%%%%%%%%%%%%%%%%%%%%%%%%%%%%%
{\bf Необходимость}. Доказательство необходимости полностью
повторяет доказательство  теоремы \ref{difconvgomogthm2}. Отметим
лишь, что ограниченность длины кривой $l_m$ неравенства
(\ref{DiffconvthreevarGomog5a}) следует из
(\ref{DiffconvthreevarGomog8}) и того факта, что проекция $l_m$ на
плоскость $\Pi$ ограничивает звездное множество. Поэтому длина
проекции $l_{1m}$ на плоскость $\Pi$ не больше, чем длина кривой,
ограничивающей $\p \tilde \varphi_{1m} (0)$. Функция $\tilde
\varphi_{1m}(0)$ определяется в (\ref{DiffconvthreevarGomog7a}).
Так как все отрезки ломаной $l_{1m}$ образуют углы с плоскостью
$\Pi$, не большие $\frac{\pi}{4}$, то  длину $l_{1m}$ можно
оценить через длину ее проекции, использовав неравенства,
аналогичные (\ref{DiffconvthreevarGomog8}).

{\bf Достаточность}. Пусть задана п.о. липшицевая функция
$\varphi(\cdot):\mathbb{R}^3 \rar \mathbb{R}$ с константой Липшица
$L$. Разобьем шар $B_1^3(0)$ на непересекающиеся по внутренности
конусы. Для этого проведем плоскости, проходящие через начало
координат и две диаметрально противоположные точки $N$ и $S$,
образующие между собой равные углы. Проведем также плоскости,
перпендикулярные отрезку $NS$ и расположенные на одинаковом
расстоянии  друг от друга. Точки пересечения окружностей,
получающихся в результате пересечения проведенных плоскостей со
сферой $ S_1^2(0)$, обозначим через $\{ q_i \}, \, i \in 1:I.$
Введем векторы, начало которых в $0$, а концы $-$ в точках $\{ q_i
\}, \, i \in 1:I.$ Эти векторы будем обозначать, как и точки,
через $\{ q_i \}$. Конусы $K_i=\mbox{con} (q_{i}, q_{i+1},
q_{i+2}) $, являющиеся конической оболочкой ближайших линейно
независимых векторов $ \{ q_{i}, q_{i+1}, q_{i+2} \} $, задают
упомянутое разбиение шара $B_1^3(0)$.

Построим по данному разбиению п.о. многогранную функцию
$\varphi_m(\cdot): \mathbb{R}^3 \rar \mathbb{R}, $  где $m-$ число
конусов. Для этого в каждом конусе $K_i$ зададим линейную функцию
$\eta_i(\cdot):\mathbb{R}^3  \rar \mathbb{R}, $ принимающую в
точках $0, q_{i}, q_{i+1}, q_{i+2}$ следующие значения:
$$
\eta_i(0)=0, \, \eta_i(q_i)=\varphi_m(q_i)= \varphi(q_i), \,\,
\eta_i(q_{i+1})=\varphi_m(q_{i+1}) = \varphi(q_{i+1}), \,
$$
\be \eta_i(q_{i+2})=\varphi_m(q_{i+2}) = \varphi(q_{i+2}).
\label{DiffconvthreevarGomog10} \ee Функция, совпадающая в каждом
конусе $K_i$ с соответствующей функцией $\eta_i(\cdot)$, и есть
п.о. многогранная функция $\varphi_m(\cdot)$.

Покажем, что $\varphi_m(\cdot) -$ липшицевая функция с константой
Липшица $L_m$, ограниченной сверху константой,  независящей от
$m$. Оценим сверху норму градиента $p_i$ функции $\eta_i(\cdot)$.
Пусть $\Pi_{i1}, \, \Pi_{i2}, \, \Pi_{i3} - $ плоскости - грани
конуса $K_i$, а $\Pi -$ плоскость, перпендикулярная $NS$ и
проходящая через начало координат.

Пусть вектор $p_i$ образует среди векторов $q_{i}, q_{i+1},
q_{i+2}$ наименьший угол $\beta_i $ с вектором $q_i$. Без
ограничения общности будем считать, что $0 \leq \beta_i \leq \pi /
4.$ В противном случае всегда можно прибавить к функции
$\varphi(\cdot)$ такую п.о. выпуклую функцию с положительными
значениями на сфере $S_1^2(0)$, чтобы это неравенство выполнялось.
Мы не уменьшаем общности рассуждений по трем причинам.  Во-первых,
для представления функции в виде разности выпуклых прибавление или
вычитание выпуклой функции роли не играет. Во-вторых, для п.о.
выпуклой функции от трех переменных неравенство, фигурирующее в
условии теоремы, для произвольной кривой $r(\cdot) \in \wp(\Pi)$
было доказано выше. В третьих, у нас есть неравенство для суммы
вариаций функций, которым мы пользовались при доказательстве
необходимости.

Из (\ref{DiffconvthreevarGomog10}) и по свойству п.о. функций
\cite{demrub} имеем

\be (p_i, q_{i})= \eta_i(q_{i})=\varphi_m(q_{i})= \varphi(q_{i})
\leq L \| q_{i} \|. \label{DiffconvthreevarGomog11} \ee Здесь
$(p_i, q_{i}) = \| p_i \| \, \| q_i \| \cos \beta_i -  $ скалярное
произведение векторов $p_i, \, q_i$, а  $\| p_i \| , \, \| q_{i}
\| -$ их длины.

Из (\ref{DiffconvthreevarGomog11}) имеем
$$
(p_i, q_{i}) = \| p_i \| \, \| q_i \| \cos \beta_i \leq L \| q_{i}
\|.
$$
Откуда
$$
\| p_i \| \leq \frac{L}{\cos \beta_i} \leq L \sqrt{2},
$$
так как $1 \ \meq \cos \beta_i \meq 1 / \sqrt{2},$ что и
требовалось доказать. Таким образом, $\varphi_m(\cdot) -$ п.о.
липшицевая многогранная функция с константой Липшица, независящей
от $m$.

Заметим, что линейные функции, определенные в соседних конусах, с
общим вектором $q_i$, имеют одну и ту же производную по
направлению $q_i$. Соединяя ломаной вершины градиентов (начальная
точка в нуле) линейных функций, построенных по соседним конусам с
общим вектором $q_i$, мы получим либо замкнутый восьмиугольник,
либо замкнутый четырехугольник в плоскости $\Pi_{q_i}$,
перпендикулярной вектору $q_i$. Восьмиугольник $A_1A_2 \dots  A_8$
не будет выпуклым и его стороны могут иметь самопересечения только
тогда, когда среди двугранных углов, субдифференциалами Кларка
\cite{clark} которых являются стороны многоугольника $A_1A_2 \dots
A_8$, есть как выпуклые. так и вогнутые двугранные углы. Под
выпуклым (вогнутым) двугранным углом будем понимать выпуклую
(вогнутую) функцию, равную максимуму (минимуму) линейных функций,
графики которых образуют данный двугранный угол. Сказанное насчет
сторон многоугольника $A_1A_2 \dots A_8$ следует из того, что
вектор $a_i \in \mathbb{R}^2$, равный проекции вектора $OA_i$ ($O
-$ начало координат) на плоскость $\Pi_{q_i}$ и записанный в
системе координат плоскости $\Pi_{q_i}$, есть градиент линейной
функции двух переменных, определенной на плоскости $\Pi_{q_i}$ и
равной $\varphi_m(q) - \varphi_m(q_i), q \in \Pi_{q_i},$ в
соответствующем секторе. Если соединить концы градиентов всех
линейных функций, построенных по соседним конусам отрезками, то
получим многогранник $G_m$, ребра и грани которого могут
пересекаться. Если функция $\varphi(\cdot)$ выпуклая, то
многогранник $G_m -$ выпуклый. В общем случае, когда $G_m$ не
выпуклый, преобразуем $G_m$ в выпуклый многогранник $G'_m$.

Рассмотрим одну из граней многогранника $G'_m$. Это многоугольник
$A_1A_2 \dots A_8$. Нашей задачей будет преобразовать
многоугольник $A_1A_2 \dots A_8$ в выпуклый многоугольник
$A'_1A'_2 \dots A'_8$ таким образом, чтобы стороны, равные
субдифференциалам выпуклых двугранных углов, не уменьшались по
длине и оставались параллельными самим себе. Это надо для того,
чтобы разность между выпуклой п.о. функцией $\varphi'_m(\cdot):
\mathbb{R}^3 \rar  \mathbb{R}$, равной
$$
\varphi'_m(q) = \max_{v \in G'_m} \, (v, q), \,\,\,\, \forall q
\in \mathbb{R}^3
$$
и функций $\varphi_m(\cdot)$ была выпуклой. Отсюда получаем, что
функция $\varphi_m(\cdot)$ есть ПРВ функция.

Многоугольник $A'_1A'_2 \dots A'_8$ после преобразования $A_1A_2
\dots A_8$ станет одной из граней выпуклого многогранника $G'_m$.
Рассмотрим далее все случаи выпуклости и вогнутости двугранных
углов, построенных по восьми конусам с общим вектором $q_i$.
Пронумеруем все эти восемь конусов по часовой стрелке (см. рис.
2). Двугранные углы, построенные по конусам  1-2 либо 5-6, мы
назовем {\em горизонтальными}, а двугранные углы 3-4 либо 7-8 $-$
{\em вертикальными} двугранными углами.

Пусть двугранный угол 1-2 вогнутый, а остальные двугранные углы
$-$ выпуклые. Восьмиугольник для этого случая изображен на рисунке
3. Отрезки $A_1A_2$ и $A_5A_6$ параллельны друг другу, так как они
перпендикулярны одной и той же вертикальной плоскости. Преобразуем
$A_1A_2 \dots A_8$ в выпуклый многоугольник $A'_1A'_2 \dots A'_8$,
сдвигая ломаную $A_2A_3A_4A_5$ параллельно самой себе так, чтобы
точки $A_1$ и $A_2$ совпали (см. рис. 3). При этом
$$
 | A'_5 A_6 | = | A_5 A_6 | +| A_1 A_2 |.
$$
Пусть теперь оба двугранных угла 1-2 и 5-6 вогнутые. Для этого
случая многоугольник $A_1A_2 \dots A_8$ изображен на рис 4.
Предположим, что длина отрезка $A_5 A_6$ больше длины отрезка $A_1
A_2$. Сдвигая параллельно самой себе ломаную $A_2A_3A_4A_5$ так,
чтобы точки $A_5$ и $A_6$ совпали и увеличив при этом длину
стороны $A_1 A_2$ на длину $| A_5 A_6 | - |A_1 A_2 |$, получим
выпуклый многоугольник $A'_1A'_2 \dots A'_8$ (см. рис. 4).
Описанное преобразование многоугольника $A_1A_2 \dots A_8$ в
$A'_1A'_2 \dots A'_8$ назовем {\em растяжением по горизонтали}.
При таком преобразовании многоугольников, построенных по конусам
из одной и той же вертикальной полосы, все их горизонтальные
двугранные углы будут выпуклыми, если эти многоугольники мы
растянем по горизонтали на наибольший по длине отрезок, являющийся
субдифференциалом Кларка вогнутого горизонтального двугранного
угла из этой вертикальной полосы.

Рассмотрим случай, когда вогнутыми будут вертикальные двугранные
углы. Пусть вогнутым будет только двугранный угол 3-4. В этом
случае многоугольник $A_1A_2 \dots A_8$ показан на рис. 5. Для
того чтобы $A_1A_2 \dots A_8$ стал выпуклым, надо ломаную $A_4A_5
A_6A_7$ перенести параллельно самой себе вниз так, чтобы точки
$A_3$ и $A_4$ совпали. При этом длины сторон $A_6 A_7$ и $A_7 A_8$
увеличим с таким расчетом, чтобы получился выпуклый замкнутый
многоугольник.

Пусть вогнутыми будут двугранные углы 3-4 и 7-8. Многоугольник для
этого случая показан на рис. 6. Ломаную $A_4A_5 A_6A_7$ перенесем
параллельно самой себе на вектор, длина которого равна наибольшей
длине отрезка, являющегося субдифференциалом Кларка вогнутого
вертикального двугранного угла. Описанное преобразование назовем
{\em растяжением по вертикали}. Для того чтобы у всех
многоугольников, построенных по конусам из одной и той же
горизонтальной полосы, после растяжения по вертикали не было
вогнутых вертикальных двугранных углов, надо все указанные
многоугольники растянуть по вертикали на максимальный по длине
отрезок, являющийся субдифференциалом Кларка вогнутого
вертикального двугранного угла, построенного по конусам из
рассматриваемой горизонтальной полосы.

Все другие случаи выпуклых и вогнутых двугранных углов 1-2, 2-3,
$\dots$ 7-8 рассматриваются аналогично. Растяжение по горизонтали
и вертикали многоугольника $A_1A_2 \dots A_8$ оценивается сверху
по его длинам сторон. Так как число различных случаев выпуклости и
вогнутости двугранных углов 1-2, 2-3, $\dots$ 7-8 конечно, то
растяжение по горизонтали многоугольников из одной и той же
вертикальной полосы и растяжение по вертикали многоугольников из
одной и той же горизонтальной полосы таким образом, чтобы все
многоугольники стали выпуклыми, оценивается сверху по максимальным
длинам горизонтальных  или вертикальных сторон многоугольников в
каждой из полос.

Получившийся в итоге многогранник $G'_m$ будет выпуклым. Покажем,
что независимо от $m$ диаметр $G'_m$ ограничен. Предположим
противное, а именно: существует такое направление $g \in
\mathbb{R}^3,$ вдоль которого растяжение $G_m$ при $m \rar \infty$
неограничено. Докажем тогда, что на сфере $S_1^2(0)$ существует
точка $x$ такая, что при неограниченном разбиении произвольного
конуса, содержащем внутри точку $x$, на подконусы, как это было
сделано ранее, растяжение в направлении $g$ многоугольников,
построенных по этим конусам и являющихся сторонами многогранника
$G_m$, бесконечно.

Разобьем шар $B_1^3(0)$ на конечное число конусов $K_iб ш \in
1:I$. Выберем из них тот, где при неограниченном уменьшении
разбиения растяжение $G_m$ вдоль направления $g$ бесконечно.
Разобьем выбранный конус на $K_{i_1}$ на конечное число подконусов
и проделаем ту же процедуру, что и ранее. Выбранный новый конус
обозначим через $K_{i_2}$. Так как $K_{i_2} \cap S_1^2(0) \subset
K_{i_1}\cap S_1^2(0)$, то, продолжая процесс, приходим к искомой
точке $x \in S_1^2(0)$.

Возьмем теперь произвольный конус, содержащий во внутренности
вектор $x$, и разобьем ее на меньшие конусы $K_i$ с общим вектором
$x \in S_1^2(0)$. Выберем из них тот, где при неограниченном
уменьшении разбиения растяжение вдоль направления $g$ бесконечно.
Выбранный конус опять разбиваем на подконусы и так далее.
Последняя процедура позволяет выбрать семейство вложенных конусов
с общим вектором $x$, которое обозначим через $\aleph$, где
растяжение многогранника $G_m$ при $m \rar \infty$ вдоль
направления $g$ бесконечно. Поскольку функции $\varphi_m(\cdot)$
для любого $m$ являются п.о. функциями, то согласно описанной выше
процедуре поиска векторов $x$ и $g$  сумма проекций на направление
$g$ длин отрезков, являющихся субдифференциалами Кларка вогнутых
двугранных углов функции $\varphi_m(\cdot)$, наибольшая, когда
вектор $g$ перпендикулярен вектору $x$. Без ограничения общности
будем считать, что одна из граней всех конусов семейства $\aleph$
принадлежит плоскости, которую обозначим через $\Pi$ и $x, g \in
\Pi$ и $g \perp x. $

Выберем в $\mathbb{R}^3 $ систему координат. две оси которой OX и
OY принадлежат плоскости $\Pi$, а ось OZ перпендикулярна плоскости
$\Pi$. Будем разбивать конусы $K_i$ на подконусы. Для этого строим
вертикальные плоскости, проходящие через OZ, а также плоскости,
перпендикулярные оси OZ. По построенному разбиению строим, как это
делали ранее, многогранник $\tilde{G}_m$.   Ясно, что существует
такое разбиение $K_i$, что растяжение многогранника $\tilde{G}_m$
вдоль направления $g$ может быть как угодно большое при $m \rar
\infty$. В каждой вертикальной полосе, пересекающей конусы $K_i$,
отметим вогнутый двугранный угол, проекция субдифференциала Кларка
которого на направление $g$ максимальна. Далее берем конус
$K_{i_1} \subset K_i$, принадлежащий семейству $\aleph$.
Аналогично предыдущему разбиваем конус $K_{i_1}$ на меньшие по
включению конусы. Опять в каждой вертикальной полосе отмечаем
вогнутый двугранный угол с максимальной проекцией субдифференциала
Кларка на направление $g$ и так далее. Все отмеченные отрезки мы
можем соединить кривой $r(\cdot) \in \wp(\Pi)$. Кроме того, путем
замены некоторых сторон многоугольников $A_1A_2 \dots A_8$ другими
сторонами, можно добиться, чтобы общая сумма длин не уменьшалась.
Это можно сделать, так как грани многогранника $G_m$ пересекаются
и каждая сторона многоугольника $A_1A_2 \dots A_8$ есть сторона
другого многоугольника, построенного по соседним конусам. Так,
например, длина стороны $A_5 A_6$ многоугольника $A_1A_2 \dots
A_8$, изображенного на рис. 4, не превосходит длину стороны $A_4
A_5$. Как только мы нашли вектор $x \in S_1^2(0)$, в произвольной
окрестности которого функция $\varphi(\cdot)$ не является ПРВ
функцией, а также направление $g \perp x$, вдоль которого сумма
отрезков, являющихся субдифференциалами Кларка вогнутых двугранных
углов функции $\varphi_m(\cdot)$ из этой окрестности, бесконечна
при $m \rar \infty$, то мы можем выбирать подмножество кривой
$r(\cdot) \in \wp(\Pi)$, для которой \be \sum_i \vee(r'; T_i,
T_{i+1}) < \infty. \label{DiffconvthreevarGomog11a} \ee и в то же
время \be \sum_i\vee(\Phi';T_i, T_{i+1}) = \infty,
\label{DiffconvthreevarGomog11b} \ee

Поскольку сумма длин отрезков вдоль направления $g$ бесконечна, то
всегда можно выбрать кривую $r(\cdot) \in \wp(\Pi)$, обходящую
выбранные отрезки в каждой вертикальной полосе. Так как длину
каждого отрезка можно оценить сверху вариацией функции
$\Phi'(t)=\varphi'_t(r(t))$ (см. \cite{lupikovdissertation}) по
соответствующему отрезку значений параметра $t$, то верны
неравенства (\ref{DiffconvthreevarGomog11a} ) и
(\ref{DiffconvthreevarGomog11b})

Поэтому константы $C$, о которой говорится в условии теоремы, не
существует.

Без ограничения общности будем считать, что
$$
\lim_{m \rar \infty} \, \rho_H(G'_m, G_1) = 0,
$$
где $\rho_H - $   метрика Хаусдорфа \cite{demvas}. Определим п.о.
выпуклую функцию $\varphi_{1m}(\cdot) : \mathbb{R}^3 \rar
\mathbb{R}$
$$
\varphi_{1m}(q)=\max_{v \in G_{1m}} \, (v,q).
$$
Тогда $\varphi_{1m}(\cdot)$ равномерно на $B_1^3(0)$ сходится к
п.о. выпуклой функции  $\varphi_{1}(\cdot) : \mathbb{R}^3 \rar
\mathbb{R}$, которая по определению есть
$$
\varphi_{1}(q)=\max_{v \in G_{1}} \, (v,q).
$$
Покажем, что разность
$$
\varphi_{2m}(q) = \varphi_{1m}(q) - \varphi_{m}(q) \,\,\, \forall
q \in \mathbb{R}^3
$$
есть также выпуклая п.о. первой степени функция. Действительно,
поскольку многогранник $G'_m$ был получен из $G_m$ растяжением по
вертикали и горизонтали с параллельным переносом граней и при этом
длины отрезков, являющихся субдифференциалами выпуклых двугранных
углов, не уменьшаются, то все двугранные углы функции
$\varphi_{2m}(\cdot) -$ выпуклые, а следовательно, и сама функция
$\varphi_{2m}(\cdot) -$ выпуклая и липшицевая с константой,
независящей от $m$. Отсюда следует, что из последовательности
функций $\{ \varphi_{2m}(\cdot) \}$ можно выбрать равномерный
предел на множестве $B_1^3(0).$

Пусть
$$
\varphi_{2}(q) = \lim_{m \rar \infty} \, \varphi_{2m}(q) \,\,\,
\forall q \in B_1^3(0).
$$
Очевидно, что $ \varphi_{2}(\cdot) -$  выпуклая п.о. функция и для
$\varphi(\cdot), \, \varphi_1(\cdot), \, \varphi_2(\cdot)  $ верно
равенство
$$
\varphi(q) = \varphi_1(q) - \varphi_2(q)  \,\,\,\,  \forall q \in
B_1^3(0),
$$
что и требовалось доказать. Достаточность и теорема доказаны.
$\Box$

\vspace{0.5cm}

\subsection{Геометрическая интерпретация теоремы для положительно однородной функции
первой степени  трех переменных}

\vspace{0.5cm}

Перефразируем теорему \ref{difconvgomogthm3}, придав им более
геометрический характер. Введем понятие поворота кривой $r(\cdot)$
на графике $\Gamma_{\varphi } = \{ (x, z) \in \mathbb{R}^4 \mid z
= \varphi(x), \, x \in \mathbb{R}^3 \}.$

Рассмотрим на $\Gamma_{\varphi }$ кривую $ R(t)=(r(t),{\varphi
}(r(t))), t \in [0, T(r)],$ где $r(\cdot) \in \wp(\Pi)$ для
некоторой плоскости $\Pi$, проходящей через начало координат. Так
как функция ${\varphi }(\cdot)$ есть липшицевая, то п.в. на
$[0,T(r)]$ существует производная $R'(\cdot),$ которую обозначим
через $\tau(\cdot)=R'(\cdot),$  а множество точек, где она
существует, $-$ через $N_{\varphi }$.

Как  и раньше, для $r(\cdot) \in \wp(\Pi) $ выделим те ее участки
$[T_i, T_{i+1}] \subset [0,T(r)]$, где производная $r'(t), t \in
[0, T(r)],$ составляет с плоскостью $\Pi$ угол $\delta $: $0 <
\delta \leq \frac{\pi }{4} $. Обозначим объединение таких участков
через $\Lambda(r)$.

\begin{defi}  Поворотом кривой $R(\cdot)$ на многообразии $\Gamma_f$
назовем величину
$$
sup_{ \{t_i\} \subset N_{\varphi }} \,\, \sum_i \Vert \tau(t_i)/
\Vert \tau (t_i) \Vert -  \tau(t_{i-1})/ \Vert \tau (t_{i-1})
\Vert \Vert = O_{\varphi }.
$$
\end{defi}

Таким образом, поворот $O_{\varphi }$ кривой $R(\cdot)$ есть
верхняя грань суммы углов между касательными $\tau(t)$  для $t \in
[0,T(r)].$

\begin{thm} Для того, чтобы произвольная липшицевая п.о. степени 1
функция $x \rightarrow {\varphi }(x) :\mathbb{R}^3 \rightarrow
\mathbb{R}$ была ПРВ функцией, необходимо и достаточно, чтобы для
любой плоскости $\Pi$, $0 \in \Pi$, и любой $r(\cdot) \in
\wp(\Pi)$ параметризованной естественным образом: $t \in [0,
T(r)]$, а также любого подмножества параметров $[T_i, T_{i+1} ]
\subset \Lambda(r) $ существовали константы $c_4=c_4({\varphi
})>0$, $c_5=c_5({\varphi})>0$ такие, что сумма поворотов
$O_{\varphi, i}$ кривой $R(\cdot)$ на $[T_i, T_{i+1}] $ для  всех
$i$ ограничена сверху неравенством
\begin{equation}
\sum_1^m O_{\varphi, i } \leq   c_4 + c_5( \sum_1^m\vee (r';
T_i,T_{i+1}) \,\,\,\,\forall  r(\cdot) \in \wp(\Pi)
\label{DiffconvthreevarGomog14}
\end{equation}
\label{difconvgomogthm4}
\end{thm}
{\bf Доказательство. }{\bf Необходимость}. Пусть $ {\varphi
}(\cdot)$ есть ПРВ функция. Покажем, что тогда справедливо
неравенство (\ref{DiffconvthreevarGomog14}). Воспользуемся
неравенством, вытекающим из неравенства треугольника,
$$
\Vert \tau(t_i) / \Vert \tau(t_i) \Vert  - \tau (t_{i-1}) / \Vert
\tau(t_{i-1} \Vert \Vert \leq
$$
$$
\leq \Vert r'(t_i) / \sqrt{ 1+\varphi'^2_t (r(t_i))}  -
r'(t_{i-1}) / \sqrt{ 1+\varphi'^2_t (r(t_{i-1}))} \Vert +
$$
$$
\mid {\varphi }'_t(r(t_i)) /  \sqrt{ 1+{\varphi }'^2_t (r(t_i))}-
{\varphi }'_t(r(t_{i-1})) / \sqrt{ 1+{\varphi }'^2_t (r(t_{i-1}))}
\mid .
$$
Так как $1 \leq  \sqrt{1+{\varphi }'^2_t (r(t_i))} \leq
\sqrt{1+L^2}$ для всех $ t_i \in [0, T(r)]$ , то очевидно,
существует такое $c_5 >1,$ для которого верно неравенство

\begin{equation} \Vert r'(t_i) / \sqrt{1+{\varphi }'^2_t (r(t_i))}-
r'(t_{i-1}) / \sqrt{1+{\varphi }'^2_t (r(t_{ i-1}))} \Vert \leq
c_5 \Vert r'(t_i) - r'(t_{i-1}) \Vert.
\label{DiffconvthreevarGomog15a}
\end{equation}

Из свойств функции $\theta(x)= x / \sqrt{ 1+x^2}$ следует
неравенство
$$
\mid {\varphi }'_t (r(t_i)) / \sqrt{ 1+{\varphi
}'^2_t (r(t_i))} - {\varphi }'_t (r(t_{i-1})) / \sqrt{1+{\varphi
}'^2_t (r(t_{i-1}))} \mid \leq
$$
\begin{equation}
\leq \mid {\varphi }'_t (r(t_i)) - {\varphi }'_t (r(t_{i-1})) \mid
. \label{DiffconvthreevarGomog16a}
\end{equation}

Из (\ref{DiffconvthreevarGomog15a}) и
(\ref{DiffconvthreevarGomog16a}) имеем
$$
\sup_{ t_j \in [T_i, T_{i+1}] } \,\, \sum_j \, \Vert \tau(t_j) /
\Vert \tau(t_j) \Vert - \tau (t_{j-1}) / \Vert \tau(t_{j-1}) \Vert
\Vert \leq
$$
\begin{equation}
\leq c_5 (\vee (  r' ; T_i ,T_{i+1}) + \vee (\Phi' ; T_i,T_{i+1})
). \label{DiffconvthreevarGomog17a}
\end{equation}

Так как по условию $ \varphi(\cdot)-$  ПРВ функция, то согласно
неравенству теоремы \ref{difconvgomogthm2a} для $\vee (\Phi';
T_i,T_{i+1})$  с учетом (\ref{DiffconvthreevarGomog17a}) следует
неравенство (\ref{DiffconvthreevarGomog14}). Необходимость
доказана.

{\bf Достаточность}.  Пусть справедливо неравенство
(\ref{DiffconvthreevarGomog14}). Покажем, что ${\varphi }(\cdot)$
- ПРВ функция. Воспользуемся неравенством
$$
\Vert \tau(t_i) / \Vert \tau( t_i) \Vert  - \tau(t_{i-1}) /  \Vert
\tau (t_{i-1}) \Vert \Vert \geq \mid  {\varphi }'_t (r(t_i)) /
\sqrt{1+{\varphi }'_t (r(t_i))} -
$$
\begin{equation} - {\varphi }'_t (r(t_{i-1})) / \sqrt {1+{\varphi }'^2_t
(r(t_{i-1}))} \label{DiffconvthreevarGomog18a} \end{equation} Из
свойств функции $\theta(x) = x / \sqrt{1+x^2}$ и из $\Vert
{\varphi }'(z) \Vert \leq L$ для всех $z \in S_1^2(0),$ где
производная существует, следует существование константы $ c_6 >
0,$ для которой
$$
\mid {\varphi }'_t (r(t_i)) / \sqrt{1+{\varphi }'^2_t (r(t_i))} -
{\varphi }'_t (r(t_{i-1})) / \sqrt{1+{\varphi }'^2_t(r(t_{i-1}))}
\geq
$$
$$
\geq c_6 \mid {\varphi }'_t (r(t_i)) - {\varphi }'_t (r(t_{i-1}))
\mid,
$$
откуда с учетом (\ref{DiffconvthreevarGomog18a}) имеем
$$
\sup_{  t_j \in [T_i,T_{i+1}] }  \sum_j \Vert \tau(t_j) / \Vert
\tau( t_j) \Vert  - \tau(t_{j-1}) / \Vert \tau (t_{j-1}) \Vert
\Vert \geq
$$
$$
\geq c_6 \vee (\Phi';  T_i,T_{i+1})).
$$
Поэтому для суммы
$$
\sum_i \vee (\Phi';  T_i,T_{i+1})).
$$
выполняется неравенство теоремы  \ref{difconvgomogthm2a}.
Следовательно, согласно теореме  \ref{difconvgomogthm2a}
${\varphi }(\cdot)$ - ПРВ функция. Достаточность и теорема
доказаны.  $\Box$

%%%%%%%%%%%%%%%%%%%%%%%%%%%%%%%%%%%%%%%%%%%%%%%%%%%%%%%%%%%%%%%%%%%%%%%%%%%%%%%%%
%%%% Добавление: необходимые условия в R^3 п.о. 1 степени. Проекция кривых на пл П -выпукл множ

\vspace{0.5cm}

\subsection{Следствия их теорем \ref{difconvgomogthm2} и
\ref{difconvgomogthm2a}}

\vspace{0.5cm}

Определим класс $\hat \wp (\Pi)$ кривых $r(t)$, $t \in [0, T(r)]$,
на сфере $S^2_1(0)$, параметризованных естественным образом,
проекции которых на некоторую плоскость $ \Pi \in \mathbb{R}^3$,
проходящую через начало координат, ограничивают на $\Pi$ выпуклые
компактные множества. Плоскость $\Pi$ выбирается произвольным
образом. Здесь $T(r)-$ длина кривой $r(\cdot)$.

Отметим, что классы кривых $\hat \wp(\Pi) $ принадлежат классу
кривых $\wp(\Pi)$ на сфере $S^2_1(0)$, так как выпуклые множества
являются звездными множествами. Но  класса кривых $\hat \wp(\Pi) $
будет недостаточно, чтобы сформулировать необходимые и достаточные
условия представимости функции в виде разности выпуклых.

Отметим, что из неравенства
$$
    \mid r(t_2) - r(t_1) \mid \leq \mid t_2  - t_1 \mid
$$
следует, что вектор-функция $r(\cdot) \in \hat \wp(\Pi)$ п.в. на
$[0, T(r)] $ имеет касательную $r'(\cdot)$. Потребуем, чтобы
$r'(\cdot)$ образовывала с плоскостью $\Pi$ угол $\al \in [0, \pi
/ 4]$ для всех $t \in [0, T(r)]$.

Возьмем произвольную кривую $r(\cdot) \in \hat \wp(\Pi)$.
Определим функцию
$$
\Phi(t)=\varphi(r(t)) \,\,\,\, \forall t \in [0, T(r)].
$$
Покажем, что $\Phi(\cdot)-$ липшицевая на отрезке $[0, T(r)]$. Для
любых $t_1, t_2 \in [0, T(r)]$
$$
| \Phi(t_1) - \Phi(t_2) | = |\varphi(r(t_1)) - \varphi(r(t_2)|
\leq L \| r(t_1) - r(t_2) \|.
$$
Из очевидного неравенства $\| r(t_1) - r(t_2) \| \leq |t_1 - t_2
|$ имеем
$$
| \Phi(t_1) - \Phi(t_2) | \leq L |t_1 - t_2 |.
$$
Откуда следует, что $\Phi(\cdot)-$ п.в. дифференцируема на $[0,
T(r)]$. Также, как и ранее, определяем вариацию
$\vee(\Phi';0,{T(r)}) $ функции $\Phi'(\cdot)$ на отрезке $
[0,T(r)]$.

Из теорем \ref{difconvgomogthm2} и \ref{difconvgomogthm2a} можно
получить часто используемую теорему.

\begin{thm}
Для того чтобы п.о. первой степени липшицевая функция
$\varphi(\cdot):\mathbb{R}^3 \rar \mathbb{R} $ с константой
Липшица $L$ была ПРВ функцией, необходимо, чтобы для всех кривых
$r(\cdot) \in \hat \wp(\Pi)$, параметризованных естественным
образом с параметром $t \in [0,T], \, T=T(r),$ где $\Pi-$
плоскость, проходящая через начало координат, нашлась константа
$C_6=C_6(\varphi)>0$ такая, что \be \vee(\Phi';0,T) \leq C_6
\,\,\,\,\,\, \forall r(\cdot) \in \hat \wp(\Pi). \ee
\label{difconvgomogthm5}
\end{thm}

{\bf Доказательство}. {\bf Необходимость}. Доказательство
необходимости следует из доказательства необходимости теоремы
\ref{difconvgomogthm2a}. Действительно, поскольку $\hat \wp(\Pi)
\subset \wp(\Pi)$, то доказательство необходимости теоремы
\ref{difconvgomogthm5} будет повторять доказательство
необходимости теорем \ref{difconvgomogthm2a} и
\ref{difconvgomogthm2}.

\newpage
\vspace{1cm}

\section{ПРЕДСТАВИМОСТЬ ПОЛОЖИТЕЛЬНО ОДНОРОДНОЙ
ФУНКЦИИ $M -$ ОЙ СТЕПЕНИ ОТ $N$ ПЕРЕМЕННЫХ В ВИДЕ РАЗНОСТИ
ВЫПУКЛЫХ ФУНКЦИЙ}

\vspace{0.5cm}

В данной главе приведены необходимые и достаточные условия
представимости произвольной положительно однородной функции $m -$
ого порядка от произвольного количества переменных в виде разности
выпуклых функций. Дана также геометрическая интерпретация этих
условий. Приведен алгоритм такого представления, результатом
которого есть последовательность равномерно сходящихся на
единичном шаре выпуклых положительно однородных функций степени $m
\geq 1$.

\noindent

\vspace{0.5cm}

\subsection{Введение}

\vspace{0.5cm}

Обозначим через $\mathbb{R}^n - n$-мерное евклидово пространство
со скалярным произведением $(a,b)$ векторов $a$ и $b$. Пусть
задана липшицевая, дифференцируемая по направлениям функция
$f:\mathbb{R}^n \rar \mathbb{R}$. Обозначим через $ f'(x,q) $
производную по направлению $q \in \mathbb{R}^n$ функции $f(\cdot)$
в точке $x \in \mathbb{R}^n$.

В предыдущих параграфах были получены необходимые и достаточные
условия представимости положительно однородной функции (п.о.)
первой степени от двух и трех переменных в виде разности выпуклых.
Функции, представимые в виде разности выпуклых, ради сокращения
называют ПРВ (DC) функциями.

Такие функции  наряду с выпуклыми функциями играют важную роль в
оптимизации и теории управления. Функции, у которых производные по
направлениям $g \in \mathbb{R}^n$ в точке $x \in \mathbb{R}^n$,
рассматриваемые как функции от этого направления $g$,
представляются в виде разности двух выпуклых п.о. первой степени
функций от $g$, называются квазидифференцируемыми (КВД функциями)
\cite{demvas} в точке $x$. КВД функции являются расширением
множества выпуклых функций. Развиты методы оптимизации таких
функций.

Доказанные теоремы для п.о. $m-$ ой степени функций от двух
переменных и п.о. функций первой степени от трех переменных тесно
связаны с геометрией. Показано, что условия представимости в виде
разности выпуклых эквивалентны условию равномерной ограниченности
поворота кривых из определенного класса на графиках исследуемых
функций.

Ниже мы рассматриваем произвольные п.о. $m -$ ой степени функции
от произвольного количества переменных, $m -$ натуральное.
Приводятся необходимые и достаточные условия представимости таких
функций в виде разности выпуклых. Дается также геометрическая
интерпретация таких условий.

\vspace{0.5cm}

\subsection{Необходимые условия представимости п.о. функции
первой степени в виде разности выпуклых  в $n-$ мерном случае }

\vspace{0.5cm}

Вначале рассмотрим п.о. липшицевую функцию первой степени от $n$
переменных. Пусть $\varphi(\cdot): \mathbb{R}^n \rar \mathbb{R} $
такая функция с константой Липшица $L$. Докажем теоремы,
аналогичные доказанным для трехмерного случая.

В дальнейшем нам понадобится понятие координатной плоскости. Под
{\em координатной плоскостью } будем понимать плоскость,
образованную двумя произвольными осями координат.

Обозначим через $\wp $ множество кривых $r(\cdot)$ на сфере
$S_1^{n-1}(0) $, проекции которых на координатные плоскости $\Pi$
есть кривые, ограничивающие выпуклые компактные множества на
$\Pi$. Все кривые $r(\cdot) \in  \wp $ параметризуем естественным
образом. Параметризованную кривую обозначим через $r(t)$, $t \in
[0,T]$, где $T=T(r) -$ длина кривой $r(\cdot)$. Отметим, что
отсюда следует существование п.в. касательной $r'(\cdot)$ у кривой
$r(\cdot)$  на $[0, T]$.

Построим проекцию кривой $r(\cdot) \in  \wp$ на координатную
плоскость $\Pi$. Обозначим эту проекцию через $r_{\Pi}(\cdot)$.

Возьмем произвольную  $r(\cdot) \in \wp$. Введем функцию
$$
\Phi(t)=\varphi(r(t)) \,\,\, \forall t \in [0,T].
$$
Обозначим через $L$ константу Липшица функции $\varphi(\cdot)$.
Так как для любых $t_1, t_2 \in [0,T]$
$$
| \Phi(t_1) - \Phi(t_2) | = | \varphi(r(t_1)) - \varphi(r(t_2)) |
\leq L \| r(t_1) - r(t_2) \| \leq L |t_1 - t_2 |,
$$
то $\Phi(\cdot)-$ липшицевая с константой $L$ и, следовательно,
п.в. дифференцируемая на $[0,T]$.

Докажем следующую теорему.

\begin{thm}
Для того чтобы п.о. первой степени липшицевая функция
$\varphi(\cdot):\mathbb{R}^n \rar \mathbb{R} $ с константой
Липшица $L$ была ПРВ функцией, необходимо, чтобы для всех кривых
$r(\cdot) \in \wp$, параметризованных естественным образом с
параметром $t \in [0,T], \, T=T(r),$ нашлась константа $C =
C(\varphi)>0$ такая, что \be \vee(\Phi';0,T) \leq C \,\,\,\,\,\,
\forall r(\cdot) \in \wp. \label{diffconvGomFuncMDegNVer1} \ee
\label{diffconvGomFuncMDegNVerThm1}
\end{thm}
{\bf Доказательство. Необходимость.} Пусть функция
$\varphi(\cdot)$ представима в виде разности п.о. выпуклых функций
$\varphi_i(\cdot):\mathbb{R}^n \rar \mathbb{R}, \, i=1,2,$ с
константами Липшица $L_i, i=1,2,$ соответственно. Зафиксируем
произвольную кривую $r(\cdot) \in \wp$. Пусть $\Pi$ одна из
координатных плоскостей. Параметризуем $ r(\cdot)$ естественным
образом и через $[0,T]$ обозначим отрезок значений параметра $t$,
где $T=T(r)$. Для функций $\varphi(\cdot), \varphi_i(\cdot), \,
i=1,2,$ определим соответственно функции $\Phi(\cdot),
\Phi_i(\cdot), i=1,2,$ как это делали ранее.

Поскольку \cite{kolmogorovfomin}
$$
\vee(\Phi'; 0,T) \leq \vee(\Phi_1'; 0,T)+\vee(\Phi_2'; 0,T),
$$
то достаточно доказать неравенство \be \vee(\Phi_1';0,{T}) \leq
C_1. \label{diffconvGomFuncMDegNVer2} \ee Первоначально будем
считать, что $ \varphi_1(\cdot), r(\cdot) -$   дифференцируемые по
своим переменным функции. К общему случаю перейдем позднее.

Пусть
$$
 \varphi_{1} (r(t)) = \max _{v \in
 \partial \varphi_{1}(0)} (v, r(t))=(v(t),r(t)), \;\;
$$
где $\partial \varphi_{1}(0)$ - субдифференциал функции
$\varphi_{1}(\cdot)$ в нуле, $v(t)-$ крайняя точка множества
$\partial \varphi_{1}(0)$ с нормалью $r(t)$. Поскольку функция $
\varphi_1(\cdot)-$ дифференцируемая, то $v(t)$ единственный
крайний вектор с нормалью $r(t)$.

Как показано выше $\Phi_{1}(\cdot)-$ липшицевая, а поэтому п.в.
дифференцируемая на $[0,T]$. Для точек дифференцируемости  $t$
$$
\Phi_{1}'(t)=(v'(t),r(t))+(v(t),r'(t)).
$$
Но поскольку $r(t)$ есть нормальный вектор к множеству $\partial
\varphi_{1}(0)$, то $(v'(t),r(t))=0.$ В нашем случае производная
$v'(t)$ будет существовать, так как граница выпуклого множества
$\partial \varphi_{1}(0)$ для рассматриваемого случая имеет
касательную гиперплоскость в каждой крайней точке.

Поскольку кривая $r(\cdot)$ параметризована естественным образом,
то $\Vert r'(t) \Vert =1$ для любых $t \in [0,T]$. Нетрудно
проверить следующую цепочку неравенств
$$
\mid \Phi_{1}' (r(t_1)) - \Phi_{1}' (r(t_2)) \mid = \mid
(v(t_1),r'(t_1)) - (v(t_2),r'(t_2)) \mid =
$$
$$
\mid (v(t_1)- v(t_2),r'(t_1)) +(v(t_2),r'(t_1))- (v(t_2),r'(t_2))
\mid \leq
$$
$$
\leq  \Vert v(t_1) - v(t_2) \Vert \; \Vert r'(t_1) \Vert + \Vert
r'(t_1) - r'(t_2) \Vert \; \Vert v(t_2) \Vert \leq
$$
\begin{equation}
\leq \Vert v(t_1) - v(t_2) \Vert  + L_{1} \Vert r'(t_1) - r'(t_2)
\Vert. \label{diffconvGomFuncMDegNVer3}
\end{equation}
Из (\ref{diffconvGomFuncMDegNVer3}) следует, что
$$
   \vee (\Phi_{1}'; 0,T) \leq \mbox{длина кривой}\,\,\,
   v(t)\,\,\,
   \mbox{для } \,\,\, t \in [0, T]
   + L_{1} \cdot \vee (r'; 0,T)=
$$
\be
   = l_{1} + L_{1} \cdot \vee (r'; 0,T),
\label{diffconvGomFuncMDegNVer4} \ee где $l_1-$ длина кривой
$v(t), t \in [0, T]$,  вариация $ \vee (r'; 0,T) $ по определению
есть
$$
\vee(r';0,{T}) = \sup_{{t_i}, t_i \in [0,T]} \sum_i \| r'(t_i) -
r'(t_{i+1}) \|.
$$
Покажем, что для всех $ r(\cdot) \in  \wp $ верно неравенство \be
  \vee(r'; 0, T) < C_2,
\label{diffconvGomFuncMDegNVer5} \ee где $C_2 -$ некоторая
константа. Очевидно
$$
\| r'(t_i) - r'(t_{i+1}) \| \leq  \sum_{j} \, | r_j'(t_i) -
r_j'(t_{i+1}) | ,
$$
где $r_j(\cdot), i \in 1:n, -$ координаты вектор-функции
$r(\cdot).$ Так как кривая $r(\cdot)$ параметризована естественным
образом, то $\| r'(\cdot) \| =1, $ а следовательно, $| r'_j
(\cdot) | \leq 1 $ для всех $j \in 1:n$. Поэтому пары координат
вектор-функции $r'(\cdot)$ образуют вектор на координатной
плоскости $\Pi$, являющийся касательной к проекции кривой $v(t), t
\in [0, T],$ на плоскость $\Pi$. Проекция $v_{\Pi}(\cdot) $ кривой
$v(t), t \in [0, T],$ на координатную плоскость $\Pi$  есть
кривая, ограничивающая выпуклое компактное множество, поскольку
касательными к проекции этой кривой являются проекции касательной
к кривой $r(\cdot)$. Проекция $r_{\Pi}(\cdot) $ самой кривой
$r(\cdot)$ ограничивает по предположению выпуклое компактное
множество. По свойству касательных к кривым, ограничивающих
выпуклые компактные множества, следует, что проекция кривой $v(t),
t \in [0, T],$ на плоскость $\Pi$ ограничивает выпуклое компактное
множество, периметр которого равномерно ограничен для всех кривых
$r(\cdot) \in  \wp $, поскольку ограничены диаметры выпуклых
множеств, являющихся проекциями субдифференциала $\p \varphi_1(0)
$ на координатные плоскости. Отсюда следует ограниченность
вариации проекции $r'(\cdot)$  на произвольную координатную
плоскость, а следовательно и вариации кривой $r'(\cdot)$. А
поэтому верно неравенство (\ref{diffconvGomFuncMDegNVer5}) для
всех $ r(\cdot) \in \wp $.

Докажем теперь, что длины $l_1$ равномерно ограничены для всех $
r(\cdot) \in \wp. $ Длина  $l_1$ для любой $ r(\cdot)$ не
превосходит суммы длин проекций кривой $v(t), t \in [0, T],$ на
координатные плоскости. Проекция кривой $v(t), t \in [0, T],$ на
координатную плоскость $\Pi$  есть кривая, ограничивающая выпуклое
компактное множество, о чем говорилось ранее. А поэтому длина
проекции кривой $v(t), t \in [0, T],$ на произвольную координатную
плоскость $\Pi$  не превосходит периметра проекции множества
$\partial \varphi_{1}(0)$ на плоскость $\Pi$. Периметр проекции
$\partial \varphi_{1}(0)$ на любую плоскость $\Pi$, очевидно,
ограничен. Отсюда следует ограниченность длин  $l_1$ равномерно
для всех $ r(\cdot) \in \wp. $

%%%%%%%%%%%%%%%%%%%%%%%%%%%%%%%%%%%%%%%%%%%%%%%%%%%%%%%%%%%%%

Из сказанного и неравенств (\ref{diffconvGomFuncMDegNVer4}) и
(\ref{diffconvGomFuncMDegNVer5}) следует необходимость для случая
гладких функций $\varphi(\cdot), r(\cdot)$. Перейдем к общему
случаю.

Мы всегда можем приблизить выпуклую функцию $\varphi_1(\cdot) $ и
кривую $ r(\cdot) \in \wp$ гладкой выпуклой функцией
$\varphi_{1k}(\cdot) $ и гладкой кривой $ r_k(\cdot) \in \wp$ так,
чтобы длина кривой $l_{1k}$, построенная для $\varphi_{1k}(\cdot)
$, и вариация $\vee(r'_{k};0,T) $ как угодно мало отличались от
длины кривой $l_{1}$ и вариации $\vee(r';0,T) $.

Отсюда следует справедливость неравенства
(\ref{diffconvGomFuncMDegNVer1}) для произвольной выпуклой функции
$\varphi_1(\cdot) $ и любой кривой $ r(\cdot) \in \wp$.
Необходимость доказана.

\vspace{0.5cm}

\subsection{Необходимые условия представимости положительно однородной функции $m-$
степени в виде разности выпуклых  в $n-$ мерном случае }

\vspace{0.5cm}

Рассмотрим теперь п.о. $m-$ ой степени ($m-$ натуральное число,
большее 1) липшицевую функцию $\varphi(\cdot): \mathbb{R}^n \rar
\mathbb{R}$. Нас будут интересовать условия ее представимости в
виде разности выпуклых.

Для произвольной кривой $r(\cdot) \in \wp,$ параметризованной
естественным образом $t \in [0, T(r)]$, определим функцию
$$
  \Phi(t)=\varphi(r(t)) \;\;\;\; \forall t \in [0, T(r)].
$$

\begin{thm}
Для того чтобы п.о. степени $m$ липшицевая  функция
$\varphi(\cdot)$ была ПРВ функцией, необходимо, чтобы для всех
$r(\cdot) \in \wp $ нашлась константа $c=c(\varphi)>0$ такая, что
\be \vee(\Phi'; 0,T) < c \,\,\,\,\,\, \forall r(\cdot) \in \wp.
\label{diffconvGomFuncMDegNVer6} \ee
\label{diffconvGomFuncMDegNVerThm2}
\end{thm}
{\bf Доказательство. Необходимость.} Пусть $\varphi(\cdot)$ есть
ПРВ функция, т.е.
$$
\varphi(q) = \varphi_1(q) - \varphi_2(q) \,\,\,\, \forall q \in
\mathbb{R}^n,
$$
где $\varphi_i(\cdot), \, i=1,2, - $ выпуклые п.о. степени $m$
функции. Построим п.о. степени 1 функцию $\psi(\cdot):\mathbb{R}^n
\rar \mathbb{R}$, принимающую те же значения на $S_1^{n-1}(0)$,
что и функция $\varphi(\cdot)$, т.е. \be \psi(q) = \| q \| \,\,
\varphi(\frac{q}{\| q \|}) \,\,\,\, \forall q \in \mathbb{R}^n.
\label{diffconvGomFuncMDegNVer7} \ee Так же, как и для двумерного
случая \cite{proudconvex2} для выпуклых функций $\varphi_i(\cdot),
\, i=1,2, $ строим  п.о. первой степени функции $\psi_i(\cdot), \,
i=1,2,$ которые мы назовем {\em соответствующими} для функций
$\varphi_i(\cdot), \, i=1,2.$ По построению функции
$\psi_i(\cdot), \, i=1,2,$ принимают те же значения на единичной
сфере с центром в нуле $S_1^{n-1}(0)$, что и функции
$\varphi_i(\cdot), \, i=1,2. $ В \cite{proudconvex2} для
двумерного пространства было доказано, что $\psi_i(\cdot), \,
i=1,2,-$ выпуклые функции.

Для доказательства выпуклости функций $\psi_i(\cdot), i=1,2,$ в
$n-$ мерном пространстве воспользуемся отличительным свойством
выпуклых функций.

\begin{thm}\cite{pshenichnyi} Для того чтобы функция $z \rar \theta(z):
\mathbb{R}^n \rar \mathbb{R}$ была выпуклой, необходимо и
достаточно, чтобы она была дифференцируема по направлениям и для
любой точки $z \in \mathbb{R}^n $ и любого направления $p \in
\mathbb{R}^n $ функция $\alpha \rar h( \alpha ): \mathbb{R}^+ \rar
\mathbb{R}:$
$$
h(\alpha) = \frac{\p \theta (z + \alpha p)}{\p p }
$$
была неубывающей по $\alpha >0.$
\label{diffconvGomFuncMDegNVerThm3}
\end{thm}

Из п. о.  следует, что для доказательства  достаточно рассмотреть
$q \in S_1^{n-1}(0) $ и вектор $p$, перпендикулярный $q$.
Зафиксируем $q$ и $p$. Проведем плоскость $\Pi$, содержащую
векторы $q$ и $p$ и проходящую через начало координат. Теперь мы
приходим к двумерному случаю, где в качестве плоскости XOY служит
плоскость $\Pi$, а в качестве окружности $S_1^1(0)- $ пересечение
$S_1^{n-1}(0) \cap \Pi$. Для двумерного случая уже было доказано
\cite{proudconvex2}, что функции $\alpha \rar h_i( \alpha ):
\mathbb{R}^+ \rar \mathbb{R}:$
$$
h_i(\alpha) = \frac{\p \psi_i (q + \alpha p)}{\p p }, \,\,\,
i=1,2,
$$
являются неубывающими по $\alpha >0.$ Поэтому функции
$\psi_i(\cdot), \, i=1,2,-$ выпуклые. Тогда
$$
\psi(q)=\psi_1(q) - \psi_2(q) \,\,\,\,\, \forall q \in
\mathbb{R}^n,
$$
т.е. $\psi(\cdot) -  $ ПРВ функция.

По доказанной теореме \ref{diffconvGomFuncMDegNVerThm1} существует
константа $c>0$ такая, что \be \vee(\Psi';0,T) < c,
\label{diffconvGomFuncMDegNVer8} \ee где $\Psi(t)=\psi(r(t))$ для
$t \in [0,T]$. Поскольку
$$
\vee(\Psi';0,T)=\vee(\Phi';0,T),
$$
то из (\ref{diffconvGomFuncMDegNVer6}) имеем
$$
\vee(\Phi';0,T) < c.
$$
Необходимость доказана.

\vspace{0.5cm}

\subsection{Геометрическая интерпретация теоремы, дающей необходимые
условия представимости положительно однородной функции $m-$ ой
степени  от $n$ переменных в виде разности выпуклых }

\vspace{0.5cm}

Перефразируем теорему \ref{diffconvGomFuncMDegNVerThm2}, придав им
более геометрический характер. Введем понятие поворота кривой
$r(\cdot)$ на графике $\Gamma_{\varphi } = \{ (x,z) \in
\mathbb{R}^{n+1} \mid z = \varphi(x), \, x \in \mathbb{R}^n \}.$

Рассмотрим на $\Gamma_{\varphi }$ кривую $ R(t)=(r(t),{\varphi
}(r(t))), t \in [0, T(r)],$ где $r(\cdot) \in \wp).$ Так как
функция ${\varphi }(\cdot)$ есть липшицевая, то п.в. на $[0,T(r)]$
существует производная $R'(\cdot),$ которую обозначим через
$\tau(\cdot)=R'(\cdot),$  а множество точек, где она существует,
$-$ через $N_{\varphi }$.

\begin{defi}  Поворотом кривой $R(\cdot)$ на многообразии $\Gamma_f$
назовем величину
$$
sup_{ \{t_i\} \subset N_{\varphi }} \,\, \sum_i \Vert \tau(t_i)/
\Vert \tau (t_i) \Vert -  \tau(t_{i-1})/ \Vert \tau (t_{i-1})
\Vert \Vert = O_{\varphi }.
$$
\end{defi}

Таким образом, поворот $O_{\varphi }$ кривой $R(\cdot)$ есть
верхняя грань суммы углов между касательными $\tau(t)$  для $t \in
[0,T(r)]$.

\begin{thm} Для того, чтобы произвольная липшицевая п.о. степени
$m$ функция $x \rightarrow {\varphi }(x) :\mathbb{R}^n \rightarrow
\mathbb{R}$ была ПРВ функцией, необходимо, чтобы для всех
$r(\cdot) \in \wp $ существовала константа $c_4({\varphi })>0$
такая, что поворот кривой $R(\cdot)$ на $\Gamma_{\varphi }$
ограничен сверху константой $c_4({\varphi })>0,$ т.е.
\begin{equation} O_{\varphi } \leq c_4({\varphi}) \,\,\,\,\,
\forall r(\cdot) \in \wp. \label{diffconvGomFuncMDegNVer10}
\end{equation}
\label{diffconvGomFuncMDegNVerThm4}
\end{thm}
{\bf Доказательство. }{\bf Необходимость}. Пусть $ {\varphi
}(\cdot)$ есть ПРВ функция. Покажем, что тогда справедливо
неравенство (\ref{diffconvGomFuncMDegNVer10}). Воспользуемся
неравенством, вытекающим из неравенства треугольника,
$$
\Vert \tau(t_i) / \Vert \tau(t_i) \Vert  - \tau (t_{i-1}) / \Vert
\tau(t_{i-1} \Vert \Vert \leq
$$
$$
\leq \Vert r'(t_i) / \sqrt{ 1+\varphi'^2_t (r(t_i))}  -
r'(t_{i-1}) / \sqrt{ 1+\varphi'^2_t (r(t_{i-1}))} \Vert +
$$
$$
\mid {\varphi }'_t(r(t_i)) /  \sqrt{ 1+{\varphi }'^2_t (r(t_i))}-
{\varphi }'_t(r(t_{i-1})) / \sqrt{ 1+{\varphi }'^2_t (r(t_{i-1}))}
\mid .
$$
Так как $1 \leq  \sqrt{1+{\varphi }'^2_t (r(t_i))} \leq
\sqrt{1+L^2}$ для всех $ t_i \in [0, T(r)]$ , то очевидно,
существует такое $c_5 >1,$ для которого верно неравенство

\begin{equation} \Vert r'(t_i) / \sqrt{1+{\varphi }'^2_t (r(t_i))}-
r'(t_{i-1}) / \sqrt{1+{\varphi }'^2_t (r(t_{ i-1}))} \Vert \leq
c_5 \Vert r'(t_i) - r'(t_{i-1}) \Vert.
\label{diffconvGomFuncMDegNVer11}
\end{equation}

Из свойств функции $\theta(x)= x / \sqrt{ 1+x^2}$ следует
неравенство
$$
\mid {\varphi }'_t (r(t_i)) / \sqrt{ 1+{\varphi }'^2_t (r(t_i))} -
{\varphi }'_t (r(t_{i-1})) / \sqrt{1+{\varphi }'^2_t (r(t_{i-1}))}
\mid \leq
$$
\begin{equation}
\mid {\varphi }'_t (r(t_i)) - {\varphi }'_t (r(t_{i-1})) \mid .
\label{diffconvGomFuncMDegNVer12}
\end{equation}

Из (\ref{diffconvGomFuncMDegNVer11}) и
(\ref{diffconvGomFuncMDegNVer12}) имеем
$$
 \sup_{\{t_i \} \in N_{\varphi } } \,\, \sum_i \,
\Vert \tau(t_i) / \Vert \tau(t_i) \Vert - \tau (t_{i-1}) / \Vert
\tau(t_{i-1}) \Vert \Vert \leq
$$
\begin{equation}
\leq c_5 (\vee (  r' ; 0 ,T(r)) + \vee (\Phi' ; 0,T(r)) ).
\label{diffconvGomFuncMDegNVer13}
\end{equation}

Так как по условию $ \varphi(\cdot)-$  ПРВ функция, то согласно
теореме \ref{diffconvGomFuncMDegNVerThm2}
$$
\vee (\Phi'; 0,T(r)) \leq  c_3(\varphi),
$$
откуда с учетом (\ref{diffconvGomFuncMDegNVer13}) и ограниченности
вариации $\vee (  r' ; 0 ,T(r)) $ следует неравенство
(\ref{diffconvGomFuncMDegNVer10}). Необходимость доказана.

%%%%%%%%%%%%%%%%%%%%%%%%%%%%%%%%%%%%%%%%%%%%%%%%%%%%%%%%%%%%%%%%%%

\vspace{1cm}

\subsection{Определение более широкого множества кривых на $S^{n-1}_1(0)$}

\vspace{0.5cm}

Наша цель $-$ определить такое множество кривых на $S^{n-1}_1(0)$,
с помощью которого можно было бы сформулировать необходимые и
достаточные условия представимости положительно-однородной функции
$\varphi(\cdot)$ степени $m$ от $n$ аргументов в виде разности
выпуклых.

Дадим определение звездной области.

\begin{defi}.  Область ${D}$ называется звездной, если существует
точка $O \in \mbox{int } {D} $, для которой отрезок, соединяющий
точку $O$ с любой точкой области ${D}$, целиком принадлежит
области $D$.
\end{defi}

Введем множество $\wp(\Pi) $  кривых $r(\cdot) $  на единичной
сфере $S^{n-1}_1(0)$, параметризованных естественным образом с
параметром $t \in [0, T]$, проекция которых на плоскость $\Pi$
есть кривая, ограничивающая звездное множество на $\Pi$. Заметим,
$\wp \subset \wp(\Pi)$.

Кривая $r(\cdot)$ п.в. на отрезке $[0,T]$ имеет касательную
$r'(\cdot)$.

Для $n=3$ на сфере $S^{2}_1(0)$  в предыдущем параграфе был введен
класс кривых $ \Re $, который, как нетрудно видеть, также
принадлежит множеству $\wp(\Pi) $.

\vspace{0.5cm}

\subsection{Необходимые и достаточные условия представимости
положительно однородной функции первой степени в виде разности
выпуклых  в $n-$ мерном случае }

\vspace{0.5cm}

Сперва рассмотрим п.о. функцию первой степени от $n$ аргументов
$\psi(\cdot):\mathbb{R}^n \rar \mathbb{R} $. Возьмем произвольную
плоскость $\Pi$, $0 \in \Pi$, и произвольную $r(\cdot) \in
\wp(\Pi)$. Введем функцию
$$
\Phi(t)=\varphi(r(t)) \,\,\, \forall t \in [0,T].
$$
Функции $\Phi(\cdot), r(\cdot) $ п.в. дифференцируемые на $[0,
T]$, так как $\Phi(\cdot)- $ липшицевая функция, а кривая
$r(\cdot) $ параметризована естественным образом.

Для $r(\cdot) \in \wp(\Pi) $ выделим те ее участки с параметрами
$[T_i, T_{i+1}] \subset [0,T(r)]$, где производная $r'(t), t \in
[0, T(r)],$ составляет с плоскостью $\Pi$ угол $\delta $: $0 <
\delta \leq \frac{\pi }{4} $. Обозначим объединение таких участков
через $\Lambda(r)$.

\begin{thm}
Для того чтобы п.о. первой степени липшицевая функция
$\Phi(\cdot):\mathbb{R}^n \rar \mathbb{R} $ с константой Липшица
$L$ была ПРВ функцией, необходимо и достаточно, чтобы для любой
плоскости $\Pi$ и любой кривой $r(\cdot) \in  \wp(\Pi)$,
параметризованной естественным образом с параметром $t \in [0,T],
\, T=T(r),$ и любого подмножества параметров $[T_i, T_{i+1} ]
\subset \Lambda(r) $ нашлись константы $C_1 = C(\varphi)>0$, $C_2
= C_2(\varphi)>0$ такие, что \be
 \sum_{i=1}^m \vee (\Phi'; T_i,T_{i+1})<  C_1 + C_2 \sum_{i=1}^m\vee (r';
 T_i,T_{i+1}),
\label{diffconvGomFuncMDegNVer14} \ee
где $\Phi(t)=\varphi(r(t))$
и производные берутся там, где они существуют.
\label{diffconvGomFuncMDegNVerThm5}
\end{thm}
\begin{rem}
Необходимость рассмотрения системы подмножеств $\Lambda(r) $
объясняется тем, что вариация $\vee (r'; 0,T(r))$ может быть
неограниченной для всей кривой $r(\cdot)$.
\end{rem}

{\bf Необходимость.} Неравенство (\ref{diffconvGomFuncMDegNVer14}
) достаточно доказать для выпуклой п.о. первой степени непрерывно
дифференцируемой на $\mathbb{R}^n \backslash \{ 0 \}  $ функции
$\varphi_1(\cdot) $. Пусть $ r(\cdot) \in \wp$ и
$\Phi_1(t)=\varphi_1(r(t)), t \in [0, T],$ и $v(t) \in \p
\varphi_1(0) -$ крайние векторы с нормалями $r(t)$.

Неравенство \be \vee(\Phi'_1;0,T)  \leq l_{1} + L_{1} \cdot \vee
(r'; 0,T) \label{diffconvGomFuncMDegNVer15}  \ee выводится
аналогично (\ref{diffconvGomFuncMDegNVer4}), где $l_1$ длина
кривой $v(t)$ для $t \in \Lambda(r) $.

Требуется доказать равномерную ограниченность $l_1$ для всех
$r(\cdot) \in \wp(\Pi)$. Ясно, что $l_1$ не превосходит суммы длин
проекций кривой $v(t), t \in [0,T],$ на координатные плоскости.

Возьмем произвольную координатную плоскость $\Pi$. Спроектируем
кривые $v(t), t \in [0,T],$ и $r(t), t \in [0,T],$  на плоскость
$\Pi$. Проекции этих кривых на плоскость $\Pi $ обозначим
соответственно через $v_{\Pi}(t), t \in [0,T],$ и $r_{\Pi}(t), t
\in [0,T].$

Утверждается, что кривая $v_{\Pi}(t), t \in [0,T],$ ограничивает
на плоскости $\Pi$, выпуклое компактное множество.

Разобьем кривую $r(t), t \in [0,T], $ точками $r(T_i) \in
r(\cdot), i \in 1:s,$ на подмножества и включим эти точки в число
тех, которые участвуют в разбиении шара $B^{n}_1(0)$ на конусы
$K_j, j \in 1:m.$ В каждом конусе $K_j$ строим линейную функцию,
значения которой совпадают со значениями функции $\varphi(\cdot)$
в точках $r(T_i)$  на сфере $S^{n-1}_1(0)$ и в нуле, определяющими
конус $K_j$.

Как и ранее, под двугранным углом будем понимать функцию, равную
максимуму или минимуму двух линейных функций. Под выпуклым
(вогнутым) двугранным углом будем понимать выпуклую (вогнутую)
функцию, равную максимуму (минимуму) линейных функций.

По разбиению шара $B^{n}_1(0)$ на конусы строим многогранную
функцию $\varphi_{1,m}(\cdot)$. Так как по предположению
$\varphi_1(\cdot)-$ выпуклая, то все двугранные углы, часть
графика которых принадлежит графику функции
$\varphi_{1,m}(\cdot)$, являются выпуклыми. Проекции векторов
$r(T_i), i \in 1:s,$ на плоскость $\Pi$ образуют секторы на этой
плоскости. По кривой $r_{\Pi}(t)$ и проекциям точек $r(T_i) \in
r(\cdot), i \in 1:s,$ на плоскость $\Pi$, а также значениям
$\varphi_1(r(T_i)), i \in 1:s,$ строим п.о. многогранную функцию
$\varphi_{1,{\Pi,s}}(\cdot):\mathbb{R}^2 \rar \mathbb{R}$ линейную
в каждом секторе.

Все двугранные углы функции $\varphi_{1,m}(\cdot)-$ выпуклые.
Градиенты линейных функций, определенных в каждом секторе на
плоскости $\Pi$, есть проекции соответствующих градиентов линейных
функции, определенных в конусах, содержащих точки $r(T_i) \in
r(\cdot), i \in 1:s.$ Проекция касательной к кривой $r(\cdot)$ на
плоскость $\Pi$ есть касательная к кривой $r_{\Pi}(t)$ в
соответствующей точке. Поэтому очередность следования проекций
градиентов линейных функций, определенных в каждом конусе с
векторами $r(T_i), i \in 1:s,$   при построении функции
$\varphi_{1,m}(\cdot)$, и в каждом секторе на плоскости $\Pi$ при
построении функции $\varphi_{1,{\Pi,s}}(\cdot)$ при увеличении $t
\in [0, T] $ совпадают. Так как все двугранные углы функции
$\varphi_{1,m}(\cdot)-$ выпуклые, то отсюда следует, что все
двугранные углы функции $\varphi_{1,{\Pi,s}}(\cdot)$ также
выпуклые. Следовательно, $\varphi_{1,{\Pi,s}}(\cdot): \mathbb{R}^2
\rar \mathbb{R} - $ выпуклая п.о. функция, градиентами
$v_{\Pi,i}(t)$ линейных функций, определенных в каждом секторе,
являются проекции на плоскость $\Pi$ градиентов $v_{i}(t)$
линейных функций, определенных в соответствующих конусах разбиения
шара $B_1^n(0)$.

При неограниченном увеличении $m$  многогранные функции
$\varphi_{1,m}(\cdot)$  равномерно на $B^{n}_1(0)$ стремятся к
$\varphi_{1}(\cdot)$. Ломаная  $v_{\Pi,i}(\cdot)$ в пределе
стремится к кривой $v_{\Pi}(\cdot)$, ограничивающей
субдифференциал в нуле выпуклой п.о. функции
$\varphi_{1,{\Pi}}(\cdot)$. Поэтому верно утверждение: кривая
$v_{\Pi}(t), t \in [0,T],$ ограничивает на плоскости $\Pi$,
выпуклое компактное множество.

Нетрудно видеть, что длина кривой $v_{\Pi}(t), t \in [0,T],$ не
превосходит длины кривой $\hat v_{\Pi}(t), t \in [0,T],$ где $\hat
v (t) \in \p \varphi_1(0) -$  вектор с нормалью $r_{\Pi}(t) \in
\Pi$, если его рассматривать как $n-$ мерный вектор.  Поскольку
проекция множества  $\p \varphi_1(0) $ на плоскость $\Pi$ есть
выпуклое компактное множество с равномерно ограниченным периметром
для любой плоскости $\Pi$ и кривой $r(\cdot) \in \wp(\Pi)$.

Нам нужно доказать ограниченность длины кривой $v(t)$ для $t \in
[T_i, T_{i+1} ] \subset \Lambda(r) $, $i \in 1:m$. Для таких $t$
угол наклона касательной $v'(t)$ с плоскостью $\Pi$ составляет
угол $\delta $: $0 < \delta \leq \frac{\pi }{4} $. Поэтому длина
кривой $v(t)$ для таких $t$ не превосходит  $\sqrt{2} \,
l_{\Pi,1}$, где $l_{\Pi,1}-$ длина кривой $v_{\Pi}(t)$ на
плоскости $\Pi$ для тех же $t$. Равномерную ограниченность длин
кривых $v_{\Pi}(t)$ для любой плоскости $\Pi$ и кривой $r(\cdot)
\in \wp(\Pi)$ мы уже доказали. Необходимость доказана.

{\bf Достаточность.} Докажем, что если выполняется неравенство
(\ref{diffconvGomFuncMDegNVer14}) , то функция $\varphi(\cdot) $
представима в виде разности выпуклых. Разбиваем шар $B_1^n(0)$ на
конусы $K_i(q_1, q_2, \dots ,q_n), i \in 1:k$, являющиеся
конической оболочкой линейно независимых векторов $q_1, q_2, \dots
q_n$. В каждом конусе $K_i$ строим линейную функцию, значения
которой на векторах $q_1, q_2, \dots , q_n$ и в нуле  совпадают со
значениями функции $\varphi(\cdot)$. Функцию, равную построенным
линейным в каждом конусе $K_i$, обозначим по числу конусов через
$\varphi_k(\cdot)$.

Функцию, равную максимуму или минимуму  линейных функций,  назовем
{\em двугранным углом}. Если эта функция выпуклая (вогнутая), то
двугранный угол называется выпуклым (вогнутым).

Мы рассмотрим двугранные углы, построенные по векторам из соседних
конусов $K_i$ и $K_{i+1}$, имеющих общую часть гиперплоскости
$\sigma_i$. Градиенты этих линейных функций обозначим через $a_i$
и $a_{i+1}$. Поскольку $\sigma_i=K_i \cap K_{i+1}$, то проекции
векторов $a_i$ и $a_{i+1}$ на $\sigma_i$ равны. Следовательно,
вектор $a_i - a_{i+1} $ перпендикулярен $\sigma_i$. Для всех $i
\in 1:k$ соединим векторы $a_i$ и $a_{i+1}$ отрезками. Получим
многогранник $G_k$, стороны которого будут пересекаться по
внутренним точкам, если среди построенных двугранных углов есть
выпуклые и вогнутые. Если $\varphi(\cdot)-$ выпуклая функция, то
все двугранные углы, построенные по разбиению шара $B_1^n(0)$ на
конусы, - выпуклые, и многогранник $G_k -$ выпуклый.

В случае невыпуклости и невогнутости функции $\varphi(\cdot)$
стороны $b_i$ и $b_{i+1}$ многогранника $G_k$ для большого $k$
будут пересекаться по внутренним точкам. Мы также, как в
\cite{proudconvex1}, растягиваем многогранник $G_k$ по всем
направлениям, чтобы получить выпуклый многогранник $G'_k$.

Предположим, что $\varphi(\cdot)-$ не ПРВ функция. При
доказательстве достаточности теоремы
\ref{diffconvGomFuncMDegNVerThm2} из \cite{proudconvex1} было
показано, что поиск кривой $r(\cdot)$ ограничивается поиском
направления $g$ и вектора $x \in S_1^{n-1}(0)$, $g \perp x$.

Функция $ \varphi_k(\cdot)$, построенная по произвольному
разбиению шара $B_1^{n}(0)$ на конусы, является ПРВ функцией как
функция, график которой имеет конечное числом граней. Пусть
растяжение многогранника $G_k$, соответствующего этому разбиению,
вдоль направления $g$ неограниченно возрастает при $k \rar
\infty$. В противном случае  в пределе при $k \rar \infty$
многогранники $G'_k $ имели бы предел, и  мы получили бы, что
$\varphi(\cdot)-$ ПРВ функция.

Поэтому есть точка  $x \in S_1^{n-1}(0)$, в окрестности которой
вдоль направления $g \perp x$ растяжение многогранника $G_k$ при
преобразовании его в выпуклый, диаметр $G_k$ неограниченно
увеличивается.

В качестве вектора $g$ можно взять тот, вдоль которого растяжение
многогранника $G_k$ согласно описанному выше алгоритму для
преобразования его в выпуклый многогранник $G'_k$ происходит
наибольшим образом при $k \rar \infty$. В качестве плоскости $\Pi$
берем  плоскость, проходящую через начало координат и параллельную
вектору $g$. Все отрезки $b_s$, являющиеся субдифференциалами
вогнутых двугранных углов, на которые мы растягиваем многогранник
$G_k$, чтобы получить выпуклый многогранник, будут образовывать с
плоскостью $\Pi$ угол $0 < \delta < \frac{\pi}{4} $.

Поскольку сказанное выше справедливо для произвольно малого
конуса, содержащего во внутренности вектор $x$, то кривую
$r(\cdot) \in \wp(\Pi)$ можно выбирать с началом в точке $x$.

Кривая $r(\cdot) \in  \wp(\Pi)$ будет обходить, иначе говоря,
содержать отрезки $b_i$, которым соответствует отрезок параметра
$t \in [T_i, T_{i+1}] $, и иметь конечную сумму вариаций
производных $\sum_i \vee(r'; T_i, T_{i+1}) < \infty $.

Вариация функции  $\vee (\Phi'; T_i,T_{i+1})$ оценивает сверху
длину отрезка $b_i$. Поэтому если сумма длин отрезков $b_i$
бесконечна, то с учетом конечности вариации  $\vee(r'; 0, T) <
\infty $ получаем, что бесконечна сумма  $\sum_i \vee (\Phi';
T_i,T_{i+1})$. Поэтому неравенство
(\ref{diffconvGomFuncMDegNVer14}) теоремы
\ref{diffconvGomFuncMDegNVerThm5} не может выполняться ни при
каких $C_1, C_2$. Пришли к противоречию. Достаточность и теорема
доказаны. $\Box$

\vspace{0.5cm}

\subsection{Геометрическая интерпретация необходимых и достаточных
условий представимости в виде разности выпуклых для положительно
однородной функции первой степени от $n$ переменных}

\vspace{0.5cm}

Перепишем теорему \ref{diffconvGomFuncMDegNVerThm5} в
геометрической форме.

Введем понятие поворота кривой $r(\cdot)$ на графике
$\Gamma_{\varphi } = \{ (x,z) \in \mathbb{R}^{n+1} \mid z =
\varphi(x), \, x \in \mathbb{R}^n \}$, как это мы делали ранее.

Рассмотрим на $\Gamma_{\varphi }$ кривую $ R(t)=(r(t),{\varphi
}(r(t))), t \in [0, T(r)],$ где $r(\cdot) \in \wp(\Pi)$ для
некоторой плоскости $\Pi$. Так как функция ${\varphi }(\cdot)$
есть липшицевая, то п.в. на $[0,T(r)]$ существует производная
$R'(\cdot),$ которую обозначим через $\tau(\cdot)=R'(\cdot),$  а
множество точек, где она существует, $-$ через $N_{\varphi }$.

\begin{defi}  Поворотом кривой $R(\cdot)$ на многообразии $\Gamma_f$
назовем величину
$$
sup_{ \{t_i\} \subset N_{\varphi }} \,\, \sum_i \Vert \tau(t_i)/
\Vert \tau (t_i) \Vert -  \tau(t_{i-1})/ \Vert \tau (t_{i-1})
\Vert \Vert = O_{\varphi }.
$$
\end{defi}

\begin{thm} Для того, чтобы произвольная липшицевая п.о. степени 1
функция $\varphi(x) :\mathbb{R}^n \rightarrow \mathbb{R}$ была ПРВ
функцией, необходимо и достаточно, чтобы для любой плоскости
$\Pi$, $0 \in \Pi$, и всех $r(\cdot) \in \wp(\Pi)$
параметризованной естественным образом: $t \in [0, T(r)]$, а также
любого подмножества параметров $[T_i, T_{i+1} ] \subset \Lambda(r)
$ существовали константы $c_1(\varphi )>0, c_2(\varphi)>0$ такие,
что сумма поворотов $O_{\varphi, i}$ кривой $R(\cdot)$ для $t \in
[T_i, T_{i+1}] $ для  всех $i$ была ограничена сверху неравенством
\begin{equation}
\sum_1^m O_{\varphi, i } \leq   c_1(\varphi) + c_2(\varphi)
\sum_1^m\vee (r';  T_i,T_{i+1}) \,\,\,\,\forall  r(\cdot) \in
\wp(\Pi) \label{DiffconvMVarGomog16}
\end{equation}
\label{diffconvGomFuncMDegNVerThm6}
\end{thm}

 {\bf Доказательство. }{\bf Необходимость}. Пусть $ {\varphi
}(\cdot)$ есть ПРВ функция. Покажем, что тогда справедливо
неравенство (\ref{DiffconvMVarGomog16}). Воспользуемся
неравенством, вытекающим из неравенства треугольника,
$$
\Vert \tau(t_i) / \Vert \tau(t_i) \Vert  - \tau (t_{i-1}) / \Vert
\tau(t_{i-1} \Vert \Vert \leq
$$
$$
\leq \Vert r'(t_i) / \sqrt{ 1+\varphi'^2_t (r(t_i))}  -
r'(t_{i-1}) / \sqrt{ 1+\varphi'^2_t (r(t_{i-1}))} \Vert +
$$
$$
\mid {\varphi }'_t(r(t_i)) /  \sqrt{ 1+{\varphi }'^2_t (r(t_i))}-
{\varphi }'_t(r(t_{i-1})) / \sqrt{ 1+{\varphi }'^2_t (r(t_{i-1}))}
\mid .
$$
Так как $1 \leq  \sqrt{1+{\varphi }'^2_t (r(t_i))} \leq
\sqrt{1+L^2}$ для всех $ t_i \in [0, T(r)]$ , то очевидно,
существует такое $c_3 >1,$ для которого верно неравенство

\begin{equation} \Vert r'(t_i) / \sqrt{1+{\varphi }'^2_t (r(t_i))}-
r'(t_{i-1}) / \sqrt{1+{\varphi }'^2_t (r(t_{ i-1}))} \Vert \leq
c_3 \Vert r'(t_i) - r'(t_{i-1}) \Vert. \label{DiffconvMVarGomog17}
\end{equation}

Из свойств функции $\theta(x)= x / \sqrt{ 1+x^2}$ следует
неравенство
$$
\mid {\varphi }'_t (r(t_i)) / \sqrt{ 1+{\varphi }'^2_t (r(t_i))} -
{\varphi }'_t (r(t_{i-1})) / \sqrt{1+{\varphi }'^2_t (r(t_{i-1}))}
\mid \leq
$$
\begin{equation}
\leq \mid {\varphi }'_t (r(t_i)) - {\varphi }'_t (r(t_{i-1})) \mid
. \label{DiffconvMVarGomog18}
\end{equation}

Из (\ref{DiffconvMVarGomog17}) и (\ref{DiffconvMVarGomog18}) имеем
$$
\sup_{ t_j \in [T_i, T_{i+1}] } \,\, \sum_j \, \Vert \tau(t_j) /
\Vert \tau(t_j) \Vert - \tau (t_{j-1}) / \Vert \tau(t_{j-1}) \Vert
\Vert \leq
$$
\begin{equation}
\leq c_3 (\vee (  r' ; T_i ,T_{i+1}) + \vee (\Phi' ; T_i,T_{i+1})
). \label{DiffconvMVarGomog19}
\end{equation}

Так как по условию $ \varphi(\cdot)-$  ПРВ функция, то согласно
неравенству теоремы \ref{diffconvGomFuncMDegNVerThm5} для $\vee
(\Phi'; T_i,T_{i+1})$  с учетом (\ref{DiffconvMVarGomog19})
следует неравенство (\ref{DiffconvMVarGomog16}). Необходимость
доказана.

{\bf Достаточность}.  Пусть справедливо неравенство
(\ref{DiffconvMVarGomog16}). Покажем, что ${\varphi }(\cdot)$ -
ПРВ функция. Воспользуемся неравенством
$$
\Vert \tau(t_i) / \Vert \tau( t_i) \Vert  - \tau(t_{i-1}) /  \Vert
\tau (t_{i-1}) \Vert \Vert \geq \mid  {\varphi }'_t (r(t_i)) /
\sqrt{1+{\varphi }'_t (r(t_i))} -
$$
\begin{equation} - {\varphi }'_t (r(t_{i-1})) / \sqrt {1+{\varphi }'^2_t
(r(t_{i-1}))} \label{DiffconvMVarGomog20} \end{equation} Из
свойств функции $\theta(x) = x / \sqrt{1+x^2}$ и из $\Vert
{\varphi }'(z) \Vert \leq L$ для всех $z \in S_1^2(0),$ где
производная существует, следует существование константы $ c_4 >
0,$ для которой
$$
\mid {\varphi }'_t (r(t_i)) / \sqrt{1+{\varphi }'^2_t (r(t_i))} -
{\varphi }'_t (r(t_{i-1})) / \sqrt{1+{\varphi }'^2_t(r(t_{i-1}))}
\geq
$$
$$
\geq c_4 \mid {\varphi }'_t (r(t_i)) - {\varphi }'_t (r(t_{i-1}))
\mid,
$$
откуда с учетом (\ref{DiffconvMVarGomog20}) имеем
$$
\sup_{  t_j \in [T_i,T_{i+1}] }  \sum_j \Vert \tau(t_j) / \Vert
\tau( t_j) \Vert  - \tau(t_{j-1}) / \Vert \tau (t_{j-1}) \Vert
\Vert \geq
$$
$$
\geq c_6 \vee (\Phi';  T_i,T_{i+1})).
$$
Тогда для суммы
$$
\sum_i \vee (\Phi';  T_i,T_{i+1}))
$$
выполняется неравенство теоремы \ref{diffconvGomFuncMDegNVerThm5}.
Из теоремы \ref{diffconvGomFuncMDegNVerThm5} следует, что
${\varphi }(\cdot)$ - ПРВ функция. Достаточность и теорема
доказаны.  $\Box$

\begin{rem}
Поскольку любой вектор или отрезок образует угол, не больший
$\frac{\pi}{4}$ с какой то координатной плоскостью, то достаточно
рассматривать только координатные плоскости $\Pi$.
\end{rem}

\vspace{0.5cm}

\subsection{Необходимые и достаточные условия представимости
положительно однородной функции $m-$ степени в виде разности
выпуклых  в $n-$ мерном случае }

\vspace{0.5cm}

Пусть $\varphi(\cdot)$ есть п.о. $m-$ ого порядка функция от $n$
переменных, $m-$ натуральное число, большее 1. По функции
$\varphi(\cdot)$ мы можем построить, как уже делали ранее в
параграфе 2.3, { \em соответствующую } п.о. функцию 1-ого порядка
$\psi(\cdot): \mathbb{R}^n \rar \mathbb{R} $

\be
\psi(q) = \| q \| \,\, \varphi(\frac{q}{\| q \|}) \,\,\,\,
\forall q \in \mathbb{R}^n. \label{diffconvGomFuncMDegNVer21}
\ee

Было доказано в параграфе 2.3., что если $ \varphi(\cdot)- $
выпуклая, то  соответствующая ей функция $\psi(\cdot) $ тоже
выпуклая. Но поскольку  функции $\varphi(\cdot)$, $\psi(\cdot)$
полностью определяются значениями на единичной $(n-1)-$ мерной
сфере $S_1^{n-1}(0)$,  то верно и обратное утверждение. А именно:
если п.о. 1-ого порядка функция $\psi(\cdot)- $ выпуклая, то
построенная по ней п.о. $m-$ ого порядка функция $\varphi(\cdot)-
$ тоже выпуклая.

Поэтому, если есть представимость в виде разности выпуклых функции
$ \varphi(\cdot) $, то есть представимость в виде разности
выпуклых и функции $\psi(\cdot)$. Верно обратное: если функция
$\psi(\cdot)$ представима в виде разности выпуклых, то и функция
$\varphi(\cdot) $ также выпуклая.

Отметим, что та же техника доказательства позволяет установить
факт \cite{proudconvex2}, что в двумерном случае для произвольной
выпуклой функции $f(\cdot):\mathbb{R}^2 \rar \mathbb{R}$,
$f(0)=0$, п.о. функция $\psi(\cdot)$, построенная по значениям
выпуклой функции $f(\cdot)$ на произвольной кривой $r(\cdot)$,
ограничивающей выпуклое компактное множество $D$ с не пустой
внутренностью, $0 \in \mbox{int } D$, есть выпуклая функция.

Из сказанного выше следует, что написанные выше теоремы
\ref{diffconvGomFuncMDegNVerThm5},
\ref{diffconvGomFuncMDegNVerThm6} справедливы для п.о. функции
$m-$ ого порядка.

Возьмем произвольную плоскость $\Pi$, $0 \in \Pi$, и произвольную
$r(\cdot) \in \wp(\Pi)$. Введем функцию
$$
\Phi(t)=\varphi(r(t)) \,\,\, \forall t \in [0,T].
$$

Для $r(\cdot) \in \wp(\Pi) $ выделим те ее участки с параметрами
$[T_i, T_{i+1}] \subset [0,T(r)]$, где производная $r'(t), t \in
[0, T(r)],$ составляет с плоскостью $\Pi$ угол $\delta $: $0 <
\delta \leq \frac{\pi }{4} $. Обозначим объединение таких участков
через $\Lambda(r)$.

\begin{thm}
Для того чтобы п.о. $m$-ой степени липшицевая функция
$\Phi(\cdot):\mathbb{R}^n \rar \mathbb{R} $ с константой Липшица
$L$ была ПРВ функцией, необходимо и достаточно, чтобы для любой
плоскости $\Pi$ и любой кривой $r(\cdot) \in  \wp(\Pi)$,
параметризованной естественным образом с параметром $t \in [0,T],
\, T=T(r),$ и любого подмножества параметров $[T_i, T_{i+1} ]
\subset \Lambda(r) $ нашлись константы $C_1 = C(\varphi)>0$, $C_2
= C_2(\varphi)>0$ такие, что \be
 \sum_1^m \vee (\Phi'; T_i,T_{i+1})<  C_1 + C_2 \sum_1^m\vee (r';
 T_i,T_{i+1}),
\label{diffconvGomFuncMDegNVer22} \ee
где $\Phi(t)=\varphi(r(t))$
и производные берутся там, где они существуют.
\label{diffconvGomFuncMDegNVerThm7}
\end{thm}
\begin{rem}
Необходимость рассмотрения системы подмножеств $\Lambda(r) $
объясняется тем, что вариация $\vee (r'; 0,T(r))$ может быть
неограниченной для всей кривой $r(\cdot)$.
\end{rem}

Вводим понятие поворота кривой также, как это было сделано в
предыдущем параграфе.

\begin{thm} Для того, чтобы произвольная липшицевая п.о. степени
$m$ функция $\varphi(x) :\mathbb{R}^n \rightarrow \mathbb{R}$ была
ПРВ функцией, необходимо и достаточно, чтобы для любой плоскости
$\Pi$, $0 \in \Pi$, и любой $r(\cdot) \in \wp(\Pi)$
параметризованной естественным образом: $t \in [0, T(r)]$, а также
любого подмножества параметров $[T_i, T_{i+1} ] \subset \Lambda(r)
$ существовали константы $c_1(\varphi )>0$, $ c_2(\varphi)>0$
такие, что сумма поворотов $O_{\varphi, i}$ кривой $R(\cdot)$ для
$t \in [T_i, T_{i+1}] $ для  всех $i$ была ограничена сверху
неравенством
\begin{equation}
\sum_1^m O_{\varphi, i } \leq   c_1(\varphi) + c_2(\varphi)
\sum_1^m\vee (r';  T_i,T_{i+1}) \,\,\,\,\forall  r(\cdot) \in
\wp(\Pi). \label{DiffconvMVarGomog23}
\end{equation}
\label{diffconvGomFuncMDegNVerThm8}
\end{thm}

\begin{rem}
При доказательстве достаточности мы растягивали многогранник $G_m$
и превращали его в выпуклый многогранник, являющийся
субдифференциалом в нуле выпуклой п.о. первой степени функции.
Такая операция приводит к ограниченному выпуклому многограннику
только тогда, когда исходная п.о. функция представима в виде
разности выпуклых. Если исходная п.о. функция не представима в
виде разности выпуклых, то существует направление неограниченного
расширения получающихся на каждом шаге выпуклых многогранников.

Поскольку любой вектор или отрезок образует угол, не больший
$\frac{\pi}{4}$ с какой-то координатной  плоскостью $\Pi$, то
достаточно рассматривать только координатные плоскости.
\end{rem}

\newpage

\vspace{1cm}

\section{К ВОПРОСУ О ПРЕДСТАВИМОСТИ ФУНКЦИИ ДВУХ ПЕРЕМЕННЫХ В
ВИДЕ РАЗНОСТИ ВЫПУКЛЫХ ФУНКЦИЙ}

\vspace{1cm}

В данном параграфе приведены исправленные необходимые и
достаточные условия представимости произвольной липшицевой функции
двух переменных в виде разности выпуклых функций
\cite{proudconvex2}. Дана также геометрическая интерпретация этих
условий. Приведен алгоритм такого представления, результатом
которого есть последовательность равномерно сходящихся выпуклых
функций.

\noindent

\vspace{0.5cm}

\subsection{Введение}

\vspace{0.5cm}

Эта проблема была впервые сформулирована академиком
А.Д.Александровым в статье \cite{aleksandrov1} и исследована
многими российскими и зарубежными  математиками
\cite{aleksandrov1}, \cite{aleksandrov2}, \cite{zalgaller1},
\cite{lupikovdissertation}, \cite{proudconvex1},
\cite{proudconvex2}, \cite{ginchev}, \cite{pallaschke},
\cite{Hartman}, \cite{hiriarturruty}, \cite{veselyzajicek}.
Решение этой проблемы интересно как для геометров, так и для
математиков, занимающихся оптимизацией, например, для построения
квазидифференциального исчисления \cite{demvas}.

Необходимые и достаточные условия представимости функции одной
переменной в виде разности выпуклых, т.е. условия. когда функция
является ПРВ функцией, хорошо известны. Эти условия могут быть
записаны в следующем виде.

Пусть $x \rightarrow f(x): [a,b] \rightarrow \mathbb{R}$ -
произвольная липшицевая функция. Известно, что множество $N_f$,
где функция $ f(\cdot)$ дифференцируемая, есть множество полной
меры на [a,b]. Для того, чтобы функция $f(\cdot)$ была представима
в виде разности выпуклых функций, необходимо и достаточно, чтобы
выполнялось условие
$$
         \vee (f'; a,b) < \infty ,
$$
где производные вычисляются там, где они существуют. Символ $\vee$
означает вариацию функции $ f'$ на отрезке [a,b].

В той же статье \cite{aleksandrov1}  А.Д.Александров задает вопрос
о представимости функции в виде разности выпуклых, если она
является таковой для любой прямой в области определения. Ответ на
этот вопрос отрицательный (см. \cite{hiriarturruty},
\cite{veselyzajicek} ).

Согласно терминологии А.Д.Александрова под многогранной
кусочно-линейной функцией с конечным числом граней будем понимать
такую функцию, график которой состоит из конечного числа частей
плоскостей, которые называются гранями.

В статье  \cite{proudconvex1} даны необходимые и достаточные
условия представимости произвольной липшицевой положительно
однородной (п.о.) функции трех переменных в виде разности выпуклых
функций. Результат может быть распространен на положительно
однородные функции $m$-ой степени. Теперь откажемся от условия
положительной однородности и будем рассматривать произвольную
липшицевую функцию $f(\cdot)$ с константой Липшица $L$ от двух
переменных $(x,y) \rightarrow f(x,y):D \rightarrow \mathbb{R},$
где $D$ есть выпуклое открытое ограниченное множество в
$\mathbb{R}^2$, так что его замыкание $\bar{D}$ - компакт.
Приведем алгоритм такого представления и найдем необходимые и
достаточные условия сходимости построенной последовательности
функций.

Дадим определение звездной области.

\begin{defi} Область ${D}$ называется звездной, если существует точка $O
\in {D} $, для которой отрезок, соединяющий точку $O$ с любой
точкой области ${D}$, целиком принадлежит области $D$.
\end{defi}

Пусть ${\wp}(D)$ - класс кривых на плоскости  $XOY$ в множестве
$D$, ограничивающих звездные области, замыкание которых компактно.
Параметризуем кривые $r \in {\wp}(D)$ естественным образом, т.е.
параметр $\tau$ точки $M$ на кривой $r(\cdot)$ равен длине кривой
между $M$ и начальной точкой. Обозначим такую кривую как $r(t), t
\in [0,T_r]$.

Заметим, что множеству кривых ${\wp}(D)$ принадлежит класс кривых
на плоскости  $XOY$ в множестве $D$, параметризованных
естественным образом и ограничивающих выпуклые компактные
множества на множестве $D$. Обозначим этот класс кривых через
$\varrho(D)$.

С помощью кривых  $r \in \wp (D)$ необходимые и достаточные
условия представимости функции $f(\cdot)$ в виде разности выпуклых
функций  могут быть записаны в следующем виде.

\begin{thm}
Для того, чтобы липшицевая функция $z \rightarrow f(z):D
\rightarrow \mathbb{R}$ была представима в виде разности выпуклых
функций (была ПРВ функцией), необходимо и достаточно, чтобы для
любой кривой $r \in \wp (D)$ и любой ее системы подмножеств
выполнялось неравенство
$$
 (\exists c_1(D,f), c_2(D,f)>0) (\forall r \in \wp(D)) ,
$$
\begin{equation}
 \vee (\Phi'; 0,T_r)<  c_1(D,f) + c_2(D,f) \vee (r';
 0,T_r), \label{difconv3}
\end{equation} где $\Phi(t)=f(r(t)) \;\;\;\; \forall t \in
[0,T_r]$. \label{difconvexthm1}
\end{thm}

\begin{rem}
Условие теоремы означает следующее.

Существуют константы $ c_1(D,f), c_2(D,f)>0 $, что для любой $r
\in \wp(D)$ и любой системы подмножеств $[T_i, T_{i+1} ] \subset
[0, T_r] $, $i \in 1:m,$
$$
 \sum_1^m \vee (\Phi'; T_i,T_{i+1})<  c_1(D,f) + c_2(D,f) \sum_1^m\vee (r';
 T_i,T_{i+1}).
$$
Необходимость рассмотрения системы подмножеств для кривой $r \in
\wp (D)$, для которой выполняется неравенство (\ref{difconv3}),
объясняется тем, что вариация $\vee (r'; 0,T_r)$ может быть
неограниченной для всей кривой $r(\cdot)$.
\end{rem}

\begin{rem}  Впервые подобные условия для вариации функции на кривых
были приведены  в \cite{proudconvex1}.
\end{rem}

Доказательство основано на специальном алгоритме представления
функции $f(\cdot)$ в виде разности выпуклых функций. В результате
получаем конечные или бесконечные последовательности выпуклых
функций, равномерно сходящиеся на $D$  к выпуклым функциям,
разность которых есть исходная функция $f(\cdot)$, если условия
теоремы \ref{difconvexthm1} выполняются.

Ниже приведен алгоритм представления функции $f(\cdot)$ в виде
разности выпуклых функций и доказана его сходимость, если условия
теоремы \ref{difconvexthm1} выполняются.

Для представления функции $f(\cdot)$ в виде разности выпуклых
функций будем использовать две операции, в результате которых
получаем конечное или счетное число выпуклых многогранных
кусочно-линейных функций, определенных на $D$.

{\em Первая операция} - это приближение функции $f(\cdot)$
многогранной кусочно-линейной функцией $f_n(\cdot)$ с конечным
числом граней.

{\em Вторая операция} - это представление функции $f_n(\cdot)$ в
виде разности выпуклых многогранных кусочно-линейных функций
$f_{1,n}(\cdot):D \rightarrow \mathbb{R}$ и $f_{2,n}(\cdot):D
\rightarrow \mathbb{R}$ согласно алгоритму, описанному ниже.

Далее доказывается, что, если выполняются условия теоремы
\ref{difconvexthm1}, то из последовательностей $f_{1,n}(\cdot)-
c_{1,n}$ и $f_{2,n}(\cdot) - c_{1,n}$, где $c_{1,n}=f_{1,n}(a)$,
$a - $ произвольная внутренняя точка области $D$, можно выбрать
сходящиеся подпоследовательности.

Когда условия теоремы выполняются, то, как будет показано,
вариация производной вдоль любого отрезка множества $D$ выпуклых
функций $f_{1,n}(\cdot)$ и $f_{2,n}(\cdot)$ ограничена сверху
константой, зависящей от $D$ и $f$.

Метод представления конечной многогранной функции в виде разности
выпуклых подобен методу, использованному А.Д.Александровым в
\cite{aleksandrov1} при исследовании возможности представления
специального вида функций в виде разности выпуклых.

\vspace{1cm}

\subsection {Доказательство теоремы}

Начнем доказательство теоремы с описания алгоритма.

\vspace{0.5cm}

{\bf ОПИСАНИЕ АЛГОРИТМА}

\vspace{0,5cm}

1. Производим достаточно равномерную триангуляцию области $D$ и
строим по каждому треугольнику линейную функцию, значения которой
равны значениям функции $f(\cdot)$ в вершинах треугольника.
Функцию с получившимся графиком обозначим через $f_n(\cdot):D
\rightarrow \mathbb{R}$, где $n$ равно числу треугольников, на
которые мы разбиваем область $D$.

2. Представляем функцию $f_n(\cdot)$ в виде разности выпуклых
согласно алгоритму, как это описано ниже.

Предварительно введем понятие двугранного угла. Будем понимать под
двугранным углом функцию, определенную на $D$, график которой
состоит из полуплоскостей с общей граничной прямой, называемого
ребром двугранного угла.

Рассмотрим все выпуклые двугранные углы, части графиков которых
принадлежат графику функции $f_n(\cdot)$. Определяем эти
двугранные углы на всей области $D$. Просуммируем все такие
выпуклые двугранные углы. В итоге получим выпуклую многогранную
функцию $f_{1,n}(\cdot):D \rightarrow \mathbb{R}$. Доказывается
\cite{aleksandrov1}, что разность \begin{equation}
 f_{1,n}(\cdot) - f_n(\cdot) = f_{2,n}(\cdot)
\label{difconv1} \end{equation} есть также выпуклая многогранная
функция.

Действительно, для доказательства достаточно показать, что все
двугранные углы, части графиков которых принадлежат графику
функции $f_{1,n}(\cdot) - f_{n}(\cdot), - $ выпуклые. Для этого
покажем, что любая точка, лежащая на проекции ребра произвольного
двугранного угла функции $f_{1,n}(\cdot) - f_{n}(\cdot),$ имеет
малую окрестность, где функция $f_{1,n}(\cdot) - f_{n}(\cdot)$
выпукла.

Если берем точку, в малой окрестности которой функция
$f_{n}(\cdot)$ линейна, то локальная выпуклость разности
$f_{1,n}(\cdot) - f_{n}(\cdot)$ очевидна. Пусть берем точку,
лежащую на проекции на плоскость ребра выпуклого двугранного угла,
часть графика которого принадлежит графику функции $f_{n}(\cdot)$.
Поскольку согласно алгоритму этот же двугранный угол входит в
сумму выпуклых двугранных углов, образующих функцию
$f_{1,n}(\cdot)$, то опять разность $f_{1,n}(\cdot) -
f_{n}(\cdot)$ будет локально выпуклой в окрестности
рассматриваемой точки. Если же точка лежит на проекции ребра
вогнутого двугранного угла, часть графика которого принадлежит
графику функции $f_{n}(\cdot)$, то $-f_{n}(\cdot)-$ локально
выпукла в окрестности этой точки, а поэтому разность
$f_{1,n}(\cdot) - f_{n}(\cdot)$ снова локально выпукла в той же
окрестности. Из локальной выпуклости всех двугранных углов функции
$f_{1,n}(\cdot) - f_{n}(\cdot)$ следует ее выпуклость на всем
множестве $D$.

Покажем, что при выполнении теоремы \ref{difconvexthm1} из
последовательности функций $f_{1,n}(\cdot) - c_{1,n} $ можно
выделить подпоследовательность, равномерно сходящуюся  на $D$ к
выпуклой функции $f_1(\cdot)$ при $n \rightarrow +\infty$. Тогда
из (\ref{difconv1}) будет следует, что подпоследовательность
функций $f_{2,n}(\cdot)-c_{1,n}$ также равномерно сходится к
выпуклой функции $f_2(\cdot)$. Для функций $f_1(\cdot)$ и
$f_2(\cdot)$ верно равенство \begin{equation}
 f_{1}(\cdot) - f_2(\cdot) = f(\cdot).
\label{difconv2} \end{equation}

Начнем доказательство с одномерного случая, когда $D=[a,b] \subset
\mathbb{R}$.

Приблизим функцию $f(\cdot)$ кусочно-линейной функцией
$f_n(\cdot)$  с любой степенью точности. На первом шаге выделяем
все выпуклые двугранные углы, части графиков которых принадлежат
графику функции $f_n(\cdot)$. Распространяем их на весь отрезок
$[a,b]$ и просуммируем. В итоге получим выпуклую кусочно-линейную
функцию $f_{1,n}(\cdot):[a,b] \rightarrow \mathbb{R}$. Согласно
сказанному выше разность $f_{1,n}(\cdot) - f_{n}(\cdot)$ есть
снова выпуклая кусочно-линейная функция на $[a,b]$.

Покажем, что вариация производных функций $f_{1,n}(\cdot)$ и
$f_{2,n}(\cdot)$ на отрезке $[a,b]$ ограничена сверху той же
константой $c$, что вариация производной функции $f(\cdot)$, т..е.
$$
\vee (f_{1,n}' ; a,b) \leq  c.
$$
Последнее следует из цепочки неравенств
$$
\vee (f_{1,n}' ; a,b) \leq \vee (f_n' ; a,b) \leq \vee (f' ; a,b)
\leq c.
$$
Но тогда из $f_{1,n}(\cdot)$ можно вычесть константу
$c_{1,n}=f_{1,n}(a)$, $ a \in \mbox{int} \, D$, чтобы функции
$f_{1,n}(\cdot)-c_{1,n}$ были ограниченными на отрезке $[a,b]$ в
совокупности по $n$, т.е. равностепенно ограниченными. Из оценки
для вариации производной, не зависящей от $n$, следует
равностепенная непрерывность функций $f_{1,n}(\cdot)-c_{1,n}$. Из
теоремы Арцела  получим, что из последовательности выпуклых
функций $f_{1,n}(\cdot) - c_{1,n}$ можно выделить
подпоследовательность функций $f_{1,n_k}(\cdot) - c_{1,n_k}$,
которая сходится равномерно при $n_k \rightarrow \infty$ к
некоторой выпуклой на $[a,b]$ функции $f_{1}(\cdot)$.
Соответственно, последовательность функций $f_{2,n_k}(\cdot) -
c_{1,n_k}$ также равномерно на $[a,b]$ сходится при $n_k
\rightarrow \infty$ к некоторой выпуклой на $[a,b]$ функции
$f_{2}(\cdot)$. В итоге будем иметь
$$
f(\cdot)=f_1(\cdot)- f_2(\cdot).
$$

Перейдем к двумерному случаю и покажем, что тот же алгоритм
приводит к паре выпуклых функций на $D$, разность которых есть
исходная функция $f(\cdot)$.

Возьмем произвольную кривую $r(\cdot) \in \wp(D).$Пусть
$$
    \Phi(t)=f(r(t)) \;\;\; \forall t \in [0,T_r].
$$
Покажем, что $\Phi(\cdot)$ - липшицевая с константой $L$.
Действительно, для любых $t_1,t_2 \in [0,T_r]$ имеем
$$
 \mid \Phi (t_1) - \Phi (t_2) \mid = \mid f(r(t_1)) - f(r(t_2)) \mid \le
 L \Vert  r(t_1) - r(t_2) \Vert  \le  L \mid t_1 - t_2 \mid.
$$
Поэтому \cite{kolmogorovfomin} $\Phi (\cdot)$ почти всюду (п.в.)
дифференцируемая на $[0,T_r].$ Множество точек дифференцируемости
функции $\Phi (\cdot)$ на $[0,T_r]$ обозначим через $N_r.$

Докажем, что если существуют константы $c_1(D, f), c_2(D,f) > 0$
такие, что для произвольной кривой $r(\cdot) \in \wp(D)$
выполняется условие теоремы 1, то из последовательностей функций
$f_{1,n}(\cdot)-c_{1,n}$, $f_{2,n}(\cdot)-c_{1,n}$,  можно выбрать
подпоследовательности, равномерно на $D$ сходящиеся к выпуклым
функциям $f_{1}(\cdot)$, $f_{2}(\cdot)$ соответственно, для
которых верно равенство (\ref{difconv2}).

Доказательство будем основывать на леммах, приведенных ниже.

\begin{lem} Для любой выпуклой п.о. степени 1 функции $q
\rightarrow \psi (q): \mathbb{R}^2 \rightarrow \mathbb{R}$ и
кривой $ r(\cdot) \in \wp(D)$, а также любого ее подмножества
верно неравенство
$$
    \vee (\Theta '; 0,T_r)  <  c_1(D, \psi) + c_2(D, \psi)
    \vee(r'; 0,T_r ),
$$
где $\Theta(t) = \psi(r(t))$  для всех $t \in [0,T_r], c_1(D,
\psi), c_2(D, \psi)$ - некоторые положительные константы.
\label{1lemdifconv}
\end{lem}

{\bf Доказательство.}  Без ограничения общности будем считать, что
$\psi(\cdot)$ есть гладкая функция на $\mathbb{R}^2 \backslash
\{0\}.$ Пусть
$$
 \psi (r(t)) = \mbox{max} \,_{v \in \partial \psi(0)} (v, r(t))=(v(t),r(t)), \;\;
v(t) \in \partial \psi(0),
$$
где $\partial \psi (0)$ - субдифференциал функции $\psi(\cdot)$ в
нуле. Будем также считать,  что $r(\cdot)$ - дифференцируемая
кривая по $t \in [0,T_r]$ .

Очевидно, что
$$
\psi'(r(t))=(v'(t),r(t))+(v(t),r'(t)).
$$
Так как $r(t)$ есть нормаль к границе множества $\partial \psi(0)$
в точке $v(t)$, то векторы $v'(t)$ и $r(t)$ перпендикулярны друг к
другу, а следовательно, $(v'(t),r(t))=0.$  Поскольку кривая
$r(\cdot)$ параметризована естественным образом, то $\Vert r'(t)
\Vert =1$ для любых $t \in [0,T_r]$.

Нетрудно проверить следующую цепочку неравенств
$$
\mid \psi' (r(t_1)) - \psi' (r(t_2)) \mid= \mid (v(t_1),r'(t_1)) -
(v(t_2),r'(t_2)) \mid = \mid (v(t_1)- v(t_2),r'(t_1)) +
$$ $$
(v(t_2),r'(t_1))- (v(t_2),r'(t_2)) \mid \leq  \Vert v(t_1) -
v(t_2) \Vert \; \Vert r'(t_1) \Vert + \Vert r'(t_1) - r'(t_2)
\Vert \; \Vert v(t_2) \Vert
$$
Отсюда следует, что
$$
   \vee (\Theta '; 0,T_r) < 2P(\partial \psi (0)) + L(D) \vee(r'; 0; T_r ),
$$
где $P(\psi (0))$- длина кривой, ограничивающей выпуклое
компактное множество $\partial \psi(0) \subset \mathbb{R}^2$ и
$L(D)$ - константа Липшица функции $\psi(\cdot).$

Пусть теперь $\psi (\cdot)$ - произвольная выпуклая ПО функция. С
любой степенью точности ее можно приблизить  на единичном круге
выпуклой ПО дифференцируемой на $\mathbb{R}^2 \backslash \{0\}$
функцией $\hat{\psi}(\cdot)$ так, чтобы в метрике Хаусдорфа
субдифференциалы в нуле этих функций отличались друг от друга как
угодно мало. Но тогда и длины кривых, ограничивающих их
субдифференциалы, будут отличаться друг от друга как угодно мало.
Также кривую $r(\cdot)$ можно приблизить дифференцируемой кривой
таким образом, чтобы их производные по $t$ в точках
дифференцируемости кривой $r(\cdot)$ отличались друг от друга по
норме на произвольно малую величину. Таким образом, любые конечные
суммы, используемые при вычислении вариаций  функций
$\Theta'(\cdot)$ и $\hat{\Theta}'(\cdot)$  для негладкого и
гладкого случая, могут быть сделаны за счет приближения как угодно
близкими друг к другу. Но поскольку вариацию функции
$\hat{\Theta}'(\cdot)$ можно ограничить сверху величиной,
зависящей только от множества $D$ и некоторых констант, то лемма 1
доказана. $\Box$

На основе этой леммы докажем утверждение.

\begin{lem}
Пусть $(x,y) \rightarrow f_1(x,y): \mathbb{R}^2 \rightarrow
\mathbb{R}-$ непрерывная выпуклая функция и $r(\cdot) \in \wp(D),
t \in [0,T_r].$ Тогда существуют константы $c_1(D, f_1),
c_2(D,f_1)>0,$ что для кривой $r(\cdot)$ и любого ее подмножества
\begin{equation}
   \vee (\Phi_1'; 0,T_r)  <  c_1(D, f_1) + c_2(D, f_1) \vee(r'; 0,T_r ),
\label{difconv4} \end{equation} где $\Phi_1(t) = f_1(r(t)), t \in
[0,T_r].$ \label{2lemdifconv}
\end{lem}

{\bf  Доказательство.} На начальном этапе будем считать, что
$f_1(\cdot,\cdot)$ дважды непрерывно дифференцируемая функция на
$D$, которая принимает неотрицательные значения и начало координат
$-$ ее точка минимума, а также, что $0=(0,0)$ принадлежит
внутренности выпуклой области на $\mathbb{R}^2$ с границей
$r(\cdot).$

Построим для функции $f_1(\cdot,\cdot)$ п.о. степени 1 функцию
$\psi(\cdot),$  которая на $r(\cdot)$ принимает значения, равные
$f_1(r(\cdot)).$ Покажем, что $\psi(\cdot)$ - выпуклая.

Рассмотрим функцию
$$
f_{\varepsilon}(x,y)=f_1(x,y)+\varepsilon ( \mid \mid x \mid
\mid^2 + \mid \mid y \mid \mid^2  ), \,\,\, \varepsilon>0.
$$
Разобьем отрезок $[ 0,T_r] $ точками  $\{t_i\} , i \in 1:J ,$ на
равные отрезки. Построим плоскости $\pi_i$ в $\mathbb{R}^3,$
проходящие соответственно через точки $(0,0,0),
(r(t_i),f_{\varepsilon}(r(t_i))),
(r(t_{i+1}),f_{\varepsilon}(r(t_{i+1})), i \in 1:J $. Части
плоскостей $\pi_i ,i \in 1:J$ , определенных в секторах,
образуемых векторами $(0,0), r(t_i), r(t_{i+1})$, определяют
график п.о. степени 1 многогранной функцию $(\psi_{\varepsilon})_J
(r(\cdot)).$ Будем понимать под двугранным углом функцию, график
которой состоит из полуплоскостей с общей граничной прямой,
включающих плоскости $\pi_i,$ построенные в соседних секторах.
Покажем, что все двугранные углы функции $(\psi_{\varepsilon})_J
(r(\cdot)),$ образуемые смежными плоскостями $\pi_i, i \in J,-$
выпуклые.

Поскольку всегда любую кривую $r(\cdot) \in \wp(D)$ можно
приблизить  с любой степенью точности гладкой кривой из $\wp(D),$
то без ограничения общности будем считать, что $r(\cdot)$ -
гладкая дифференцируемая кривая с производной $r'(\cdot).$

Под градиентом плоскости $\pi_i$ будем понимать градиент линейной
функции, график которой совпадает с плоскостью $\pi_i$. Обозначим
градиенты плоскостей $\pi_i$ и $\pi_{i+1}$ через $\nabla \pi_i$ и
$\nabla \pi_{i+1}$ соответственно.  Воспользуемся теоремой о
средней точке, согласно которой существует такая точка $t_m \in
[t_i, t_{i+1}],$ что
$$
           \partial f_{\varepsilon}(r(t_m))/ \partial e_i =
           (\nabla \pi_i, e_i),
$$
где
$$
e_i=(r(t_{i+1})-r(t_i))/ \mid \mid r(t_{i+1})-r(t_i) \mid \mid .
$$
Аналогично для плоскости $\pi_{i+1}$ и некоторой точки $t_c \in
[t_{i+1},  t_{i+2}]$ имеем
$$
           \partial f_{\varepsilon}(r(t_c))/ \partial e_{i+1} =
           (\nabla \pi_{i+1}, e_{i+1}),
$$
где
$$
e_{i+1}=(r(t_{i+2})-r(t_{i+1}))/ \mid \mid r(t_{i+2})-r(t_{i+1})
\mid \mid .
$$
Функция $f_{\varepsilon}(\cdot)$ сильно выпуклая, так как ее
матрица вторых частных производных положительно определенная.
Любая выпуклая функция имеет неубывающую производную по
направлению вдоль произвольного луча. Но для сильно выпуклой
функции производная по касательному направлению к кривой вида
$r(x_0,\tau, g)= x_0+\tau g +o_{\varepsilon}(\tau)$, $g \in
\mathbb{R}^n$,$ \tau >0$ в малой окрестности точки $x_0$ есть
возрастающая функция вдоль этой кривой. Поэтому для достаточно
большом $J$ и равномерном разбиении кривой $r(\cdot)$ точками
$t_i$ имеем
$$
\partial f_{\varepsilon}(r(t_m))/ \partial e_i <
\partial f_{\varepsilon}(r(t_c))/ \partial e_{i+1},
$$
или
$$
(\nabla \pi_i, e_i) < (\nabla \pi_{i+1}, e_{i+1}).
$$
Учтем также, что разность $\nabla \pi_{i+1} - \nabla \pi_i$
перпендикулярна вектору $r(t_{i+1}).$ Отсюда и из неравенства выше
следует, что двугранный угол $\pi_i, \pi_{i+1}$ - выпуклый. При $J
\rightarrow \infty$
$$ (\psi_{\varepsilon})_J(\cdot) \Rightarrow (\psi_{\varepsilon})(\cdot).$$
Так как точечный предел для выпуклых функций равносилен
равномерному пределу, то $\psi_{\varepsilon}(\cdot)$ - выпуклая
функция. Также $\psi_{\varepsilon}(\cdot) \Rightarrow \psi
(\cdot)$ при $\varepsilon \rightarrow +0,$ т.е. $\psi(\cdot)-$
выпуклая, что и требовалось доказать.

Очевидно, что градиенты линейных функций, графики которых есть
$\pi_i, i \in J,$ ограничены константой, зависящей только от
множества $D$ и самой функции $ f_1(\cdot,\cdot).$ Верно равенство
$$
    \psi (r(t)) = f_1(r(t)) \;\; \forall t \in [0,T_r].
$$
Ясно, что $\psi (\cdot,\cdot)$ строится однозначно по функции
$f_1(\cdot,\cdot)$ и выбранной кривой $r(\cdot).$ Из сказанного
выше следует, что функция $\psi (\cdot,\cdot)$ есть липшицевая с
константой $L(D,f)$.

Пусть
$$
 \Psi_1(t) = \psi (r(t)) \;\; \forall t \in [0,T_r].
$$
Поскольку
$$
  \vee (\Phi_1 '; 0,T_r)  =  \vee( \Psi_1 ' ; 0,T_r) ,
$$
то из леммы 1 следует, что для некоторых констант $ c_1(D, f_1),
c_2(D, f_1) > 0 $
$$
\vee (\Phi_1 '; 0,T_r) \leq   c_1(D, f_1) + c_2(D, f_1) \vee(r';
0,T_r ).
$$
Если функция $ f_1(\cdot,\cdot)$ не есть дважды непрерывно
дифференцируемая, то ее можно приблизить выпуклой дважды
непрерывно дифференцируемой функцией $ \tilde{f}_1(\cdot,\cdot)$ и
построить соответствующую ей функцию $\tilde{\psi} (\cdot,\cdot)$
так, чтобы значения  функций $\psi (\cdot,\cdot)$, $\tilde{\psi}
(\cdot,\cdot)$ и их производных там, где они существуют, как
угодно мало отличались друг от друга. Но тогда аналогичное будет
верно для функций $\Psi_1(\cdot)$, $\tilde{\Psi}_1(\cdot)$,
построенных по $\psi (\cdot,\cdot), \tilde{\psi} (\cdot,\cdot)$
соответственно, и их производных. Значит написанное выше
неравенство для вариации производных функции $\Psi_1(\cdot)$ верно
для общего случая. Лемма \ref{2lemdifconv} доказана. $\Box$

Из леммы \ref{2lemdifconv} следует, что если $f(\cdot,\cdot)$
представима  в виде разности выпуклых функций, т.е.
$$ f(z) = f_1(z) - f_2(z) \;\;\;\; \forall z \in D,   $$
где $f_i(\cdot,\cdot), i=1,2,$ - выпуклые, то условие
(\ref{difconv3})  c необходимостью выполняется. Действительно, для
произвольной $r(\cdot) \in \wp(D)$ введем обозначения
$$
\Psi_1(t) = f_1(r(t)), \Psi_2(t) = f_2(r(t)) \;\;\; \forall t \in
[0,T_r].
$$
Поскольку \cite{kolmogorovfomin}
$$
\vee (\Phi ';0,T_r) \leq \vee (\Phi_1 '; 0,T_r) + \vee (\Phi_2
';0,T_r)
$$
то, принимая во внимание неравенство (\ref{difconv4}), неравенство
(\ref{difconv3}) с необходимостью выполняется.

Докажем достаточность условия (\ref{difconv3}) для представления
функции $f(\cdot)$ в виде разности выпуклых функций.

Прежде всего покажем, что  для любого $r(\cdot) \in \wp(D)$ и
констант $c_1(D,f), c_2(D,f)>0 $ верно неравенство
$$
\vee (\Phi_n ';0,T_r) \leq c_1(D,f) + c_2(D,f) \vee(r'; 0, T_r),
$$
где $\Phi_n(t) = f_n(r(t))$. Действительно, для любой триангуляции
области $D$ градиенты в точках $r(t_k) \in r(\cdot), t_k \in
[0,T_r],$ линейных функций, графики которых есть грани функции
$f_n(\cdot)$, будут с любой степенью точности $\varepsilon_n$, где
$\varepsilon_n \rightarrow +0$, близки к обобщенным градиентам
функции $f(\cdot)$. Поэтому произвольная конечная сумма
$$
\sum_{i=1}^{N} \mid \Phi_n'(t_i) - \Phi_n'(t_{i+1}) \mid
$$
для больших $n$ будет как угодно мало отличаться от суммы
$$
\sum_{i=1}^{N} \mid \Phi'(t_i) - \Phi'(t_{i+1}) \mid .
$$
А поскольку вариация функции $\Phi_n'(\cdot)$  может только
возрастать при вложенности триангуляций области $D$ при увеличении
$n$, то отсюда и из сказанного выше следует, что
\begin{equation} \vee (\Phi_n';0,T_r) \leq \vee
(\Phi';0,T_r)+\delta(n) \leq c_1(D,f) + c_2(D,f) \vee(r'; 0, T_r),
\label{difconv5}
\end{equation} где $\delta(n) \rightarrow +0 $ при $n \rightarrow
\infty$.

Вариация производных по направлению вдоль произвольного отрезка
суммы выпуклых функций равна сумме вариаций производных этих
выпуклых функций по тому же отрезку. Если будет доказано, что
сумма вариаций  производных всех выпуклых двугранных углов функции
$f_n(\cdot)$ вдоль любого отрезка области $D$  ограничена сверху
константой, независящей от $n$, то отсюда будет следовать, что
ограничена сверху той же константой вариация производной функции
$f_{1,n}(\cdot)$ вдоль произвольного отрезка области $D$. Отсюда
следует равностепенная ограниченность и непрерывность функций
$f_{1,n}(\cdot) - c_{1,n}$. Но тогда по теореме Арцела из
последовательности $f_{1,n}(\cdot) - c_{1,n}$ можно выбрать
подпоследовательность, равномерно сходящуюся на $D$ к выпуклой
функции $f_{1}(\cdot)$. Соответствующая подпоследовательность
последовательности  $f_{2,n}(\cdot) - c_{1,n}$ будет стремиться к
выпуклой функции $f_2(\cdot)$, что означает, что $f(\cdot)$ есть
ПРВ функция.

Пусть условия теоремы выполняются, но $f(\cdot)$ не есть ПРВ
функция. Проделаем следующую процедуру. Путем разбиения множества
$D$ на выпуклые подобласти можно выделить ту подобласть, где
функции $f_{1,n}(\cdot)$ имеют предельное бесконечное значение
вариации производной вдоль некоторых отрезков этой подобласти при
$n \rightarrow \infty$. Действительно, в противном случае из
последовательности функции $f_{1,n}(\cdot) - c_{1,n}$ можно было
бы выбрать сходящуюся подпоследовательность, и $f(\cdot)$ была бы
ПРВ функцией.

Далее разбиваем выделенную подобласть на меньшие области и опять
выделяем ту, где вариация производной функций $f_{1,n}(\cdot)$
вдоль некоторых отрезков неограничена при $n \rightarrow \infty$.
В итоге определяем точку $M$, в произвольной окрестности которой
вариация производной функций $f_{1,n}(\cdot)$ вдоль некоторых
отрезков неограничена при $n \rightarrow \infty$. Без ограничения
общности можно считать, что $M-$ внутренняя точка множества
$\bar{D},$ так как все получаемые в процессе применения алгоритма
функции $-$ равномерно липшицевы и могут быть распространены во
вне множества $\bar{D},$

Берем произвольную окрестность точки $M$ и разбиваем ее на
конечное число секторов. Выбираем произвольный из них, где
вариация производной функций $f_{1,n}(\cdot)$ вдоль некоторых
отрезков неограничена при $n \rightarrow \infty$. Далее выбранный
сектор разбиваем на конечное число секторов и опять выбираем тот
из низ, где вариация производной функций $f_{1,n}(\cdot)$ вдоль
некоторых отрезков неограничена при $n \rightarrow \infty$ и т.д.
Множество выбранных секторов стягивается к некоторому направлению,
определяемому единичным вектором $l$ с вершиной в точке $M$.
Очевидно, что в произвольном секторе $K$ с вершиной с точке $M$,
содержащем вектор $\alpha l$ в $\mbox{int} \, K$, $\alpha
>0$, вариация производной функций $f_{1,n}(\cdot)$ вдоль некоторых
отрезков неограничена при $n \rightarrow \infty$.

Возможны два случая:

a) вариация производных функций $f_{1,n}(\cdot)$ по направлению
$l$ неограничена при $n \rightarrow \infty$;

б) вариация производных функций $f_{1,n}(\cdot)$ по направлению
$\eta,$ перпендикулярном направлению $l$, неограничена при $n
\rightarrow \infty$ .

Сказанное можно перефразировать следующим образом, а именно: сумма
вариаций производных выпуклых двугранных углов функции
$f_{n}(\cdot)$ вдоль указанного направления неограничена при $n
\rightarrow \infty$.

Рассмотрим случай а). Возьмем произвольный сектор $K$, содержащий
вектор $\alpha l$ в $\mbox{int} \, K$, $\alpha >0$. Будем
рассматривать выпуклые двугранные углы функций $f_{1,n}(\cdot)$ из
конуса $K$ для всех $n$.

За счет равномерной липшицевости по $n$ всех двугранных углов
функций $f_{n}(\cdot)$ вариации производных по направлению этих
двугранных углов равномерно непрерывны относительно направления и
$n$.

Для каждого выпуклого $k-$ ого  двугранного  функции
$f_{n}(\cdot)$ выделим отрезок $v_{k,n}$, вариация производной
вдоль которого для $k-$ ого двугранного угла максимальна и равна
$a_{k,n}$. Ясно, что отрезок $v_{k,n}$ должен быть перпендикулярен
проекции на плоскость $XOY$ линии раздела двух граней $k-$ ого
двугранного угла.

Пусть угол наклона отрезков $v_{k,n}$ с направлением $l$ не
превосходит $\pi / 2- \delta$ для некоторого $ \delta>0.$

Путем разбиения сектора $K$ на меньшие секторы, стягивающиеся к
вектору $\alpha l$ и точку $M$, и рассмотрения в каждом из них
своей группы отрезков $v_{k,n}$ для всех значений $k$ и $n$, можно
выделить одну или несколько групп указанных отрезков, каждую из
которых можно пересечь кривой $r(\cdot) \in \wp(D),$ образующей в
точке пересечения с отрезками $v_{k,n}$ угол, не превосходящий
$\pi /2 -\delta, \delta>0.$ Поскольку сектор $K$ произвольный,
содержащий вектор $\alpha l$, то можно рассматривать такие кривые,
для которых $r'(t) \rightarrow -l,$ когда $r(t) \rightarrow M.$
Сама кривая $r(\cdot)$ будет включать в себя отрезки, близкие к
отрезкам $v_{k,n}$.

Если для рассматриваемого случая подгруппа отрезков $\{ v_{k,n}
\}$ существует только одна, то вдоль найденной кривой $ r(\cdot)
\in \wp(D)$ вариация производной суммы выпуклых двугранных углов
  стремится к бесконечности при $n \rightarrow +\infty$.

Кривая $r(\cdot)$, как упоминалось, строится таким образом, чтобы
она включала отрезки, близкие к отрезкам $\{ v_{k,n} \}$. Так как
при выполнении неравенства (\ref{difconv3}) выполняется
неравенство (\ref{difconv5}), а мы нашли кривую $r(\cdot)$, вдоль
которой сумма вариаций производных двугранных углов бесконечна, то
из (\ref{difconv5}) следует, что неограничена вдоль $r(\cdot)$
вариация производной функции $\Phi(\cdot)$. Приходим к
противоречию.

Кроме того, возможен случай, когда у нас есть несколько групп
отрезков $\{ v_{k,n} \}_i,$ для каждой из которых найдется кривая
$r_i(\cdot) \in \wp(D)$, что
$$
\vee (\Phi_n';0,T_{r_i})= c_i, \,\,\,  r_i' (t) \rightarrow_{t
\rightarrow T_{r_i}} -l,
$$
где $T_{r_i}-$  параметр кривой $r_i(\cdot)$ при естественной
параметризации в точке $M$,  а также
$$
 \sum_i \, c_i = \infty.
$$
Тогда кривую $r(\cdot) \in \wp(D),$ вдоль которой сумма вариаций
производных двугранных углов  стремится к бесконечности при $n
\rightarrow +\infty$, будем строить следующим образом.

Кривая $r(\cdot)$ должна содержать достаточное количество $k_i$
отрезков из каждой группы отрезков $\{ v_{k,n} \}_i,$ (либо
близких к ним), чтобы
$$
\vee (\Phi'_n ;t_{r_{i}},t_{r_{i+1}}) =  c_i - \mu_i,
$$
где $t_{r_i}>0-$ значения параметра $t$ для $i$-ой группы отрезков
при естественной параметризации кривой $r_i(\cdot)$, $\mu_i < c_i
-$ малые положительные числа, для которых
$$
\sum_i \, \mu_i < \infty.
$$

Нетрудно видеть, что всегда такую кривую $r(\cdot)$ построить
можно. Она будет состоять из набора  кривых $r_i(\cdot)$.   Для
этого надо осуществить плавный переход от одной кривой
$r_i(\cdot)$ к кривой $r_{i+1}(\cdot),$ не выходя из множества
$\wp(D)$. Поскольку $r'_i(t) \rightarrow -l$ при $t \rightarrow
T_{r_i}$ для вех $i$, то подобная процедура осуществима всегда.

Но тогда
$$
\vee (\Phi'_n ; 0,T_r) \geq \sum_i  \vee
(\Phi'_n;t_{r_{i}},t_{r_{i+1}}) =
$$
$$
=\sum_i(c_i-\mu_i)= \sum_i c_i -\sum_i \mu_i  =\infty.
$$

Но тогда, как следует из (\ref{difconv5}), нарушается неравенство
(\ref{difconv3}), которое по предположению достаточности условия
теоремы является верным. Опять приходим к противоречию.

Если вариация производной суммы выпуклых двугранных углов функции
$f_n(\cdot)$ конечна вдоль направления, определяемого вектором
$l$, при любом $n$, то для случая неограниченности при $n
\rightarrow \infty$ вариации производной функции $f_{1,n}(\cdot)$
в произвольно малом секторе с вершиной $M$, содержащем вектор
$\alpha l$, $\alpha >0$, следует, что вариация производной суммы
выпуклых двугранных углов функции $f_n(\cdot)$ бесконечна при $n
\rightarrow \infty$  вдоль направления $\eta $,

Случай б). Все отрезки $v_{k}$ можно разбить на такие группы $\{ m
\}$ отрезков, которые можно пересечь кривой $r_{m}(\cdot) \in
\wp(D),$ для которой
$$
        r'_{m} (\tau) \rightarrow_{\tau \rightarrow T_{r_{m}}} -l,
$$
где $T_{r_{m}}-$ есть параметр кривой $r_{m}(\cdot)$ при
естественной параметризации в точке $M,$  и кривизна кривой $r_{m}
(\cdot)$ стремится к бесконечности при $\tau \rightarrow
T_{r_{m}}.$ Кривая $r_{m}(\cdot)$ пересекает свою группу отрезков
под острыми углами $\alpha_{k_m}$ в точках $\tau_{k_m}$, причем
$\alpha_{k_m} \rightarrow \pi / 2$ при $\tau_{k_m} \rightarrow
T_{r_{m}}$. Ясно, что сказанное всегда выполнимо путем разбиения
множества всех отрезков $v_{k}$ на подмножества с требуемыми
свойствами.

Кроме того, углы $\alpha_{k_m},$ кривые $r_{m}(\cdot)$ и группы
отрезков $\{ v_{k} \}_m$ можно выбрать такими, чтобы предел по $m$
суммы вариаций производных функций $\Phi'_m(\cdot)$ вдоль отрезков
кривых $r_{m}(\cdot)$  был равен бесконечности. В противном случае
функции $f_n(\cdot)$ имели бы ограниченную вариацию вдоль
направления $\eta$ при $n \rightarrow +\infty$ (см. замечание).

Построение кривых $r_{m}(\cdot)$ с неограниченно увеличивающейся
кривизной в точке $M$, для которой
$$
   \lim_{m \rightarrow \infty} \sum_{m}
   \vee (\Phi'_m ; T_{k_m},T_{{k_m+1}})=\infty,
$$
осуществляется аналогичным способом, как и в случае a). Для этого
надо построить кривую $r_{m}(\cdot) \in \wp(D)$ с описанными выше
свойствами, состоящую из достаточного  количества $k_{m}$ отрезков
$\{v_{k}\}$ (либо близких к ним), чтобы
$$
\sum_{k_m} \vee (\Phi'_m ;T_{k_{m}},T_{k_{m}+1}) = c_{m},
$$
и
$$
\sum_m  c_{m} =\infty,
$$
$[ T_{k_{m}}, T_{k_{m}+1} ]$- значение параметра $t$ для отрезка
кривой $r_m (\cdot)$ при ее естественной параметризации. Такие
кривые $r_{m}(\cdot)$ всегда можно построить. При $m \rightarrow
\infty$ кривые $r_{m}$ будут пересекать под острыми углами все
большее число указанных отрезков из произвольно малого сектора,
содержащем вектор $\alpha l$, с вершиной в точке $M$. Кривизны
кривых $r_{m}$ вблизи точки $M$ неограниченно увеличиваются при $m
\rightarrow \infty$.

Но тогда из последовательности кривых $r_m(\cdot)$ можно составить
кривую $r(\cdot) \in \wp(D) $ с конечной суммой вариаций
производной на отрезках $[ T_{k_m},T_{k_m+1} ]$, что для
$\Phi(t)=f(r(t)) $
$$
\sum_{k_m} \vee (\Phi'; T_{k_m},T_{k_m+1}) = \infty.
$$
Отсюда и из (\ref{difconv5}) приходим к противоречию с
(\ref{difconv3}).

Итак, доказано, что при выполнении условия теоремы, сумма вариаций
производных выпуклых двугранных углов функции $f_n(\cdot)$ вдоль
любого отрезка области $D$ при $n \rightarrow \infty$ ограничена
сверху константой, независящей от $n$. Отсюда, как отмечалось
выше, следует, что $f(\cdot)-$ ПРВ функция.

Итак, теорема \ref{difconvexthm1} доказана. $\Box$

\begin{rem} Рассуждения с выбором углов $\alpha_{k_m}$ и кривых
$r_m(\cdot)$ аналогичны следующим.

Пусть имеем расходящийся ряд
$$
\sum_i \, a_i = \infty, \,\,\,\, a_i >0 \,\,\,\,\, \forall i.
$$
Всегда можно выбрать монотонно убывающую по $i$ последовательность
$\{ \beta_i \},  \beta_i \rightarrow_{i \rightarrow \infty} 0,$
чтобы
$$
\lim_{m \rightarrow \infty} \sum_{i=1}^{m} \beta_i \, a_i =
\infty.
$$
Здесь $a_i$ является аналогом вариации производной двугранного
угла вдоль отрезка $v_i$, а $\beta_i$ - аналог косинуса угла,
образуемого кривой $r_i$ с этим отрезком в точке пересечения.
\end{rem}

\vspace{1cm}

\subsection{Геометрическая интерпретация теоремы 1}

\vspace{0.5cm}

Перефразируем теорему \ref{difconvexthm1}, придав ей более
геометрический характер. Для этого введем понятие поворота кривой
$r(\cdot)$  на графике $\Gamma_f = \{(x,y,z) \in \mathbb{R}^3 \mid
z = f(x , y)\}.$

Рассмотрим на $\Gamma_f$ кривую $ R(t)=(r(t),f(r(t))),$ где
$r(\cdot) \in \wp(D).$  Так как функция $f(\cdot,\cdot)$ есть
липшицевая, то п.в. на $[0,T_r]$ существует производная
$R'(\cdot),$ которую обозначим через $\tau(\cdot)=R'(\cdot).$

{\bf Определение 4.4.1.} {\em  Поворотом кривой $R(\cdot)$ на
многообразии $\Gamma_f$ назовем величину
$$
sup_{ \{t_i\} \subset N_r} \,\, \sum_i \Vert \tau(t_i)/ \Vert \tau
(t_i) \Vert -  \tau(t_{i-1})/ \Vert \tau (t_{i-1}) \Vert \Vert =
O_r.
$$
}

Таким образом, поворот $O_r$ кривой $R(\cdot)$ есть верхняя грань
суммы углов между касательными $\tau(t)$  для $t \in [0,T_r].$
Нетрудно видеть, что для плоской гладкой кривой, параметризованной
естественным образом, величина $ O_r$ равна интегралу
$$
\int^{T_r}_0 \mid k(s) \mid ds,
$$
где $k(s)$ - кривизна рассматриваемой кривой $ r(\cdot)$ в точке
$s \in [0,T_r],$ т.е. совпадает с обычным определением поворота
кривой в точке \cite{pogorelov1} .

\begin{thm}
Для того, чтобы произвольная липшицевая функция \\
$z \rightarrow f(z) :D \rightarrow \mathbb{R}$ была ПРВ функцией
на выпуклом компактном множестве $ D \in \mathbb{R}^2,$ необходимо
и достаточно, чтобы для любой кривой $r(\cdot) \in \wp(D)$ и любой
ее системы подмножеств существовали константы $c_3(D,f), c_4(D,f)
>0$ такие, что поворот кривой $R(\cdot)$ на $\Gamma_f$ ограничен
сверху неравенством, т.е.
\begin{equation} O_r \leq c_3(D,f) + c_4(D,f) \vee(r';0,T_r)  \;\;
\forall r \in \wp(D). \label{difconv6} \end{equation}
\end{thm}

\begin{rem} Условие теоремы означает следующее. Сумма поворотов любых
участков кривой $R(\cdot)$ на $\Gamma_f$ ограничена сверху
неравенством (\ref{difconv6}). Здесь также, как и выше, приходится
оговаривать правило подсчета поворота кривой $R(\cdot)$, так как
поворот всей кривой $R(t)$, $t \in [0, T_r] $, может оказаться
бесконечной.
\end{rem}

{\bf Доказательство. }{\bf Необходимость}. Пусть $ f(\cdot,\cdot)$
есть ПРВ функция. Покажем, что тогда справедливо неравенство
(\ref{difconv6}). Воспользуемся неравенством, вытекающим из
неравенства треугольника,
$$
\Vert \tau(t_i) / \Vert \tau(t_i) \Vert  - \tau (t_{i-1}) / \Vert
\tau(t_{i-1} \Vert \Vert \leq
$$
$$
\leq \Vert r'(t_i) / \sqrt{ 1+f'^2_t (r(t_i))}  - r'(t_{i-1}) /
\sqrt{ 1+f'^2_t (r(t_{i-1}))} \Vert +
$$
$$
+ \mid f'_t(r(t_i)) /  \sqrt{ 1+f'^2_t (r(t_i))}- f'_t(r(t_{i-1}))
/ \sqrt{ 1+f'^2_t (r(t_{i-1}))} \mid .
$$
Так как $1 \leq  \sqrt{1+f'^2_t (r(t_i))} \leq \sqrt{1+L^2}$ для
всех $ t_i \in [0,T_r]$ , то очевидно, существует такое $c_3 >1,$
для которого верно неравенство

\begin{equation} \Vert r'(t_i) / \sqrt{1+f'^2_t (r(t_i))}-
r'(t_{i-1}) / \sqrt{1+f'^2_t (r(t_{ i-1}))} \Vert \leq c_3 \Vert
r'(t_i) - r'(t_{i-1}) \Vert. \label{difconv7} \end{equation}

Из свойств функции $\theta(x)= x / \sqrt{ 1+x^2}$ следует
неравенство
$$
\mid f'_t (r(t_i)) / \sqrt{ 1+f'^2_t (r(t_i))} - f'_t (r(t_{i-1}))
/ \sqrt{1+f'^2_t (r(t_{i-1}))} \mid \leq
$$
\begin{equation}
\leq \mid f'_t (r(t_i)) - f'_t (r(t_{i-1})) \mid .
\label{difconv8}
\end{equation}

Из (\ref{difconv7}) и (\ref{difconv8}) имеем \begin{equation}
\sup_{\{t_i \} \in N_r } \,\, \sum_i \, \Vert \tau(t_i) / \Vert
\tau(t_i) \Vert - \tau (t_{i-1}) / \Vert \tau(t_{i-1}) \Vert \Vert
\leq c_3 (\vee (  r'; 0 ,T_r) + \vee (\Phi' ;  0,T_r) ).
\label{difconv9} \end{equation}

Так как по условию $ f(\cdot,\cdot)-$  ПРВ функция, то согласно
теореме \ref{difconvexthm1}
$$
\vee (\Phi'; 0,T_r) \leq  c_1(D,f) + c_2(D,f) \vee(r';0,T_r),
$$
откуда с учетом (\ref{difconv9}) следует неравенство
(\ref{difconv6}). Необходимость доказана.

{\bf Достаточность}.  Пусть справедливо неравенство
(\ref{difconv6}). Покажем, что $f(\cdot,\cdot)$ - ПРВ функция.
Воспользуемся неравенством
$$
\Vert \tau(t_i) / \Vert \tau( t_i) \Vert  - \tau(t_{i-1}) /  \Vert
\tau (t_{i-1}) \Vert \Vert \geq \mid  f'_t (r(t_i)) / \sqrt{1+f'_t
(r(t_i))} -
$$
\begin{equation} - f'_t (r(t_{i-1})) / \sqrt {1+f'^2_t
(r(t_{i-1}))} \label{difconv10} \end{equation} Из свойств функции
$\theta(x) = x / \sqrt{1+x^2}$ и из $\Vert f'(z) \Vert \leq L$ для
всех $z \in D,$ где производная существует, следует существование
константы $ c_4(L) > 0,$ для которой
$$
\mid f'_t (r(t_i)) / \sqrt{1+f'^2_t (r(t_i))} - f'_t (r(t_{i-1}))
/ \sqrt{1+f'^2_t(r(t_{i-1}))} \geq
$$
$$
\geq c_4 \mid f'_t (r(t_i)) - f'_t (r(t_{i-1})) \mid,
$$
откуда с учетом (\ref{difconv10}) имеем
$$
c_2(D,f) + c_3(D,f) \vee(r';0,T_r) \geq  \sup_{ \{ t_i \} \subset
N_r} \sum_i \Vert \tau(t_i) / \Vert \tau( t_i) \Vert  -
\tau(t_{i-1}) / \Vert \tau (t_{i-1}) \Vert \Vert \geq
$$
$$
\geq c_4 \vee (\Phi';  0,T_r) .
$$
Из теоремы \ref{difconvexthm1} следует, что $f(\cdot)$ - ПРВ
функция. Достаточность доказана.  $\Box$

Возьмем произвольную кривую $r(\cdot) \in \varrho (D)$. Поскольку
вариация производной $\vee (r'; 0,T_r) $ ограничена сверху для
любой кривой $r(\cdot) \in \varrho (D)$, то из Теоремы 1 и Теоремы
2 вытекают следствия.

\begin{cor}
Для того, чтобы липшицевая функция $z \rightarrow f(z):D
\rightarrow \mathbb{R}$ была представима в виде разности выпуклых
функций, необходимо, чтобы для любой кривой $r \in \varrho (D)$ и
для некоторой константы $c_5(D,f)>0$ выполнялось неравенство
$$
 (\exists c(D,f),) (\forall r \in \varrho(D, f)) \,\,
 \vee (\Phi'; 0,T_r) <  c_5 (D, f),
$$
где $\Phi(t)=f(r(t)) \;\;\;\; \forall t \in [0,T_r]$.
\label{difconvexconl1}
\end{cor}

\begin{cor}
Для того, чтобы липшицевая функция $z \rightarrow f(z):D
\rightarrow \mathbb{R}$ была представима в виде разности выпуклых
функций, необходимо, чтобы для любой кривой $r \in \varrho (D)$ и
для некоторой константы $c_6(D,f)>0$ поворот кривой $R(\cdot)$ на
$\Gamma_f$ ограничен сверху неравенством, т.е. $$ O_r \leq
c_6(D,f)\;\; \forall r \in \varrho(D). $$
\end{cor}

В статье \cite{veselyzajicek} приведен пример, подтверждающий
несправедливость достаточности утверждений Следствия 1 и Следствия
2 \cite{proudconvex2}. Фактически авторы статьи
\cite{veselyzajicek} привели пример, показывающий, что кривых
множества $ \varrho(D) $ недостаточно, чтобы проверить
представимость функции в виде разности выпуклых функций. Класс
кривых $\wp(D) $ значительно шире класса $\varrho(D)$.

\newpage
\vspace{1cm}

\section{НЕОБХОДИМЫЕ И ДОСТАТОЧНЫЕ УСЛОВИЯ ПРЕДСТАВИМОСТИ
ФУНКЦИИ МНОГИХ ПЕРЕМЕННЫХ В ВИДЕ РАЗНОСТИ ВЫПУКЛЫХ ФУНКЦИЙ}

\vspace{0.5cm}

Результаты, полученные в этом параграфе,  являются продолжением
работ автора \cite{proudconvex1}, \cite{proudconvex2}, где
приведены необходимые и достаточные условия представимости
произвольной функции двух переменных в виде разности выпуклых
функций. Результат распространяется на функции от произвольного
количества аргументов, определенные в конечномерном евклидовом
пространстве. Геометрическая интерпретация этих условий, как и в
работе \cite{proudconvex2}, также приведена. Описан алгоритм
такого представления, применение которого есть последовательность
равномерно сходящихся выпуклых функций.

\noindent

\vspace{0.5cm}

\subsection{Введение}

\vspace{0.5cm}

Автор пришел к проблеме об условиях представимости функции в виде
разности выпуклых в процессе изучения теории квазидифференцируемых
функций, развитой специалистами по негладкой оптимизации
\cite{demvas}, \cite{clark}.

Эта проблема была впервые сформулирована академиком
А.Д.Александровым в статье \cite{aleksandrov1} и исследована
многими российскими и зарубежными математиками
\cite{aleksandrov1}, \cite{aleksandrov2}, \cite{zalgaller1},
\cite{proudconvex1}, \cite{proudconvex2}, \cite{ginchev},
\cite{Hartman}, \cite{hiriarturruty}, \cite{veselyzajicek}.
Решение этой проблемы важно для геометров и математиков,
занимающихся теорией управления и оптимизацией \cite{strecal2},
например, для построения квазидифференциального исчисления
\cite{demvas}, \cite{demrub}. Проблема оказалась довольно не
простой.  Решение ее не было получено с 40-ых годов 20-ого
столетия, хотя многие предпринимали попытку ее решить. В работе
\cite{proudconvex2} дана предыстория вопроса.

Согласно терминологии А.Д.Александрова под многогранной
кусочно-линейной функцией с конечным числом граней будем понимать
такую функцию, график которой состоит из конечного числа частей
плоскостей (гиперплоскостей), которые называются гранями.

В статье  \cite{proudconvex1} даны необходимые и достаточные
условия представимости произвольной липшицевой положительно
однородной (п.о.) функции трех переменных в виде разности выпуклых
функций.

В статье \cite{proudconvex2} рассматривается произвольная
липшицевая функцию $f(\cdot)$ с константой Липшица $L$ от двух
переменных $(x,y) \rightarrow f(x,y):D \rightarrow \mathbb{R},$
где $D$ есть выпуклое открытое ограниченное множество в
$\mathbb{R}^2$ с непустой внутренностью, так что его замыкание
$\bar{D}$ - компакт. Там же приводится алгоритм такого
представления в виде последовательности выпуклых многогранных
функций, а также находятся необходимые и достаточные условия
поточечной сходимости построенной последовательности функций.

Анализ \cite{veselyzajicek} результатов, полученных в
\cite{proudconvex1}, \cite{proudconvex2}, показал, что неравенства
для вариаций производной исходной  функции, определенной на кривых
из класса кривых, указанных в \cite{proudconvex1},
\cite{proudconvex2}, задают только необходимые, но не достаточные
условия представимости функции в виде разности выпуклых. Но если
расширить класс кривых, то, практически не меняя доказательства
теорем статей \cite{proudconvex1}, \cite{proudconvex2}, можно
сформулировать необходимые и достаточные условия представимости
функции в виде разности выпуклых функций.

Пусть $\wp(D)$ - класс всех кривых на плоскости  $XOY$ на
множестве $D \in \mathbb{R}^2$, ограничивающих звездные области.

\begin{defi} Звёздная область $D$ относительно
точки $O -$ это подмножество евклидова пространства $\mathbb{R}^n$
такое, что отрезок, соединяющий любую точку области $D$  с точкой
$O$, целиком принадлежит этой области. Множество называется
звёздной областью, если существует точка, относительно которой это
подмножество звёздное.
\end{defi}

Вопрос заключается в том, хватит ли таких кривых для ответа на
вопрос о представимости функции в виде разности выпуклых?
Показывается, что таких кривых хватает, если изменить правило
подсчета вариации производной функции одной переменной,
определенной на таких кривых.

Параметризуем кривые $r \in \wp(D)$ естественным, или натуральным
образом, т.е. параметр $\tau$ точки $M$ на кривой $r(\cdot)$ равен
длине кривой между $M$ и начальной точкой. Обозначим такую кривую
как $r(t)$, $t \in [0, T_r]$. Кривая $r(\cdot) $ п.в. на $[0,
T_r]$ имеет касательную $r'(\cdot)$.

С помощью кривых  $r \in \wp (D)$ необходимые и достаточные
условия представимости функции $f: D \longrightarrow \mathbb{R}$ в
виде разности выпуклых функций записываются в следующем виде.

С помощью кривых  $r \in \wp (D)$ необходимые и достаточные
условия представимости функции $f(\cdot)$ в виде разности выпуклых
функций  могут быть записаны в следующем виде.

\begin{thm}
 Для того, чтобы в двумерном случае липшицевая функция $z \rightarrow f(z):D
\rightarrow \mathbb{R}$ была представима в виде разности выпуклых
функций (была ПРВ функцией), необходимо и достаточно, чтобы для
любой кривой $r \in \wp (D)$ и любой ее системы подмножеств $[T_i,
T_{i+1} ] \subset [0, T_r] $, $i \in 1:m$, выполнялось неравенство
$$(\exists c_1(D,f), c_2(D,f)>0) (\forall r \in \wp(D)),$$
\begin{equation}
 \sum_1^m \vee (\Phi'; T_i,T_{i+1})<  c_1(D,f) + c_2(D,f) \sum_1^m\vee (r';
 T_i,T_{i+1}),
 \label{difconv3}
\end{equation} где $\Phi(t)=f(r(t)) \;\;\;\; \forall t \in
[0,T_r]$. \label{difconvexthm1}
\end{thm}

\begin{rem}
Впервые подобные условия для вариации функции на кривых были
приведены  в \cite{proudconvex1}.
\end{rem}

Применяется  специальный алгоритм представления функции $f(\cdot)$
в виде разности выпуклых функций. В результате применения этого
алгоритма получается конечная или бесконечная последовательность
выпуклых функций, равномерно сходящаяся на $D$  к выпуклым
функциям, разность которых есть исходная функция $f(\cdot)$, если
условия теоремы 1 выполняются.

Для представления функции $f(\cdot)$ в виде разности выпуклых
функций используются две операции, в результате которых получается
конечное или счетное число выпуклых многогранных кусочно-линейных
функций, определенных на $D$.

{\em Первая операция} - это приближение функции $f(\cdot)$
многогранной кусочно-линейной функцией $f_N(\cdot)$ с конечным
числом граней. График функции $f_N(\cdot)$ состоит из конечного
числа частей плоскостей, которые строятся по разбиению области
$D_N \subset D$ на подобласти в виде треугольников с непустыми
внутренностями. Диаметры подобластей равномерно стремятся к нулю,
$\rho_H (D,D_N) \rar_N 0$ при $N \rar \infty$, где $\rho_H - $
метрика Хаусдорфа \cite{demvas}.

{\em Вторая операция} - это представление функции $f_N(\cdot)$ в
виде разности выпуклых многогранных кусочно-линейных функций
$f_{1,N}(\cdot):D \rightarrow \mathbb{R}$ и $f_{2,N}(\cdot):D
\rightarrow \mathbb{R}$ согласно алгоритму, описанному ниже. За
счет равномерной зависимости от $N$ коэффициента липшицевости
функции $f_{1,N}(\cdot), f_{2,N}(\cdot)$ всегда можно
распространить на всю область $D$.

Доказывается, что если выполняются условия теоремы 1, то из
последовательностей $f_{1,N}(\cdot)- c_{1,N}$ и $f_{2,N}(\cdot) -
c_{1,N}$, где $c_{1,N}=f_{1,N}(a)$, $a - $ произвольная внутренняя
точка области $D$, можно выбрать сходящиеся подпоследовательности.

Метод представления многогранной функции в виде разности выпуклых
подобен методу, предложенному  А.Д.Александровым в
\cite{aleksandrov1} при исследовании возможности представления
многогранной функции с конечным числом граней в виде разности
выпуклых. В нашем случае в процессе применения алгоритма число
граней неограниченно увеличивается.

Также в \cite{proudconvex2} дана геометрическая интерпретация
полученного результата через поворот кривой $R(t)=(r(t),f(r(t))),$
где $r(\cdot) \in \wp(D).$ Оказывается, что  липшицевая функция
$f(\cdot)$ от двух переменных представима в виде разности выпуклых
тогда и только тогда, когда для поворота кривых
$R(t)=(r(t),f(r(t)))$ справедлива оценка, аналогичная неравенству
(\ref{difconv3}) для всех $r(\cdot) \in \wp(D)$.

Вопрос об условиях представления функции в виде разности выпуклых
интересен для специалистов многих специальностей. Авторы статей на
эту тему стараются получить необходимые и достаточные условия
представимости, которые можно легко проверить. Некоторые такие
условия представлены ниже.

\vspace{1cm}

\subsection{Многомерный случай}

Рассмотрим многомерный случай $x \in \mathbb{R}^n$, $n>2$. Пусть
$D-$ выпуклое открытое множество в $\mathbb{R}^n$ с непустой
внутренностью, замыкание которого есть компакт. Нас будут
интересовать необходимые и достаточные условия представимости
функции $f(\cdot): D \rightarrow \mathbb{R}$ в виде разности
выпуклых.

В работе \cite{ginchev} приводятся необходимые и достаточные
условия представимости функции в виде разности выпуклых на
произвольном выпуклом множестве $\Omega$ из линейного векторного
пространства $X$.

\begin{thm} Пусть $X$ - линейное пространство, $\Omega$ -
выпуклое множество в $X$ и $f : \Omega \rar R$ - произвольная
функция. $f$ является разностью двух выпуклых функций, если и
только если найдётся (конечное или бесконечное) множество индексов
$I$ и множество выпуклых функций $h_i : \Omega \rar R$, $i \in I$,
таких, что $\sum_{i \in I} h_i (x)$ существует и конечно в
$\Omega$ и) для любой пары точек $a,b \in  \Omega$ существует
множество $J \subset I$ такое, что $f + \sum_{i \in J} h_i$
является выпуклой функцией на сегменте $[a,b]$.
\end{thm}

В приведенной теореме в отличие от условий представимости функции
в виде разности выпуклых, указанных в данной статье, требуется
существование выпуклых функций $h_i, i \in I,$ сумма которых
конечна на $\Omega $. Кроме того, для каждой пары точек $a, b  \in
\Omega$ выбирается подмножество $J \subset I$, зависящее от $a,b$,
такая что сумма $f + \sum_{i \in J} h_i$ также выпуклая. В
теоремах данной статьи не требуется существование выпуклых
функций, а наоборот, доказывается, что при выполнении условий
теорем выпуклые функции, разность которых есть исходная функция,
существуют. Приводится алгоритм их построения. В этом заключается
существенная разница между статьей автора и статьей
\cite{ginchev}. Статьи \cite{ginchev}, \cite{veselyzajicek}
подтверждают большой интерес к этой тематике ученых со всего мира
вплоть до наших дней.

Введем множество замкнутых кривых  $r(\cdot) \in \tilde{\wp}(D)$,
принадлежащих $D \subset \mathbb{R}^n$. Кривую $r(\cdot)$
параметризуем естественным образом, т. е.  $t-$ натуральный
параметр, равный длине кривой $r$ от начальной точки с параметром
$t=0$  до рассматриваемой точки c параметром $t$. Отрезок значений
параметра $t$ обозначим через $[0,T_r].$ Поскольку выполняется
неравенство
$$
\Vert  r(t_1) - r(t_2) \Vert  \le  \mid t_1 - t_2 \mid ,
$$
то кривая $r(\cdot)$ почти всюду дифференцируема на $[0,T_r]$.
Множество точек дифференцируемости кривой $r(\cdot)$ обозначим
через $N_r$.

Множество $\tilde{\wp}(D)$ будет состоять из всех  кривых $r(t),t
\in [0, T_r],$ которые имеют взаимно однозначные проекции на одну
из плоскостей $\Pi$  (для каждой кривой своя плоскость $\Pi$),
проходящую через начало координат,  в виде кривых из множества
${\wp}(Pr(D))$, где $Pr(D)-$ проекция области $D$ на ту плоскость
$\Pi$, на которую проектируется кривая $r(\cdot)$. Определение
множества кривых ${\wp}(Pr(D))$ было дано выше для двумерного
случая в предыдущем параграфе.

Для $r(\cdot) \in \tilde\wp(D) $ выделим те ее отрезки параметров
$[T_i, T_{i+1}] \subset [0,T(r)]$, где производная $r'(t)$  и
вектор $r(t)$ составляют с плоскостью $\Pi$ угол $\delta $: $0 <
\delta \leq \frac{\pi }{4} $. Обозначим объединение таких участков
через $\Lambda_{\Pi}(r)$.

Если, например, $D=B^3_R(0)=\{ x \in \mathbb{R}^3 \mid \parallel x
\parallel \leq R \}-$ шар радиуса $R$ с центром в начале
координат трехмерного пространства, то кривые, получающиеся в
результате пересечения произвольных плоскостей с концентрическими
сферами $S^2_{\varepsilon}(0)=\{ x \in \mathbb{R}^3 \mid \parallel
x \parallel=\varepsilon \}$ радиуса $\varepsilon \leq R$, будут
принадлежать множеству $\tilde{\wp}(D)$.

С помощью кривых  $r \in \tilde{\wp} (D)$ необходимые и
достаточные условия представимости функции $f: D \longrightarrow
\mathbb{R}$ в виде разности выпуклых функций записываются в
следующем виде.

\begin{thm} Для того, чтобы липшицевая функция $z \rightarrow
f(z):\mathbb{R}^n \rightarrow \mathbb{R}$ была представима в виде
разности выпуклых функций (была ПРВ функцией) на $D$, необходимо и
достаточно,  чтобы для любой плоскости $\Pi$, $0 \in \Pi$, и любой
кривой $r(\cdot) \in  \tilde\wp(D)$, параметризованной
естественным образом с параметром $t \in [0,T], \, T=T(r),$ и
любого подмножества параметров $[T_i, T_{i+1} ] \subset
\Lambda_{\Pi}(r) $
$$
 (\exists c(D,f)>0) (\forall r \in \tilde{\wp}(D)) \;\;\;\;
$$
\be
 \sum_1^m \vee (\Phi'; T_i,T_{i+1})<  c(D,f) (1+  \sum_1^m\vee (r';
 T_i,T_{i+1})),
\label{difconvexthm1a} \ee где $\Phi(t)=f(r(t)) \;\;\;\; \forall t
\in [0,T_r]$. \label{difconvmultexthm2}
\end{thm}

\begin{rem}
Из теоремы \ref{difconvmultexthm2} следует теорема для двумерного
случая, поскольку в двумерном случае  векторы $r(\cdot),
r'(\cdot)$ лежат на плоскости $\mathbb{R}^2$, которой принадлежит
множество $D$. Вариация производной $r'(\cdot)$ определяется
также, как длина кривой (см. \cite{pogorelov1}), и ее точное
определение дается ниже.
\end{rem}

\begin{rem}
В следующем параграфе будут рассматриваться не произвольные
плоскости  $\Pi$, а только координатные плоскости $\Pi_i$,
образованные осями координат.
\end{rem}

\begin{rem}
В дальнейшем неравенство (\ref{difconvexthm1a}) будем записывать в
в виде
$$
 \vee (\Phi'; t \in \Lambda_{\Pi}(r))<  c(D,f) (1+  \sum_1^m\vee (r';
 t \in \Lambda_{\Pi}(r)))
$$
\end{rem}

\vspace{0.5cm}

{\bf Необходимость}. Доказательство необходимости начнем с
описания алгоритма представимости функции в виде разности
выпуклых.

\vspace{0,5cm}

{\bf ОПИСАНИЕ АЛГОРИТМА}

1. Произведем разбиение выпуклого компактного множества $D_N, D_N
\subset D, \, \mbox{int } D_N \neq \varnothing, \, \rho_H (D, D_N)
\rar_N 0,$ на выпуклые $n-$ мерные многогранники $G_k,$ $ k=1,2,
\dots, N,$ с непустыми внутренностями и $n+1$ вершинами, диаметры
которых равномерно стремятся к нулю при $N \rar \infty$. Строим по
каждому многограннику $G_k$ гиперплоскость $\pi_k(\cdot)$,
являющуюся графиком линейной функции, определенной на $D$,
значения которой равны значениям функции $f(\cdot)$ в вершинах
многогранника $G_k$. Функцию, график которой внутри каждого
многогранника $G_k$ совпадает с гиперплоскостью $\pi_k(\cdot)$, $
k \in \overline{1, N},$ обозначим через $f_N(\cdot):D \rightarrow
\mathbb{R}$. Назовем $f_N(\cdot)$ многогранной функцией.

2. Представляем функцию $f_N(\cdot)$ в виде разности выпуклых
согласно алгоритму, как это описано ниже.

Предварительно введем понятие {\em двугранного угла}. Будем
понимать под двугранным углом функцию, определенную на $D$, равную
максимуму или минимуму линейных функций, графиками которых
являются гиперплоскости $\pi_k(\cdot)$ и $\pi_{l}(\cdot)$,
построенные по соседним многогранникам $G_k$, имеющим общие грани
размерности $n-1$.

Рассмотрим все выпуклые двугранные углы, построение которых
описано выше. Доопределим эти двугранные углы на всю область $D$.
Просуммируем все такие выпуклые двугранные углы. В итоге получим
выпуклую многогранную функцию $f_{1,N}(\cdot):D \rightarrow
\mathbb{R}$. Доказывается \cite{aleksandrov1}, что разность
\begin{equation}
 f_{1,N}(\cdot) - f_N(\cdot) = f_{2,N}(\cdot)
\label{difconvmult1} \end{equation} есть также выпуклая
многогранная функция.

Действительно, для доказательства достаточно показать, что все
двугранные углы, части графиков которых размерности $n$
принадлежат графику функции $f_{1,N}(\cdot) - f_{N}(\cdot),$ есть
выпуклые. Для этого покажем, что любая точка, лежащая на проекции
на $D$ пересечения $\pi_{kl}(\cdot) = \pi_k(\cdot) \cap
\pi_{l}(\cdot)$ произвольных гиперплоскостей $\pi_k(\cdot)$ и
$\pi_{l}(\cdot)$, образующих график двугранного угла функции
$f_{1,N}(\cdot) - f_{N}(\cdot),$ имеет малую окрестность, где
функция $f_{1,N}(\cdot) - f_{N}(\cdot)$ выпуклая.

Если берем точку, в малой окрестности которой функция
$f_{N}(\cdot)$ линейная, то локальная выпуклость разности
$f_{1,N}(\cdot) - f_{N}(\cdot)$ очевидна. Пусть берем точку,
лежащую на проекции на $D$ множества $\pi_{kl}(\cdot)$ выпуклого
двугранного угла, часть графика которого принадлежит графику
функции $f_{N}(\cdot)$. Поскольку согласно алгоритму двугранный
угол, график которого образован гиперплоскостями $\pi_k(\cdot)$ и
$\pi_{l}(\cdot)$, входит в сумму выпуклых двугранных углов,
образующих функцию $f_{1,N}(\cdot)$, то опять разность
$f_{1,N}(\cdot) - f_{N}(\cdot)$ будет локально выпуклой в
окрестности рассматриваемой точки. Если же точка лежит на проекции
на $D$ множества $\pi_{kl}(\cdot)$ вогнутого двугранного угла,
часть графика которого принадлежит графику функции $f_{N}(\cdot)$,
то $-f_{N}(\cdot)$ $ -$ локально выпуклая в окрестности этой
точки, а поэтому разность $f_{1,N}(\cdot) - f_{N}(\cdot)$ снова
локально выпуклая в той же окрестности. Из локальной выпуклости
всех двугранных углов функции $f_{1,N}(\cdot) - f_{N}(\cdot)$
следует ее выпуклость на всем множестве $D$.

Покажем, что при выполнении условия теоремы
\ref{difconvmultexthm2} из последовательности функций
$f_{1,N}(\cdot) - c_{1,N}, $ где $c_{1,N}=f_{1,N}(a)$, $a - $
произвольная внутренняя точка области $D$,  можно выделить
подпоследовательность, равномерно сходящуюся на $D$ к выпуклой
функции $f_1(\cdot)$ при $N \rightarrow +\infty$. Тогда из
(\ref{difconvmult1}) будет следует, что подпоследовательность
функций $f_{2,N}(\cdot)-c_{1,N}$ также равномерно сходится к
выпуклой функции $f_2(\cdot)$. Для функций $f_1(\cdot)$ и
$f_2(\cdot)$ верно равенство \begin{equation}
 f_{1}(\cdot) - f_2(\cdot) = f(\cdot).
\label{difconvmult2} \end{equation}

В статье \cite{proudconvex2} показано, что данный алгоритм
приводит к равномерно сходящейся последовательности выпуклых
функций от одной переменной. В двумерном случае, как доказано там
же, при выполнении условий теоремы этот алгоритм также проводит к
равномерно сходящейся на $D$ последовательности выпуклых функций.

Перейдем к случаю $n > 2$ и покажем, что тот же алгоритм при
выполнении приведенной выше теоремы \ref{difconvmultexthm2} также
приводит к паре выпуклых функций на $D$, разность которых есть
исходная функция $f(\cdot)$.

Возьмем произвольную кривую $r(\cdot) \in \tilde{\wp}(D)$. Пусть
$$
    \Phi(t)=f(r(t)) \;\;\; \forall t \in [0,T_r].
$$

Покажем, что $\Phi(\cdot)$ - липшицевая с константой $L$.
Действительно, для любых $t_1,t_2 \in [0,T_r]$ имеем
$$
 \mid \Phi (t_1) - \Phi (t_2) \mid = \mid f(r(t_1)) - f(r(t_2)) \mid \le
 L \Vert  r(t_1) - r(t_2) \Vert  \le  L \mid t_1 - t_2 \mid.
$$
Поэтому \cite{kolmogorovfomin} $\Phi (\cdot)$ почти всюду (п.в.)
дифференцируемая на $[0,T_r]$. Кривая $r(\cdot)$,
параметризованная естественным образом, имеет п.в.  на $[0,T_r]$
касательную $r'(\cdot)$.

Множество точек дифференцируемости функции $\Phi (\cdot)$ и
$r(\cdot)$ на $[0,T_r]$  обозначим  через $N_r.$ Докажем, что если
существует константа $c(D) > 0$ такая, что для произвольной кривой
$r(\cdot) \in \tilde{\wp}(D)$ и произвольного подмножества из
$\Lambda_{\Pi}(r)$
\begin{equation}
    \bigvee (\Phi '; t \in \Lambda_{\Pi}(r))  <  c(D) (1+ \bigvee (r'; t \in \Lambda_{\Pi}(r)))        ,
\label {difconvmult3} \end{equation} то из последовательностей
функций $f_{1,N}(\cdot)-c_{1,N}$, $f_{2,N}(\cdot)-c_{1,N}$,  можно
выбрать подпоследовательности, равномерно на $D$ сходящиеся к
выпуклым функциям $f_{1}(\cdot)$, $f_{2}(\cdot)$ соответственно,
для которых верно равенство (\ref{difconvmult2}).

Доказательство будет основываться на леммах, приведенных ниже.

{\bf Определение 1.} {\em  Под вариацией кривой $r'(\cdot) \in
\tilde{\wp}(D)$ на множестве $\Lambda_{\Pi}(r)$ будем понимать
величину
$$
\bigvee (r';  t \in \Lambda_{\Pi}(r)) =   sup_{  t_i \in
\Lambda_{\Pi}(r)} \,\, \sum_i \Vert r'(t_i) -  r'(t_{i-1}) \Vert.
$$
}

\begin{lem} Для любой выпуклой п.о. степени 1 функции $q \rightarrow
\psi (q): \mathbb{R}^n \rightarrow \mathbb{R}$, любой кривой $
r(\cdot) \in \tilde{\wp}(D)$   и  для произвольного подмножества
из  $\Lambda_{\Pi}(r) $ верно неравенство
$$
    \bigvee (\Theta '; t \in \Lambda_{\Pi}(r)) <  c_1(D, \psi) ( 1 + \bigvee (r';  t \in \Lambda_{\Pi}(r)),
$$
где $\Theta(t) = \psi(r(t))$  для всех $t \in [0,T_r], c_1(D,
\psi)$ - некоторая константа. \label{1lemdifconvmult}
\end{lem}

{\bf Доказательство.}  Рассмотрим сперва случай, когда
$\psi(\cdot)$ есть гладкая функция на $\mathbb{R}^n \backslash
\{0\}.$ Пусть
$$
 \psi (r(t)) = \mbox{max} \,_{v \in \partial \psi(0)} (v, r(t))=(v(t),r(t)), \;\;
v(t) \in \partial \psi(0),
$$
где $\partial \psi (0)$ - субдифференциал функции $\psi(\cdot)$ в
нуле. Будем также считать,  что $r(\cdot)-$  дифференцируемая
кривая по $t \in [0,T_r]$ .

Очевидно, что
$$
\psi'(r(t))=(v'(t),r(t))+(v(t),r'(t)).
$$
Так как $r(t)$ есть нормаль к границе множества $\partial \psi(0)$
в точке $v(t)$, то векторы $v'(t)$ и $r(t)$ перпендикулярны друг к
другу, а следовательно, $(v'(t),r(t))=0.$  Поскольку кривая
$r(\cdot)$ параметризована естественным образом, то $\Vert r'(t)
\Vert =1$ для любых $t \in [0,T_r]$.

Нетрудно проверить следующую цепочку неравенств
\begin{equation}
\mid \psi' (r(t_1)) - \psi' (r(t_2)) \mid= \mid (v(t_1),r'(t_1)) -
(v(t_2),r'(t_2)) \mid = \mid (v(t_1)- v(t_2),r'(t_1)) +
$$ $$
(v(t_2),r'(t_1))- (v(t_2),r'(t_2)) \mid \leq  \Vert v(t_1) -
v(t_2) \Vert \; \Vert r'(t_1) \Vert + \Vert r'(t_1) - r'(t_2)
\Vert \; \Vert v(t_2) \Vert \leq
$$ $$
\Vert v(t_1) - v(t_2) \Vert  + L(D) \Vert r'(t_1) - r'(t_2) \Vert.
\label{difconvmult3a}
\end{equation}
Отсюда следует, что
$$
   \bigvee (\Theta ';  t \in \Lambda_{\Pi}(r)) \leq \mbox{длина кривой v(t) для }  t \in \Lambda_{\Pi}(r)
   + L(D) \cdot \bigvee (r';  t \in \Lambda_{\Pi}(r)),
$$
где $v(t)$ - граничные векторы множества $\partial \psi(0)$ с
нормалями $r(t)$ и $L(D)$ - константа Липшица функции
$\psi(\cdot).$ Если будет показано, что длина кривой $v(t), t \in
\Lambda_{\Pi}(r),$  ограничена сверху одной и той же константой
для всех кривых $ r(\cdot) \in \tilde{\wp}(D)$, то лемма
\ref{1lemdifconvmult} для рассматриваемого случая будет доказана.

Докажем равномерную ограниченность длины кривой $v(t),  t \in
\Lambda_{\Pi}(r),$ независимо от $n$. Для $n=2$ утверждение верно,
что доказано в \cite{proudconvex1}.

Рассмотрим проекцию $Pr(r(\cdot))$ кривой $ r(\cdot) \in
\tilde{\wp}(D)$ на одну из плоскостей $\Pi$, с которой векторы
$r(t),\, r'(t),  t \in \Lambda_{\Pi}(r),$ образуют угол, не
больший $ \pi / 4 $, а также  $Pr(r(\cdot))$ принадлежат
$\wp(Pr(D))$, где $Pr(D)-$ проекция множества $D$ на $\Pi$.

Так как нас интересуют крайние векторы $v(\cdot)$ и нормали
$r(\cdot)$ к множеству $\partial \psi(0)$, то без ограничения
общности будем считать, что $0 \in int \,\, \mbox{co}
Pr(r(\cdot))$.

Спроектируем кривую $v(\cdot)$ на плоскость $\Pi$. Обозначим
получившуюся кривую через $Pr(v(t)),  t \in [0,T_r]$. Кривая
$Pr(v(\cdot))$ принадлежит границе выпуклого компактного множества
$V_r \subset \mathbb{R}^2$ на плоскости $\Pi$. Действительно,
нормалями к кривой $Pr(v(\cdot))$  являются векторы
$Pr(r(\cdot))$, а кривая $Pr(r(\cdot))$ принадлежит множеству
$\wp(Pr(D))$. Это также следует из того, что при приближении
функции $\psi(\cdot)$ многогранной $\psi_m (\cdot)$ получаем
разбиение кривой $r(\cdot)$ на подмножества точками $t_i$, $i \in
1:m$. Все двугранные углы, построенные по секторам на плоскости
$\Pi$ разбиения кривой $Pr(r(\cdot))$ точками $t_i$, $i \in 1:m$,
также будут выпуклыми, так как обход кривых $r(\cdot)$ и  $
Pr(r(\cdot))$ происходит в одном направлении, поскольку $
Pr(r(\cdot)) \in \wp(Pr(D))$. Отсюда следует сказанное выше.

Докажем, что для всех $r(\cdot) \in \tilde{\wp}(D)$ множества
$V_r$ равномерно ограничены в совокупности. Построим по кривой
$Pr(v(\cdot))$ п.о. выпуклую функцию $\eta(\cdot): \mathbb{R}^2
\rightarrow \mathbb{R}$. Положим по определению
$$
\eta(q)=\max_{y \in V_r} (y,q) \,\,\,\,\, \forall q \in S^2_1(0).
$$
Множество $V_r$ является субдифференциалом в нуле функции
$\eta(\cdot)$, т.е. $V_r=\partial \eta (0)$. Функция $\eta(\cdot)$
является липшицевой с константой $L(D)$, так как все ее обобщенные
градиенты ограничены по норме той же константой, какой ограничены
по норме обобщенные градиенты функции $\psi(\cdot)$, т.е. $L(D)$.

Из сказанного выше следует, что длины кривых $Pr(v(\cdot))$
ограничены в совокупности для всех $r(\cdot) \in \tilde{\wp}(D)$.
Векторы $Pr(v(t)), t \in [0, T_r],$ являются проекциями векторов
$v(t), t \in [0, T_r],$ на $\Pi$, но поскольку векторы $r'(t), t
\in \Lambda_{\Pi}(r),$ образуют с $\Pi$ угол, не больший $ \pi / 4
$, то длина вектора $v(t_1)- v(t_2)$ при малом $| t_1 - t_2 | $
для  $t_1, t_2 \in \Lambda_{\Pi}(r) $ оценивается сверху величиной
$S_1(D,\psi)  \cdot \| Pr(v(t_1)) - Pr(v(t_2)) \|$, где
$S_1(D,\psi) - $ константа, определяемая рассматриваемым классом
кривых $\tilde{\wp}(D)$, а именно: максимальным углом, образуемым
векторами $r(t), r'(t), t \in \Lambda_{\Pi}(r),$ с плоскостью $\Pi
$, множеством $D$ и самой функцией $\psi$. Отсюда можно
утверждать, что длины кривых $v(\cdot)$ равномерно ограничены
сверху для всех $r(\cdot) \in \tilde{\wp}(D)$.

Из равномерной ограниченности длин кривых  $v(\cdot)$ для всех
$r(\cdot) \in \tilde{\wp}(D)$  и из неравенства
(\ref{difconvmult3a}) следует утверждение леммы
\ref{1lemdifconvmult} при сделанном предположении.

Пусть теперь $\psi (\cdot)-$ произвольная выпуклая п.о. функция. С
любой степенью точности ее можно приблизить  на $B^{n}_1(0)
\backslash {0}$ выпуклой п.о. дифференцируемой функцией
$\hat{\psi}(\cdot)$ так, чтобы в метрике Хаусдорфа
субдифференциалы в нуле этих функций отличались друг от друга как
угодно мало. Но тогда длины кривых $v(\cdot)$ для любых $r(\cdot)
\in \tilde{\wp}(D)$, построенных для функций $\hat{\psi}(\cdot)$ и
$\psi (\cdot)$, также будут отличаться друг от друга как угодно
мало. Кривую $r(\cdot)$ можно приблизить дифференцируемой кривой
таким образом, чтобы производные по $t$ этих кривых в точках
дифференцируемости  отличались друг от друга по норме на
произвольно малую величину. Таким образом, любые конечные суммы,
используемые при вычислении вариаций функций $\Theta'(\cdot)$ и
$\hat{\Theta}'(\cdot)$  для негладкого и гладкого случая, могут
быть сделаны как угодно близкими друг к другу. Но поскольку
вариацию функции $\hat{\Theta}'(\cdot)$ можно ограничить сверху
величиной, зависящей только от множества $D$, кривой $r'(\cdot)$ и
самой функции $\psi(\cdot)$, а также некоторых констант, то лемма
\ref{1lemdifconvmult} доказана. $\Box$

На основе этой леммы докажем утверждение, аналогичное в
\cite{proudconvex2}.

\begin{lem} Пусть $f_1(\cdot): \mathbb{R}^n \rightarrow \mathbb{R}-$
непрерывная выпуклая функция и $r(\cdot) \in \tilde{\wp}(D).$
Тогда существует константа $c_2(D,f_1)>0,$ что для произвольного
подмножества отрезка $[0, T_r]$ из  $\Lambda_{\Pi}(r) $
\begin{equation}
   \bigvee (\Phi_1' ; t \in \Lambda_{\Pi}(r) \leq  c_2(D,f_1),( 1 + \bigvee (r';
   t \in \Lambda_{\Pi}(r)),
\label{difconvmult4} \end{equation} где $\Phi_1(t) = f_1(r(t)), t
\in [0,T_r].$ \label{2lemdifconvmult}
\end{lem}

{\bf  Доказательство.} На начальном этапе будем считать, что
$f_1(\cdot)$ дважды непрерывно дифференцируемая функция на $D$,
которая принимает неотрицательные значения и начало координат $-$
ее точка минимума, где $f_1(0)=0$. Обозначим константу Липшица
функции $f_1(\cdot)$ на $D$   через $L_1(D)$.

Считаем,  что $0$ принадлежит внутренности выпуклой области на
$\Pi$ с границей $Pr(r(\cdot))-$ проекцией кривой $r(\cdot)$ на
плоскость $\Pi$, с которой векторы $r(t), r'(t),$ $ t \in [0,
T_r],$ образуют углы не более $\pi /4$ и $Pr(r(\cdot)) \in
\wp(Pr(D)) $.

Построим для функции $f_1(\cdot)$ выпуклую п.о. степени 1 функцию
$\eta(\cdot): \mathbb{R}^2 \rightarrow \mathbb{R}$,  которая на
$Pr(r(t)) $ принимает значения, равные $f_1(r(\cdot)),$ а в начале
координат $-$ нуль. В данном случае под $Pr(r(t)) $ будем понимать
двумерные векторы координатной плоскости $\Pi$.

Положим по определению
$$
\eta(Pr(r(t)))=f_1(r(t)) \,\,\,\, \forall t \in [0, T_r]
$$
и для любого $\lambda >0$
$$
\eta(\lambda  Pr(r(t)))=\lambda \eta(  Pr(r(t))) \,\,\,\, \forall
t \in [0, T_r].
$$
Ясно, что $\eta (\cdot)$ строится однозначно по функции
$f_1(\cdot)$ и выбранной кривой $r(\cdot).$

Функция $\eta(\cdot) $ будет п.о, так как для любого $\lambda
>0$ и $z=\mu Pr(r(t)) \in \Pi , \mu >0,$
$$
\eta(\lambda z)= \eta (\lambda  \mu  Pr(r(t))) = \lambda  \mu (
\eta (Pr( r(t)))) =\lambda  \eta (\mu Pr(r(t)))=\lambda \eta (z).
$$
Функция $\eta(\cdot) $ липшицева с константой $\sqrt{2} L_1(D)$
для $t \in  \Lambda_{\Pi}(r)$, так как
$$
| \eta(Pr(r(t))) | = | f_1 (r(t)) | \leq L_1(D) \| r(t) \| \leq
\sqrt{2} L_1(D) \| Pr(r(t)) \| \,\, \forall t \in
\Lambda_{\Pi}(r).
$$
Функция $\eta(\cdot) $ будет выпуклой. Покажем это.

Рассмотрим функцию
$$
f_{\varepsilon}(x)=f_1(x)+\varepsilon ( \mid \mid x \mid \mid^2),
\,\,\,\, \varepsilon>0, \,\, x \in \mathbb{R}^n.
$$
Разобьем отрезок $[ 0,T_r] $ точками  $\{t_i\} , i \in 1:J ,$ на
равные отрезки. Построим плоскости $\pi_i$ в $\mathbb{R}^{3},$
проходящие соответственно через точки $0,
(Pr(r(t_i)),f_{\varepsilon}(r(t_i))),
(Pr(r(t_{i+1})),f_{\varepsilon}(r(t_{i+1})), i \in 1:J $. Части
плоскостей $\pi_i ,i \in 1:J$ , определенных в секторах,
образованных векторами $0, Pr(r(t_i)), Pr(r(t_{i+1}))$, определяют
график п.о. степени 1 многогранной функцию $(\eta_{\varepsilon})_J
(Pr(r(\cdot))).$

Будем понимать под двугранным углом функцию, равную максимуму или
минимуму линейных функций, графики которых совпадают с
гиперплоскостями ( плоскостями) $\pi_i, \pi_{i+1}$ построенными по
соседним секторам, как это описано выше. Покажем, что все
двугранные углы функции $(\eta_{\varepsilon})_J (Pr(r(\cdot)))$
$-$ выпуклые.

Поскольку всегда любую кривую $r(\cdot) \in \tilde{\wp}(D)$ можно
приблизить  с любой степенью точности гладкой кривой из
$\tilde{\wp}(D),$ то без ограничения общности будем считать, что
$r(\cdot)$ - гладкая дифференцируемая кривая с производной
$r'(\cdot).$

Под градиентом гиперплоскости (плоскости) $\pi_i$ будем понимать
градиент линейной функции, график которой совпадает с плоскостью
(гиперплоскостью) $\pi_i$. Обозначим градиенты плоскостей $\pi_i$
и $\pi_{i+1}$ через $\nabla \pi_i$ и $\nabla \pi_{i+1}$
соответственно. Воспользуемся теоремой о средней точке, согласно
которой существует такая точка $t_m \in [t_i, t_{i+1}],$ что
$$
           \partial f_{\varepsilon}(r(t_m))/ \partial {e}_i =
           (\nabla \pi_i, e_i),
$$
где
$$
e_i=(Pr(r(t_{i+1}))-Pr(r(t_i)))/ \mid \mid
Pr(r(t_{i+1}))-Pr(r(t_i)) \mid \mid.
$$
Аналогично для плоскости $\pi_{i+1}$ и некоторой точки $t_c \in
[t_{i+1},  t_{i+2}]$ имеем
$$
           \partial f_{\varepsilon}(r(t_c))/ \partial {e}_{i+1} =
           (\nabla \pi_{i+1}, e_{i+1}),
$$
где
$$
e_{i+1}=(Pr(r(t_{i+2}))-Pr(r(t_{i+1})))/ \mid \mid
Pr(r(t_{i+2}))-Pr(r(t_{i+1})) \mid \mid.
$$
Здесь под векторами $e_i, e_{i+1}$ надо понимать либо двумерные,
либо $n$- мерные векторы в зависимости от контекста. Так, если
векторы $e_i, e_{i+1}$ относятся к производной по направлению
функции $f_{\varepsilon}(\cdot)$, то это $n$- мерные векторы, если
они относятся к производной по направлению функции $\pi(\cdot)$,
то это двумерные векторы.

Функция $f_{\varepsilon}(\cdot)$ сильно выпуклая, так как ее
матрица вторых частных производных положительно определенная.
Любая выпуклая функция имеет неубывающую производную по
направлению вдоль произвольного луча. Но для сильно выпуклой
функции производная по касательному направлению к проекции на
$\Pi$ кривой $r(x_0,\tau, g)= x_0+\tau g +o_{\varepsilon}(\tau)$,
$g \in \mathbb{R}^n$, $ \tau >0,$ есть возрастающая функция вдоль
этой кривой в малой окрестности точки $x_0$. Поэтому для
достаточно большом $J$ и равномерном разбиении кривой $r(\cdot)$
точками $t_i$ на подмножества, длины которых стремятся к нулю при
$J \rar \infty$, имеем
$$
\partial f_{\varepsilon}(r(t_m))/ \partial {e}_i < \partial
f_{\varepsilon}(r(t_c))/ \partial {e}_{i+1},
$$
или
$$
(\nabla \pi_i, e_i) < (\nabla \pi_{i+1}, e_{i+1}),
$$
поскольку значения функции $\eta_{\varepsilon}(\cdot)$ на кривой
$r(\cdot)$ совпадают со значениями функции
$f_{\varepsilon}(\cdot)$ согласно построению.

Учтем также, что разность $\nabla \pi_{i+1} - \nabla \pi_i$
перпендикулярна вектору $Pr(r(t_{i+1})).$ Отсюда и из неравенства
выше следует, что двугранный угол $\pi_i, \pi_{i+1}$ - выпуклый.
При $J \rightarrow \infty$
$$
(\eta_{\varepsilon})_J(\cdot) \Rightarrow
(\eta_{\varepsilon})(\cdot).
$$
Так как точечный предел для выпуклых функций равносилен
равномерному пределу, то $\eta_{\varepsilon}(\cdot)$ - выпуклая
функция. Также $\eta_{\varepsilon}(\cdot) \Rightarrow \eta
(\cdot)$ при $\varepsilon \rightarrow +0,$ т.е. $\eta(\cdot)-$
выпуклая функция, что и требовалось доказать.

Очевидно, что градиенты линейных функций, графики которых есть
$\pi_i, i \in J,$ ограничены константой, зависящей только от
множества $D$ и самой функции $ f_1(\cdot),$ поскольку функции
$\pi_i, i \in J,$ строятся по функции $\eta(\cdot)$, которая
липшицевая с константой $\sqrt{2} L_1(D)$, где $ L_1(D)$ константа
Липшица функции $ f_1(\cdot).$

Пусть
$$
 \Theta(t) = \eta (Pr(r(t)))=f_1(r(t))=\Phi_1(t) \;\; \forall t \in [0,T_r].
$$
Поскольку
$$
  \bigvee (\Theta '; 0,T_r)  =  \bigvee( \Phi_1 ' ; 0,T_r)
$$
и
$$
\bigvee (Pr(r');0,T_r) \leq \bigvee (r';0,T_r),
$$
то из леммы 1 следует, что для некоторой константы $c_2(D, f_1)>0$
верно неравенство
$$
\bigvee (\Theta '; t \in \Lambda_{\Pi}(r)) \leq  c_2(D, f_1)(1+
\bigvee (r';t \in \Lambda_{\Pi}(r))).
$$
Следовательно, для вариации производной функции $\Phi_1(\cdot)$
также верно неравенство
$$
\bigvee (\Phi_1'; t \in \Lambda_{\Pi}(r)) \leq  c_2(D, f_1)(1+
\bigvee (r';t \in \Lambda_{\Pi}(r))).
$$
Если функция $ f_1(\cdot)$ не есть дважды непрерывно
дифференцируемая, то ее можно приблизить выпуклой дважды
непрерывно дифференцируемой функцией $ \tilde{f}_1(\cdot)$ и
построить соответствующую ей функцию $\tilde{\eta} (\cdot)$ так,
чтобы значения  функций $\eta (\cdot)$,  $\tilde{\eta} (\cdot)$ и
их производных там, где они существуют, как угодно мало отличались
друг от друга. Но тогда написанные выше неравенства  будут верны
для функций $\Theta(\cdot)$, $\tilde{\Theta}_1(\cdot)$,
построенных по $\eta (\cdot), \tilde{\eta} (\cdot)$
соответственно, и их производных. Значит неравенство для вариации
производных функции $\Phi_1(\cdot)$ верно для общего случая. Лемма
\ref{2lemdifconvmult} доказана. $\Box$

Из леммы \ref{2lemdifconvmult} следует, что если $f(\cdot,\cdot)$
представима  в виде разности выпуклых функций, т.е.
$$ f(z) = f_1(z) - f_2(z) \;\;\;\; \forall z \in D,   $$
где $f_i(\cdot,\cdot), i=1,2,$ - выпуклые, то условие
(\ref{difconvmult3})  c необходимостью выполняется. Действительно,
для произвольной $r(\cdot) \in \wp(D)$ введем обозначения
$$
\Phi_1(t) = f_1(r(t)), \Phi_2(t) = f_2(r(t)) \;\;\; \forall t \in
[0,T_r].
$$
Поскольку \cite{kolmogorovfomin}
$$
\bigvee (\Phi ';0,T_r) \leq \bigvee (\Phi_1 '; 0,T_r) + \bigvee
(\Phi_2 ';0,T_r)
$$
то, принимая во внимание неравенство (\ref{difconvmult4}),
неравенство (\ref{difconvmult3}) с необходимостью выполняется.

{\bf Достаточность}.   Докажем достаточность условия
(\ref{difconvmult3}) для представления функции $f(\cdot)$ в виде
разности выпуклых функций.

Прежде всего покажем, что  для любого $r(\cdot) \in
\tilde{\wp}(D)$ верно неравенство
$$
\bigvee (\Phi_N ';t \in \Lambda_{\Pi}(r)) \leq c (1+ \bigvee (r';t
\in \Lambda_{\Pi}(r)),
$$
где $\Phi_N(t) = f_N(r(t))$, и множества кривых $\tilde{\wp}(D)$
достаточно для ответа на вопрос о представимости функции в виде
разности выпуклых.

Действительно, для достаточно равномерного разбиения области $D_N
\subset D$, $\rho_H (D_N, D) \rar 0$, на многогранники $G_k, k \in
\overline{1, N},$  с непустыми внутренностями функции $\Phi_N, \,
\Phi \,$  и их производные $\Phi'_N, \, \Phi' \,$, вычисленные в
точках,  где они существуют, близки друг к другу с любой степенью
точности $\varepsilon_N$, где $\varepsilon_N \rightarrow +0$ при
$N \rar \infty$.

Поэтому произвольная конечная сумма
$$
\sum_{i=1}^{N} \mid \Phi_N'(t_i) - \Phi_N'(t_{i+1}) \mid
$$
для больших $N$ будет как угодно мало отличаться от суммы
$$
\sum_{i=1}^{N} \mid \Phi'(t_i) - \Phi'(t_{i+1}) \mid .
$$
А поскольку вариация функции $\Phi_N'(\cdot)$  может только
возрастать при вложенности разбиений области $D_N$ при увеличении
$N$, то отсюда и из сказанного выше следует, что для произвольного
подмножества из  $\Lambda_{\Pi}(r)$
\begin{equation} \bigvee (\Phi_N';t \in \Lambda_{\Pi}(r)) \leq \bigvee
(\Phi';t \in \Lambda_{\Pi}(r))+\delta(N) \leq c (1+ \bigvee (r';t
\in \Lambda_{\Pi}(r))) , \label{difconvmult5}
\end{equation} где $\delta(N) \rightarrow +0 $ при $N \rightarrow
\infty$, $c$- константа.

Вариация производных по направлению вдоль произвольного отрезка
суммы выпуклых функций равна сумме вариаций производных этих
выпуклых функций по тому же отрезку. Если будет доказано, что
сумма вариаций  производных всех выпуклых двугранных углов функции
$f_N(\cdot)$ вдоль любого отрезка области $D$  ограничена сверху
константой, независящей от $N$, то отсюда будет следовать, что
ограничена сверху той же константой вариация производной функции
$f_{1,N}(\cdot)$ вдоль произвольного отрезка области $D$. Отсюда
следует равномерная ограниченность и равностепенная непрерывность
функций $f_{1,N}(\cdot) - c_{1,N}$, где $c_{1,N}$ некоторые
константы. Но тогда по теореме Арцела - Асколи  из
последовательности $f_{1,N}(\cdot) - c_{1,N}$ можно выбрать
подпоследовательность, равномерно сходящуюся на $D$ к выпуклой
функции $f_{1}(\cdot)$. Соответствующая  последовательность
$f_{2,N}(\cdot) - c_{1,N}$ будет стремиться к выпуклой функции
$f_2(\cdot)$. Переходя для любого $x \in \mbox{int} D_N$ в
равенстве
$$
f_N(x) = (f_{1,N}(x) - c_{1,N}) - (f_{2,N}(x) - c_{1,N})
$$
к пределу по $N \rar \infty$, получим представление функции $f$ в
виде разности выпуклых, т.е. $f(\cdot)$ есть ПРВ функция.

Пусть условия теоремы выполняются, но $f(\cdot)$ не есть ПРВ
функция. Проделаем следующую процедуру. Путем разбиения множества
$D$ на выпуклые непересекающиеся подобласти можно выделить ту
подобласть, где функции $f_{1,N}(\cdot)$ имеют предельное
бесконечное значение вариации производной вдоль некоторых отрезков
этой подобласти при $N \rightarrow \infty$. Действительно, в
противном случае из последовательности функции $f_{1,N}(\cdot) -
c_{1,N}$ можно было бы выбрать сходящуюся подпоследовательность в
каждой подобласти, а значит,  $f(\cdot)$ была бы ПРВ функцией на
всем множестве $D$.

Далее разбиваем выделенную подобласть на меньшие области и опять
выделяем ту, где вариация производной функций $f_{1,N}(\cdot)$
вдоль некоторых отрезков неограничена при $N \rightarrow \infty$.
В итоге определяем точку $M$, в произвольной окрестности которой
вариация производной функций $f_{1,N}(\cdot)$ вдоль некоторых
отрезков неограничена при $N \rightarrow \infty$. Без ограничения
общности можно считать, что $M-$ внутренняя точка множества
$\bar{D},$ так как все получаемые в процессе применения алгоритма
функции $-$ равномерно липшицевые и могут быть распространены во
вне множества $\bar{D},$

Берем произвольную окрестность $S$ точки $M$. Рассмотрим конечное
число конусов $K$ с общей вершиной $M$ и с непересекающимися
внутренностями. Разбиваем $S$ на подмножества путем пересечения с
конусами $K$. Выбираем одно из таких подмножеств $K \cap S$, где
вариация производной функций $f_{1,N}(\cdot)$ вдоль некоторых
отрезков неограничена при $N \rightarrow \infty$. Далее выбранное
множество $K \cap  S$ разбиваем на конечное число подмножеств
непересекающимися по внутренности конусами с вершиной в точке $M$.
Выбираем из них такое,  где вариация производной функций
$f_{1,N}(\cdot)$ вдоль некоторых отрезков выбранного подмножества
неограничена при $N \rightarrow \infty$. Процесс повторяем до
бесконечности.

Множество выбранных конусов стягивается к некоторому направлению,
определяемому единичным вектором $l$ с вершиной в точке $M$.
Очевидно, что в произвольном конусе $K$ с вершиной с точке $M$,
содержащем вектор $\alpha l$ в $\mbox{int} \, K$, $\alpha >0$,
вариация производной функций $f_{1,N}(\cdot)$ вдоль некоторых
отрезков окрестности $S$ неограничена при $N \rightarrow \infty$.

Векторов, подобных вектору $l$, может быть не один. Для найденного
вектора $l$ возможны два случая:

a) вариация производных функций $f_{1,N}(\cdot)$ по направлению
$l$ неограничена при $N \rightarrow \infty$;

б) вариация производных функций $f_{1,N}(\cdot)$ по направлению
$\zeta,$ перпендикулярному направлению $l$, неограничена при $N
\rightarrow \infty$. Вектор $\zeta$ находится аналогично вектору
$l$, о чем будет сказано ниже.

Но если вариация производной функции $f_{1,N}(\cdot)$ неограничена
вдоль некоторых отрезков множества $K \cap  S$  при $N \rightarrow
\infty$ , то и сумма вариаций производных вдоль этих отрезков
выпуклых двугранных углов, построенных по многогранникам разбиений
множества $K \cap S $ при построении функций $f_{N}(\cdot)$, также
неограничена при $N \rightarrow \infty$.

Рассмотрим случай а). Возьмем произвольный конус $K$, содержащий
вектор $\alpha l$ в $\mbox{int} \, K$, $\alpha >0$. Будем
рассматривать выпуклые двугранные углы функций $f_{1,N}(\cdot)$,
построенные по многогранникам из $K \cap S $.

%За счет равномерной липшицевости по $N$ всех двугранных углов
%функций $f_{N}(\cdot)$ вариации производных по направлению этих
%двугранных углов равномерно непрерывны относительно направления и
%$N$.

Для каждого выпуклого $k-$ ого  двугранного  функции
$f_{N}(\cdot)$ выделим отрезок $v_{k,N}$, вариация производной
вдоль которого для $k$- ого двугранного угла максимальна и равна
$a_{k,N}$. Ясно, что отрезок $v_{k,N}$ должен быть, как и в
двумерном случае, перпендикулярен проекции на $\mathbb{R}^n$ грани
$k$-ого двугранного угла, совпадающей с пересечением
гиперплоскостей $\pi_i$ и $\pi_j$, образующих этот двугранный
угол. Очевидно, что грань произвольного двугранного угла имеет
размерность $n-1$ и $a_{k,N}= || \nabla \pi_i - \nabla \pi_j || $,
где $\nabla \pi_i, \nabla \pi_j-$ градиенты гиперплоскостей
$\pi_i$ и $\pi_j$ соответственно.

Пусть угол наклона отрезков $v_{k,N}$ с направлением $l$ не
превосходит $\pi / 2- \delta$ для некоторого $ \delta>0.$

Путем уменьшения окрестности $S$ и разбиения $K \cap S $ на
меньшие подконусы, стягивающихся к вектору $\alpha l, \alpha>0,$ и
к точке $M$, и рассмотрения в каждом из них своей группы отрезков
$\{v_{k,N} \}$ для всех значений $k$ и $N \rightarrow + \infty$,
можно выделить одну или несколько групп указанных отрезков, каждую
из которых можно пересечь кривой $r(\cdot) \in \tilde{\wp}(D),$
образующей в точке пересечения с отрезками $v_{k,N}$ угол, не
превосходящий $\pi /2 -\delta_1, \delta_1>0.$ Поскольку конус $K$,
содержащий вектор $\alpha l, \alpha>0,$  окрестность $S$ точки $M
- $ произвольные, то можно рассматривать такие кривые, для которых
$r'(t) \rightarrow -l,$ когда $r(t) \rightarrow M$, а также $
r(\cdot) \in \tilde{\wp}(D)$.  Сама кривая $r(\cdot)$ будет
включать в себя отрезки, близкие к отрезкам $v_{k,N}$.

Если для рассматриваемого случая подгруппа отрезков $\{ v_{k,N}
\}$ существует только одна, то вдоль найденной кривой $ r(\cdot)
\in \tilde{\wp}(D)$  с конечной вариацией $\bigvee (r';t \in
\Lambda_{\Pi}(r))$ сумма вариаций производных выпуклых двугранных
углов равна бесконечности. Так как при выполнении неравенства
(\ref{difconvmult3}) выполняется неравенство (\ref{difconvmult5}),
а мы нашли кривую $r(\cdot)$, вдоль которой сумма вариаций
производных двугранных углов бесконечна, то из
(\ref{difconvmult5}) следует, что  вдоль $r(\cdot)$ не ограничена
вариация производной функции $\Phi(\cdot)$ для $t \in
\Lambda_{\Pi}(r)$. Приходим к противоречию насчет справедливости
неравенства в условии теоремы. Противоречие получилось, так как мы
предположили, что $f$ не ПРВ функция.

Кроме того, возможен случай, когда у нас есть несколько групп
отрезков $\{ v_{k,N} \}_i,$ для каждой из которых найдется кривая
$r_i(\cdot) \in \tilde{\wp}(D)$, что
$$
\bigvee (\Phi_N';t_{1,r_i},t_{2, r_i})= c_i, \,\,\,  r_i' (t)
\rightarrow_{t \rightarrow t_{r_i}} -l,
$$
где $t_{1,r_i},t_{2, r_i}$ - параметры кривой $r_i(\cdot)$ при
натуральной параметризации,  а также
$$
 \sum_i \, c_i = \infty.
$$
Тогда кривую $r(\cdot) \in \tilde{\wp}(D),$ для которой сумма
вариаций производных вдоль  отрезков этой   кривой стремится к
бесконечности при $N \rightarrow +\infty$, будем строить следующим
образом.

Кривая $r(\cdot)$ должна содержать достаточное количество $k_i$
отрезков из каждой группы отрезков $\{ v_{k,N} \}_i,$ (либо
близких к ним), чтобы
$$
\bigvee (\Phi'_N ;t_{1, r_{i}},t_{2, r_{i}}) =  c_i - \mu_i,
$$
где $t_{1,r_{i}}, t_{2,r_{i}}   >0-$ значения параметра $t$ для
$i$-ой группы отрезков при естественной параметризации кривой
$r_i(\cdot)$, $\mu_i < c_i -$ малые положительные числа, для
которых
$$
\sum_i \, \mu_i < \infty.
$$

Нетрудно видеть, что всегда такую кривую $r(\cdot)$ построить
можно. Она будет состоять из частей  кривых $r_i(\cdot)$. Для
этого надо осуществить плавный переход от одной кривой
$r_i(\cdot)$ к кривой $r_{i+1}(\cdot),$ не выходя из множества
$\tilde{\wp}(D).$ Поскольку $r'_i(t) \rightarrow -l$ при $t
\rightarrow T_{r_i}$ для всех $i$, то подобная процедура
осуществима всегда, так как конус $K$ можно разбить на меньшие
подконусы с общим вектором $\alpha l, \alpha >0$, и рассматривать
тот, где $f(\cdot)$ не ПРВ функция. Причем данная процедура
приводит к кривой  $r(\cdot) \in \tilde{\wp}(D),$ содержащей или
пересекающей под острыми углами отрезки $\{ v_{k,N} \}_i,$   с
бесконечной суммой длин.

Но тогда
$$
\sum_i  \bigvee (\Phi';t_{1, r_{i}},t_{2 ,r_{i}}) \geq  \sum_i
\bigvee (\Phi'_N;t_{1, r_{i}},t_{2 ,r_{i}}) = \sum_i(c_i-\mu_i)=
\sum_i c_i -\sum_i \mu_i  =\infty.
$$
Нетрудно видеть, что согласно алгоритму всегда можно построить
кривую $r(\cdot)$ с конечной суммой вариаций  $\sum_i \bigvee
(r';t_{1,r_i}, t_{2,r_{i+1}}))$. Но как следует из
(\ref{difconvmult5}), неравенство (\ref{difconvmult3}), которое
должно быть верно по предположению достаточности условия теоремы,
нарушается. Приходим к противоречию с условием теоремы, так как
предположили, что $f$ не ПРВ функция.

Случай б). Если сумма выпуклых двугранных углов функции
$f_N(\cdot)$ в произвольном конусе с вершиной $M$, содержащем во
внутренности вектор $\alpha l$, $\alpha >0$, бесконечна, а ее
вариация производной вдоль направления $l$ конечна при $N \rar
\infty $, то отсюда следует, что вариация производной суммы
выпуклых двугранных углов функции $f_N(\cdot)$ бесконечна при $N
\rightarrow \infty$ вдоль некоторого направления $\zeta $,
перпендикулярного направлению $l$.

Найти направление $\zeta$ можно следующим образом. Возьмем
произвольную окрестность $S$ точки $M$. Разобьем $S$ на
подмножества, совпадающие с пересечениями $S$ и конусов $V$ с
вершиной в точке $M$, образованных вектором $l$ и векторами из
ортогонального к $l$ подпространства $L^\bot$. Пересечения $V \cap
L^\bot$ также разбивают  $S \cap L^\bot$ на подмножества в виде
конусов из $L^\bot$ с вершиной в точке $M$. Далее для построенных
конусов из $S \cap L^\bot$ повторяем рассуждения аналогично тому,
как это делали при нахождении вектора $l$ до тех пор, пока не
придем к направлению $\zeta \in L^\bot$.

В произвольном конусе $K$ с вершиной $M,$  содержащем во
внутренности векторы $\alpha l$, $\alpha >0$, сумма вариаций
производных выпуклых двугранных углов функции $f_N(\cdot)$ вдоль
направления $\zeta $   бесконечна при $N \rightarrow \infty$.

Все отрезки $v_{k,N}$ можно разбить на такие группы $\{ m \}$
отрезков, которые можно пересечь гладкой кривой $r_{m,N}(\cdot)
\in \tilde{\wp}(D),$ для которой
$$
        r'_{m,N} (\tau) \rightarrow_{\tau \rightarrow T_{r_{m,N}}} -l,
$$
где $T_{r_{m,N}}-$ есть параметр кривой $r_{m,N}(\cdot)$ при
натуральной параметризации в точке $M$.  Кривизны кривых $r_{m,N}
(\cdot)$ стремятся к бесконечности при $\tau \rightarrow
T_{r_{m,N}}$ и $m,N \rar \infty$. Кривая $r_{m,N}(\cdot)\in
\tilde{\wp}(D)$ содержит свою группу отрезков  $\{v_{k,N}\}_m$.
Ясно, что сказанное всегда выполнимо путем разбиения конуса $K$ на
подконусы с общим вектором $\alpha l, \alpha>0,$ а также множества
всех отрезков $v_{k,N}$ на подмножества с требуемыми свойствами.

Кроме того,  кривые $r_{m,N}(\cdot) \in \tilde{\wp}(D)$ и группы
отрезков $\{ v_{k,N} \}_m$ можно выбрать такими, чтобы предел по
$m, N$ суммы вариаций производных функций
$\Phi'_{m,N}(\cdot)=f'_N(r_{m,N}(\cdot))$ вдоль этих отрезков
кривых $r_{m,N}(\cdot)$ был равен бесконечности. Как это сделать,
будет описано ниже. В противном случае функции $f_N(\cdot)$ имели
бы ограниченную вариацию вдоль направления $\zeta$ при $N
\rightarrow +\infty$.

Построение кривых $r_{m,N}(\cdot)$ с неограниченно увеличивающейся
кривизной в точке $M$, для которой
$$
  \lim_{m,N \rar \infty} \sum_i \bigvee (\Phi'_{m,N} ;  t_{i,
r_{m,N}}, t_{i+1, r_{m,N}})=\infty,
$$
осуществляется аналогичным способом, как и в случае a). Для этого
надо построить кривую $r_{m,N}(\cdot) \in \tilde{\wp}(D)$ с
описанными выше свойствами, состоящую из достаточного  количества
$k_{m,N}$  групп отрезков $\{v_{k_,N}\}_m$ (либо близких к ним),
чтобы
$$
\sum_i \bigvee (\Phi'_{m,N} ; t_{i, r_{m,N}}, t_{i+1, r_{m,N}}) =
c_{m,N},
$$
и
$$
  \lim_{m,N \rar \infty} c_{m,N} =\infty,
$$
$ k_{m,N} \rightarrow \infty$ при $m,N \rightarrow \infty$, $[
t_{i, r_{m,N}}, t_{i+1, r_{m,N}} ]$- значение параметра $t$ для
$m$-ой группы отрезков $\{v_{k_,N}\}_m$  при естественной
параметризации кривой $r_{m,N}(\cdot)$. Такие кривые
$r_{m,N}(\cdot) \in \tilde{ \wp}(D) $ всегда можно построить путем
разбиения конуса $K$ на подконусы с общим вектором $\alpha l,
\alpha >0,$ а также плавного перехода от одной группы отрезков
$\{v_{k_,N}\}_m$  к другой. Последнее возможно, так как сумма
вариаций производных выпуклых двугранных углов, построенных по
многогранникам разбиений области $S \cap K$ для любого конуса $K$,
содержащего $l$ в $\mbox{int} K$, вдоль направления $\zeta$
бесконечна при $N \rightarrow \infty$. При увеличении $m,N$ кривые
$r_{m,N}$ будут содержать все большее число отрезков
$\{v_{k_,N}\}_m$ из множества $ S \cap K$, содержащего вектор
$\alpha l$, $\alpha >0$. Кривизны кривых $r_{m,N}$ вблизи точки
$M$ неограниченно увеличиваются при $m,N \rightarrow \infty$.
Причем данная процедура приводит к последовательности кривых $r_m
(\cdot) \in \tilde{\wp}(D)$ с конечной суммой вариаций $\sum_{i,m}
\bigvee (r_m';t_{i, r_{m,N}}, t_{i+1, r_{m,N}})$.

Но тогда
$$
\sum_i \bigvee (\Phi'; t_{i, r_{m,N}}, t_{i+1, r_{m,N}}) \geqslant
\sum_i \bigvee (\Phi'_{m,N}; t_{i, r_{m,N}}, t_{i+1, r_{m,N}}) =
c_{m,N},
$$
где
$$
\lim_{m,N \rar \infty} c_{m,N} =\infty.
$$

Отсюда приходим к противоречию с (\ref{difconvmult3}), так как из
(\ref{difconvmult3}) следует (\ref{difconvmult5}). В то же время
неравенство (\ref{difconvmult5}) не выполняется для всех $m,N$ для
одной и той же константы $c$, участвующей в неравенстве.
Противоречие с нарушением неравенства в условии теоремы получилось
в связи с тем, что мы предположили, что $f$ не ПРВ функция.

Итак, доказано, что при выполнении условия теоремы, сумма вариаций
производных выпуклых двугранных углов функции $f_N(\cdot)$ вдоль
любого отрезка области $D$ при $N \rightarrow \infty$ ограничена
сверху константой, независящей от $N$. Отсюда, как отмечалось
выше, следует, что $f(\cdot)-$ ПРВ функция.

Итак, теорема \ref{difconvmultexthm2} доказана. $\Box$

\vspace{0.5cm}

\subsection{Геометрическая интерпретация теоремы 2}

\vspace{0.5cm}

Перефразируем теорему \ref{difconvmultexthm2}, придав ей более
геометрический характер. Для этого введем понятие поворота кривой
$r(\cdot)$ на графике $\Gamma_f = \{(x,y) \in \mathbb{R}^{n+1}
\mid y = f(x), \, x \in \mathbb{R}^n \}.$

Рассмотрим на $\Gamma_f$ кривую $ R(t)=(r(t),f(r(t))),$ где
$r(\cdot) \in \tilde{\wp}(D).$  Так как функция $f(\cdot)$ есть
липшицевая, то п.в. на $[0,T_r]$ существует производная
$R'(\cdot),$ которую обозначим через $\tau(\cdot)=R'(\cdot).$
Множество $t \in [0, T_r]$, где существует $R'(\cdot),$ обозначим
через $N_R$.

{\bf Определение 2.} {\em  Поворотом кривой $R(\cdot)$ на
многообразии $\Gamma_f$ назовем величину
$$
sup_{ \{t_i\} \subset N_R} \,\, \sum_i \Vert \tau(t_i)/ \Vert \tau
(t_i) \Vert -  \tau(t_{i-1})/ \Vert \tau (t_{i-1}) \Vert \Vert =
O_R.
$$
}

Таким образом, поворот $O_R$ кривой $R(\cdot)$ есть верхняя грань
суммы углов между касательными $\tau(t)$  для $t \in [0,T_r].$
Обозначим через $O_{R,i}$ поворот кривой $R(\cdot)$ на отрезке $
[t_i, t_{i+1}]$.

Нетрудно видеть, что для плоской гладкой кривой, параметризованной
естественным образом, величина $ O_R$ равна интегралу
$$
\int^{T_r}_0 \mid k(s) \mid ds,
$$
где $k(s)$ - кривизна рассматриваемой кривой $ r(\cdot)$ в точке
$s \in [0,T_r],$ т.е. совпадает с обычным определением поворота
кривой в точке \cite{pogorelov1} .

\begin{thm} Для того, чтобы произвольная липшицевая функция
$z \rightarrow f(z) :\mathbb{R}^n \rightarrow \mathbb{R}$ была ПРВ
функцией, необходимо и достаточно, чтобы для любой плоскости
$\Pi$, $0 \in \Pi$, и любой $r(\cdot) \in \tilde \wp(D)$,
параметризованной естественным образом: $t \in [0, T(r)]$, а также
любого подмножества параметров $[T_i, T_{i+1} ] \subset
\Lambda_{\Pi}(r) $ существовали константы $c_{21}(D,f),
c_{22}(D,f) >0$  такие, что сумма поворотов $O_{R, i}$ кривой
$R(\cdot)$ для $t \in [T_i, T_{i+1}] $ для  всех $i$ на $\Gamma_f$
была ограничена сверху неравенством
\begin{equation}
\sum_1^m O_{R, i } \leq   c_{21}(D,f) + c_{22}(D,f) \sum_1^m\vee
(r'; T_i,T_{i+1}) \,\,\,\,\forall  r(\cdot) \in \tilde \wp(D)
\label{difconvmult6}
\end{equation}
\label{difconvmultexthm3}
\end{thm}

\begin{rem}

Известно (см. Теорему \ref{difconvmultexthm6}), что вопрос о
представимости функции в виде разности выпуклых сводится к вопросу
о ее локальной представимости. Поскольку плоскость, определяемая
векторами $r(t), r'(t)$, образует угол, не больший $\frac{\pi}{4}$
с какой-то координатной плоскостью, то в данной  и дальнейших
Теоремах вместо произвольной плоскости $\Pi$ достаточно
рассматривать только координатные плоскости $\Pi_i$.

\end{rem}

{\bf Доказательство. }{\bf Необходимость}. Пусть $ f(\cdot)$ есть
ПРВ функция. Покажем, что тогда справедливо неравенство
(\ref{difconvmult6}). Воспользуемся неравенством, вытекающим из
неравенства треугольника,
$$
\Vert \tau(t_i) / \Vert \tau(t_i) \Vert  - \tau (t_{i-1}) / \Vert
\tau(t_{i-1} \Vert \Vert \leq
$$
$$
\Vert r'(t_i) / \sqrt{ 1+f'^2_t (r(t_i))}  - r'(t_{i-1}) / \sqrt{
1+f'^2_t (r(t_{i-1}))} \Vert +
$$
$$
\mid f'_t(r(t_i)) /  \sqrt{ 1+f'^2_t (r(t_i))}- f'_t(r(t_{i-1})) /
\sqrt{ 1+f'^2_t (r(t_{i-1}))} \mid .
$$
Так как $1 \leq  \sqrt{1+f'^2_t (r(t_i))} \leq \sqrt{1+L^2}$ для
всех $ t_i \in [0,T_r]$ , то очевидно, существует такое $c_3 >1,$
для которого верно неравенство

\begin{equation} \Vert r'(t_i) / \sqrt{1+f'^2_t (r(t_i))}-
r'(t_{i-1}) / \sqrt{1+f'^2_t (r(t_{ i-1}))} \Vert \leq c_3 \Vert
r'(t_i) - r'(t_{i-1}) \Vert. \label{difconvmult7} \end{equation}

Из свойств функции $\theta(x)= x / \sqrt{ 1+x^2}$ следует
неравенство
$$
\mid f'_t (r(t_i)) / \sqrt{ 1+f'^2_t (r(t_i))} - f'_t (r(t_{i-1}))
/ \sqrt{1+f'^2_t (r(t_{i-1}))} \mid \leq
$$
\begin{equation}
 \leq \mid f'_t (r(t_i)) - f'_t (r(t_{i-1})) \mid .
\label{difconvmult8} \end{equation}

Из (\ref{difconvmult7}) и (\ref{difconvmult8}) имеем
$$
\sup_{\{t_i \} \in N_R \cap [T_i,T_{i+1}]} \,\, \sum_i \, \Vert
\tau(t_i) / \Vert \tau(t_i) \Vert - \tau (t_{i-1}) / \Vert
\tau(t_{i-1}) \Vert \Vert \leq
$$
\begin{equation}
\leq c_3 (\bigvee (  r'; T_i,T_{i+1}) + \bigvee (\Phi' ;
T_i,T_{i+1}) ). \label{difconvmult9} \end{equation}

Так как по условию $ f(\cdot)-$  ПРВ функция, то согласно теореме
\ref{difconvmultexthm2}
$$
\bigvee (\Phi'; T_i,T_{i+1}) \leq c(D,f) (1+ \bigvee
(r';T_i,T_{i+1})) ,
$$
откуда с учетом (\ref{difconvmult9}) следует неравенство
(\ref{difconvmult6}). Необходимость доказана.

{\bf Достаточность}.  Пусть справедливо неравенство
(\ref{difconvmult6}). Покажем, что $f(\cdot,\cdot)$ - ПРВ функция.
Воспользуемся неравенством
$$
\Vert \tau(t_i) / \Vert \tau( t_i) \Vert  - \tau(t_{i-1}) /  \Vert
\tau (t_{i-1}) \Vert \Vert \geq \mid  f'_t (r(t_i)) / \sqrt{1+f'_t
(r(t_i))} -
$$
\begin{equation} - f'_t (r(t_{i-1})) / \sqrt {1+f'^2_t
(r(t_{i-1}))} \mid \label{difconvmult10}
\end{equation} Из свойств функции
$\theta(x) = x / \sqrt{1+x^2}$ и из $\Vert f'(z) \Vert \leq L$ для
всех $z \in D,$ где производная существует, следует существование
константы $ c_4(L) > 0,$ для которой
$$
\mid f'_t (r(t_i)) / \sqrt{1+f'^2_t (r(t_i))} - f'_t (r(t_{i-1}))
/ \sqrt{1+f'^2_t(r(t_{i-1}))}  \mid \geq
$$
$$
\geq c_4 \mid f'_t (r(t_i)) - f'_t (r(t_{i-1})) \mid,
$$
откуда с учетом (\ref{difconvmult10}) имеем
$$
c_{21}(D,f)+c_{22}(D,f) \bigvee (  r'; T_i,T_{i+1}) \geq
$$
\begin{equation}
\geq \sup_{ \{ t_i \} \subset N_R \cap [T_i,T_{i+1}]} \sum_i \Vert
\tau(t_i) / \Vert \tau( t_i) \Vert  - \tau(t_{i-1}) / \Vert \tau
(t_{i-1}) \Vert \Vert \geq c_4 \bigvee (\Phi';  T_i,T_{i+1}) .
\label{difconvmult11}
\end{equation}
Из (\ref{difconvmult11}) следует неравенство
$$
\frac{c_{21}(D,f)}{c_4} + \frac{c_{22}(D,f)}{c_4} \bigvee (  r';
T_i,T_{i+1}) \geq \bigvee (\Phi'; T_i,T_{i+1}).
$$
Из теоремы \ref{difconvmultexthm2} следует, что $f(\cdot)$ - ПРВ
функция. Достаточность доказана. $\Box$

Из доказанных теорем следует ответ на один вопрос, интересный для
математиков разных специальностей. Будет ли  из локальной
представимости функции $f(\cdot)$ в виде разности выпуклых
следовать ее глобальная представимость? Иначе говоря, если функция
$f(\cdot)$ представима в виде разности выпуклых в некоторой
окрестности произвольной точки $x$ выпуклого компактного множества
$D$, то является ли $f(\cdot)$ ПРВ функцией на всем $D$? Ответ на
этот вопрос дает следующая теорема.

\begin{thm}
Если функция $f(\cdot)$ представима в виде разности выпуклых в
окрестности произвольной точки $x \in \bint Q, $ где $Q -$
выпуклое компактное множество в $\mathbb{R}^n$ с непустой
внутренностью, то она будет представима в виде разности выпуклых
на всем множестве $Q$. \label{difconvmultexthm6}
\end{thm}
{\bf Доказательство.} Покроем множество $Q$ окрестностями для
каждой точки $x \in D$, где по условию функция $f(\cdot)$
представима в виде разности выпуклых. В связи с компактностью
множества $Q$ из такого покрытия можно выбрать конечное
подпокрытие окрестностями $Q_i, i \in 1:N,$ $Q \subset \cup_{i=1}
^{N} Q_i$. Возьмем произвольную кривую $r \in \tilde{\wp}(D)$.
Кривая $r(\cdot)$ пересекает конечное число окрестностей
выбранного подпокрытия. Обозначим эти окрестности через $\{ Q_j
\}, j \in J$. Для каждого $j \in J$ пересечение $r(\cdot) \cap
Q_j$ представляет собой часть кривой $r(\cdot)$, принадлежащей
множеству кривых $r \in \tilde{\wp}(Q_j)$, а поэтому  для
некоторой константы $c_j(D_j,f)$ верно неравенство
$$
\bigvee (\Phi'; T_{j_1},T_{j_2})<  c_j(D_j,f)(1 +
\vee(r',T_{j_1},T_{j_2})),
$$
где $T_{j_1},T_{j_2}$- значение параметра $t$ при натуральной
параметризации части кривой $r(\cdot)$, принадлежащей окрестности
$Q_j$, $\Phi(t)=f(r(t))$. Всю кривую $r(\cdot)$ можно покрыть
конечным числом пересечений $r \cap Q_j$, каждое из которых
принадлежит $\tilde{\wp}(Q_j), j \in J$. Следовательно, верно
неравенство
$$
\sum_{j=1}^N \bigvee (\Phi'; T_{j_1},T_{j_2})<  \sum_{j =1}^{N}
c_j(Q_j,f)(1 + \vee(r',T_{j_1},T_{j_2})) <
$$
$$
< C   (1 + \sum_{j =1}^{N} \vee(r',T_{j_1},T_{j_2})),
$$
для константы $C$, независящей от выбранной кривой $r$, что
согласно Теореме \ref{difconvmultexthm4} и означает, что функция
$f(\cdot)$ является ПРВ функцией на всем множестве $Q$. Теорема
доказана. $\Box$

Отсюда следует вывод, что если функция $f(\cdot) $  не является
ПРВ (DC) функцией на большом компактном  множестве $Q$, то
обязательно найдется точка $x \in \mbox{ int } Q$ с произвольно
малой окрестностью $\Delta $, где функция  $f(\cdot) $  не
является ПРВ функцией.

\vspace{0.5cm}

\subsection{Поиск более узкого класса кривых $r(\cdot)$,
характеризующих ПРВ функции от $n$ переменных}

\vspace{0.5cm}

В этом разделе мы укажем более узкий класс кривых $r(\cdot)$, с
помощью которого можно сформулировать необходимые и достаточные
условия представимости произвольной функции от $n$ переменных в
виде разности выпуклых. Все обозначения, введенные ранее, остаются
в  силе.

Введем класс всех непрерывных замкнутых кривых $\hat{\wp}(D)$,
параметризованных естественным образом, проекция которых на  одну
из координатных плоскостей $\Pi_i, \, i \in 1:I,$ образованную
двумя векторами-ортами, есть кривая из множества $\wp(Pr(D))$, где
$Pr(D)- $ проекция множества $D$ на ту же координатную плоскость
$\Pi_i$. Последнее означает, что проекции кривых $r(\cdot)$ на
координатную плоскость $\Pi_i$ ограничивают на $\Pi_i$ звездные
множества.

Ясно, что $\hat{\wp}(D) \subset \tilde{\wp}(D)$, так как в
определении множества $ \tilde{\wp}(D)$ участвуют произвольные
плоскости $\Pi$, проходящие через начало координат, а в
определении $\hat{\wp}(D)$ - только координатные плоскости.

Как и ранее для кривой $r(\cdot) \in \hat\wp(D) $ выделим те ее
отрезки параметров $[T_i, T_{i+1}] \subset [0,T(r)]$, где векторы
$r(t), r'(t)$, $ t \in [0, T_r],$ образуют с какой-то координатной
плоскостью $\Pi_i, i \in 1:I,$ угол $\delta $: $0 < \delta \leq
\frac{\pi }{4} $. Обозначим объединение таких участков через $\hat
\Lambda_i(r)$.

Наша задача доказать, что кривых множества $ \hat{\wp}(D)$ и
системы отрезков $\hat \Lambda(r) $ достаточно для ответа на
вопрос о представимости функции в виде разности выпуклых.

\begin{thm} Для того, чтобы липшицевая функция $z \rightarrow
f(z):\mathbb{R}^n \rightarrow \mathbb{R}$ была представима в виде
разности выпуклых функций (была ПРВ функцией) на $D$, необходимо и
достаточно,  чтобы для любой координатной плоскости $\Pi_i$ и
любой кривой $r(\cdot) \in  \hat\wp(D)$, параметризованной
естественным образом с параметром $t \in [0,T], \, T=T(r),$ и
любого подмножества параметров $[T_i, T_{i+1} ] \subset
\hat\Lambda_i(r) $
$$
 (\exists c(D,f)>0) (\forall r \in \hat{\wp}(D)) \;\;\;\;
$$
$$
 \sum_1^m \vee (\Phi'; T_i,T_{i+1})<  c(D,f) (1+  \sum_1^m\vee (r';
 T_i,T_{i+1}))
$$
где $\Phi(t)=f(r(t)) \;\;\;\; \forall t \in [0,T_r]$.
\label{difconvmultexthm4}
\end{thm}

{\bf Доказательство.} Доказательство теоремы будет полностью
повторять доказательство Теоремы \ref{difconvmultexthm2}, если
показать, что любая кривая $r(\cdot) \in \tilde \wp(D)$ также
принадлежит $\hat \wp(D)$ и  для некоторой координатной плоскости
$\Pi_i, i \in 1:I,$ проекции кривых $r(\cdot)$ на $\Pi_i $
принадлежат $\wp(Pr(D))$. Единственно что надо показать, так это
то, что координатных плоскостей $\Pi_i, i \in 1:I,$ и множества
отрезков из $\hat\Lambda_i(r)$ хватает для ответа на вопрос:
представима ли липшицевая функция $f(\cdot): \mathbb{R}^n \rar
\mathbb{R} $ в виде разности выпуклых функций.

Возьмем произвольную кривую $r(\cdot) \in \tilde\wp(D)$ и
плоскость $\Pi, 0 \in \Pi$, проекция на которую кривой $r(\cdot)$
принадлежит $\wp(Pr(D))$. По свойству локальности в вопросе о
представимости функции в виде разности выпуклых, а также из
доказательства Теоремы \ref{difconvmultexthm6} следует, что можно
взять плоскость $\Pi$, образованную векторами
$$
a=\lim_{t_k \rar t_0 }r(t_k)
$$
и
$$
b=\lim_{t_k \rar t_0 }r'(t_k)
$$
для некоторой последовательности точек $t_k \in [0, T(r)]$. Также
будем считать, что все векторы $r(t_k), \, r'(t_k) $ для $t_k$ из
некоторой малой окрестности $\Upsilon(t_0)$ точки $t_0 \in
[0,T(r)]$ также принадлежат плоскости $\Pi$.

Ясно, что плоскость $\Pi$ образует с некоторой координатной
плоскостью $\Pi_i$ угол $\gamma$: $0 < \gamma \leq \frac{\pi}{4}$.
Но  тогда проекция части кривой $r(\cdot)$ и без ограничения
общности можно считать и вся кривая ограничивают на координатной
плоскости $\Pi_i$ также звездное множество.

Поскольку угол $\gamma$ удовлетворяет написанному выше
неравенству, то, как нетрудно видеть, векторы $r(t_k)$ и $ r'(t_k)
$ для указанных $t \in \Upsilon(t_0)$  также образуют углы из
отрезка $[0, \frac{\pi }{4}]$ с координатной плоскостью $\Pi_i$, о
которой речь шла выше. Кроме того, проекция $r(\cdot)$ для $t \in
\Upsilon(t_0)$  на координатную плоскость $\Pi_i$ принадлежит
$\wp(Pr(D))$.  Отсюда следует, что можно брать кривую $r(\cdot)$
из множества $\hat \wp(D)$ и достаточно проверять неравенство
теоремы только для таких кривых.

Верны и обратное утверждение. А именно: возьмем  кривую $r(\cdot)
\in \hat\wp(D)$ и координатную плоскость $\Pi_i$, проекция на
которую кривой $r(\cdot)$ принадлежит $\wp(Pr(D))$.

По Теореме \ref{difconvmultexthm6}, согласно которой вопрос о
глобальной представимости функции в виде разности выпуклых
сводится к вопросу о локальной ее представимости. Из
доказательства теоремы \ref{difconvmultexthm2} следует, что без
ограничения общности можно считать, что все векторы $r(t), \,
r'(t) $ для $t$ из малой окрестности $\Upsilon(t_0)$ точки $t_0
\in [0,T(r)]$ принадлежат некоторой плоскости $\Pi$.

Ясно, что найдется координатная плоскость $\Pi_i$, образующая угол
$\gamma$: $0 < \gamma \leq \frac{\pi }{4}$ с плоскостью $\Pi$,
проекция на которую кривой $r(\cdot) $ будет из множества кривых
$\wp(Pr(D))$.  Но тогда получаем, что кривая $r(\cdot)$
принадлежит множеству $\tilde \wp(D)$, поскольку из звездности
проекции множества на плоскость $\Pi_i $ следует звездность
проекции множества на плоскость $\Pi$.

Повторив доказательство Теоремы \ref{difconvmultexthm2}, мы
получим доказательство Теоремы \ref{difconvmultexthm4}. $\Box $

%%%%%%%%%%%%%%%%%%%%%%%%%%%%%%%%%%%%%%%%%%%%%%%%%%%%%%%%%%%%%%%

Теперь можно воспользоваться предыдущими  результатами. Из Теоремы
\ref{difconvmultexthm3} и Теоремы \ref{difconvmultexthm4} следует

\begin{thm} Для того, чтобы произвольная липшицевая функция
$z \rightarrow f(z) :\mathbb{R}^n \rightarrow \mathbb{R}$ была ПРВ
функцией, необходимо и достаточно, чтобы для любой координатной
плоскости $\Pi_i$ и любой $r(\cdot) \in \hat\wp(D)$,
параметризованной естественным образом: $t \in [0, T(r)]$, а также
любого подмножества параметров $[T_i, T_{i+1} ] \subset
\hat\Lambda_i(r) $ существовали константы $c_{31}(D,f),
c_{32}(D,f)
>0$ такие, что сумма поворотов $O_{R, i}$ кривой $R(\cdot)$ для $t
\in [T_i, T_{i+1}] $ для  всех $i$ на $\Gamma_f$ была ограничена
сверху неравенством
$$
\sum_1^m O_{R, i } \leq c_{21}(D,f) + c_{22}(D,f) \sum_1^m\vee
(r'; T_i,T_{i+1}) \,\,\,\,\forall r(\cdot) \in \hat \wp(\Pi)
$$
\label{difconvmultexthm5}
\end{thm}

Приведем, наконец, пример липшицевой функции $f(\cdot): \mathbb{R}
\rar \mathbb{R}$, которую легко распространить на пространства
большей размерности, не представимой в виде разности выпуклых
функций.

\begin{ex}\hspace{-2mm}. Граф функции $f:\mathbb{R} \rar \mathbb{R}$
состоит из сегментов с тангенсами углов наклона $\pm 1,$
расположенных между кривыми $-x^2, +x^2$. Функция $f(\cdot)$ не
представима в виде разности двух выпуклых функций в окрестности
точки $0,$ поскольку $\bigvee ^a _0 f' = \infty $ для
произвольного $a>0.$
 \label{ex1322}
\end{ex}

\subsection{Заключение}

Проблема, которую решает автор, является далеко не простой. Ею
интересовались и интересуются крупные ученые  в России и за
рубежом еще с начала 20-ого столетия с момента введения выпуклых
функций. Тогда же были получены необходимые и достаточные условия
представимости в виде разности выпуклых функций для одномерного
случая. Никаких существенных публикаций для многомерного случая с
тех пор не было.  В 40-ых годах проблема о представимости в виде
разности выпуклых была четко сформулирована и опубликована в
\cite{aleksandrov1}. Выпуклые и ПРВ функции находят широкое
применение в геометрии и в оптимизации \cite{strekal},
\cite{strecal2} благодаря своим хорошим свойствам.

Если для липшицевой функции от одной переменной необходимым и
достаточным условием ее представимости на отрезке $[a,b]$  в виде
разности выпуклых функций является ограниченность вариации ее
производной на этом отрезке, то для многомерного случая это только
необходимое условие, но не достаточное. Автор упорно искал класс
кривых, по которым можно было бы ответить на вопрос, является ли
функция ПРВ или нет. Первыми кандидатами для функции от двух
переменных стали кривые, ограничивающие выпуклые компактные
множества. В работе \cite{proudconvex2} используются такие кривые.

Авторы статьи \cite{veselyzajicek} дали пример функции
$f(\cdot):\mathbb{R}^2 \rar \mathbb{R}$ и фактически доказали в
этой статье, что класса кривых $r(\cdot)$ на плоскости $XOY$,
ограничивающих выпуклые компактные множества в области $D$, о
которых говорилось в \cite{proudconvex2}, недостаточно для ответа
на вопрос о представимости функции в виде разности выпуклых. При
доказательстве достаточности Теоремы 1 в случае б) статьи
\cite{proudconvex2} найти кривую (или кривые) $r(\cdot)$,
ограничивающую выпуклое компактное множество и состоящую из
отрезков $v_{k,n}$, сумма вариаций вдоль которых производной
функции $\Phi(t)=f(r(t))$ бесконечна для не ПРВ функции
$f(\cdot)$, не всегда, как оказалось, удается. Но тогда возникает
вопрос: "А какой класс кривых надо взять?"

Автору пришла идея рассмотреть кривые на плоскости, ограничивающие
звездные множества. Но для таких кривых $r(\cdot)$ вариация
производной $r'(\cdot)$ может оказать неограниченной сверху.
Поэтому вариация $\Phi'(\cdot)$ тоже может оказаться
неограниченной.  Но если мы расширим класс кривых и переопределим
правило подсчета вариации $\Phi'(\cdot)$  вдоль кривой $r(\cdot)
\in \hat\wp(D)$, то можно сформулировать необходимые и достаточные
условия представимости функции в виде разности выпуклых (см.
Теоремы \ref{difconvmultexthm4}, \ref{difconvmultexthm5}).

Так как вопрос о представимости функции в виде разности выпуклых
сводится к локальной ее представимости и векторы $r(t), r'(t)$ для
любого $t$ образуют плоскость $\Pi$, составляющую с какой- то
координатной плоскостью угол, не больший $\pi / 4$, то достаточно
рассматривать только координатные плоскости, образованные
координатными осями. Все сказанное отражено в формулировках Теорем
\ref{difconvmultexthm4}, \ref{difconvmultexthm5}.

Поскольку координатные плоскости $\Pi_i$ для двумерного случая
совпадают с $\mathbb{R}^2$, то условие представимости функции от
двух переменных в виде разности выпуклых можно сформулировать в
виде Теоремы \ref{difconvexthm1}

\newpage

\vspace{1cm}

\section{АЛГОРИТМЫ ПРЕДСТАВЛЕНИЯ ФУНКЦИИ В ВИДЕ РАЗНОСТИ
ВЫПУКЛЫХ}

\vspace{0.5cm}

Вопрос о методах представления функции в виде разности выпуклых
интересовал специалистов разных специальностей давно. Объясняется
это хорошими свойствами выпуклых функций для их  практического
применения. Первые результаты на эту тему были написаны в работе
\cite{aleksandrov1}. На этом все и закончилось. Как только мы
имеет дело с негладкой функцией, вопрос о способах представления
этой функции в виде разности выпуклых остается открытым. Ниже
будут описаны наиболее часто встречаемые в практике методы.

\subsection{Дважды непрерывно дифференцируемая функция}

Начнем с наиболее распространенного случая, когда функция
$f(\cdot): \mathbb{R}^n \rar \mathbb{R} $ дважды непрерывно
дифференцируемая.

Пусть $D-$ выпуклая область в $ \mathbb{R}^n$, замыкание которой
компактно, т.е. $\bar D-$ компакт.

Обозначим матрицу вторых смешанных производных функции $f(\cdot)$,
определенной в $D$, через $f''(\cdot) $, или $\nabla^2 f(\cdot)$.
Пусть  в области $D$ функция $f(\cdot) $ дважды непрерывно
дифференцируема с матрицей вторых смешанных производных $
f''(\cdot)$, или $\nabla^2 f( \cdot) $. Пусть
$$
\| f''(x) \| \leq L  \,\,\,\,\forall x \in D.
$$
Представим $f(\cdot)$ в виде разности выпуклых следующим образом.
В качестве одной из выпуклых функций возьмем функцию $f_1(\cdot) =
L \| x \|^2. $ Покажем, что разность
$$
  f_2(x)= L \| x \|^2 - f(x)
$$
есть выпуклая функция в $D$.

Для этого найдем матрицу вторых смешанных производных функции
$f_2(\cdot)$. Имеем
$$
  f_2''(x)=2L-f''(x).
$$
Оценим значения квадратичной  функции
$$
  (f_2''(x)g,g)=2L-(f''(x)g,g) \geq 0 \,\,\,\,\,\forall g \in
  S^{n-1}_1(0)= \{ g \in \mathbb{R}^n \vl \| g \| \leq 1 \}.
$$
поскольку по теореме Гершгорина \cite{voevodin} все собственные
значения матрицы $f_2''(x)$ для любого $x \in D$ положительные.
Отсюда следует, что матрица $f_2''(x)$ положительно определенная в
$D$. Следовательно, функция $f_2(\cdot)-$ выпуклая в $D$, и верно
разложение
$$
  f(x)=f_1(x)-f_2(x) \,\,\,\, \forall x \in D,
$$
где $f_1(\cdot), \, f_2(\cdot)-$ выпуклые, что и требовалось. То
есть $f(\cdot-)$ ПРВ функция.

\vspace{0.5cm}

\subsection{Функции с липшицевой первой производной}

\vspace{0.5cm}

Докажем теорему о представлении функции в виде разности выпуклых в
более общем случае.
\begin{thm} \cite{aleksandrov1}
Если функция $f(\cdot):\mathbb{R}^n \rar \mathbb{R}$ имеет
липшицевую первую производную  в области $D$ с константой Липшица
$L$, то в области $D$ функция $f(\cdot)$ представима в виде
разности выпуклых функций. \label{diffconvAlgorithmThm1}
\end{thm}
{\bf Доказательство}. Пусть функция $f(\cdot)$ удовлетворяет
условиям теоремы в области $D$. Рассмотрим разность
$$
  f_2(x)=L \, \| x \|^2 - f(x).
$$
Если мы докажем, что $f_2(\cdot)-$ выпуклая, то отсюда будет
следовать, что $f(\cdot)$ представима в виде разности выпуклых
функций $f_1(x)= L \, \| x \|^2$ и $f_2(x)$.

Для доказательства воспользуемся характерным свойством выпуклой
функции.

\begin{thm}\cite{pshenichnyi} Для того чтобы функция $z \rar \theta(z):
\mathbb{R}^n \rar \mathbb{R}$ была выпуклой, необходимо и
достаточно, чтобы она была дифференцируема по направлениям и для
любой точки $z \in \mathbb{R}^n $ и любого направления $p \in
\mathbb{R}^n $ функция $\alpha \rar h( \alpha ): \mathbb{R}^+ \rar
\mathbb{R}:$
$$
h(\alpha) = \frac{\p \theta (z + \alpha p)}{\p p }
$$
была неубывающей по $\alpha >0.$
\label{diffconvAlgorithmThm2}
\end{thm}

 Возьмем произвольное направление $g \in \mathbb{R}^n$. Определим
 функцию $f_1(x)=L \, \| x \|^2$. Рассмотрим
 разность
 $$
   f_2(x) = f_1(x) - f(x).
 $$
 Покажем, что $f_2(\cdot)-$ выпуклая функция. Используем Теорему
 \ref{diffconvAlgorithmThm2}. Докажем, что производная по
 направлению $g \in \mathbb{R}^n$ функции $f_2(x)$ монотонно
 возрастающая вдоль этого направления. .

 Покажем, что для $\al >0$  разность производных по направлению
 $$
  (f'_2(x+\al g), g)- (f'_2(x), g).
 $$
 положительная, где круглые скобки означают скалярное произведение векторов.
 Имеем
 $$
   \mid  (f'_2(x+\al g), g)- (f'_2(x), g) \mid =
 $$
 $$
  = \mid 2 \, L (x+\al g, g)- (f'(x+\al g), g) -
  2 \, L (x, g) + (f'(x), g) \mid.
 $$
 Отсюда делаем оценку для разности производных по направлению $g$
 вдоль направления  $g$
 $$
   \mid 2 L \, \al \| g \|^2 - (f'(x+\al g),g)  + (f'(x), g) \mid \geq
   \mid 2 L \, \al \| g \| ^2 - L \, \al \| g \|^2 > 0,
 $$
т.е. производная по направлению $g$ функции $f_2(\cdot)$ монотонно
возрастает вдоль направления $g$. По Теореме
\ref{diffconvAlgorithmThm2} получаем, что функция $f_2(\cdot)-$
выпуклая в $D$. Итак, получили, что функция $f(\cdot)$ представима
в виде разности двух выпуклых функций $f_1(\cdot)$ и $f_2(\cdot)$.
Теорема доказана. $\Box$

\vspace{0.5cm}

\subsection{Представление недифференцируемых функции в виде разности
выпуклых }

\vspace{0.5cm}

Алгоритма представления недифференцируемой функции в виде разности
выпуклых в общем случае не существует. Даже в случае, когда у
функции $f(\cdot)$ есть только одна точка недифференцируемости, мы
не знаем, как ее представить в виде разности выпуклых. Есть только
процессы, дающие в пределе (при выполнении написанных выше
условий) функции, разность которых есть исходная функция. 

Рассмотрим несколько частных случаев. Пусть $f_N(\cdot): D \rar
\mathbb{R}^n -$ многогранная функция, т.е. функция, график которой
состоит из конечного число $N$ частей гиперплоскостей. Мы
воспользуемся алгоритмом, описанным в \cite{aleksandrov1}.

Согласно терминологии А.Д.Александрова под многогранной
кусочно-линейной функцией с конечным числом граней будем понимать
такую функцию, график которой состоит из конечного числа частей
плоскостей (гиперплоскостей), которые называются гранями.

Введем понятие {\em двугранного угла}. Будем понимать под
двугранным углом функцию, определенную на $D$, равную максимуму
или минимуму линейных функций, графиками которых являются
гиперплоскости $\pi_k(\cdot)$ и $\pi_{l}(\cdot)$,  имеющие общие
грани размерности $n-1$. Если такая функция выпуклая, то
двугранный угол назовем выпуклым, если эта функция вогнутая, то
двугранный угол назовем вогнутым.

Рассмотрим все выпуклые двугранные углы, части графиков которых
принадлежат графику функции $f_N(\cdot)$. Доопределим эти
двугранные углы на всю область $D$. Просуммируем все такие
выпуклые двугранные углы. В итоге получим выпуклую многогранную
функцию $f_{1,N}(\cdot):D \rightarrow \mathbb{R}$. Докажем, что
разность
\begin{equation}
 f_{1,N}(\cdot) - f_N(\cdot) = f_{2,N}(\cdot)
\label{difconvmult1} \end{equation} есть также выпуклая
многогранная функция. Для полноты изложения я повторю
доказательство, приведенное ранее в предыдущей главе.

Для доказательства достаточно показать, что все двугранные углы,
части графиков которых размерности $n$ принадлежат графику функции
$f_{1,N}(\cdot) - f_{N}(\cdot),$ есть выпуклые. Для этого покажем,
что любая точка, лежащая на проекции на $D$ пересечения
$\pi_{kl}(\cdot) = \pi_k(\cdot) \cap \pi_{l}(\cdot)$ произвольных
гиперплоскостей $\pi_k(\cdot)$ и $\pi_{l}(\cdot)$, образующих
график двугранного угла функции $f_{1,N}(\cdot) - f_{N}(\cdot),$
имеет малую окрестность, где функция $f_{1,N}(\cdot) -
f_{N}(\cdot)$ выпуклая.

Если берем точку, в малой окрестности которой функция
$f_{N}(\cdot)$ линейная, то локальная выпуклость разности
$f_{1,N}(\cdot) - f_{N}(\cdot)$ очевидна. Пусть берем точку,
лежащую на проекции на $D$ множества $\pi_{kl}(\cdot)$ выпуклого
двугранного угла, часть графика которого принадлежит графику
функции $f_{N}(\cdot)$. Поскольку согласно алгоритму двугранный
угол, график которого образован гиперплоскостями $\pi_k(\cdot)$ и
$\pi_{l}(\cdot)$, входит в сумму выпуклых двугранных углов,
образующих функцию $f_{1,N}(\cdot)$, то опять разность
$f_{1,N}(\cdot) - f_{N}(\cdot)$ будет локально выпуклой в
окрестности рассматриваемой точки. Если же точка лежит на проекции
на $D$ множества $\pi_{kl}(\cdot)$ вогнутого двугранного угла,
часть графика которого принадлежит графику функции $f_{N}(\cdot)$,
то $-f_{N}(\cdot)$ $ -$ локально выпуклая в окрестности этой
точки, а поэтому разность $f_{1,N}(\cdot) - f_{N}(\cdot)$ снова
локально выпуклая в той же окрестности. Из локальной выпуклости
всех двугранных углов функции $f_{1,N}(\cdot) - f_{N}(\cdot)$
следует ее выпуклость на всем множестве $D$, что и требовалось
доказать.

Рассмотрим второй случай, когда  $f(\cdot): \mathbb{R} \rar
\mathbb{R}^n$ есть липшицевая дифференцируемая во всех точках
функция за исключением лишь одной точки $x_0$. В точке $x_0$
функция $f(x)$ дифференцируема по направлениям.

Построим конус $K(x_0, f(x_0)) $ касательных направлений к графику
$\Gamma_f $ функции $f(\cdot)$ в точке в точке $(x_0, f(x_0))$.
Обозначим функцию, графиком которой  является конус $K(x_0,
f(x_0)) $, через $\psi(\cdot)$. Ясно, что  функция
$$\varphi(x)=\psi(x+x_0)-f(x_0)$$
будет положительно-однородной функцией (п.о.): $\varphi(0)=0$,
$\varphi(\lambda g) = \lambda \varphi(g)$, $\lambda >0$ и
$$
  \varphi(g)=\frac{\p f(x_0)}{\p g}.
$$
Для представления функции $f(\cdot)$ в виде разности выпуклых на
{\em первом шаге} надо проверить, представима ли функция
$\varphi(\cdot)$ в виде разности выпуклых, т.е. является ли она
ПРВ функцией. Это надо сделать с помощью теорем, сформулированных
в предыдущих главах. Если $\varphi(\cdot)$ не является ПРВ
функцией, то функция $ f(\cdot)$  также не будет ПРВ функцией. Это
происходит потому, что в случае представимости функции $f(\cdot)$
в виде разности выпуклых $f_1(\cdot)$ и $f_2(\cdot)$, будет верно
равенство
$$
   \varphi(g)= \varphi_1(g) -  \varphi_2(g),
$$
где функции  $\varphi_!(\cdot)$  и $ \varphi_2(\cdot)$ для функций
$f_1(\cdot)$ и $f_2(\cdot)$ соответственно определяются аналогично
тому, как мы определяли $\varphi(\cdot)$ для функции $f(\cdot)$.
Но функции  $\varphi_!(\cdot)$  и $ \varphi_2(\cdot)$, построенные
для выпуклых функций $f_1(\cdot)$ и $f_2(\cdot)$, являются
выпуклыми.

На {\em втором шаге} надо построить функцию $\hat f(\cdot) =
f(\cdot) -\psi(\cdot)$. Если $f(\cdot)-$ ПРВ функция, то функция
$\hat f(\cdot)$ будет непрерывно дифференцируемой в точке $x_0$.
Действительно, если
$$
   f(x)=f_1(x)-f_2(x),
$$
где  $f_1(\cdot)$ и $f_2(\cdot)-$ выпуклые, то $\psi_1(\cdot)$,
$\psi_2(\cdot)-$ также выпуклые и
$$
 \hat f(x)=f(x)-\psi(x)=f_1(x)-f_2(x)-\psi_1(x)+\psi_2(x)=
$$
$$
 = (f_1(x)-\psi_1(x))-(f_2(x)-\psi_2(x)).
$$
Но функции $f_1(x)-\psi_1(x) $, $f_2(x)-\psi_2(x) -$ выпуклые,
производные по направлениям которых в точке $x_0$ равны нулю, а
значит субдифференциалы этих функций в точке $x_0$ по свойствам
выпуклых функций \cite{demvas}  состоят из одной точки $0$. Но
если субдифференциал выпуклой функции в некоторой точке состоит из
одного вектора, то в этой точке выпуклая функция непрерывно
дифференцируемая. Это следует из того, что для выпуклой функции
$\tilde f(\cdot) $ субдифференциал в точке $x_0$ есть
субдифференциал Кларка, который в свою очередь равен  предельным
значениям всех градиентов функции $\tilde f(\cdot) $, ыычисленных
в точках $x$, когда $x \rar x_0$ \cite{demvas}. А значит выпуклые
функции $f_1(x)-\psi_1(x) $ и $f_2(x)-\psi_2(x) $ непрерывно
дифференцируемы в точке $x_0$, и производные в точке $x_0$ равны
нулю.

На {\em третьем шаге} для ответа на вопрос о представимости
функции $f(\cdot)$ в виде разности выпуклых надо ответить на
вопрос о представимости в виде разности выпуклых функции $\hat
f(\cdot)$.

Если  функция $f(\cdot)$ имеет липшицевую с константой $L$
производную на множестве $\mathbb{R}^n \backslash
B_{\delta}(x_0)$, где $B_{\delta}(x_0)-$ $\delta$ окрестность
точки $x_0$, $\delta-$ произвольно малое положительное число, то
функция $\hat f(\cdot)$ также будет иметь липшицевую производную
на $\mathbb{R}^n \backslash B_{\delta}(x_0)$, так как липшицевую
производную на $\mathbb{R}^n \backslash B_{\delta}(x_0)$ будет
иметь функция $\psi(\cdot)$. Докажем сказанное.

Сначала покажем, что $\varphi(\cdot)$ и $\psi(\cdot)-$ липшицевые
с константой $L$.  Действительно,
$$
 \vl \varphi(g_1)-\varphi(g_2) \vl = \vl \lim_{\al \rar 0^+}
 \frac{f(x_0+\al g_1)-f(x_0)}{\al} - \lim_{\al \rar 0^+}
 \frac{f(x_0+\al g_2)-f(x_0)}{\al} \vl=
$$
$$
= \vl \lim_{\al \rar 0^+}  \frac{f(x_0+\al g_1)-f(x_0+\al
g_2)}{\al} \vl \leq  L \| g_1 - g_2 \|,
$$
т.е. функция $\varphi(\cdot)$, а значит и $\psi(\cdot)-$
липшицеввые с константой $L$.

Докажем, что функция $\varphi(\cdot)$, а значит и $\psi(\cdot)$
представимы  в виде разности выпуклых функций, если  функция
$f(\cdot)$ имеет липшицевую с константой $L$ производную на
множестве $\mathbb{R}^n \backslash B_{\delta}(x_0)$.

Определим п.о. функцию $\varphi_{\delta}(\cdot): \mathbb{R}^n \rar
\mathbb{R}$, равную нулю в нуле и имеющей те же значения на сфере
$\| g \| = \delta$, что и функция $f(\cdot) - f(x_0) $ на сфере
$\| x-x_0 \| = \delta $. Здесь $\delta>0$, как ранее,  произвольно
малое число, $g-$ переменная функции $\varphi_{\eps}(\cdot)$, а
$x-$ переменная функции $f(\cdot)$. Функция $\varphi_{\eps}(\cdot)
$ будет иметь липшицевую производную с константой $L$ на множестве
$\mathbb{R}^n \backslash B_{\delta}(0)$, как и функция $f(\cdot)$
на множестве $\mathbb{R}^n \backslash B_{\delta}(x_0)$.

Известно, что выпуклость или вогнутость любой п.о. функции
определяется ее значениями на сфере произвольно малого радиуса.

%%%%%%%%%%%%%%%%%%%%%%%%%%%%%%%%%%%%%%%%%%%%%%%%%%%%%%%%%%

Согласно Теореме \ref{diffconvAlgorithmThm1} функция
\be
  \bar \varphi_{\delta} (g)= L \| g \|^2 -  \varphi_{\delta}(g)
  \,\,\, \forall g
  \in S^{n-1}_1(0) = \{ v \in \mathbb{R}^n \vl \| v \| =1 \}.
  \label{difconvmult2}
\ee является выпуклой на множестве  $\mathbb{R}^n \backslash
B_{\delta}(0)$. По функции $\bar \varphi_{\delta }(\cdot)$
построим п.о. первой степени функцию $\varphi_{2,\delta}(\cdot)$,
равную функции  $\bar \varphi_{\delta}(\cdot)$ на сфере
$S^{n-1}_{1}(0)$. Ранее, когда изучали задачу о представимости
п.о. степени $m$ функции в виде разности выпуклых, мы уже
проделывали подобную операцию и доказывали, что получившаяся
функция будет выпуклой. Поэтому $\varphi_{2,\delta}(\cdot)-$
выпуклая в шаре $B^n_1(0)$. Аналогично, по функции $L \| g \|^2$
также построим п.о. степени 1 функцию $ \varphi_1(\cdot)$, которая
тоже будет выпуклой в шаре $B^n_1(0)$ и значения которой равны
значениям функции $L \| g \|^2$ на сфере $S^{n-1}_1(0)$. Из
равенства (\ref{difconvmult2}) и из п.о. функции
$\varphi_{\delta}(\cdot)$ и  всех построенных функций
$\varphi_1(\cdot), \, \varphi_{2, \delta}(\cdot) $ получаем
представление функции $\varphi(\cdot)$ в виде разности выпуклых
п.о. функций \be
  \varphi_{\delta}(g)=\varphi_1(g) - \varphi_{2, \delta}(g) \,\,\,\,
  \forall g \in S^{n-1}_1(0).
  \label{difconvmult3}
\ee

Перейдем в равенстве (\ref{difconvmult3}) к пределу по $\delta
\rar +0$. Выпуклая п.о. функция $\varphi_1(\cdot)$, построенная по
$L \| g \|^2 $, постоянная и не зависит от $\delta$.  Функция
$\varphi_{\delta}(\cdot)$ равномерно по $\delta \rar +0$ на $B_1^n
(0) $ сходится к $\varphi(\cdot)$. Поэтому  выпуклые п.о. функции
$\varphi_{2,\delta}(\cdot)$ равномерно по $\delta \rar +0 $ на
единичном шаре $B_1^n(0)$ сходятся к выпуклой п.о. функции
$\varphi_{2}(\cdot)$. Из равенства (\ref{difconvmult3}) при
переходе к пределу по $\delta \rar +0$ следует равенство \be
  \varphi(g)=\varphi_1(g) - \varphi_{2}(g) \,\,\,\,
  \forall g \in S^{n-1}_1(0),
  \label{difconvmult4}
\ee где   $\varphi_1(\cdot), \varphi_{2}(\cdot)- $ выпуклые п.о.
функции, что и требовалось показать.

Производные $\varphi'_{2,\delta}(\cdot)$ равномерно ограничены на
множестве $B_1^n(0) \backslash B_{\delta}(x_0)$ и равностепенно
непрерывны по $\delta$, так как такими свойствами обладают
производные $ \varphi'_{\delta}(\cdot), \varphi'_1(\cdot)$  на том
же множестве. Поэтому в пределе по $\delta \rar 0^+ $ мы получим
функции $\varphi(\cdot), \varphi_1(\cdot), \varphi_{2}(\cdot) $,
которые будут иметь липшицевые с константой $L$  производные на
множестве $B_1^n(0) \backslash B_{\delta}(0)$.

%%%%%%%%%%%%%%%%%%%%%%%%%%%%%%%%%%%%%%%%%%%%%%%%%%%%%%

Если исходная функция $f(\cdot)$ имеет липшицевую производную с
константой $L$ на множестве $\mathbb{R}^n \backslash \{ x_0 \}$,
то функция $\hat f(\cdot)$, как разность двух функций с липшицевой
производной с константой $L$ на множестве $\mathbb{R}^n \backslash
\{x_0\}$ также будет иметь липшицевую производную с константой $L$
на всем $\mathbb{R}^n$, так как по доказанному выше функция $\hat
f(\cdot)$ непрерывно дифференцируема в точке $x_0$ и $\hat
f'(x_0)=0$.

По теореме, приведенной выше, мы можем сделать вывод, что функция
$\hat f(\cdot)-  $ ПРВ функция. Алгоритм ее представления в виде
разности выпуклых мы знаем (см. Теорему
\ref{diffconvAlgorithmThm1}). Согласно описанному алгоритму можно
представить функцию $f(\cdot)$ в виде разности выпуклых.

Итак, мы разобрали несколько примеров представления негладкой
функции в виде разности выпуклых. Рассмотренные случаи можно
распространить на более сложные с их различными комбинациями.

\newpage
\vspace{1cm}

\section{УСЛОВИЯ СХОДИМОСТИ ПОСЛЕДОВАТЕЛЬНОСТИ
КВАЗИДИФФЕРЕНЦИРУЕМЫХ ФУНКЦИЙ К КВАЗИДИФФЕРЕНЦИРУЕМОЙ ФУНКЦИИ}

\vspace{1cm}

В этом параграфе мы рассмотрим одно из многочисленных применений
доказанных теорем. Будут приведены достаточные условия, при
которых последовательность квазидифференцируемых  функций $\{
f_k(\cdot) \}, k = 1,2, \dots$ от двух переменных сходится
равномерно на компакте $D$ с $\mbox{int } D \neq \emptyset $ к
квазидифференцируемой функции $f(\cdot)$. Представляет интерес
распространить эти результаты на функции от произвольного
количества аргументов.

\vspace{0.5cm}

\subsection{Последовательность одномерных функций}

\vspace{0.5cm}

Пусть задана последовательность квазидифференцируемых функций $\{
f_k(\cdot) \}, k = 1,2, \dots$, определенных на компактном
множестве $D \subset \mathbb{R}^2$. Обозначим через $L_k$
константу Липшица функций $f_k(\cdot)$  на множестве $D$, а через
$f'_k(x,g)-$ ее производную в точке  $x$  по направлению $g \in
\mathbb{R}^2$. Предположим, что в любой точке $x \in \mbox{int} D$
функции $f_k(\cdot)$ квазидифференцируемые (КВД функции), т.е.
верны равенства
$$
f'_k(x,g)=\max_{v \in \underline{\p} f_k(x) } (v,g) + \min_{w \in
\overline{\p} f_k(x) } (w,g) \,\,\,\,\, \forall g \in
\mathbb{R}^2,
$$
где  $\underline{\p} f_k(x), \overline{\p} f_k(x) - $
субдифференциал и суппердифференциал функции $f_k(\cdot)$ в точке
$x$  соответственно, которые по определению являются выпуклыми
компактами на плоскости $\mathbb{R}^2. $

Пусть последовательность $\{ f_k(\cdot) \}, k = 1,2, \dots$
сходится равномерно на множестве $D$ к функции $f(\cdot)$.
Зафиксируем точку $x \in \mbox{int} D$. Спрашивается, при каких
условиях, наложенных на функции $\{ f_k(\cdot) \}$, функция
$f(\cdot)$ будет КВД в точке $x$?

Через $r(t)$  обозначим радиус-вектор единичной окружности
$S_1^1(0)=\{ z \in \mathbb{R}^2 | \| z \| =1 \}$, где $t \in [0,
2\pi]- $ угол поворота вектора $r(\cdot)$. Положим по определению
$$
\Phi_k (t)= f'_k(x,r(t)), \,\,\,\, \Phi (t)= f'(x,r(t)) \,\,\,\,\,
\forall t \in [0,2 \pi].
$$
Здесь $f'_k(x,r(t)), \, f'(x,r(t)) - $ производные по направлению
$r(t)$ функций $f_k(\cdot),\, f(\cdot) $ в точке $x$
соответственно.

Легко показать, что функции $ \Phi_k(\cdot)- $  липшицевые на
отрезке $[0,2\pi]$ с константой Липшица $L_k$. Поэтому для каждой
из них существует множество  ${\cal N}_k$   точек
дифференцируемости, которое всюду плотное, полной меры на отрезке
$[0, 2 \pi]$. Положим по определению \be \vee(\Phi'_k; 0,
2\pi)=\sup_{ \{ t_i \} \in {\cal N}_k } \sum_{i=1}^k |
\Phi'_k(t_{i+1}) - \Phi'_k(t_{i}) | \label{DiffConvSupp1} \ee
Рассмотрим теперь произвольную последовательность одномерных
функций $\Psi_k:\mathbb{R} \rar \mathbb{R}$, определенных на
отрезке $[a,b]$  и равномерно сходящихся на этом отрезке к функции
$\Psi(\cdot)$. Обозначим через $\aleph_k- $   множество точек
дифференцируемости функции $\Psi_k(\cdot)$ на отрезке $[a,b]$.
Предположим, что множества $\aleph_k $  для любого $k -$ всюду
плотные, полной меры на $[a,b]$. Аналогично (\ref{DiffConvSupp1})
определим $\vee(\Psi'_k; a, b) $.
\begin{thm}
Если существует такая константа $c>0$, что для любого $k$   верно
неравенство \be \vee(\Psi'_k; a, b) \leq c, \label{DiffConvSupp2}
\ee то множество точек дифференцируемости  $\aleph  $   функции
$\Psi(\cdot) $ полной меры, всюду плотное на отрезке $[a,b]$   и
верно неравенство
$$
\vee(\Psi'; a, b) \leq 2 c.
$$
\label{DiffConvSuppThm1}
\end{thm}
{\bf Доказательство}. Зафиксируем произвольное $k$. Известно
\cite{kolmogorovfomin}, что из условия  (\ref{DiffConvSupp2})
следует, что функция $\Psi_k(\cdot)  $ представима на отрезке
$[a,b]$ в виде разности выпуклых функций $\Psi_{1k}(\cdot)$ и
$\Psi_{2k}(\cdot)$ , т.е. \be \Psi_k(t) = \Psi_{1k}(t) -
\Psi_{2k}(t) \,\,\,\,\, \forall t \in [a,b]. \label{DiffConvSupp3}
\ee Представление (\ref{DiffConvSupp3})  неединственно. Однако,
существует такое представление, для которого верны неравенства \be
\vee(\Psi'_{1k}; a, b) \leq  c, \,\, \vee(\Psi'_{2k}; a, b) \leq
c. \label{DiffConvSupp4} \ee Например, можно положить
$$
\Psi'_{1k}(t) = \vee(\Psi'_k;a,t), \,\,\, \Psi_{1k}(t)= \int_a^t
\Psi'_{1k}(\eta) d \eta, \,\,\, \Psi_{2k}(t)= \Psi_{1k}(t) -
\Psi_{k}(t).
$$
Без ограничения общности будем считать, что функции
$\Psi_{1k}(\cdot)$ и  $\Psi_{2k}(\cdot)$  на  $[a,b]$  сходятся к
выпуклым функциям $\Psi_{1}(\cdot)$   и  $\Psi_{2}(\cdot)$
соответственно. Из (\ref{DiffConvSupp4})   и свойств выпуклых
функций  \cite{pshenichnyi}  следует, что
$$
\vee(\Psi'_1; a, b) \leq  c, \,\,\,\, \vee(\Psi'_2; a, b) \leq  c.
$$
В равенстве
$$
\Psi_{k}(t)= \Psi_{1k}(t) - \Psi_{2k}(t).
$$
перейдем к пределу по  $k \rar \infty$. Получим
$$
\Psi(t)= \Psi_{1}(t) - \Psi_{2}(t) \,\,\,\,\, \forall t \in [a,b].
$$
Отсюда следует, что функция  $\Psi(\cdot)$  почти всюду (п.в.)
дифференцируема на $[a,b]$, а также верно неравенство
\cite{kolmogorovfomin}
$$
\vee(\Psi'; a, b) \leq \vee(\Psi'_1; a, b)+ \vee(\Psi'_2; a, b)
\leq 2c,
$$
что и требовалось доказать. Теорема доказана. $\Box$

\vspace{0.5cm}

\subsection{Последовательность квазидифференцируемых функций}

\vspace{0.5cm}

Применим теперь Теорему \ref{DiffConvSuppThm1} для вывода условий
квазидифференцируемости функции  $f(\cdot)$  в точке $x$.
Воспользуемся следующей теоремой, доказанной нами ранее в
предыдущих параграфах.
\begin{thm} {\em \cite{proudconvex1}, \cite{lupikovdissertation}}
Для того, чтобы липшицевая функция $f(\cdot): \mathbb{R}^2 \rar
\mathbb{R}$ была КВД функцией в точке $x$, необходимо и
достаточно, чтобы
$$
 (\exists c(f)>0):  \;\;\; \vee (\Phi'; 0,2\pi)<  c(f),
$$
где $\Phi(t)=f'(x, r(t)), \,\, r(t)=(\cos t, \sin t) \in
\mathbb{R}^2, \,\, t \in [0, 2\pi] $. \label{DifconvGomogThm1}
\end{thm}

Воспользовавшись этой теоремой, можно получить следующий
результат.
\begin{thm}
Если последовательность функций  $\{ f_{k}(\cdot) \}$ равномерно
на $D \subset \mathbb{R}^2$, $x \in \mbox{int  } D,$ сходится к
функции $f(\cdot)$ и имеют место условия

а) функции  $ \Phi_{k}(\cdot),  $ где $\Phi_{k}(t)= f'_{k}(x,
r(t)), $ равномерно при $k \rar \infty$  на отрезке $[0, 2\pi]$
стремятся к функции $\Phi(t)$,

б) существует такая константа  $c>0$, что
$$
\vee(\Phi'_{k}; 0, 2 \pi) < c   \,\,\,\,\,\, \forall k,
$$
то предельная функция $f(\cdot)$  является КВД в точке $x$.
\label{DiffConvSuppThm3}
\end{thm}

{\bf Доказательство.} Возьмем в качестве функции $\Psi_n(\cdot) $
, фигурирующей в Теореме \ref{DiffConvSuppThm1}, функцию
$\Phi_k(\cdot)$. Полагаем $a=0, \, b= 2\pi$. Нетрудно видеть, что
условия Теоремы \ref{DiffConvSuppThm1} для последовательности
$\Phi_k(\cdot)$ выполнены. Поэтому
$$
\vee(\Phi';0,2\pi) \leq c.
$$
Откуда из теоремы  \ref{DifconvGomogThm1} следует утверждение
теоремы. Теорема доказана. $\Box$

Условие  а) Теоремы  \ref{DiffConvSuppThm3}   не всегда
выполняется. приведем соответствующий пример.
\begin{ex}
Рассмотрим последовательность функций $\{ f_k(\cdot) \} $,
$f_k:\mathbb{R}^2 \rar \mathbb{R}$. Положим по определению
\begin{equation*}
f_k(x)  = \left\{ \begin{aligned}
            0, \mbox{   если} \,\,\, \| x \| \leq \frac{1}{k},  \\
            \| x \| - \frac{1}{k}, \mbox{   если} \,\,\,
             \| x \| > \frac{1}{k}.
\end{aligned} \right.
\end{equation*}
где $x \in \mathbb{R}^2$. Очевидно, что функции $\{ f_k(\cdot) \}$
равномерно на единичном круге $B_1^1(0)$ стремятся к функции
$f(x)=\| x \| $ при $k \rar \infty$. Однако $\Phi_k(t)=0$, для
всех $k$ и $\Phi(t)=1$ для $t \in [0, 2\pi]$. Видно, что не
существует последовательности $\{ \Phi_k(\cdot) \}$, для которой
выполняется условие а) теоремы \ref{DiffConvSuppThm3}.
\end{ex}

Попытаемся изменить условие а) теоремы  \ref{DiffConvSuppThm3}.

Рассмотрим теперь последовательность функций
$$
\varphi_{m}(g)=\frac{f_m(x+\al_m g) - f_m(x)}{\al_m} \,\,\,\,
\forall g \in S_1^2(0).
$$
Для других $g$ доопределяем функцию $\varphi_{m}(\cdot)$
положительно однородным образом. Ясно, что функции $
\varphi_{m}(\cdot) - $ липшицевые по своему аргументу. Положим по
определению
$$
\hat{\Phi}_{m}(t)=\varphi_{m}(r(t)) \,\,\,\,\, \forall t \in [0,
2\pi].
$$
Очевидно, что из последовательности функций $
\hat{\Phi}_{m}(\cdot) $    всегда можно выбрать
подпоследовательность  $\hat{\Phi}_{{m}_i}(\cdot) $  такую, что
функции $\hat{\Phi}_{{m}_i}(\cdot) $ равномерно сходятся к
$\Phi(t) $ на отрезке $[0, 2 \pi]$   к функции  $\Phi(\cdot)$.
Поэтому верна следующая теорема.
\begin{thm}
Если для всех индексов $m$    найдется константа  $c>0$ такая, что
$$
\vee(\hat{\Phi}'_{m}; 0, 2\pi) \leq c,
$$
то функция $f(\cdot)$, равная равномерному пределу
последовательности $\{ f_m (\cdot) \} $ на $D$,  является КВД в
точке $x$.
\end{thm}

Рассмотрим многомерный случай.

Пусть имеется последовательность функций КВЛ функций $f_m(\cdot)$
в точке $x$, равномерно сходящуюся в некоторой окрестности точки
$x$ к функции $f(\cdot)$. Возникает вопрос: при каких условиях
предельная функция $f(\cdot)$ будет КВД в точке $x$?

Введем, как и выше, последовательность п.о. функций $
\varphi_{m}(\cdot): \mathbb{R}^n \rar \mathbb{R}$
$$
\hat \varphi_{m}(g)=\frac{f_m(x+\al_m g) - f_m(x)}{\al_m} \,\,\,\,
\forall g \in S^{n-1}_1(0), \,\,
$$
а для других $g$ функция $\varphi_{m}(\cdot)$ определяется
положительно однородным образом
$$
  \varphi_m (g) = \| g \| \hat \varphi_m(\frac{g}{\| g \|}).
$$

Введем множество $\wp(\Pi) $  кривых $r(\cdot) $  на единичной
сфере $S^{n-1}_1(0)$, параметризованных естественным образом с
параметром $t \in [0, T]$, проекция которых на плоскость $\Pi$,
проходящую через начало координат, есть кривая, ограничивающая
звездное множество на $\Pi$.

Кривая $r(\cdot)$ параметризована естественным образом и п.в. на
отрезке $[0,T]$ имеет касательную $r'(\cdot)$.

Возьмем произвольную плоскость $\Pi$, $0 \in \Pi$, и произвольную
$r(\cdot) \in \wp(\Pi)$. Введем функцию
$$
\Phi_m(t)=\varphi_m(r(t)) \,\,\, \forall t \in [0,T].
$$
Функции $\Phi_m(\cdot), r(\cdot) $ п.в. дифференцируемы на $[0,
T]$ как липшицевая функция.

Для $r(\cdot) \in \wp(\Pi) $ выделим те ее участки с параметрами
$[T_i, T_{i+1}] \subset [0,T(r)]$, где производная $r'(t), t \in
[0, T(r)],$ составляет с плоскостью $\Pi$ угол $\delta $: $0 <
\delta \leq \frac{\pi }{4} $. Обозначим объединение таких участков
через $\Lambda(r)$.

\begin{thm}
Для того, чтобы последовательность  п.о. первой степени липшицевых
функций $\varphi_m (\cdot):\mathbb{R}^n \rar \mathbb{R} $ с
константой Липшица $L$ сходилась к ПРВ функции, необходимо и
достаточно, чтобы для любой плоскости $\Pi$ и любой кривой
$r(\cdot) \in \wp(\Pi)$, параметризованной естественным образом с
параметром $t \in [0,T], \, T=T(r),$ и любого подмножества
параметров $[T_i, T_{i+1}] \subset \Lambda(r)$, $i \in 1:k$,
нашлись константы $C_1>0$, $ C_2>0$ такие, что \be
 \sum_1^k \vee (\Phi_m'; T_i,T_{i+1})<  C_1 + C_2 \sum_1^m\vee (r';
 T_i,T_{i+1}),
\label{DiffConvSupp5} \ee для всех $m$, где $\Phi_m(t)=\varphi_m
(r(t))$, и производные берутся там, где они существуют.
\label{DiffConvSuppThm4}
\end{thm}

{\bf Доказательство.} Доказательство основывается на  ходе
доказательства Теоремы \ref{diffconvGomFuncMDegNVerThm5}. Мы
приближали функцию  $\varphi_m(\cdot)$ многогранной функцией и
выделяли выпуклые двугранные углы. При доказательстве этой теоремы
константа $C_1$ ограничивала сверху сумму вариаций всех выпуклых
двугранных углов, которые пересекает кривая   $r(\cdot) \in
\wp(\Pi)$. Но кривых из множества $\wp(\Pi) $, как доказано в
Теореме \ref{diffconvGomFuncMDegNVerThm5}, достаточно, чтобы
проверить представимость функции в виде разности выпуклых.
Константа $C_2$ ограничивает сверху норму градиентов линейных
функций, образующих двугранные углы. При выполнении условий
Теоремы \ref{diffconvGomFuncMDegNVerThm5} функции
$\varphi_m(\cdot)$ могут быть представлены как разность выпуклых
п.о. функций, из которых можно выбрать равномерно сходящуюся на
единичном шаре. Предельная функция $\varphi(\cdot)$   будет
выпуклой п.о. функцией, представимой в виде разности выпуклых п.о.
функций. Таким образом, функция $f(\cdot)$ является КВД функцией в
точке $x$. Теорема доказана. $\Box$

\newpage
\vspace{1cm}

\section{Заключение}

\vspace{0.5cm}

Мы рассмотрели различные классы кривых для липшицевых функций. Это
был длинный путь к истине. В итоге мы пришли к классами кривых,
ограничивающих выпуклые компактные множества на сфере
$S_1^{n-1}(0)$ для п.о. функций $m-$ ой степени от $n$ переменных
(теоремы \ref{difconvgomogthm1}, \ref{difconvgomogthm3},
\ref{difconvgomogMDegreethm3}, \ref{diffconvGomFuncMDegNVerThm1})
либо выпуклые компактные множества на плоскости $\mathbb{R}^2$ для
произвольной липшицевой функции от двух переменных (теорема
\ref{difconvexthm1}). Для произвольных функций от многих
переменных для формулировки условия представимости их в виде
разности выпуклых достаточно рассмотреть кривые в $\mathbb{R}^n$,
проекции которых на координатные плоскости, образованные двумя
осями координат, ограничивают выпуклые компактные множества на
этих плоскостях. Сформулирована теорема \ref{difconvmultexthm4},
дающая необходимые и достаточные условия представимости таких
функций в виде разности выпуклых. В каждом параграфе дана
геометрическая трактовка полученных результатов через поворот
кривых из рассматриваемых классов на графиках функций (теоремы
\ref{difconvgomogthm4}, \ref{difconvgomogMDegreethm4},
\ref{diffconvGomFuncMDegNVerThm4}, \ref{difconvmultexthm5}).

Рассматривается одно из применений полученных результатов.
Приводятся достаточные условия сходимости липшицевых
квазидифференцируемых (КВД) функций от двух переменных к КВД
функции. Эти результаты могут быть распространены на функции от
произвольного количества аргументов.

Полученные результаты интересны не только для геометров,
занимающихся построением внутренней геометрии поверхностей (см.
\cite{aleksandrov0}, \cite{aleksandrov2}), но для специалистов в
области оптимизации (см. \cite{demvas}, \cite{demrub},
\cite{strekal}), поскольку выпуклые, а также ПРВ функции нашли
широкое применение в разных областях математики и техники из-за
хороших свойств выпуклых функций.

%\vspace{0.2cm}

%\input{DiffConvPicture.tex}


\newpage
\addcontentsline{toc}{section}{Список литературы}


\newpage

\begin{thebibliography}{21}

\bibitem{aleksandrov0} Александров А. Д. Внутренняя геометрия выпуклых поверхностей. -
Гостехиздат, 1948.

\bibitem{aleksandrov1} Александров А.Д. О поверхностях, представимых
в виде разности выпуклых функций //  Изв. АН Каз.ССР 1949. N 3. С.
3-20.

\bibitem{aleksandrov2}  Александров А.Д.  Поверхности, представимые
разностями выпуклых функций //  Докл. АН СССР. 1950. Т. 72, №4.
С.613-616.

\bibitem{alexandrovSecondDeriv} Александров А.Д.  Существование почти всюду
второго дифференциала выпуклой функции и некоторые свойства
выпуклых поверхностей, связанные с этим свойством, ЛГУ. Ученые
записки. Матем. Сер.6 (1939), 3-35.

\bibitem{zalgaller1} Залгаллер В.А. О представимости функции двух
переменных в виде разности выпуклых функций //  Вестник ЛГУ, N 1,
1963, С.44-45.

\bibitem{demvas} Демьянов В.Ф., Васильев Л.В. Недифференцируемая
оптимизация -  М,: Наука, 1981. 384 С.

\bibitem{demrub} Демьянов В.Ф., Рубинов А.М.  Основы негладкого
анализа и квазидифференциальное исчисление. - М.:Наука, 1990. 432
с.

\bibitem{karpelsadovskii}  Карпелевич Ф.И., Садовский Л.Е.
Элементы линейной алгебры и линейного программирования.- М.:Наука,
1967.

\bibitem{clark}  Кларк Ф. Оптимизация и негладкий анализ. М.:Наука,
1988. 280 с.

\bibitem{kolmogorovfomin} Колмогоров А.Н., Фомин С.В. Елементы
теории функций и функционального анализа -  М.: Наука, 1976. 544
С.

\bibitem{kutateladzerubinov} Кутателадзе С.С., Рубинов А.М.
Двойственность Минковского и ее приложения. Новосибирск: Наука.
1976. 254 с.

\bibitem{lupikovdissertation} Лупиков И.М.   Многозначные
отображения, их описание и применение к оптимизации. - Л.:
Автореферат кандидатской диссертации. - ЛГУ. 1985. 16 с.

\bibitem{proudconvex1}  Прудников И.М.  Необходимые и достаточные
условия представимости положительно однородной функции трех
переменных в виде разности выпуклых  функций //  Известия АН РАН
Т.59. N 5, 1992, С.1116-1128.

\bibitem{proudconvex2} Прудников И.М. К вопросу о  представимости
функции двух переменных в виде разности выпуклых функций //
Сибирский математический журнал РАН, 55(6), 2014. С. 1368-1380.

\bibitem{pshenichnyi} Пшеничный Б.Н. Выпуклый анализ и экстремальные
задачи. - М.:Наука, 1980, 320с.

\bibitem{pogorelov1} Погорелов А.В. Дифференциальная геометрия -
M,: Наука, 1974, 176 С.

\bibitem{rocafellar} Рокафеллар Р. Выпуклый анализ. М.: Мир, 1973.

\bibitem{strekal} Стрекаловский А.С. Элементы невыпуклой
оптимизации - Новосибирск: Наука, 2003.

\bibitem{strecal2} Стрекаловский А.С., Янулевич М.В. Глобальный поиск
оптимального управления с терминальным целевым функционалом,
представленным разностью выпуклых функций// Ж. вычисл. матем. и
матем. Физики. 2008. Т. 48. № 7. С. 1187-1201.

\bibitem{aleksandrov3} Aleksandrov A.D., Reshetnyak Yu.,G. General Theory of
Irregular Curves, Kluwer, 1989.

\bibitem{ginchev} Ginchev I., Gintcheva D.  Characterization and
recognition of d.c. functions. J. Glob. Optim. 57, No. 3, 633-647
(2013).

\bibitem{pallaschke} Grzybowski J., Pallaschke D., Urbanski
R.  Characterization of Differences of Sublinear Functions on the
Plane // Thesises of the international conference "Nonlinear
analisys and related questions", pp. 21-24. 2017.


\bibitem{Hartman} Hartman P.  On functions representable as a difference of convex
functions, Pacific J.Math., 9 (1959), pp.707-713.

\bibitem{hiriarturruty}  Hiriart-Urruty J.B. Generalized
differentiability, duality and optimization for problems dealing
with differences of convex functions. In: Convexity and Duality in
Optimization (Groningen, 1984), Lecture Notes in Econom. and Math.
Systems, 256, Springer, Berlin-New York, 1985, pp. 37-70.


\bibitem{veselyzajicek} L.Vesely, L.Zajicek. A non-DC function which is DC
along all convex curves // J. Math. Anal. Appl. 463(2018),
167-175.


\bibitem{voevodin} Воеводин.В. В. Вычислительные основы линейной алгебры. - М.:
Наука, 1977.

\end{thebibliography}
\end{document}